\documentclass{amsart}

\usepackage{amsmath, amssymb, amsthm, amsfonts}

\usepackage{tikz}
\usepackage{tikz-cd}
\usetikzlibrary{arrows,matrix,positioning, calc}
\usetikzlibrary{patterns,decorations.pathreplacing, decorations.markings, decorations.pathmorphing}
\usetikzlibrary{angles, calc, cd, intersections, math, quotes, through}

\definecolor{cof}{RGB}{219,144,71}
\definecolor{pur}{RGB}{186,146,162}
\definecolor{greeo}{RGB}{91,173,69}
\definecolor{greet}{RGB}{52,111,72}

\usepackage{mathtools}

\usepackage{mathrsfs}

\input xy
\xyoption{all}

\usepackage{hyperref}

\swapnumbers
\numberwithin{equation}{section}

\theoremstyle{plain}
\newtheorem{theorem}[subsubsection]{Theorem}
\newtheorem{lemma}[subsubsection]{Lemma}
\newtheorem{prop}[subsubsection]{Proposition}
\newtheorem{cor}[subsubsection]{Corollary}
\newtheorem{conj}[subsubsection]{Conjecture}

\theoremstyle{definition}
\newtheorem{defn}[subsubsection]{Definition}

\newtheorem{remark}[subsubsection]{Remark}
\newtheorem{exam}[subsubsection]{Example}

\newtheorem{defnlemma}[subsubsection]{Definition/Lemma}
\newtheorem{warn}[subsubsection]{Warning}
\newtheorem{set}[subsubsection]{Setting}
\newtheorem{claim}[subsubsection]{Claim}

\def\mathbi#1{\textup{\textbf{\emph #1}}} 

\def\AA{\mathbb{A}}
\def\BB{\mathbb{B}}
\def\CC{\mathbb{C}}

\def\FF{\mathbb{F}}
\def\GG{\mathbb{G}}

\def\II{\mathbb{I}}
\def\JJ{\mathbb{J}}
\def\KK{\mathbb{K}}
\def\LL{\mathbb{L}}

\def\NN{\mathbb{N}}

\def\PP{\mathbb{P}}
\def\QQ{\mathbb{Q}}
\def\RR{\mathbb{R}}
\def\SS{\mathbb{S}}

\def\XX{\mathbb{X}}

\def\ZZ{\mathbb{Z}}

\newcommand\cA{\mathcal{A}}
\newcommand\cB{\mathcal{B}}
\newcommand\cC{\mathcal{C}}
\newcommand\cD{\mathcal{D}}
\newcommand\cE{\mathcal{E}}
\newcommand\cF{\mathcal{F}}
\newcommand\cG{\mathcal{G}}
\newcommand\cH{\mathcal{H}}
\newcommand\cI{\mathcal{I}}
\newcommand\cJ{\mathcal{J}}
\newcommand\cK{\mathcal{K}}
\newcommand\cL{\mathcal{L}}
\newcommand\cM{\mathcal{M}}
\newcommand\cN{\mathcal{N}}
\newcommand\cO{\mathcal{O}}
\newcommand\cP{\mathcal{P}}
\newcommand\cQ{\mathcal{Q}}
\newcommand\cR{\mathcal{R}}
\newcommand\cS{\mathcal{S}}
\newcommand\cT{\mathcal{T}}
\newcommand\cU{\mathcal{U}}
\newcommand\cV{\mathcal{V}}
\newcommand\cW{\mathcal{W}}
\newcommand\cX{\mathcal{X}}
\newcommand\cY{\mathcal{Y}}

\def\bI{\mathbf{I}}
\def\bJ{\mathbf{J}}

\def\bL{\mathbf{L}}

\def\bP{\mathbf{P}}

\def\bR{\mathbf{R}}

\def\bT{\mathbf{T}}

\def\bl{\textbf{\em{l}}}

\def\bm{\mathbf{m}}

\newcommand\frD{\mathfrak{D}}

\newcommand\frF{\mathfrak{F}}
\newcommand\frG{\mathfrak{G}}

\newcommand\frI{\mathfrak{I}}

\newcommand\frK{\mathfrak{K}}

\newcommand\frP{\mathfrak{P}}
\newcommand\frQ{\mathfrak{Q}}

\newcommand\frS{\mathfrak{S}}
\newcommand\frT{\mathfrak{T}}
\newcommand\frU{\mathfrak{U}}

\newcommand\frW{\mathfrak{W}}

\newcommand\frZ{\mathfrak{Z}}

\newcommand\fra{\mathfrak{a}}
\newcommand\frb{\mathfrak{b}}
\newcommand\frc{\mathfrak{c}}

\newcommand\frf{\mathfrak{f}}
\newcommand\frg{\mathfrak{g}}
\newcommand\frh{\mathfrak{h}}

\newcommand\frj{\mathfrak{j}}

\newcommand\frl{\mathfrak{l}}
\newcommand\fm{\mathfrak{m}}
\newcommand\frn{\mathfrak{n}}

\newcommand\frq{\mathfrak{q}}

\newcommand\frs{\mathfrak{s}}
\newcommand\frt{\mathfrak{t}}

\newcommand\frz{\mathfrak{z}}

\newcommand\tilW{\widetilde{W}}

\newcommand{\dG}{G^{\vee}}
\newcommand{\dB}{B^{\vee}}

\newcommand{\dT}{T^{\vee}}
\newcommand{\dH}{H^{\vee}}
\newcommand{\dL}{L^{\vee}}

\newcommand\aff{\textup{aff}}

\newcommand\alg{\textup{alg}}
\newcommand\Alg{\textup{Alg}}
\newcommand\Alt{\textup{Alt}}
\newcommand\AS{\textup{AS}}

\newcommand{\Bun}{\textup{Bun}}

\newcommand{\ch}{\textup{char}}

\newcommand{\codim}{\textup{codim}}
\newcommand{\Coh}{\textup{Coh}}
\newcommand{\coker}{\textup{coker}}

\newcommand{\Corr}{\textup{Corr}}

\newcommand\ev{\textup{ev}}

\newcommand\Gal{\textup{Gal}}

\newcommand{\Gr}{\textup{Gr}}

\newcommand\id{\textup{id}}
\renewcommand{\Im}{\textup{Im}}
\newcommand{\Ind}{\textup{Ind}}

\newcommand\inv{\textup{inv}}

\newcommand{\inin}{\textup{in}}
\newcommand\Jac{\textup{Jac}}
\newcommand\Lie{\textup{Lie}\ }
\newcommand\Loc{\textup{Loc}}

\newcommand\Mod{\textup{Mod}}

\newcommand{\Nm}{\textup{Nm}}

\newcommand{\out}{\textup{out}}

\newcommand{\Perf}{\textup{Perf}}

\newcommand{\Pic}{\textup{Pic}}
\newcommand\pr{\textup{pr}}

\newcommand\pt{\textup{pt}}

\newcommand{\red}{\textup{red}}
\newcommand{\reg}{\textup{reg}}
\newcommand\Rep{\textup{Rep}}

\newcommand\res{\textup{res}}

\newcommand\rs{\textup{rs}}

\newcommand\Span{\textup{Span}}
\newcommand\Spec{\textup{Spec}\ }

\newcommand\st{\textup{st}}
\newcommand\Stab{\textup{Stab}}

\newcommand\Sym{\textup{Sym}}
\newcommand{\tors}{\textup{tors}}
\newcommand\Tot{\textup{Tot}}

\newcommand\triv{\textup{triv}}

\newcommand{\Vect}{\textup{Vect}}

\newcommand\Aut{\textup{Aut}}
\newcommand\Hom{\textup{Hom}}
\newcommand\End{\textup{End}}

\newcommand\uAut{\underline{\Aut}}

\newcommand{\proj}{\textup{proj}}
\newcommand{\Sh}{\textup{Sh}}
\newcommand{\Real}{\textup{Re}}
\newcommand{\Imagine}{\textup{Im}}

\newcommand{\sfC}{\mathsf{C}}
\newcommand{\sfc}{\mathsf{c}}
\newcommand{\sfH}{\mathsf{H}}
\newcommand{\sfS}{\mathsf{S}}
\newcommand{\sfZ}{\mathsf{Z}}

\newcommand{\scU}{\mathscr{U}}

\newcommand\GL{\textup{GL}}

\newcommand\PGL{\textup{PGL}}

\newcommand\SL{\textup{SL}}
\renewcommand\sl{\mathfrak{sl}}

\newcommand\Sp{\textup{Sp}}

\newcommand{\Gm}{\GG_m}
\def\Ga{\GG_a}

\newcommand{\ad}{\textup{ad}}
\newcommand{\Ad}{\textup{Ad}}
\renewcommand\sc{\textup{sc}}
\newcommand{\der}{\textup{der}}

\newcommand\xch{\mathbb{X}^*}
\newcommand\xcoch{\mathbb{X}_*}

\newcommand{\incl}{\hookrightarrow}
\newcommand{\isom}{\stackrel{\sim}{\to}}

\newcommand{\surj}{\twoheadrightarrow}

\newcommand{\Qlbar}{\overline{\QQ}_\ell}

\newcommand{\twtimes}[1]{\stackrel{#1}{\times}}

\renewcommand{\j}[1]{\langle{#1}\rangle}
\newcommand{\wt}[1]{\widetilde{#1}}
\newcommand{\wh}[1]{\widehat{#1}}
\newcommand\quash[1]{}
\newcommand\mat[4]{\left(\begin{array}{cc} #1 & #2 \\ #3 & #4 \end{array}\right)}  % 2-by-2 matrix
\newcommand\un{\underline}

\newcommand{\ov}{\overline}
\newcommand{\bs}{\backslash}
\newcommand{\tl}[1]{[\![#1]\!]}
\newcommand{\lr}[1]{(\!(#1)\!)}
\newcommand\sss{\subsubsection}
\newcommand\xr{\xrightarrow}
\newcommand\op{\oplus}
\newcommand\ot{\otimes}

\newcommand{\sslash}{\mathbin{/\mkern-6mu/}}
\renewcommand\c\circ
\newcommand\sne{\subsetneq}
\newcommand\vn{\varnothing}

\newcommand\Nb{\mathrm{Nb}}

\newcommand\lng{\langle}
\newcommand\rng{\rangle}

\newcommand\trg{\triangle}

\newcommand{\cohog}[2]{\textup{H}^{#1}({#2})}     % plain group
\newcommand{\cohoc}[2]{\textup{H}_{c}^{#1}({#2})}     % compact support
\newcommand\upH{\textup{H}}

\newcommand{\oll}[1]{\overleftarrow{#1}}
\newcommand{\orr}[1]{\overrightarrow{#1}}

\newcommand{\ovl}[1]{\overline{#1}}

\renewcommand\a\alpha
\renewcommand\b\beta
\newcommand\g\gamma
\newcommand\G\Gamma
\renewcommand\d\delta
\newcommand\D\Delta
\newcommand{\e}{\epsilon}
\newcommand{\io}{\iota}
\renewcommand{\k}{\kappa}
\renewcommand{\th}{\theta}

\newcommand{\ph}{\varphi}
\renewcommand\r\rho
\newcommand{\s}{\sigma}
\newcommand{\Sig}{\Sigma}
\renewcommand{\t}{\tau}

\newcommand{\y}{\eta}
\newcommand{\z}{\zeta}
\newcommand{\ep}{\epsilon}

\renewcommand{\l}{\lambda}

\newcommand{\om}{\omega}
\newcommand{\Om}{\Omega}

\newcommand\dm{\diamondsuit}
\newcommand\hs{\heartsuit}

\newcommand\sh{\sharp}
\newcommand\da{\dagger}

\newcommand{\hol}{\mathrm{hol}}
\newcommand{\Graph}{\mathrm{Graph}}

\newcommand{\Fun}{\textup{Fun}}
\newcommand{\QCoh}{\textup{QCoh}}
\newcommand{\Tor}{\textup{Tor}}
\newcommand{\Cone}{\textup{Cone}}

\newcommand{\Cat}{\textup{Cat}}
\newcommand{\Gpd}{\textup{Gpd}}

\newcommand{\sfh}{\mathsf{h}}

\newcommand{\sfL}{\mathsf{L}}

\newcommand{\sfN}{\mathsf{N}}
\newcommand{\sfI}{\mathsf{I}}

\newcommand{\fin}{\mathsf{fin}}

\newcommand{\CAlg}{\mathrm{CAlg}}

\newcommand{\CatEx}{\Cat_\infty^{\mathrm{st}}}

\newcommand{\cores}{\mathrm{co}\text{-}\mathrm{res}}

\newcommand\Wa{W_{\aff}}
\newcommand\colim{\mathrm{colim}}

\newcommand\vk{\varkappa}
\newcommand\ab{\mathrm{ab}}

\newcommand{\fF}{\mathfrak{F}}

\newcommand{\cpt}{\mathrm{cpt}}
\newcommand{\bq}{\mathbf{q}}
\newcommand{\bp}{\mathbf{p}}

\newcommand{\dagg}{\dagger}
\newcommand{\dda}{\ddagger}

\newcommand{\fz}{\mathfrak{z}}

\newcommand{\std}{\mathrm{std}}

\newcommand{\Kos}{\textup{Kos}}
\newcommand{\mt}{\mapsto}
\newcommand{\lmod}{\textup{-mod}}
\newcommand{\Grot}{\mathbb{G}^\textup{rot}_m}
\newcommand{\ani}{\textup{ani}}
\newcommand{\Iu}{\bI^+_0, \bI^+_\infty}
\newcommand{\bt}{\boxtimes}
\newcommand{\frev}{\mathfrak{ev}}
\newcommand{\aut}{\mathfrak{aut}}
\newcommand{\coev}{\textup{coev}}
\newcommand{\vs}{\textup{vs}}
\newcommand{\cod}{cod}

\newcommand{\dGad}{G^{\vee,\ad}}

\newcommand{\dGreg}{G^{\vee,\reg}}

\newcommand{\dBad}{B^{\vee,\ad}}

\newcommand{\sft}{\subset_{ft}}

\newcommand{\geqsl}{\geqslant}
\newcommand{\leqsl}{\leqslant}

\newcommand{\CI}{\cC\cI}
\newcommand{\Hsm}{H^{sm}}

\newcommand{\bRI}{b_{\RR}^I}
\newcommand{\PfI}{\cP_{ft}(\wt{I})}
\newcommand{\PdI}{\cP^\dagg(I)}
\newcommand{\PdIH}{\cP^\dagg(I_H)}

\newcommand{\rigid}{\textup{rig}}

\newcommand{\Spc}{\textup{Spc}}
\newcommand{\cPI}{\cP\cI}
\newcommand{\Maps}{\textup{Maps}}

\newcommand{\LamJ}{\Lambda_{L_J}^{U_J}}
\newcommand{\LamJflat}{\Lambda_{L_J^\flat}^{U_J}}
\newcommand{\hP}{\mathrm{h}\cP}
\newcommand{\ovs}{\overset}
\newcommand{\lrar}{\longrightarrow}

\newcommand{\hrar}{\hookrightarrow}
\newcommand{\car}{\circlearrowright}
\newcommand{\cal}{\circlearrowleft}

\newcommand{\cSHfv}{\cS_{H^{\flat, \vee}}}

\newcommand{\fflat}{{\flat\flat}}
\newcommand{\LJf}{L_J^\flat}

\newcommand{\yL}{{y_L}}
\newcommand{\yH}{{y_H}}
\newcommand{\yf}{{y^\flat}}

\newcommand{\yJf}{{y_J^\flat}}

\newcommand{\zfO}{z_{1}^\flat}
\newcommand{\zfT}{z_{2}^\flat}
\newcommand{\zfi}{z_{i}^\flat}
\newcommand{\zf}{z^\flat}
\newcommand{\zfv}{z^\flat_\vn}
\newcommand{\zvnf}{z_\vn^\flat}
\newcommand{\zvn}{z_\vn}

\newcommand{\yHf}{y_H^\flat}
\newcommand{\yLf}{y_L^\flat}
\newcommand{\yLpf}{y_L^{\prime, \flat}}
\newcommand{\yHpf}{y_H^{\prime, \flat}}

\newcommand{\yfO}{y_{1}^\flat}
\newcommand{\yfT}{y_{2}^\flat}
\newcommand{\yfi}{y_{i}^\flat}

\newcommand{\yfv}{y^{\flat}_\vn}
\newcommand{\yffv}{y^{\fflat}_\vn}
\newcommand{\yvnf}{y^{\flat}_\vn}
\newcommand{\yfvn}{y^{\flat}_\vn}
\newcommand{\yff}{y^{\fflat}}

\newcommand{\yvnJ}{{y^{J}_\vn}}
\newcommand{\yvnJf}{{y^{J,\flat}_\vn}}

\newcommand{\gaOT}{\gamma_{1,2}}
\newcommand{\gaTO}{\gamma_{2,1}}

\newcommand{\gafTO}{\gamma^\flat_{2,1}}
\newcommand{\gafOT}{\gamma^\flat_{1,2}}

\newcommand{\Hf}{{H^\flat}}
\newcommand{\frPH}{\frP_{H}}
\newcommand{\frPHf}{\frP_{\Hf}}

\newcommand{\UHSc}{U_H^{\SS, \cpt}}
\newcommand{\UHfS}{U_{H^\flat}^{\SS}}
\newcommand{\UHfSc}{U_{H^\flat}^{\SS, \cpt}}
\newcommand{\UfS}{U^{\flat, \SS}}

\newcommand{\UJfS}{U_J^{\flat, \SS}}

\newcommand{\wUpyf}{\wt{\Upsilon}_{\yf}}
\newcommand{\wUpyfO}{\wt{\Upsilon}_{\yfO}}
\newcommand{\wUpyfT}{\wt{\Upsilon}_{\yfT}}
\newcommand{\Upyf}{\Upsilon_{\yf}}
\newcommand{\UpyfO}{\Upsilon_{\yfO}}
\newcommand{\UpyfT}{\Upsilon_{\yfT}}
\newcommand{\wUpy}{\wt{\Upsilon}_{y}}
\newcommand{\wUpyO}{\wt{\Upsilon}_{y_1}}
\newcommand{\wUpyT}{\wt{\Upsilon}_{y_2}}
\newcommand{\UpPf}{\Upsilon_{\frP^{op}, \flat}}
\newcommand{\wUpPf}{\wt{\Upsilon}_{\frP^{op}, \flat}}

\newcommand{\cTH}{\cT_H}
\newcommand{\THf}{\cT_{H^\flat}}
\newcommand{\WTH}{\cW(\cT_H)}
\newcommand{\WTHp}{\cW(\cT'_H)}
\newcommand{\WTHf}{\cW(\cT_{H^\flat})}

\newcommand{\WTLJf}{\cW(\cT_{L_J^\flat})}
\newcommand{\WP}{\cW_{\frP}}
\newcommand{\WPHf}{\cW_{\frP, \Hf}}
\newcommand{\WPH}{\cW_{\frP, H}}
\newcommand{\WPL}{\cW_{\frP, L}}
\newcommand{\WPHL}{\cW_{\frP, H, L}}

\newcommand{\wPhiHP}{\wt{\Phi}_{H, \frP^{op}}}
\newcommand{\PhiHP}{\Phi_{H, \frP^{op}}}
\newcommand{\IWf}{\cI\cW_\flat}
\newcommand{\Qf}{\cQ_\flat}
\newcommand{\Qft}{\cQ_\flat^t}

\newcommand{\QPf}{\cQ_{\frP^{op}}^{\flat}}
\newcommand{\wQPf}{\wt{\cQ}_{\frP^{op}}^{\flat}}
\newcommand{\pf}{p^\flat}
\newcommand{\Tf}{T^{\flat}}
\newcommand{\Tfv}{T^{\flat, \vee}}
\newcommand{\Hfv}{H^{\flat, \vee}}
\newcommand{\Lfvn}{L^{\flat}_{\vn}}
\newcommand{\Lffvn}{L^{\fflat}_{\vn}}
\newcommand{\bia}{\mathbi{a}}
\newcommand{\uds}{\underset}

\newcommand{\btau}{\boldsymbol{\tau}}
\newcommand{\CCline}{\CC\text{-lines}}

\newcommand{\chiO}{{\chi_1}}

\newcommand{\Lf}{L^\flat}
\newcommand{\Hff}{H^{\fflat}}

\newcommand{\Kdu}{K^{*}}

\newcommand{\frGK}{\frG^{K^*}}
\newcommand{\frGKt}{\frG^{\Kdu}}
\newcommand{\Gpds}{\mathrm{Gpds}}
\newcommand{\Ham}{\mathrm{Ham}}

\newcommand{\bil}{\mathbi{l}}
\newcommand{\tsfa}{\textsf{a}}

\newcommand{\AHft}{A_{\Hf}^t}
\newcommand{\Avnt}{A^{\vn, t}}
\newcommand{\AHf}{A_{\Hf}}
\newcommand{\ALf}{A_{\Lf}}
\newcommand{\Avn}{A^{\vn}}

\newcommand{\iotvn}{\iota_{\vn}^t}
\newcommand{\iotHf}{\iota_{\Hf}^t}
\newcommand{\dmd}{\diamond}

\newcommand{\CF}{\mathrm{CF}}
\newcommand{\sfF}{\mathsf{F}}
\newcommand{\ZH}{Z(H)}
\newcommand{\ZHc}{Z(H)_\cpt}
\newcommand{\ZHfc}{Z(H^\flat)_\cpt}
\newcommand{\ZHf}{Z(H^\flat)}
\newcommand{\SHv}{\cS_{H^{\vee}}}
\newcommand{\SHfv}{\cS_{H^{\flat, \vee}}}

\newcommand{\frWP}{\frW_{\frP}}
\newcommand{\frIW}{\frI\frW}

\newcommand{\hor}{\mathrm{hor}}

\newcommand{\Isom}{\mathrm{Isom}}

\newcommand{\Isog}{\mathrm{Isog}}

\newcommand{\pathz}{\gamma}

\newcommand{\IndW}{\Ind\cW}

\newcommand{\CRez}{\CC_{\Real z\geqsl 0}}
\newcommand{\bd}{\mathrm{bd}}
\newcommand{\tame}{\textup{tame}}

\title{Mirror symmetry of the affine Toda systems}

\dedicatory{}
\author{Xin Jin}
\address{Department of Mathematics, Boston College, Chestnut Hill, MA 02467}
\email{xin.jin@bc.edu}
\author{Zhiwei Yun}
\thanks{Supported partially by the Simons Investigatorship.}
\address{Department of Mathematics, Massachusetts Institute of Technology, 77 Massachusetts Ave, Cambridge, MA 02139}
\email{zyun@mit.edu}
\date{}

\begin{document}

\begin{abstract}
For a complex reductive group $G$, we prove a homological mirror symmetry between the wrapped Fukaya category of the affine Toda system for $G$ and coherent sheaves on the regular centralizer group scheme for the Langlands dual group $\dG$. This can be interpreted as a geometric Langlands equivalence for $\PP^1$ with mildest wild ramification at $0$ and $\infty$. 
\end{abstract}

\maketitle

\tableofcontents

\section{Introduction}

\subsection{Toda system and its constructible quantization}
The classical Toda system for a complex reductive group $G$, as formulated by Kostant \cite{Kos}, is the group scheme of regular centralizers $\pi: \cJ_G\to\frc:=\frg\sslash G$. More precisely,  the fiber of $\pi$ over a point $a\in \frc$ is the centralizer $C_G(\wt a)$ where $\wt a\in \frg$ is any {\em regular} element with image $a$ in $\frc$. 

Let $U$ be a maximal unipotent subgroup of $G$ and $\psi: U\to \Ga$ is a homomorphism that is nontrivial on each simple root subgroup. We call the ``double coset'' 
\begin{equation}\label{biWh}
(U,\psi)\bs G/(U,\psi)
\end{equation}
the {\em bi-Whittaker ``space''} of $G$. With the presence of $\psi$ it is merely a symbol at this point, but the category of sheaves on the bi-Whittaker space can be made sense of in various contexts. The relation between $\cJ_G$ and the bi-Whittaker space is clarified by the isomorphism
\begin{equation*}
\cJ_G\cong T^*((U,\psi)\bs G/(U,\psi)):=(T^*G)//_{(\psi,-\psi)}U\times U.
\end{equation*}
Here the right side means the Hamiltonian reduction of $T^*G$ under the left and right action of $U$ by taking the preimage of $(d\psi,-d\psi)\in \frn^*\times \frn^*$  under the moment map and quotienting by $U\times U$. 

Therefore, a (constructible) quantization of $\cJ_G$ should be the category of sheaves on the bi-Whittaker space \eqref{biWh}.
In the $D$-module context, it means the derived category of $D$-modules on $G$ with left and right equivariant structures with respect to $U$ against pullback of the exponential $D$-module along $\psi$. If one works with $\ell$-adic sheaves over schemes in characteristic $p$, the exponential $D$-module gets replaced by the Artin-Schreier sheaf on $\Ga$ and the similar left and right $(U,\psi)$-equivariant derived category on $G$ makes sense. Alternatively, one can work with the wrapped Fukaya category $\cW(\cJ_G)$ of $\cJ_G$ as a model for microsheaves, which gives yet another constructible quantization of $\cJ_G$. We call these constructible quantizations of $\cJ_G$ {\em bi-Whittaker categories} of $G$.

In all three contexts ($D$-modules, $\ell$-adic and Fukaya), various authors have proved equivalences between the constructible quantization of $\cJ_G$ and coherent sheaves on spaces that are very close to $\dT\sslash W$, where $\dT$ is the Langlands dual of a maximal torus $T$ of $G$. In the $D$-module context, an analogous result was proved by Ginzburg \cite{Gin}, Lonergan \cite{Lon-biWh}, Ben-Zvi
 and Gunningham \cite{BZG} and Gannon \cite{Gan}. In the $\ell$-adic setting, an analogous result was proved by Bezrukavnikov and Deshpande \cite{BD}. Most relevant to the current paper is the following mirror description of the bi-Whittaker category in the Fukaya context, proved by the first-named author. We state a simplified version here.

\begin{theorem}[{\cite[Theorem 7.7]{J}}]\label{th:intro biWh}
Assume $G$ has connected center, then we have a canonical equivalence
\begin{align*}
\cW(\cJ_G)\simeq \Coh(T^\vee\sslash W)
\end{align*}
sending the unit section of $\cJ_G$ to the structure sheaf $\cO_{T^\vee\sslash W}$. 
\end{theorem}

\subsection{Main constructions, results and conjectures}
The purpose of this paper is two-fold: we study the geometry and cohomology of the affine Toda system $\cM_G$, and we prove a homological mirror symmetry theorem for $\cM_G$.

\sss{Affine Toda system} For any split reductive group $G$ of rank $r$ over a field $k$, we introduce a symplectic algebraic variety $\cM_G$ of dimension $2r$, the {\em affine Toda system}. The affine Toda system is birationally equivalent to the {\em periodic Toda system} defined in \cite{Bog} that has been studied extensively in the mathematical physics literature. The space $\cM_G$ is constructed as the moduli space of $G$-Higgs bundles over $\PP^1$ with specific singularities at $0$ and $\infty$. We show that it has the following features:

\begin{itemize}
\item $\cM_G$ is glued from (torsors of) Toda systems of $G$ as well as its pseudo-Levi subgroups.
\item There is a proper map $f_G: \cM_G\to \cA_G$, the Hitchin map, to an affine space $\cA_G$ that realizes $\cM_G$ as a completely integrable system.  
General fibers of $f_G$ are finite unions of abelian varieties of dimension $r$. This integrable system is birationally equivalent to the periodic Toda system in type $A$, whose solutions are expressed in terms of theta functions on the Jacobians of spectral curves, which are the fibers of the Hitchin map.

\item $\cM_G$ admits another affine map $\cM_G\to \PP$ to a weighted projective space $\PP$ of dimension $r$ that is also a completely integrable system. This map records the isomorphism types of the underlying $G$-bundles (with level structure) of the Higgs bundles.

\item Combining the two integrable systems gives a finite flat map $\nu: \cM_G\to \cA_G\times \PP$. It has degree $|W|$ on each connected component of $G$.
\end{itemize}

\sss{Mirror symmetry for $\cM^\c_G$} We now work with the base field $k=\CC$. The mirror dual to the neutral component $\cM^\c_G$ of $\cM_G$ turns out to be the group-version of the regular centralizer group scheme of $\dG$, the Langlands dual group of $G$. When the derived group of $\dG$ is simply-connected, $J^{\dGad}_{\dG}\to \dG\sslash \dG$ is the family of centralizers in $\dGad$ along a Steinberg section $\dG\sslash \dG\to \dGreg$. For general $\dG$, the definition of $J^{\dGad}_{\dG}$ involves a stacky version of the base $\dG\sslash \dG$, which we elaborate in \S\ref{ss:JG}.

\begin{theorem}[see Theorem \ref{th:mirror}]\label{th:intro} There is an equivalence of dg-categories over $\CC$: 
\begin{equation}\label{mirror intro}
\cW(\cM^\circ_G)\cong \Coh(J^{\dGad}_{\dG}).
\end{equation}
Here, $\cW(\cM^\circ_G)$ denotes wrapped Fukaya category of the neutral component $\cM^\circ_G$ of $\cM_G$, to be defined precisely in \S\ref{s:Fukaya}.  
\end{theorem}

\sss{Cohomology and P=W conjecture for $\cM^\c_G$}

We also compute the cohomology of $\cM^\c_G$ (the neutral component of $\cM_G$) explicitly in Corollary \ref{c:coho MG}. In particular we show that when $G$ is semisimple, $\cohog{*}{\cM_G}$ is concentrated in even degrees and pure. When $G$ is simple adjoint, we show that the Poincar\'e polynomial of $\cM^\c_G$ is
\begin{equation*}
    \sum_{d\in \NN, \wt I_d\ne \vn}t^{2(r+1-r_d)}[r_d]_{t^2}\ph(d).
\end{equation*}
Here $r_d=|\wt I_d|$, where $\wt I_d$ is the set of affine simple roots of $G$ whose Dynkin labeling is divisible by $d$; $[n]_{t^2}=1+t^2+\cdots+t^{2(n-1)}$; $\ph$ is the Euler function.

On the other hand, consider the regular centralizer group scheme $J^G_{G^\sc}$ (now on the $G$-side; $G^\sc$ being the simply-connected cover of $G$). The point-counting of $J^G_{G^\sc}$ over a finite field was done by Lusztig \cite{L-Coxeter}, and the result shows striking similarities with the point-counting of $\cM^\c_G$.  Motivated by the Simpson correspondence \cite{Simpson} and the P=W conjecture of de Cataldo--Hausel--Migliorini \cite{dCHM}, we make the following conjecture to explain the similarities.

\begin{conj}[see Conjectures \ref{c:HK} and \ref{c:P=W}]\label{conjIntrohyperkahler} Let $G$ be a semisimple group over $\CC$ and $G^\sc$ be its simply-connected cover.
\begin{enumerate}
    \item There exists a hyperK\"ahler structure on $\cM^\c_G$ such that $J^{G}_{G^\sc}$ as a complex manifold is obtained from $\cM^\c_G$ by  hyperK\"ahler rotation. In particular, there is a canonical isomorphism on cohomology
    \begin{equation*}
        \cohog{*}{\cM^\c_G,\QQ}\isom \cohog{*}{J^G_{G^\sc}, \QQ}.
    \end{equation*}
    \item The above isomorphism sends $P_{i}\cohog{*}{\cM^\c_G,\QQ}$, the $i^{\textup{th}}$ piece of the perverse filtration associated to the Hitchin fibration $\cM^\c_G\to \cA_G$, isomorphically to the $2i^{\textup{th}}$ piece of the weight filtration $W_{2i}\cohog{*}{J^G_{G^\sc},\QQ}$, for $i=0,1,\cdots, 2r$.
\end{enumerate}
\end{conj}

\subsection{Connection with ramified geometric Langlands correspondence}
The equivalence in Theorem \ref{th:intro} can be viewed as a ramified geometric Langlands correspondence, where the left side of \eqref{mirror intro} is a microlocal incarnation of an automorphic category, and on the right side, $J^{\dGad}_{\dG}$ can be interpreted as a (wild) character variety.

More precisely, we consider on the one hand the moduli stack $\Bun_G(\bI^+_0, \bI^+_\infty)$ of $G$-bundles on $\PP^1$ with level structures at $0$ and $\infty$ given by the pro-unipotent radicals $\bI^+_0$ and $\bI^+_\infty$ of the respective Iwahori subgroups. The stack $\Bun_G(\bI^+_0, \bI^+_\infty)$ can be thought of as an affine analogue of $U\bs G/U$, or rather $U^-\bs G/U$, where $U^-$ is opposite to $U$. Let $\psi_0: \bI^+_0\to \Ga$ and $\psi_\infty: \bI^+_\infty\to \Ga$ be homomorphisms that are nontrivial on each affine root subgroup. In analogy with the bi-Whittaker ``space'' \eqref{biWh}, we formally write 
\begin{equation}\label{intro BunG}
\Bun_G((\bI^+_0,\psi_0), (\bI^+_\infty,\psi_\infty))
\end{equation}
whose precise meaning can only be given when we consider sheaves on it. The space $\cM_G$ can be identified with the cotangent stack of \eqref{intro BunG}, hence $\cW(\cM_G)$ is a microlocal model for sheaves on \eqref{intro BunG}. Similarly one can consider $D$-modules or $\ell$-adic sheaves on \eqref{intro BunG} in characteristic $p$. This justifies calling the left side of \eqref{mirror intro} the automorphic side \footnote{One may wonder about a purely affine analog of the bi-Whittaker space, namely $(\bI^+,\psi)\bs G\lr{t}/(\bI^+,\psi)$. However, sheaves on this space are not interesting enough: when $G$ is simply-connected, all sheaves are supported on the double cosets of central elements of $G$.}.
As further numerical justification, in \S\ref{ss:auto} we consider the situation over a finite field $\FF_q$ and relate the point-counting on $\cM_G$ to the counting of automorphic forms with specific local conditions. 

On the other hand, $J^{\dGad}_{\dG}$ can be identified with the fiber product
\begin{equation}\label{intro wild char}
\frac{\dB w\dB}{\dBad}\times_{\frac{\dG}{\dGad}} \frac{\dB w\dB}{\dBad}
\end{equation}
where $\dB\subset \dG$ is a Borel subgroup, and $w\in W$ is any Coxeter element. This follows from the fact that $\dB w\dB$ is contained in the regular locus of $\dG$, and intersects each regular conjugacy class in exactly one $\dBad$-orbit which is also a free orbit. Now $\frac{\dB w\dB}{\dBad}$ may be interpreted as the moduli stack of Stokes data associated with a specific irregular type of local $G$-connections (called {\em isoclinic} of slope $1/h$ in the terminology of \cite{JY}, where $h$ is the Coxeter number of $G$). The fiber product  \eqref{intro wild char} then classifies $\dG$-local systems on $\PP^1-\{0,\infty\}$ together with Stokes data at $0$ and $\infty$ specified by the Coxeter element $w$. This justifies calling the right side of \eqref{mirror intro} as the spectral side. 

The matching of the automorphic and spectral sides is dictated by the local matching: the $(\bI_0^+,\psi_0)$-level structure on the automorphic side corresponds to the local wild character variety $\frac{\dB w\dB}{\dB}$. In \cite{BBAMY}, the authors propose a ramified geometric Langlands correspondence also over $\PP^1$ but with tame ramification at $0$ and wild ramification at $\infty$. The local matching at $\infty$ of the automorphic side and spectral side in \cite{BBAMY} includes the current situation as a special case.

\subsection{Geometric features of the homological mirror symmetry}
We comment on some of the geometric features of our HMS results and some connections to the literature. 

\sss{HMS of multiplicative Coulomb branches}
When $G$ is semisimple, the B side of the mirror symmetry $J^{G^{\vee,\ad}}_{G^\vee}$ is a (twisted) multiplicative Coulomb branch (MCB). As in Conjecture \ref{conjIntrohyperkahler}, $\cM_G^\circ$ is expected to be a hyperK\"ahler rotation of $J_{G^{\sc}}^{G}$. In particular, by viewing the real part of the holomorphic symplectic form on $\cM_G^\circ$ as the K\"ahler form on $J_{G^{\sc}}^{G}$, we can state our main result as HMS between MCB as follows. 

\begin{theorem}
Assume Conjecture \ref{conjIntrohyperkahler} (i) holds. For $G$ semisimple, there is an equivalence of categories 
\begin{align*}
\cW(J_{G^{\sc}}^G)\simeq \Coh(J_{G^{\vee}}^{G^{\sc,\vee}}). 
\end{align*}
\end{theorem}

We briefly mention some of the recent developments on HMS of MCB. 

First, in the abelian case, by the work of Gammage, McBreen and Webster \cite{MW, GMW}, HMS of multiplicative hypertoric varieties is established. Their approach  is based on the calculation of the Lagrangian skeleton and the category of microlocal sheaves on it.

Second, in the work of Aganagic, Danilenko, Li, Shende and Zhou \cite{ADLSZ}, HMS of MCB from quiver gauge theory, equipped with a superpotential, is established. 
They construct an embedding of a certain cylindrical KLRW category into the corresponding Fukaya-Seidel category, and the B-model is a resolved additive Coulomb branch. 
Their approach is based on counting $J$-holomorphic disks in symmetric products of surfaces. When $G$ is of type $A$ with matter $0$, there should be a close relation between their result and ours. However, our approach is quite different (see \S\ref{subsecPfsketch} below), which essentially avoids counting holomorphic disks.

\sss{SYZ mirror symmetry}
Strominger--Yau--Zaslow (SYZ) fibrations and $T$-duality is 
one of the guiding principles in the study of mirror symmetry. 
The Hitchin map $f^\c_G: \cM^\c_G\to \cA_G$ is an example of a SYZ fibration.
We comment on how to understand our HMS results conceptually from SYZ mirror symmetry.

Assume $G$ is adjoint, e.g. $\PGL_n$, there is a canonical Hitchin section $S_I$ to the Hitchin fibration that generates the wrapped Fukaya category of $\cM^\c_G$.  
The calculation of the mirror reduces to wrapping the Hitchin section and the calculation of the endomorphism algebra $\End(S_I)$.
In general, it is a challenging task to calculate the endomorphism algebra directly. Our main theorem can be stated as an isomorphism  $\End(S_I)\cong \cO(J^{\dGad}_{G^\vee})$, and our calculation is indeed indirect.

\sss{Lagrangian skeleta: Example of $G=\PGL_2$}

Lastly, we explain how to get a Lagrangian skeleton of $\cM_G^\circ$ in the case of $G=\PGL_2$, where the geometry is quite explicit. 
In this case, $\cM_G^\circ$ has real dimension $4$, and the SYZ fibration has exactly two singular fibers, each with a nodal singularity. 
There are two ways to construct a Lagrangian skeleton. 

First, using the SYZ fibration, we can construct a Lagrangian skeleton from a smooth fiber with two attaching Lagrangian thimbles. The skeleton is homeomorphic to $D^2\cup_{c_+} T^2\cup_{c_-} D^2$, where $c_+$ and $c_-$ are two primitive 1-cycles on $T^2$ with $c_+\cap c_-=\pm 2$. 

We can also construct the skeleton from a ``special" Weinstein structure (cf. Remark \ref{remSpecialZflatHandle}). The skeleton is homeomorphic to $\RR\PP^2\cup_{\RR\PP^1} \RR\PP^2$, where $\RR\PP^1$ is the equator in each $\RR\PP^2=S^2/(\ZZ/2\ZZ)$. 

These two skeleta, though look quite different, are both homotopy equivalent to $\RR\PP^2\vee S^2$. In general, it is not easy to describe Lagrangian skeleta from the SYZ fibrations.  On the other hand, the second way of constructing skeleta generalizes, in which one can explicitly describe the (generalized) Weinstein handles (cf. Theorem \ref{thmMGWeinstein} and Remark \ref{remSpecialZflatHandle}).

\subsection{Proof sketch}\label{subsecPfsketch}
The general idea to prove the equivalence \eqref{mirror intro} is to express both sides as limits of categories featuring in Theorem \ref{th:intro biWh}, with $G$ replaced by its pseudo-Levi subgroups. For simplicity we assume $G$ is almost simple and simply-connected. For each proper subset $J$ of the set $\wt I$ of affine simple roots of the loop group $LG$, we have the pseudo-Levi subgroup $L_J$ of $G$ that is isomorphic to the reductive quotient of the standard parahoric subgroup $\bP_J$ of $LG$. 

For each $L_J$, the general version of Theorem \ref{th:intro biWh} gives an equivalence
\begin{equation}\label{Wh TLJ}
    \cW(\cT_{L_J})\cong \Coh(\cS_{L^\vee_J}).
\end{equation}
Here $\cS_{L^\vee_J}$ is a partially stacky version of the GIT quotient $L_J^\vee\sslash L^\vee_J$, defined for any reductive group $H$ as follows.  Choose any central isogeny $\nu: H_1 \to H$  such that $H_1^{\der}$ is
simply-connected, then $\cS_H$ is defined to be the stack $(H_1\sslash H_1)/\ker(\nu)$. It is independent of the choice of $H_1$ up to canonical isomorphisms. More precisely, to deal with compatibility with sectorial inclusions mentioned below, we establish a canonical mirror equivalence for $\cT_{L_J}$ in \S\ref{secProofMirror}, where $\cW(\cT_L)$ is replaced by $\frP_{L_J}\times^{\pi Z(L_J)}\cW(\cT_{L_J})$ for a canonically defined $\pi Z(L_J)$-torsor $\frP_{L_J}$ (cf. Corollary \ref{corPhiHySy} for the precise statement).

On the A side, the variety $\cM_G$ is covered by open subsets $\cT_{L_J}$, which is a torsor for the regular centralizer group scheme $\cJ_{L_J}$ over $\frc_{L_J}$. In \S\ref{s:Fukaya}, we show that the wrapped Fukaya category $\cW(\cM_G)$ can be canonically written as a colimit
\begin{equation}\label{eqWMGintro}
    \cW(\cM_G)\cong \colim_{J\sne \wt I} \cW(\cT_{L_J})
\end{equation}
of wrapped Fukaya categories of $\cT_{L_J}$ (as open subsectors), with transition functors given by co-restrictions. Passing to ind-completions and right adjoints, we get an equivalence
\begin{equation}\label{lim indW}
    \Ind\cW(\cM_G)\cong \lim_{J\sne \wt I} \Ind\cW(\cT_{L_J})
\end{equation}
where the transition functors are given by restrictions. 

More precisely, to set up the wrapped Fukaya category of $\cM_G$, we construct suitable \emph{real} Liouville 1-forms $\vartheta$ on $\cM_G$ such that $(\cM_G, \vartheta)$ is a Liouville manifold (and Weinstein up to deformations). In fact, all such constructed Liouville manifold structures are canonically equivalent
(up to a contractible space of deformations).  We also construct convenient Weinstein sectorial coverings $\{\cT_{L_J}^{\ovl{\cU}_J}\}_{J\sne\wt{I}}$ of $\cM_G$ so that by sectorial descent \cite{GPS2}, we establish \eqref{eqWMGintro}. 
The lack of a  holomorphic Liouville 1-form (cf. Example \ref{exam: omegaSL2}) makes the construction of real Liouville manifold structures into a technical task. 
In Appendix \ref{AppSecTropical}, we develop several frameworks that are particularly convenient for us to construct Liouville manifold and Weinstein manifold structures on $\cM_G$. They are also substantially used in the construction of Weinstein sectorial coverings of $\cM_G$.   
The canonical property of the colimit \eqref{eqWMGintro} is established based on results on deformations of sector structures on $\cT_{L_J}$
reviewed and developed in Appendix \ref{secAppSectors}.

On the B side, we relate $\QCoh(J_{\dG})$ and the limit of $\QCoh(\cS_{L^\vee_J})$ over $J\sne \wt I$ with transition functors given by pullback along canonical maps $\cS_{L^\vee_J}\to \cS_{L^\vee_{J'}}$ for $J\subset J'$. The strategy can be summarized in the following diagram
\begin{equation}\label{intro qcoh}
    \xymatrix{\QCoh(J_{\dG}) \ar[r]^-{j_*} \ar@{-->}[d] & \QCoh(\left(R_{\dT/\cS_{\dG}}(\dT\times \dT)\right)^W) \ar[r]^-{} & \QCoh(\dT\times\dT)^W\ar[d]^{\cong}\\
    \lim_{J\sne \wt I}\QCoh(\cS_{L^\vee_J})\ar[r]^-{\lim p^*_J} &  \lim_{J\sne \wt I}\QCoh(\dT)^{W_J} & \QCoh(\dT)^{\Wa}\ar[l]_-{\sim}^-{r}}
\end{equation}
Let us explain the meaning of the functors. On the bottom row, each $p_J^*$ is the pullback along the canonical projection $p_J: \dT/W_J\to \cS_{L^\vee_J}$. Each $p_J^*$ is fully faithful and its essential image can be characterized using the behavior on quasi-coherent sheaves along root subtori, which is similar to the Lie algebra analog treated in \cite{Lon} and \cite{Gan}. The equivalence $r$ follows from the fact that $\BB W_\aff$ is the homotopy colimit of $\BB W_J$ for $J\sne \wt I$, proved by Penghui Li \cite{PL}.

The right vertical equivalence is the following simple observation: for any scheme $X$ over $\CC$, direct image along the projection $X\times \dT\to X$ identifies $\QCoh(X\times \dT)$ with the category of $\xch(\dT)=\xcoch(T)$-equivariant quasi-coherent sheaves on $X$.

On the top row, the functor $j_*$ is the direct image functor of a canonical map $j: J_{\dG}\to \left(R_{\dT/\cS_{\dG}}(\dT\times \dT)\right)^W$, which is a group-theoretic analog of the construction of Donagi--Gaitsgory \cite{DG}, and is close to be an isomorphism. The functor $\Phi$ is given by pull-push along a correspondence which makes sense for more general Weil restrictions. In \S\ref{s:Weil res} we prove that the functor $\Phi$ in greater generality is always fully faithful and characterize its essential image (Theorem \ref{th:ff}). In Theorem \ref{th:main QCoh} we show that the composition of the functors in the top row is fully faithful and its essential image is the same as the essential image of the full embedding on the bottom row. Thus the dotted arrow is defined and is an equivalence. 

Combining \eqref{lim indW}, the dotted equivalence in \eqref{intro qcoh} and the equivalences \eqref{Wh TLJ} for various $J$ (whose compatibility is encoded in a choice of compatible system of Kostant sections, cf. Definition \ref{defncompatibleKostantsec}), we conclude with an equivalence
\begin{equation*}
    \Ind\cW(\cM_G)\cong \QCoh(J_{\dG}).
\end{equation*}

%%% thanks %%%

\subsection*{Acknowledgement} 
The authors would like to thank P. Etingof, T. Gannon, A. Goncharov, Penghui Li, G. Lusztig, and J. Pardon for helpful discussions.

%%%%%%%%%%%%%%%%%%%
%%%%%%%%%%%%%%%%%%%
\section{The space $\cM_{G}$ as a Higgs moduli space}\label{secMGHiggs}

In this section we introduce the main geometric player of the paper, the affine Toda system $\cM_G$, and study its algebro-geometric properties. The space $\cM_G$ is a partial compactification of the periodic Toda system. We show that it carries two different completely integrable system structures.

\subsection{Definition of $\cM_G$}\label{ss:def MG}

\sss{The curve} Let $k$ be a field. Consider the curve $X=\PP^1$ over $k$, with affine coordinate $t$. Denote by $0$ and $\infty$ the $k$-points of $\PP^1$ defined by $t=0$ and $t=\infty$ respectively. Let $\t=t^{-1}$ be the affine coordinate at $\infty$. 

\sss{Group-theoretic data} Let $G$ be a split connected reductive group over $k$. We assume that $\ch(k)$ is either zero, or larger than twice the Coxeter number of each simple factor of $G$.

Choose a split maximal torus $T$ of $G$ and a pair of opposite Borel subgroups $B$ and $B^-$ containing $T$. Let $N$ and $N^-$ be the unipotent radicals of $B$ and $B^-$, whose Lie algebras we denote by $\frn$ and $\frn^-$ respectively.

Let $G\lr{t}$ and $G\lr{\t}=G\lr{t^{-1}}$ be the loop groups of $G$ at $0$ and $\infty$, with arc groups denoted by $G\tl{t}$ and $G\tl{\t}=G\tl{t^{-1}}$. 

Let $\bI_0\subset G\tl{t}$ be the Iwahori subgroup obtained as the preimage of $B$ under the reduction mod $t$ map $G\tl{t}\to G$.  Similarly, let $\bI_\infty\subset G\tl{\t}$ be the Iwahori subgroup obtained as the preimage of $B^-$ under the reduction mod $\t$ map $G\tl{\t}\to G$.  Let $\bI^+_0$ and $\bI^+_\infty$ denote the pro-unipotent radicals of $\bI_0$ and $\bI_\infty$ (so they are the preimages of $N$ and $N^-$ under the respective reduction maps).

Denote the Lie algebras of $G\lr{t}, G\tl{t}, \bI_0$ and $\bI^+_0$ by $\frg\lr{t}, \frg\tl{t}, \frI_0$ and $\frI^+_0$ respectively. Similarly define the Lie algebras $\frg\lr{\t}, \frg\tl{\t}, \frI_\infty$ and $\frI^+_\infty$.

\sss{Bundles} 
We will consider the moduli stack $\Bun_G(\bI_0^+, \bI_\infty^+)$ that classifies $G$-bundles on $\PP^1$ with level structures at $0$ and $\infty$ given by the subgroups $\bI_0^+$ and $ \bI_\infty^+$. More precisely, for a test scheme $S$ over $k$, the $S$-points of $\Bun_G(\bI_0^+, \bI_\infty^+)$ is the groupoid of triples $(\cE, \cF_0^+, \cF_\infty^+)$, where
\begin{itemize}
    \item $\cE$ is a $G$-bundles over $\PP^1\times S$;
    \item $\cF^+_0$ is an $N$-reduction of the $G$-bundle $\cE_0:=\cE|_{\{0\}\times S}$  over $\{0\}\times S\cong S$; 
    \item and $\cF^+_\infty$ is an $N^-$-reduction of the $G$-bundle $\cE_\infty=\cE|_{\{\infty\}\times S}$ over $\{\infty\}\times S\cong S$.
\end{itemize}

For a triple $(\cE,\cF^+_0, \cF^+_\infty)$ as above, we denote the induced $B$-bundle from $\cF^+_0$ by $\cF_0$; we denote the induced $B^-$-bundle from $\cF^+_\infty$ by $\cF_\infty$. Then $\cF_0$ and $\cF_\infty$ are Borel reductions of the $G$-bundles $\cE_0$ and $\cE_\infty$ over $S$.

\sss{Generic linear forms}
Let $\wt I$ be the index set of simple affine roots of $G\lr{t}$ with respect to the Iwahori $\bI_0$.  For $i\in \wt I$ we denote the corresponding affine simple root by $\a_i$, viewed as elements of $\xch(T\times \Grot)$, where $\Grot$ is the one-dimensional rotation torus scaling $t$. Let $\frg\lr{t}_{\a_i}$ be the corresponding root space in the loop Lie algebra $\frg\lr{t}$. There is a canonical homomorphism of pro-algebraic groups over $k$
\begin{equation*}
    \e_0: \bI^+_0\to \prod_{i\in \wt I} \frg\lr{t}_{\a_i}=:V_0
\end{equation*}
where the right side is viewed as a vector group over $k$. A linear form $\psi_0\in V_0^*$ is called {\em generic} if it is nonzero on each affine simple root space $\frg\lr{t}_{\a_i}$.

Similarly, we have the loop Lie algebra $\frg\lr{\t}$ of $G\lr{\t}$. The affine simple roots of $G\lr{\t}$ with respect to $\bI_\infty$ are $\{-\a_i\}_{i\in \wt I}$, again as elements in $\xch(T\times \Grot)$. We have a canonical homomorphism of pro-algebraic groups over $k$
\begin{equation*}
    \e_\infty: \bI^+_\infty\to \prod_{i\in \wt I} \frg\lr{\t}_{-\a_i}=:V_\infty.
\end{equation*}
We have an analogous notion of a generic linear form $\psi_\infty\in V_\infty^*$.

\sss{}\label{psi as one form} The space of continuous $k$-linear forms on $\frg\lr{t}$ (with respect to the $t$-adic topology) can be identified with $\frg^*\lr{t}dt$ by the residue pairing. Therefore the dual $V^*_0$ of $V_0$ can be identified with the subspace of $\frg^*\lr{t}dt$ that is the direct sum of weight spaces under $T\times \Grot$ with weights $-\a_i, i\in \wt I$. Similarly, $V_\infty^*$ can be identified with a subspace of $\frg^*\lr{\t}d\t$.

\begin{remark}
    Assume $G$ is almost simple. 
    Let $\{\a_i\}_{i\in I}$ be the simple roots of $G$ with respect to $(B,T)$, then $\wt I=I\cup \{0\}$, with $\a_0=-\th+\d\in \xch(T\times\Grot)$ (where $\th$ is the highest root of $G$, and $\d\in \xch(\Grot)$ is the tautological character). 
    Then 
    \begin{eqnarray*}
        V_0=\frg_{-\th}\cdot t\oplus (\oplus_{i\in I}\frg_{\a_i}),\\
        V_\infty=\frg_{\th}\cdot \t\oplus (\oplus_{i\in I}\frg_{-\a_i}).
    \end{eqnarray*}
    Dually, 
    \begin{eqnarray}
        V_0^*=(\frg_{-\th})^*\cdot \frac{dt}{t^2}\oplus \left(\oplus_{i\in I}(\frg_{\a_i})^*\frac{dt}{t}\right),\\
        V^*_\infty=(\frg_{\th})^*\cdot \frac{d\t}{\t^2}\oplus \left(\oplus_{i\in I}(\frg_{-\a_i})^*\frac{d\t}{\t}\right).
    \end{eqnarray}
    Therefore, the datum of a generic $\psi_{0}$ is the same as a nonzero element in $(\frg_{\a_i})^*$ for each $i\in I$, and a nonzero element in $(\frg_{-\th})^*$. Similarly, a generic $\psi_\infty $ is the same datum as a tuple of nonzero elements in $(\frg_{-\a_i})^*$ for each $i\in I$, and a nonzero element in $(\frg_{\th})^*$.

\end{remark}

\sss{Higgs moduli stack} Fix choices of generic linear forms $\psi_0\in V_0^*$ and $\psi_\infty\in V_\infty^*$.

We define $\cM_{G}=\cM_G(\psi_0,\psi_\infty)$ as the stack over $k$ whose $S$-points is the groupoid of quadruples $(\cE, \cF^+_{0}, \cF^+_{\infty},\ph)$ where
\begin{itemize}
    \item $(\cE, \cF^+_{0}, \cF^+_{\infty})$ is an $S$-point of $\Bun_G(\bI_0^+, \bI_\infty^+)$;
    \item $\ph$ is a Higgs field for $\cE|_{\Gm\times S}$, i.e., it is a section
    \begin{equation*}
        \ph\in \G(\Gm\times S, \Ad^{*}(\cE)\ot\om_{\Gm}).
    \end{equation*}
    Here $\Ad^*(\cE)$ is the vector bundle associated to the coadjoint representation $\frg^*$ of $G$ and the $G$-bundle $\cE$;
\end{itemize}
that satisfy the following conditions:
\begin{enumerate}
    \item Under some trivialization of $\cE$ in the formal neighborhood of $\{0\}\times S$ so that $\cF^+_0$ corresponds to $N$, the Laurent expansion of $\ph$ near $0$ takes the form
    \begin{equation}\label{ph near 0}
        \ph\in \left(\psi_{0}+\frn^\bot\frac{dt}{t}+\frg^*\tl{t}dt\right)\ot \cO_S.
    \end{equation}
    Here we identify $\psi_0\in V_0^*$ with an element in $\frg^*\lr{t}dt$ as discussed in \S\ref{psi as one form}. Note that the right side above is invariant under the coadjoint action of $\bI^+_0$, hence the above condition is independent of the choice of trivializations. 
    \item Under some (equivalently any) trivialization of $\cE$ in the formal neighborhood of $\{\infty\}\times S$ so that $\cF^+_\infty$ corresponds to $N^-$, the Laurent expansion of $\ph$ near $\infty$ takes the form
    \begin{equation}\label{ph near infty}
        \ph\in \left(\psi_{\infty}+(\frn^{-})^{\bot}\frac{d\t}{\t}+\frg^*\tl{\t}d\t\right)\ot \cO_S.
    \end{equation}
\end{enumerate}

\begin{remark}\label{r:MG as Ham red}
    Let $\bI^{++}_0=\ker(\ep_0)\subset \bI_0^+$ and $\bI^{++}_\infty=\ker(\ep_\infty)\subset \bI_\infty^+$.  We may alternatively define $\cM_G$ as the Hamiltonian reduction of $T^*\Bun_G(\bI^{++}_0,\bI^{++}_\infty)$. More precisely, we have the moment map with respect to the $V_0\times V_\infty$-action
    \begin{equation*}
        \mu: T^*\Bun_G(\bI^{++}_0,\bI^{++}_\infty)\to V_0^*\times V_\infty^*.
    \end{equation*}
    Then there is a canonical isomorphism
    \begin{equation*}
        \cM_G\cong \mu^{-1}(\psi_0,\psi_\infty)/(V_0\times V_\infty).
    \end{equation*}
\end{remark}

\begin{exam}\label{ex:torus} When $G=A$ is a torus, we necessarily have $\psi_0=0$ and $\psi_\infty=0$. In this case, $\cM_A$ is isomorphic to the cotangent bundle of $\Bun_A(0,\infty)$, the moduli space of $A$-bundles on $\PP^1$ with trivializations at $0$ and $\infty$. Note that $\Bun_A(0,\infty)\cong \xcoch(A)\times A$, hence $\cM_A\cong \xcoch(A)\times A\times \fra^*$, where $\fra=\Lie A$.  
\end{exam}

\sss{Functoriality}\label{sss:fun MG}
Let $\pi: G\to G'$ be a smooth surjective homomorphism of reductive groups. Let $B'$ and $B'^-$ be the images of $B$ and $B^-$, so they are opposite Borel subgroups od $G'$. The simple affine root spaces of $\frg'\lr{t}$ can be identified with some of the simple affine root spaces of $\frg\lr{t}$, thus the generic linear form $\psi_0$ gives rise to a generic linear form $\psi'_0$ on $\bI'^+_0$ by restriction. Similarly,  $\psi_\infty$ gives rise to $\psi'_\infty$. The Lie algebra map $d\pi: \frg\to \frg'$ admits a unique section $\io: \frg'\to \frg$ whose image is an ideal. Thus we have a map $\io^*: \frg^*\to \frg'^*$ equivariant under the coadjoint actions of $G$ and $G'$. The maps $\pi$ and $\io^*$ induce a map of moduli stacks
\begin{equation}\label{MG func}
     \Pi: \cM_G(\psi_0,\psi_\infty)\to \cM_{G'}(\psi'_0,\psi'_\infty).
\end{equation}

If $\pi$ is an \'etale central isogeny, then $\Pi$ is a $\ker(\pi)$-torsor over its image, which is open and closed in $\cM_{G'}(\psi'_0,\psi'_\infty)$. In particular, $\Pi$ is finite \'etale in this case.

When $G$ is a product group $G=\prod_j G_j$, then
\begin{equation*}
    \cM_G(\psi_0,\psi_\infty)=\prod_j \cM_{G_j}(\psi^j_0,\psi^j_\infty)
\end{equation*}
where $\psi^j_0,\psi^j_\infty$ are the restrictions of $\psi_0, \psi_\infty$ to the $j$th factor.

\subsection{The Hitchin base and Hitchin map}\label{ss:AG}
Let $\frc^*=\frg^*\sslash G\cong \frt^*\sslash W$. 

Let $\cA_G=\cA_G(\psi_0,\psi_\infty)=\frc^*$. We shall define a Hitchin map 
\begin{equation}\label{Hitchin map}
    f_G: \cM_G=\cM_G(\psi_0,\psi_\infty)\to \cA_G(\psi_0,\psi_\infty)=\cA_G.
\end{equation}

\sss{Almost simple case}\label{sss:AG almost simple} First consider the case where $G$ is almost simple. 

Let $f_{1},\cdots, f_{r}$ be homogeneous generators of $\Sym(\frg)^{G}$, of degrees $2=d_{1}\le \cdots\le d_{r}=h$, the Coxeter number of $G$. Using $f_1,\cdots, f_r$ as coordinates we have an isomorphism $\frc^*\cong \AA^r$.

Each $f_i$ gives a polynomial function on $\frg^*$. Extend $k\lr{t}$-linearly, it gives a polynomial function on $\frg^*\lr{t}dt$ valued in $k\lr{t}(dt)^{d_i}$. In particular we have $f_i(\psi_0)\in k\lr{t}(dt)^{d_i}$ and $f_i(\psi_\infty)\in k\lr{\t}(d\t)^{d_i}$.

\begin{lemma}
    We have $f_i(\psi_0)=0$ for $1\le i<r$, and $f_r(\psi_0)=a_0t^{-1}\left(\frac{dt}{t}\right)^{h}$ for some $a_0\in k^\times$. 
    
    Similarly, $f_i(\psi_\infty)=0$ for $1\le i<r$, and $f_r(\psi_\infty)=a_\infty t\left(\frac{dt}{t}\right)^{h}$ for some $a_\infty\in k^\times$.
\end{lemma}
\begin{proof}
    Let $\r^\vee:\Gm\to G^{\ad}$ be half the sum of positive roots with respect to $(B,T)$. Let $h=\j{\r^\vee,\th}+1$ be the Coxeter number of $G$. Consider the element $t^{\r^\vee/h}\in G\lr{t^{1/h}}$ that acts on each simple affine root space $\frg\lr{t}_{\a_i}$ by the scalar $t^{1/h}$. View $\psi_0$ as an element of $\frg^*\lr{t^{1/h}}$, we have
    \begin{equation*}
        \Ad^*(t^{\rho^\vee/h})\psi_0=t^{-1/h}\psi_0,
    \end{equation*}
    which implies that $f_i(\psi_0)\in k\cdot t^{-d_i/h}(dt/t)^{d_i}$. Since $f_i(\psi_0)\in k\lr{t}(dt/t)^{d_i}$, we see that $f_i(\psi_0)=0$ unless $h|d_i$, which only happens for $f_r$, and $f_r\in k\cdot t^{-1}(dt/t)^{d_i}$. Since $\psi_0$ is not nilpotent, $f_r(\psi_0)\ne 0$, hence $f_r(\psi_0)$ takes the desired form. The statement about $\psi_\infty$ is proved similarly, where $t^{-1}$ now plays the role of the uniformizer at $\infty$. 
\end{proof}

For $(\cE,\cF_0^+, \cF_\infty^+, \ph)\in \cM_{G}(S)$, $f_{i}(\ph)$ is a rational section of $\om(0+\infty)^{\ot d_{i}}\cong \cO_{\PP^1\times S}$ (via $\frac{dt}{t}$).

\begin{lemma}\label{l:fi phi}
    For $1\le i<r$, we have
    $$f_i(\ph)\in \cO_S \left(\frac{dt}{t}\right)^{d_i}.$$
    On the other hand, for $i=r$, we have
    $$ f_r(\ph)\in \left(a_0t^{-1}+\cO_S+a_\infty t\right) \left(\frac{dt}{t}\right)^{h}.$$
\end{lemma}
\begin{proof}
    Write $f_i(\ph)=b_i\left(\frac{dt}{t}\right)^{d_i}$ for some $b_i\in \G(\Gm\times S, \cO)$. Then $b_i$ has a pole of order $\le d_{i}/h$ at $0$ and $\infty$. Indeed, at $0$, we may assume \eqref{ph near 0} holds. Now embed $\frg^*\lr{t}$ into $\frg^*\lr{t^{1/h}}$, we have $\psi_0$ as an element in $\frg^*\lr{t^{1/h}}$, we have
    \begin{equation*}
        \Ad^*(t^{\rho^\vee/h})\ph\in t^{-1/h}\frg^*\tl{t^{1/h}}\frac{dt}{t}
    \end{equation*}
    which implies that $b_i=f_i(\ph)\left(\frac{dt}{t}\right)^{-d_i}$ has pole order $\le d_i/h$. Therefore $b_i\in \cO_S$ for $1\le i\le r-1$ (since $d_i<h$), and $b_r$ has a simple pole at $0$ and $\infty$. 
\end{proof}

We then define the map \eqref{Hitchin map} on the level of $S$-points by
\begin{equation*}
    f_G(\cE,\cF_0^+, \cF_\infty^+, \ph)=(b_1,b_2,\cdots, b_r)\in (\cO_S)^r=\AA^r(S)\cong \frc^*(S)
\end{equation*}
where $b_i$ are determined by the equations (after Lemma \ref{l:fi phi})
\begin{equation}\label{fi ph}
    f_i(\ph)=\begin{cases} b_i\left(\frac{dt}{t}\right)^{d_i}, & 1\le i<r;\\
    \left(a_0t^{-1}+b_r+a_\infty t\right) \left(\frac{dt}{t}\right)^{h}, & i=r.
    \end{cases}
\end{equation}

\sss{General case}\label{sss:fG gen}
When $G=A$ is a torus,  $\cA_A=\fra^*$ (where $\fra=\Lie A$). Using the description of $\cM_A$ in Example \ref{ex:torus},  the map $f_A: \cM_A=\xcoch(A)\times A\times \fra^*\to \fra^*$ is the projection.  

For general reductive $G$, we choose a central isogeny $\io:G\to  \prod_{j}G_j$ where $G_j$ is either an almost simple group or a torus. By the functoriality \eqref{MG func}, $\io$ induces a canonical map $\Pi_j: \cM_G\to \cM_{G_j}$. Also $\io$ induces a canonical isomorphism $\cA_G\cong \prod_j \cA_{G_j}$. Define $f_G: \cM_G\to \cA_G=\prod_j \cA_{G_j}$ whose $j$-component is $f_{G_j}\c\Pi_j$. 
It is easy to check that the map $f_G$ is independent of the choice of the isogeny $\io$.

\subsection{A variant of $\cM_G$}\label{ss:M flat} 

\sss{The stack $\cM^\flat_G$} Consider the moduli stack $\cM^{\flat}_G$ whose $S$-points classify $(\cE, \cF_0, \cF_\infty, \ph)$ where $\cE$ is a $G$-torsor over $\PP^1\times S$ with $B$-reduction $\cF_0$ along $\{0\}\times S$ and $B^-$-reduction $\cF_\infty$ along $\{\infty\}\times S$, and $\ph\in \cohog{0}{\Gm, \Ad^{*}(\cE)\ot \om}$ satisfying the following conditions near $0$ and $\infty$: 
\begin{enumerate}
    \item Under some (equivalently, any) trivialization of $\cE$ in the formal neighborhood of $\{0\}\times S$ so that $\cF_0$ corresponds to $B$, the Laurent expansion of $\ph$ near $0$ takes the form
    \begin{equation*}
        \ph\in \left(V^*_{0}+\frn^\bot\frac{dt}{t}+\frg^*\tl{t}dt\right)\ot \cO_S.
    \end{equation*}
    \item Under some (equivalently, any) trivialization of $\cE$ in the formal neighborhood of $\{\infty\}\times S$ so that $\cF_\infty$ corresponds to $B^-$, the Laurent expansion of $\ph$ near $\infty$ takes the form
    \begin{equation*}
        \ph\in \left(V^*_{\infty}+(\frn^{-})^{\bot}\frac{d\t}{\t}+\frg^*\tl{\t}d\t\right)\ot \cO_S.
    \end{equation*}
\end{enumerate}
The leading term of $\ph$ at $0$ in $V_0^*$ is well-defined up to the action of $T$. Therefore we have a map
\begin{equation*}
    e_0: \cM^\flat_G\to V_0^*\sslash T.
\end{equation*}
Similarly we have a map
\begin{equation*}
    e_\infty: \cM^\flat_G\to V_\infty^*\sslash T.
\end{equation*}
When $G$ is almost simple, both $V_0^*\sslash T$ and $V_\infty^*\sslash T$ are isomorphic to $\AA^1$.

\sss{The base $\cA^\flat_G$}\label{sss:A flat} The Hitchin base $\cA^\flat_G$ for $\cM^\flat_G$ is defined to be the affine space
\begin{equation*}
    \cA^\flat_G=V_0^*\sslash T\times V_\infty^*\sslash T \times \frc^*.
\end{equation*}
To define the Hitchin map 
\begin{equation*}
    f^\flat_G: \cM^\flat_G\to \cA^\flat_G,
\end{equation*}
following the same route of \S\ref{ss:AG} by considering the almost simple case first. 

When $G$ is almost simple, choose a set of homogeneous generators $f_1,\cdots, f_r$ of $\Sym(\frg)^G$ as in \S\ref{sss:AG almost simple}. Restricting $f_r$ to $V_0^*$ gives a map $V_0^*\to \AA^1 t^{-1}(\frac{dt}{t})^h$ that is $T$-invariant. Therefore it factors through an isomorphism
\begin{equation*}
    a_0: V_0^*\sslash T\isom \AA^1.
\end{equation*}
such that $f_r(\psi_0)=a_0(\psi_0)t^{-1}(\frac{dt}{t})^h$ for any $\psi_0\in V_0^*$. Similarly, there is a unique isomorphism
\begin{equation*}
    a_\infty: V_\infty^*\sslash T\isom \AA^1.
\end{equation*}
such that $f_r(\psi_\infty)=a_\infty(\psi_\infty)t(\frac{dt}{t})^h$ for any $\psi_\infty\in V_\infty^*$.  The Hitchin map 
maps $m=(\cE, \cF_0, \cF_\infty, \ph)$ to 
\begin{equation*}
    (e_0(m), e_\infty(m), f_1(\ph),\cdots, f_{r-1}(\ph), b_r(m)).
\end{equation*}
Here $b_r(m)$ is characterized by
\begin{equation*}    f_r(\ph)=\left(a_0(e_0(m))t^{-1}+b_r(m)+a_\infty(e_\infty(m))t\right)\left(\frac{dt}{t}\right)^h.
\end{equation*}
For reductive $G$, $f_G^\flat$ is defined in the same way as $f_G$ in \S\ref{sss:fG gen}.

From the definitions, there is a natural map
\begin{equation*}
    \io^\flat(\psi_0,\psi_\infty): \cM_G(\psi_0,\psi_\infty)\to \cM_G^\flat
\end{equation*}
that covers the embedding of Hitchin bases
\begin{equation*}
    \a: \cA_G(\psi_0,\psi_\infty)=\{(a_0(\psi_0),a_\infty(\psi_\infty))\}\times \frc^*\incl \cA^\flat_G=V_0^*\sslash T\times V_\infty^*\sslash T \times \frc^*.
\end{equation*}

\begin{lemma}\label{l:ZG2 torsor}  For generic linear forms $\psi_0\in V_0^*$ and $\psi_\infty\in V_\infty^*$, the map
\begin{equation*}
    \io(\psi_0,\psi_\infty): \cM_G(\psi_0,\psi_\infty)\to \cM_G^\flat(\psi_0,\psi_\infty):=\cM^\flat_G\times_{\cA_G^\flat}\cA_G(\psi_0,\psi_\infty)
\end{equation*}
induced by $\io^\flat(\psi_0,\psi_\infty)$ is a $(ZG)^2$-torsor.
\end{lemma}
\begin{proof}
    First check $ \io(\psi_0,\psi_\infty)$ is surjective. Given a geometric point $m^\flat=(\cE, \cF_0,\cF_\infty, \ph)\in \cM^\flat_G(\psi_0,\psi_\infty)$, choose any $N$-reduction $\cF^+_0$ of $\cF_0$ and $N^-$-reduction $\cF^+_\infty$  of $\cF_\infty$. Under any trivialization of $\cE$ near $0$ such that $\cF^+_0$ corresponds to $N$, $\ph$ has the form $\psi'_0+\frn^\bot \frac{dt}{t}+\frg^*\tl{t}dt$ for some $\psi'_0\in V_0^*$ whose image in $V_0^*\sslash T$ is the same as that of $\psi_0$. The fiber over a nonzero point of $V_0^*\to V_0^*\sslash T$ forms a single $T$-orbit, hence there exists $s\in T$ such that $\Ad^*(s)\psi'_0=\psi_0$. Changing $\cF^+_0$ by $s$ we get another $N$-reduction of $\cF_0$ under which $\ph$ has leading term $\psi_0$ at $0$. Similarly, one can adjust $\cF^+_\infty$ by an element in $T$ so that the leading term of $\ph$ at $\infty$ is equal to $\psi_\infty$. This proves the surjectivity of $\io(\psi_0,\psi_\infty)$.

    The group $(ZG)^2$ acts on $\cM_G(\psi_0,\psi_\infty)$ by changing the $N$ and $N^-$-reductions of $\cE_0$ and $\cE_\infty$ while keeping the $B$ and $B^-$-reductions $\cF_0$ and $\cF_\infty$. The discussion in the previous paragraph shows that the fiber of $\io(\psi_0,\psi_\infty)$ over a geometric point $m^\flat$ of $\cM^\flat_G(\psi_0,\psi_\infty)$ is a torsor under the following subgroup of $T$ 
    \begin{equation*}
        \Stab_T(\psi_0)\times \Stab_T(\psi_\infty).
    \end{equation*}
    Since $\psi_0$ and $\psi_\infty$ are generic, both stabilizers are equal to $ZG$. This shows that $\io(\psi_0,\psi_\infty)$ is a $(ZG)^2$-torsor.
\end{proof}

\sss{Neutral component}
Let $\Bun^\c_G(\bI^+_0,\bI^{+}_{\infty})\subset \Bun_G(\bI^+_0,\bI^{+}_{\infty})$ be the connected component of the unit double coset under the decomposition \eqref{Bun I unif}. Equivalently, $\Bun^\c_G(\bI^+_0,\bI^{+}_{\infty})$ is the union of $\Bun^w_G(\bI^+_0,\bI^{+}_{\infty})$ for all $w\in \Wa$. Let $\cM^\c_G\subset \cM_G$ be the open and closed substack obtained as the preimage of $\Bun^\c_G(\bI^+_0,\bI^{+}_{\infty})$ under $b$. 

It will be justified in Corollary \ref{c:MGc conn} that $\cM^\c_G$ is indeed connected, therefore deserves to be called the neutral connected component of $\cM_G$.

%%%%%%%%%%%%%%%%%%%%%%%%%
\subsection{First properties}

\begin{theorem}\label{th:smooth} The stack $\cM_{G}$ is a smooth algebraic space \footnote{We will show in Corollary \ref{c:MG scheme} that $\cM_G$ is actually a scheme} over $k$ of dimension $2r$ (where $r$ is the rank of $G$) with a canonical symplectic structure.
\end{theorem}
\begin{proof}
The tangent complex of $\cM_{G}$ at a geometric point $(\cE, \cF^{+}_{0},\cF^{+}_{\infty}, \ph)$ is $R\Gamma(\PP^{1}, \cK)$ where $\cK=[\cK^{-1}\to \cK^{0}]$ is the two-step complex (in degrees $-1$ and $0$) 
\begin{equation*}
\cK^{-1}=\Ad(\cE, \cF^{+}_{0}, \cF^{+}_{\infty})\xr{[-,\ph]} \Ad^{*}(\cE, \cF_{0}, \cF_{\infty})\ot \om_{\PP^1}(0+\infty)=\cK^{0}.
\end{equation*}
Here 
\begin{itemize}
    \item The vector bundle $\Ad(\cE, \cF^{+}_{0}, \cF^{+}_{\infty})$ is the preimage of $i_{0*}\Ad(\cF^+_0)\op i_{\infty*}\Ad(\cF^+_\infty)$ under the evaluation map
\begin{equation*}
    \Ad(\cE)\to i_{0*}\Ad(\cE_0)\op i_{\infty*}\Ad(\cE_\infty).
\end{equation*}
Here $i_0:\{0\}\incl \PP^1$ and $i_\infty: \{\infty\}\incl \PP^1$ are the inclusions. 

    \item The vector bundle $\Ad^*(\cE, \cF_{0}, \cF_{\infty})$ is the preimage of $i_{0*}\Ad(\cF^+_0)^\bot\op i_{\infty*}\Ad(\cF^+_\infty)^\bot$ under the evaluation map
\begin{equation*}
    \Ad^*(\cE)\to i_{0*}\Ad^*(\cE_0)\op i_{\infty*}\Ad^*(\cE_\infty).
\end{equation*}
    \item Taking Lie bracket with $\ph$ gives an $\cO_{\PP^1}$-linear map $\cK^0\to \cK^1$.
\end{itemize}
From the description above, we see that $\cK^{-1}$ is canonically isomorphic to the Serre dual of $\cK^{0}[1]$, i.e., $\cK^{0}\cong \un\Hom(\cK^{-1}, \om_{\PP^1})$. The map $[-,\ph]: \cK^{-1}\to \un\Hom(\cK^{-1},\om_{\PP^1})$, viewed as a pairing $\cK^{-1}\ot \cK^{-1}\to \om_{\PP^1}$, is alternating. This implies $\bR\Gamma(\PP^{1}, \cK)$ carries a perfect symplectic pairing. In particular, the classical tangent space $\cohog{0}{\PP^{1}, \cK}$ of $\cM_G$ at $m=(\cE, \cF^{+}_{0},\cF^{+}_{\infty}, \ph)$ carries a symplectic form, and $\cohog{-1}{\PP^{1}, \cK}$ (infinitesimal automorphisms of $m$)is dual to the obstruction group $\cohog{1}{\PP^{1}, \cK}$. 

If we show that $\cohog{-1}{\PP^{1}, \cK}=0$, i.e., points in $\cM_G$ do not have nonzero infinitesimal automorphisms, then the obstruction group will also vanish, hence $\cM_{G}$ is a smooth Deligne-Mumford stack  with a canonical symplectic structure. Since $\deg\cK^{0}-\deg\cK^{-1}=2r$, we see that $\chi(\PP^{1}, \cK[1])=2r$. Thus $\dim \cM_G=2r$.

It remains to show that for any geometric point $m\in \cM_{G}$, $\Aut(m)$ is trivial as an algebraic group (which implies in particular the infinitesimal automorphisms of $m$ are trivial, and that $\cM_G$ is an algebraic space).

Consider the following Hitchin moduli stack $\cH_{G}$ without level structure: it classifies pairs $(\cE,\ph)$ where $\cE$ is a $G$-torsor over $\PP^{1}$ and $\ph$ is a section of $\Ad^{*}(\cE)\ot \om\ot \cO(0+\infty)^{\ot 2}\cong\Ad^{*}(\cE)\ot\cO(2)$. We have a natural map
\begin{equation*}
    \nu: \cM_G\to \cH_G
\end{equation*}
that forgets the $N$ and $N^-$-reductions $\cF^+_0$ and $\cF^+_\infty$. The map $\nu$ is representable. For a geometric point $m=(\cE,\cF^+_0,\cF^+_\infty,\ph)\in \cM_G$, $\nu(m)$ lies in the open substack $\cH^\hs_G\subset \cH_G$ where the Higgs field is generically regular semisimple: this is because any element on the right side of \eqref{ph near 0} is regular semisimple. By \cite[Cor.4.11.3]{NgoFL}, $\Aut(\nu(m))$ is isomorphic to a subgroup scheme of $T$, hence diagonalizable. Since $\nu$ is representable, $\Aut(m)\incl \Aut(\nu(m))$, hence $\Aut(m)$ is also diagonalizable. On the other hand, $\Aut(m)$ is isomorphic to a subgroup scheme of $\bI^+_0$ by evaluating at $0$, which implies that $\Aut(m)$ is unipotent. Combining these facts we conclude that $\Aut(m)$ has to be the trivial group scheme. This completes the proof.
\end{proof}

\begin{theorem} The map $f_G: \cM_{G}\to \cA_{G}$ is proper, and is a completely integrable system, i.e., general fibers of $f_G$ are Lagrangians in $\cM_G$.
\end{theorem} 
\begin{proof}
(1) We show $f_G$ is proper. 
Recall the moduli stack $\cH_G$ introduced at the end of proof of (1). Let $h:\cH_{G}\to \cB_{G}$ be the Hitchin map for $\cH_G$ to its Hitchin base $\cB_{G}=\prod_{i=1}^{r}\cohog{0}{\PP^{1}, \cO(2d_{i})}$. Let $\wt\cH_G=\cH_{G}\times_{\Bun_G}\Bun_G(\bI_0, \bI_\infty)$. 

We have a natural embedding $\cA_{G}\incl \cB_{G}$ and a forgetful map $\nu: \cM_{G}\to \cH_{G}|_{\cA_{G}}$ compatible with the Hitchin maps. We decompose $f$ as
\begin{equation*}
f: \cM_{G}\xr{\nu} \cH_{G}|_{\cA_{G}}\xr{h|_{\cA_{G}}} \cA_{G}.
\end{equation*}

To show $f$ is proper, we show separately that $\nu$ and $h|_{\cA_{G}}$ are proper. Using the moduli stack $\cM^{\flat}_{G}$ introduced in \S\ref{ss:M flat}, the map $\nu$ further factors as
\begin{equation*}
\nu: \cM_{G}\xr{\nu_{1}} \cM^{\flat}_{G}|_{\cA_{G}}\xr{\nu_{2}|_{\cA_{G}}}\cH_{G}|_{\cA_{G}}.
\end{equation*}
Now $\nu_{1}$ is a $(ZG)^2$-torsor by Lemma \ref{l:ZG2 torsor}, hence proper. The map $\nu_{2}:  \cM^{\flat}_{G}\to \cH_{G}|_{\cA^{\flat}_{G}}$ realizes $\cM^{\flat}_{G}$ as a closed substack of $(\cH_{G}|_{\cA^{\flat}_{G}})\times_{\Bun_{G}}\Bun_{G}(\bI_{0}, \bI_{\infty})$. Since $\Bun_{G}(\bI_{0}, \bI_{\infty})\to \Bun_{G}$ is proper, so is $\nu_{2}$ and its restriction over $\cA_{G}$.

By \cite[\S4.10.5, \S5.4.7]{NgoFL}, there is an open subset $\cB^{\ani}_{G}\subset \cB_{G}$ characterized by the following property. A point $a\in \cB^{\ani}_{G}$ defines a map $\wt a: \PP^{1}\bs\{0,\infty\}\to \frc^*$, first one requires that $\ov a$ generically lands in the regular locus $\frc^{*,\rs}=\frt^{*,\rs}\sslash W$, which then defines a $W$-torsor $\y_{a}\to \y$ over the generic point $\y\in \PP^{1}$. This yields a continuous homomorphism $\pi_{a}: \Gal(\ov\y/\y)\to W$ well-defined up to conjugacy. Let $W(a)\subset W$ be the image of $\pi_{a}$,  a subgroup well-defined up to conjugacy. Then $a\in \cB^{\ani}_{G}$ if and only if $\frt^{W(a)}=0$. It is shown in \cite[Prop. 6.1.3]{NgoFL} (which essentially relies on \cite[II.4]{Faltings}) that $h|_{\cB^{\ani}_{G}}$ is proper. Therefore it remains to check that $\cA_{G}\subset \cB_{G}^{\ani}$.

For this we can reduce to the case $G$ almost simple. Consider the following $\Gm$-action on $\cM^\flat_G$: $\l\in \Gm$ scales $\ph$ by $\l$ and acts on $\PP^1$ by $t\mt \l^{-h}t$. This induces a $\Gm$-action on $\cA^\flat_G$ so that $f^\flat_G$ is $\Gm$-equivariant:
\begin{equation*}
    \l\cdot (a_0, a_\infty, b_1,\cdots, b_r)=(\l^{2h}a_0, a_\infty, \l^{d_1}b_1,\cdots, \l^{d_r}b_r).
\end{equation*}
In particular, the $\Gm$-action preserves the open subset $\cA^\flat_G\cap \cB^{\ani}_{G}$. All points in $\cA_G=\cA_G(\psi_0,\psi_\infty)$ get contracted to $(0,a_\infty(\psi_\infty), 0,\cdots, 0)$ by the $\Gm$-action. Therefore to show that $\cA_G$ is contained in $\cB^{\ani}_G$, it suffices to show that $(0,a_\infty(\psi_\infty), 0,\cdots, 0)$ is in  $\cB^{\ani}_G$, i.e., the point $b=(0,\cdots, 0, a_\infty(\psi_\infty)t)\in \cB^{\ani}_G$. Localizing at $\infty$, the point $b$ is the image of $\psi_\infty$ in $\frc^*(k\lr{\t})$. It is well-known that $\psi_\infty$ is an elliptic regular semisimple element in $\frg^*\lr{\t}$ whose centralizer is a maximal torus in $G\lr{\t}$ of type corresponding to the Coxeter conjugacy class in $W$, we conclude that $b\in \cB^{\ani}_G$, hence $\cA_G\subset \cB^{\ani}_G$.

(2) Since $\dim \cM_G=2r$ and $\dim \cA_G=r$, it suffices to show that for any geometric point $m=(\cE,\cF^{+}_{0}, \cF^{+}_{\infty}, \ph)\in \cM_{G}$ with image $a\in \cA_{G}$, the image of the differential
\begin{equation*}
(df)_{m}: T^{*}_{a}\cA_{G}\to T^{*}_{m}\cM_{G}
\end{equation*}
is isotropic. The argument is similar to that of \cite[\S2.12]{BBAMY} so we will be brief here. We use notations from the deformation calculation in the proof of Theorem \ref{th:smooth}. The tangent map of $f_G$ at $m$ is obtained from the following map of complexes by taking global sections
\begin{equation*}
    \xymatrix{\cK^{-1}\ar[r]^-{[-,\ph]} & \cK^0\ar[d]^{(df_i \textup{ at } \ph)} \\
    & \frc^*\ot \cO_{\PP^1}}
\end{equation*}
The composition of the two arrows is zero, thus we have an induced map
\begin{equation*}
    \Phi: \cH^0\cK=\coker([-,\ph])\to \frc^*\ot\cO.
\end{equation*}
Thus the tangent map of $f_G$ at $m$ is the composition
\begin{equation}\label{Tfm}
    (Tf)_m: T_m\cM_G\cong \upH^0(\PP^1,\cK)\to \upH^0(\PP^1,\cH^0\cK)\to \upH^0(\PP^1,\frc^*\ot\cO)\cong \cA_G.
\end{equation}
Note however we have a short exact sequence
\begin{equation}\label{Lag in tangent}
    0\to \upH^1(\PP^1, \cH^{-1}\cK)\to T_m\cM_G\to \upH^0(\PP^1, \cH^0\cK)\to 0.
\end{equation}
This follows from the spectral sequence calculating $\upH^*(\PP^1,\cK)$ and the vanishing of $\upH^{-1}(\PP^1,\cK)$ and $\upH^{1}(\PP^1,\cK)$ by Theorem \ref{th:smooth}. Moreover, the self-dual nature of $\cK$ shows that in \eqref{Lag in tangent}, $\upH^1(\PP^1, \cH^{-1}\cK)$ is a Lagrangian in $T_m\cM_G$. By \eqref{Tfm}, $(Tf)_m$ factors through the quotient by a Lagrangian, its adjoint map $(df)_m$ has isotropic image. This finishes the proof.

\end{proof}

\subsection{Hitchin fibers} Classically, periodic Toda systems are solved using theta functions on the Jacobians of spectral curves, which are  fibers of the Hitchin map $f_G$. We give a modern account of the spectral curves associated with $\cM_G$ in type $A$ and Hitchin fibers in general, following the paradigm of Ng\^o \cite{NgoFL}.

For $a\in \cA_G(k)$, let $\cM_a=f_G^{-1}(a)$ be the corresponding Hitchin fiber.

\begin{exam} For $G=\SL_{2}$, take $\psi_{0}=\mat{}{t^{-1}}{1}{}$, $\psi_{\infty}=\mat{}{-1}{-t}{}$.  Then $\cA_{G}\subset \G(\PP^{1}, \cO(0+\infty))$ is a one-dimensional affine space consisting of $-(t+a+t^{-1})(dt/t)^{2}$ for varying $a\in \AA^1$. The spectral curve $Y_{a}$ is a curve in the total space of $\cO(2)$ over $\PP^{1}$, with equation $y^{2}=t(t^{2}+at+1)$. We see that $p_{t}: Y_{a}\to \PP^{1}$ (sending $(y,t)\mt t$) is a double cover ramified over four points: $0,\infty$ and the two roots of $t+a+t^{-1}=0$. 

Suppose $a\ne2$, then $p_t: Y_{a}\to X=\PP^{1}$ is ramified at 4 distinct points $R\subset Y_{a}$. Hence $Y_{a}$ has genus $1$. Let $\wt 0,\wt \infty\in R$ be the unique points over $0,\infty\in \PP^1$. The pair  $(Y_a,\wt\infty)$ is an elliptic curve. Let $\s\in \Gal(Y_{a}/\PP^1)$ be the nontrivial involution attached to the double cover $p_t$. The Hitchin fiber $\cM_{a}$ classifies $(\cL, \io, v_{0}, v_{\infty})$ where 
\begin{itemize}
\item $\cL$ is a line bundle over $Y_{a}$ of degree $2$.
\item $\io: \cL\ot\s^{*}\cL\isom \cO_{Y_{a}}(R)$ is an isomorphism satisfying $\s^{*}\io=-\io$.
\item $v_{0}\in \cL_{\wt 0}$ and $v_{\infty}\in \cL_{\wt\infty}$ are nonzero vectors, such that $\io(v_{0}, yv_{0})=-1$, $\io(v_{\infty}, yv_{\infty})=1$ as sections of $\cO(R)$ at $0$ and $\infty$. Here $y$ is the canonical vertical coordinate on $\Tot(\cO(2))$. 
\end{itemize}
The pairs $(\cL,v_{\infty})$ as above are classified by the Picard variety $\Jac^{2}(Y_{a})$, which can be identified with $Y_{a}$ itself using the point $\wt\infty$. The condition $\io(v_{\infty}, yv_{\infty})=1$ then determines $\io$, and there are two choices of $v_{0}$ satisfying   $\io(v_{0}, yv_{0})=-1$. Hence $\cM_{a}\to Y_{a}$ is a $\mu_2=ZG$-torsor. Indeed, $\cM_{a}$ can be identified with the square roots of a canonical everywhere nonzero section of the line bundle $\cL\mapsto \cL_{\wt 0}^{\ot 2}$ on $\Jac^{2}(Y_{a})\cong Y_{a}$ (using that $\wt 0$ is a $2$-torsion point of $Y_{a}$).  This shows that $\cM_{a}$ itself is an elliptic curve. 

When $a=\pm2$, $Y_{a}$ is a nodal rational curve with a node over $\mp1\in \PP^{1}$. Let $R\subset Y_a$ be the divisor $\wt 0+\wt \infty+p_t^{-1}(\mp1)$. Then $\cM_{a}$ classifies  $(\cL, \io, v_{0}, v_{\infty})$ with $\cL$ allowed to be torsion-free rank 1, $\io: \cL\to \s^{*}\cL^{\vee}\ot\cO(R)$ satisfying $\s^{*}\io^{\vee}=-\io$ and $v_{0},v_{\infty}$ as above. Same argument shows that $\cM_{a}$ is a $\mu_2$-torsor over the compactified Jacobian $\ov\Jac^{2}(Y_{a})$, which is isomorphic to $Y_{a}$ itself. Since $f_G$ is proper with connected fibers for $a\ne\pm2$, and $\cA_{G}$ is normal, Zariski main theorem implies that all geometric fibers are connected.  Therefore $\cM_{a}$ for $a=\pm2$ has to be a nontrivial double cover of the nodal rational curve, i.e. a union of two $\PP^{1}$'s glued at two points.
\end{exam}

\begin{exam}\label{ex:SLn fiber} Consider the case $G=\SL_{n}$. Without loss of generality we may assume that $a_0(\psi_0)=a_\infty(\psi_\infty)=1$. Therefore $\cA_G\cong \AA^{n-1}$ so that $a=(a_2,\cdots, a_n)\in \cA_G(k)$ corresponds to the characteristic polynomial 
\begin{equation*}
    P_a(t,y)=y^{n}+a_{2}y^{n-2}-\cdots+(-1)^{n}(t+a_{n}+t^{-1}).
\end{equation*}
Consider the spectral curve $Y^\times_a\subset \Gm\times \AA^1$ defined by $P_a(t,y)=0$, where $t$ is the coordinate of the $\Gm$-factor and $y$ is the coordinate in $\AA^1$-factor. The spectral curve $Y_a$ is defined to be the unique compactification of $Y^\times_a$ that is smooth over $t=0$ and $t=\infty$.

The spectral curve $Y_a$ can be understood via the two projections $(y,t)\mt t$ and $(y,t)\mt y$.

Taking the $y$ coordinate gives a degree $2$ map $p_y: Y_a\to \PP^1$ ramified over the roots of $y^n+a_2y^{n-1}+\cdots+(-1)^n(a_n\pm 2)=0$, i.e., $2n$ points counted with multiplicities. Riemann-Hurwitz formula implies that $Y_a$ has arithmetic genus $n-1$.

On the other hand, the projection $p_t: Y_a\to X$ (taking the $t$-coordinate) has degree $n$. It is totally ramified over $0$ and $\infty$. We denote the unique points of $Y_a$ over $0$ and $\infty$ by $\wt 0$ and $\wt\infty$. The discriminant $\D(t)$ of $P_a(t,y)$ as a polynomial in $y$ is itself a polynomial of degree $n-1$ in $t+a_{n}+t^{-1}$. Therefore as a rational function in $t$, $\D(t)$ has $2(n-1)$ zeros in $\Gm$, counted with multiplicities. This means the ramification index of $Y^\times_{a}\to \Gm$ is $2(n-1)$. Riemann-Hurwitz formula again confirms that $Y_a$ has arithmetic genus $n-1$.

The Hitchin fiber $\cM_a$ can be described as the moduli stack of $(\cL, v_0,v_\infty, \d)$ where $\cL$ is a torsion-free coherent sheaf of generic rank one on $Y_a$, $\d: \det(p_{t*}\cL)\cong \cO_{\PP^1}$ (which forces $\cL$ to have degree $2n-2$), $v_0$ and $v_\infty$ are nonzero vectors of the lines $\cL_{\wt 0}$ and $\cL_{\wt \infty}$ such that $\d(v_0\wedge yv_0\wedge\cdots\wedge y^{n-1}v_0)=(-1)^{n-1}$ and $\d(v_\infty\wedge yv_\infty\wedge\cdots\wedge y^{n-1}v_\infty)=1$. From this we see that $\cM_a$ is a $\mu_n=ZG$-torsor over one component $\ov\Jac^{2n-2}(Y_a)$ of the compactified Jacobian of $Y_a$.

\end{exam}

\sss{The Picard variety}

Let $a\in \cA_{G}(k)$. Following \cite[\S4.3]{NgoFL} we have the group scheme $J_{a}^{0,\infty}$ over $\PP^{1}\bs\{0,\infty\}$ obtained by pulling back the regular centralizer group scheme $\cJ\to \frc$ via the map $\wt a: \PP^{1}\bs\{0,\infty\}\to \frc$. Let $J_{a}$ be the group scheme over $\PP^{1}$ extending $J^{0,\infty}_{a}$, whose restriction to $D_{0}$ and $D_{\infty}$ are the parahoric group schemes of the generic fibers $J_{a}|_{D^{*}_{0}}$ and $J_{a}|_{D^{*}_{\infty}}$. Note that $J_{a}|_{D^{*}_{0}}$ is a torus of Coxeter type, and hence $J_{a}(\cO_{0})$ is a pro-unipotent group.

Let $\cP_{a}$ be the moduli stack of $J_{a}$-torsors over $\PP^{1}$. This is indeed a commutative group scheme over $\CC$. As $a$ varies in $\cA_{G}$, $\cP_{a}$ form a smooth group scheme $\cP\to \cA_{G}$.

\begin{exam} If $G=\SL_{n}$, and $Y_{a}\to X$ be the spectral curve partially normalized over $0$ and $\infty$. Then $\cP_{a}$ is the moduli space of $(\cL,v_{0},v_{\infty}, \io)$ where $\cL$ is a line bundle of degree 0 on $Y_{a}$, $v_{0}\in \cL_{0}$ and $v_{\infty}\in \cL_{\infty}$ are bases and $\io$ is an isomorphism $\Nm_{Y_{a}/X}\cL\cong \cO_{\PP^{1}}$ that sends $v_{0}^{\ot n}$ and $v_{\infty}^{\ot n}$ to the canonical basis $1$ of $\cO$ at $0$ and $\infty$. Therefore $\cP_{a}$ is a $\mu_{n}$-torsor over $\Jac(Y_{a})$.
\end{exam}

\begin{exam} If $G=\PGL_{n}$, and $Y_{a}\to X$ be the spectral curve partially normalized over $0$ and $\infty$. Let $\Pic(Y_{a}; \wt 0,\wt \infty)$ be the moduli space of line bundles on $Y_{a}$ with trivializations at $\wt 0$ and $\wt \infty$. Similarly define $\Pic(\PP^{1}; 0,\infty)$. Then $\cP_{a}\cong\Pic(Y_{a};\wt 0,\wt\infty)/\Pic(\PP^{1}; 0,\infty)$. Therefore $\pi_{0}(\cP_{a})=\ZZ/n\ZZ$, and the neutral component of $\cP_{a}$ is isomorphic to $\Jac(Y_{a})$.
\end{exam}

\sss{Singularities of Hitchin fibers}
Following \cite[\S4.9]{NgoFL}, for $a\in \cA_G(k)$, define $\d_a$ to be the dimension of the affine part of the commutative $k$-group $\cP_a$. The invariant $\d_a$ is a rough measure the singularity of the Hitchin fiber $\cM_a$. 

We have a local variant $\d_{a}(t)$ for all $t\in \PP^1(k)$, nonzero for finitely many $t$, so that $\d_a=\sum_{t\in \PP^1(k)} \d_a(t)$. One can show that for $a\in \cA_{G}(k)$, the only contribution to $\d_{a}$ are from $t=\pm1$.  We denote $\d_a^+=\d_a(1)$ and $\d_a^-=\d_a(-1)$. Then 
\begin{equation}\label{da}
    \d_a=\d_a^++\d_a^-.
\end{equation}

\begin{exam}
Consider the case $G=\SL_{n}$. The only singularities of the spectral curve $Y_{a}$ are over $t=\pm1$. Let $p(y)=y^{n}+a_{2}y^{n-1}-\cdots+(-1)^{n}a_{n}$. Let $\{y_{i}\}$ be the roots of $p'(y)=0$, each with multiplicity $m_{i}$. Then $\sum_{i}m_{i}=n-1$. If $p(y_{i})\ne\pm2$, then $(y_{i},\pm1)$ does not lie on $Y_{a}$. If $p(y_{i})\pm2(-1)^{n}=0$, then $(y_{i},\pm1)\in Y_{a}$. The local equation of $Y_{a}$ at $(y_{i}, \pm1)$ is of the form $\y^{m_{i}+1}+s^{2}=0$. Therefore its contribution to the $\d$-invariant is $\d_{i}=[(m_{i}+1)/2]$. We see that
\begin{equation*}
\d_{Y_{a}}=\sum_{p'(y_{i})=0, p(y_{i})=\pm2}[\frac{m_{i}+1}{2}].
\end{equation*}
The decomposition \eqref{da} is according as whether $p(y_i)=2$ or $p(y_i)=-2$.
\end{exam}

\sss{Product formula}
When $\d_{a}=0$, the fiber $\cM_{a}:=f^{-1}(a)$ is a torsor under $\cP_{a}$. In general, we can apply Ng\^o's product formula to describe $\cM_{a}$.

Let $\Sp_{a}^{+}$ be the affine Springer fiber in the affine Grassmannian $\Gr$ at $t=1$ for the regular semisimple conjugacy class $a\in \frc(\cO_{1})$. Let $\cP^{+}_{a}$ (whose neutral component is an affine commutative group scheme) be the local version of $\cP_{a}$ at $t=1$ acting on $\Sp_{a}^{+}$. We have $\dim \Sp_a^+=\dim \cP_a^+=\d_a^+$. Similarly define $\Sp_{a}^{-}$ and $\cP_{a}^{-}$ by localizing at $t=-1$, both having dimension $\d_a^-$. The same argument as in \cite[Prop.4.15.1]{NgoFL} and \cite[\S2.10]{BBAMY} shows that there is a canonical universal homeomorphism
\begin{equation*}
\cP_{a}\twtimes{\cP_{a}^{+}\times\cP^{-}_{a}}(\Sp_{a}
^{+}\times\Sp^{-}_{a})\to \cM_{a}.
\end{equation*}

\begin{cor} The map $f_G:\cM_{G}\to \cA_{G}$ is flat \footnote{The flatness of $f_G$ also follows from Proposition \ref{p:nu finite flat}}.
\end{cor}
\begin{proof} By the product formula, and the fact that $\dim \cP_{a}^{\pm}=\dim\Sp_{a}^{\pm}=\d_a^{\pm}$, all geometric fibers of $f$ have the same dimension as $\cP_{a}$, which is $r$. Since $\cM_{G}$ is smooth of dimension $2r$ and $\cA_{G}$ is smooth of dimension $r$, $f$ is flat.
\end{proof}

\subsection{The map to $\Bun_{G}$} We compute the image of the map $\cM_G\to \Bun_{G}(\bI_{0}, \bI_{\infty})$ and use it to get a stratification on $\cM_G$.

\sss{Extended affine Weyl group}
The extended affine Weyl group $\tilW=N_{G\lr{t}}(T)/T\tl{t}$ is canonically isomorphic to $\xcoch(T)\rtimes W$. Let $\Om\subset \tilW$ be the image of $N_{G\lr{t}}(\bI_0)$. Let $s_i\in \tilW$ be the simple reflection for $i\in \wt I$; they generate the affine Weyl group $\Wa\subset \tilW$. We have $\tilW=\Wa\rtimes \Om$. The length function $\ell$ on $\Wa$ extends to $\tilW$ by declaring $\ell(w\om)=\ell(w)$ whenever $w\in \Wa$ and $\om\in\Om$.

For $J\subset \wt I$, let $W_J\subset \Wa$ be the subgroup generated by $\{s_j|j\in J\}$; $J$ is called {\em finite type} if $W_J$ is finite, in which case it is a finite Weyl group. We denote $J\sft \wt I$ when $J$ is a subset of $\wt I$ of finite type. In this case, denote by $w_{J}$ be the longest element in $W_{J}$. 

\sss{Pseudo-Levi subgroups}
Parahoric subgroups of $G\lr{t}$ containing $\bI_0$ are in bijection with proper subsets $J$ of $\wt I$. For $J\sft \wt I$, let $\bP_{J,0}$ be the corresponding parahoric subgroup of $G\lr{t}$, whose reductive quotient is denoted by $L_{J}$. There is a unique section $L_J\subset \bP_{J,0}$ whose image contains $T$, through which we identify $L_J$ as a subgroup of $G\lr{t}$, which in fact is contained in the polynomial loop group $G[t,t^{-1}]$. 

Similarly, parahoric subgroups of $G\lr{\t}$ containing $\bI_\infty$ are also in bijection with proper subsets $J$ of $\wt I$. For $J\sft \wt I$, we have the corresponding parahoric subgroup $\bP_{\infty,J}\supset \bI_\infty$ that also contains $L_J$ (as a subgroup of $G[t,t^{-1}]$, hence of $G\lr{\t}$) as a Levi subgroup. 

Let $B_J$ (resp. $B_J^{-}$) be the image of the projection $\bP_{J,0}\to L_J$ (resp. $\bP_{\infty,J}\to L_J$). Then $B_J$ and $B_J^{-}$ are opposite Borel subgroups of $L_J$ whose intersection is $T$. Let $N_J$ (resp. $N_J^{-}$) be the unipotent radical of $B_J$ (resp. $B_J^{-}$). Denote the Lie algebras of the groups introduced above by $\frb_J, \frb^-_J, \frn_J$ and $\frn^-_J$.

\sss{Strata in $\Bun_G$} Uniformization at $0$ gives an isomorphism
\begin{equation}\label{Bun I unif}
    \Bun_{G}(\bI_{0}, \bI_{\infty})\cong (\bI_\infty\cap G[\t])\bs G\lr{t}/\bI_0.
\end{equation}
Here $G[\t]\subset G\tl{\t}$ is the sub-ind-scheme of polynomial arcs.
By Birkhoff decomposition, the double cosets above are in canonical bijection with the extended affine Weyl group $\tilW$. For $w\in \tilW$, let $\Bun^{w}_{G}(\bI_{0}, \bI_{\infty})\subset \Bun_{G}(\bI_{0}, \bI_{\infty})$ be the locally closed substack corresponding to the double coset of $w$ under \eqref{Bun I unif}.

Consider the projection $\Bun_{G}(\bI^{+}_{0}, \bI^{+}_{\infty})\to \Bun_{G}(\bI_{0}, \bI_{\infty})$. For $w\in \tilW$, let $\Bun^{w}_{G}(\bI^{+}_{0}, \bI^{+}_{\infty})$ be the preimage of $\Bun^w_{G}(\bI_{0}, \bI_{\infty})$. 

Forgetting the Higgs fields gives a natural map
\begin{equation*}
\b: \cM_G\to \Bun_{G}(\bI^{+}_{0}, \bI^{+}_{\infty}).    
\end{equation*}

\begin{lemma}\label{l:rel wJ} The image of $\b$ is contained in the union 
\begin{equation*}
\bigsqcup_{J\sft \wt I, \om \in \Om}\Bun^{w_{J}\om}_{G}(\bI^{+}_{0}, \bI^{+}_{\infty}).
\end{equation*}
For each $J\sft \wt I$, the intersection of $\Bun^{w_{J}\om}_{G}(\bI^{+}_{0}, \bI^{+}_{\infty})$ with the image of $\b$ has coarse moduli space equal to a torsor under the center $Z(L_{J})$ of $L_J$.
\end{lemma}
\begin{proof}
Call a point of $\Bun_G(\bI^{+}_{0}, \bI^{+}_{\infty})$ {\em relevant} if it is in the image of $\b$. For a point in $\Bun^{w}_{G}(\bI^{+}_{0}, \bI^{+}_{\infty})$ represented by $\dot w h$, where $\dot w$ is a lifting of $w$ and $h\in T$, it is relevant if and only if $\psi_{0}$ and  ${}^{\dot wh}\psi_{\infty}$ restrict to the same character on $\bJ_{w}:=\bI_{0}^{+}\cap {}^{w}\bI^{+}_{\infty}$. If $\a_{i}$ is a simple affine root but $w^{-1}\a_{i}<0$, then the $\a_{i}$ appears in $\bJ_{w}$. If $-w^{-1}\a_{i}$ was not a simple affine root then ${}^{\dot wh}\psi_{\infty}$ would restrict to zero on the $\a_{i}=w(w^{-1}\a_{i})$ root space of $\bJ_{w}$, while $\psi_{0}$ restricts nontrivially, so $\dot wh$ would not be relevant. Therefore, if  $\dot w h$ is relevant, then for any $i\in \wt I$, either $w^{-1}\a_{i}>0$, or $w^{-1}\a_{i}$ is the negative of an affine simple root.

Let $A$ be the fundamental alcove in the building of $G\lr{t}$ with respect to $\bI_0$. The above condition on $w$ means that for any wall $H_{i}, i\in \wt I$ (defined by the vanishing of $\a_{i}$) of $A$, either $A$ and $wA$ lie on the same side of the wall, or they share this wall but lie on opposite sides. Let $J\sft \wt I$ be the set of $j$ such that $A$ and $wA$ lie on opposite sides of $H_{j}$.  Consider the alcoves $w_{J}A$ and $wA$. We see that $w_{J}A$ and $wA$ lie on the same side of $H_{j}$ for all $j\in J$. For $i\notin J$, they both lie on the same side of $H_{i}$ as $A$. Therefore $w_{J}A=wA$, hence $w=w_{J}\om$ for some $\om\in \Om$. 

For fixed $w=w_{J}\om$, the choices of $h\in T$ such that $\dot wh$ is relevant has to make $\psi_0$ and ${}^{\dot wh}\psi_{\infty}$ to restrict to the same linear function on the root space of $\a_{j}$ for $j\in J$. Such $h$ form a torsor under the kernel of all $\a_{j}, j\in J$, i.e., the center of $L_{J}$.
\end{proof}

\sss{Stratification of $\cM^\c_G$}\label{sss: strat}
Let $\cM_{G,J}=\b^{-1}\left(\Bun^{w_J}_G(\Iu)\right)$. Then Lemma \ref{l:rel wJ} implies a stratification
\begin{equation}\label{MG strat}
    \cM^\c_{G}=\bigsqcup_{J\sft \wt I}\cM_{G,J}.
\end{equation}
Since the closure order among $\Bun^w_G(\Iu)$ is opposite to the Bruhat order on $\tilW$, the union
\begin{equation*}
    \cM_{G,\subset J}:=\bigsqcup_{J'\subset J}\cM_{G,J}
\end{equation*}
is open in $\cM^\c_G$ while the union
\begin{equation*}
    \cM_{G,\supset J}:=\bigsqcup_{J'\supset J}\cM_{G,J}
\end{equation*}
is closed in $\cM^\c_G$.

\subsection{Open covering of $\cM_G$ by torsors of regular centralizers}

In this subsection we identify the open subspace $\cM_{G,\subset J}$ of $\cM^\c_G$ as a torsor for the Toda system of $L_J$.

Let $J\sft \wt I$. Consider the map $\Bun_{G}(\bI^+_0, \bI^+_\infty)\to \Bun_{G}(\bP_{0,J},\bP_{\infty,J})$. We have a double coset presentation of $\Bun_{G}(\bP_{0,J},\bP_{\infty,J})$ similar to \eqref{Bun I unif}, in which the unit double coset gives an open embedding $\pt/L_{J}\incl \Bun_{G}(\bP_{0,J},\bP_{\infty,J})$. Taking the preimage of this open substack in $\Bun_{G}(\bI^+_0, \bI^+_\infty)$ we get an open embedding
\begin{equation*}
    i^{\Bun}_J: N^-_{J}\bs L_{J}/N_{J}\subset \Bun_{G}(\Iu). 
\end{equation*}
We have
\begin{equation}\label{im bJ} \Im(i^{\Bun}_J)=\Bun^{\le w_J}_G(\Iu).
\end{equation}

\sss{The space $\cT_{L_J}$}
Since the simple roots of $L_J$ with respect to $(T, B_J)$ are exactly 
$\a_j$ for $j\in J$, $\psi_0$ restricts to a generic linear character
\begin{equation*}
    \psi_{0,J}:N_J\surj \prod_{j\in J}\frg\lr{t}_{\a_j}\xr{\psi_0}\Ga. 
\end{equation*}
Similarly, $\psi_0$ restricts to a generic linear character
\begin{equation*}
    \psi_{\infty,J}:N^{-}_J\surj \prod_{j\in J}\frg\lr{\t}_{-\a_j}\xr{\psi_\infty}\Ga. 
\end{equation*}
We view $\psi_{0,J}$ (resp. $\psi_{\infty,J}$), or rather its differential, as an element in $\frn_J^*$ (resp. $(\frn^-)^*$).

Consider the twisted cotangent bundle
\begin{equation}\label{eq: cT_LJ}
\cT_{L_{J}}:=T^{*}\left((N^-_{J},-\psi_{\infty,J})\bs L_{J}/ (N_{J},\psi_{0, J})\right)
\end{equation}
Here the right side is understood as the Hamiltonian reduction of $T^*L_J$ with respect to the left and right translations of $N_J\times N^-_J$ (both as right actions) at the moment map value $(\psi_{0,J}, -\psi_{\infty, J})\in \frn_J^*\op (\frn^-_J)^{*}$.

Unraveling the definition of $\cT_{L_J}$, we see
\begin{equation}\label{cT desc}
    \cT_{L_J}=N_J^{-}\bs \{(h, \th)\in L_J\times \frl^*_J|\th\in \psi_{0,J}+\frn_J^\bot, \Ad(h)^*\th \in -\psi_{\infty,J}+(\frn^-_J)^\bot\}/N_J
\end{equation}
Here $n\in N_J$ acts by $(h,\th)\mt (hn, \Ad^*(n^{-1})\th)$, and $n^-\in N^-_J$ acts by $(h,\th)\mt (n^-h, \th)$.

\begin{prop}\label{p:i_J cover}
For any $J\sft \wt I$, there is a canonical symplectomorphism  
    \begin{equation*}
i_{J}: \cT_{L_J}\cong \cM_{G,\subset J}
\end{equation*}
making the following diagram commutative
\begin{equation}\label{iJ bJ}
    \xymatrix{\cT_{L_J}\ar[d]\ar[r]^-{i_J}_-{\sim} & \cM_{G,\subset J}\ar[d]^{\b|_{\cM_{G,\subset J}}}\\
    N^-_J\bs L_J/N_J\ar[r]^-{i^\Bun_J}_-{\sim} & \Bun^{\le w_J}_G(\bI^+_0, \bI^+_\infty)}
\end{equation}

\end{prop}
\begin{proof}
The open subspace $\cM_{G,\subset J}$ is the preimage of the open point $\pt/L_J\incl \Bun_G(\bP_{0,J}, \bP_{\infty,J})$ under the projection $p_J: \cM_G\to \Bun_G(\Iu)\to \Bun_G(\bP_{0,J}, \bP_{\infty,J})$, which we now describe. 

Let $\wt \cT_{J}$ be the fiber of $p_J$ over the base point $\pt\to \Bun_G(\bP_{0,J}, \bP_{\infty,J})$ corresponding to the trivial bundle with the tautological $\bP_{0,J}$ and $\bP_{\infty, J}$-level structures. Then $\wt \cT_J$ classifies triples $(h_0N_J, h_\infty N^-_J, \ph)$ where $h_0N_J\in L_J/N_J$, $h_\infty\in L_J/N_J^-$, and $\ph$ is the Higgs field for the trivial $G$-bundle of the form
\begin{equation*}
    \ph=\ph_0^J+\ph_J \frac{dt}{t}+\ph_\infty^J.
\end{equation*}
Here $\ph_J\in \frl_J^*$, $\ph_0^J\in (V_0^J)^*$ and $\ph_\infty^J\in (V_\infty^J)^*$, where $V_0^J=\bP^+_{0,J}/[\bP^+_{0,J},\bP^+_{0,J}]$ (as an $L_J$-module) and $V_\infty^J=\bP^+_{\infty,J}/[\bP^+_{\infty,J},\bP^+_{\infty,J}]$. Let $\psi^J_0\in (V_0^J)^*$ (resp. $\psi_\infty^J\in (V_\infty^J)^*$) be the restriction of $\psi_0$ (resp. $\psi_\infty$) to those simple affine root spaces not in $J$. The conditions \eqref{ph near 0} and \eqref{ph near infty} translate to
\begin{eqnarray}\label{th conditions}
    \Ad^*(h_0^{-1})\ph_0^J=\psi_0^J, \quad \Ad^*(h_0^{-1})\ph_J\in \psi_{0,J}+\frn_J^\bot,\\
    \Ad^*(h_\infty^{-1})\ph_\infty^J=\psi_\infty^J, \quad \Ad^*(h_\infty^{-1})\ph_J\in -\psi_{\infty,J}+(\frn^-_J)^\bot.
\end{eqnarray}
The negative sign in the last condition comes from the fact that $\frac{dt}{t}=-\frac{d\t}{\t}$. Since $\ph_0^J$ and $\ph_\infty^J$ are determined by $h_0$ and $h_\infty$, $\wt\cT_J$ classifies triples $(h_0N_J, h_\infty N_J^-,\ph_J)\in L_J/N_J\times L_J/N_J^-\times \frl_J^*$ satisfying the two conditions in \eqref{th conditions} involving $\ph_J$. There is a map
\begin{equation}\label{wtT to T}
    \wt\cT_{J}\to \cT_{L_J}
\end{equation}
sending $(h_0,h_\infty, \ph_J)$ to $(h_\infty^{-1}h_0, \Ad^*(h_0^{-1})\ph_J)\in \cT_{L_J}$ under the description \eqref{cT desc}. It is easy to see that \eqref{wtT to T} is an $L_J$-torsor for the diagonal $L_J$-action on $\wt\cT_J$.  On the other hand, $\wt\cT_J\to p_J^{-1}(\pt/L_J)$ is also an $L_J$-torsor for the diagonal $L_J$-action, therefore we get a canonical isomorphism $p_J^{-1}(\pt/L_J)\cong \cT_{L_J}$. It is easy to see that this isomorphism preserves the symplectic structures.

\end{proof}

\begin{cor}\label{c:MGc conn}
    The space $\cM^\c_G$ is connected.
\end{cor}

\begin{proof} By \eqref{MG strat}, $\cM^\c_{G}$ is covered by open subspaces $\cM_{G,\subset J}$ for various $J\sft \wt I$. Each $\cM_{G,\subset J}$ contains $\cM_{G,\vn}$, which is isomorphic to $T^*T$ by Proposition \ref{p:i_J cover}, we conclude that $\cM^\c_G$ is connected.
\end{proof}

\sss{Torsor for regular centralizer}
Let $\frc_J=\frl_J^*\sslash L_J\cong \frt^*\sslash W_J$. Let
\begin{equation*}
    \cJ_{L_{J}}\to \frc_{J}=\frt^{*}\sslash W_{J}.
\end{equation*}
be the Lie algebra version of the regular centralizer group scheme of $L_J$, see \cite[\S2.1]{NgoFL}. This is the classical Toda system for $L_J$.

We construct an action of $\cJ_{L_{J}}$ on $\cT_{L_J}$ as follows.

Let $\k_{0,J}: \frc^*_J\to \psi_{0,J}+\frn^\bot_J$ be a Kostant slice through the regular nilpotent covector $\psi_{0,J}$ of $\frl_J$. It is a section of the $N_J$-torsor $\psi_{0,J}+\frn^\bot_J\to \frc^*_J$. Similarly define $\k_{\infty ,J}:  \frc^*_J\to -\psi_{\infty,J}+(\frn^-_J)^\bot$. We can rewrite \eqref{cT desc} as
\begin{equation*}
    \cT_{L_J}\cong \{(h,a)\in L_J\times \frc^*_J|\Ad^*(h)(\k_{0,J}(a))=\k_{\infty, J}(a)\}.
\end{equation*}
On the other hand, 
\begin{equation*}
    \cJ_{L_J}\cong \{(g,a)\in L_J\times \frc^*_J|\Ad^*(g)(\k_{0,J}(a))=\k_{0,J}(a)\}.
\end{equation*}
From this we get a right action of $\cJ_{L_J}$ on $\cT_{L_J}$ over $\frc^*_J$ by
\begin{equation*}
    (h,a)\cdot (g,a)=(hg,a), \textup{ where } a\in \frc^*_J, (g,a)\in \cJ_{L_J}, (h,a)\in \cT_{L_J}.
\end{equation*}

The following lemma follows immediately from the fact that for $a\in \frc^*_J$, $\k_{0,J}(a)$ and $\k_{0,\infty}(a)$ are in the same regular coadjoint orbit of $L_J$.

\begin{lemma}
    With the action of $\cJ_{L_{J}}$ on $\cT_{L_J}$ defined above, $\cT_{L_J}$ is a  $\cJ_{L_{J}}$-torsor.
\end{lemma}

\begin{exam}\label{exam: SL2coord} Consider $G=\SL_{2}$, we give a formula for $i_I: \cT_{G}\incl \cM_G$ (where $I=\{1\}$ indexes the unique simple root of $G$). 

To fix ideas, we identify $\frg$ with $\frg^*$ using the trace pairing, and take 
\begin{equation}\label{psi SL2}
    \psi_0=\mat{}{t^{-1}}{1}{}\frac{dt}{t}, \quad \psi_\infty=\mat{}{1}{t}{}\frac{dt}{t},
\end{equation}
We use Kostant sections $\k_{0}: a\mapsto \mat{}{-a}{1}{}$, and $\k_{\infty}: a\mapsto \mat{}{1}{-a}{}$. We get that $\cT_{G}=\{(\mat{x}{y}{y}{-ax},a)|y^{2}+ax^{2}+1=0\}$.  Let $g(a,x,y)=\mat{x}{y}{y}{-ax}\in \SL_{2}$.  The map $i_{I}: \cT_{G}\incl \cM_{G}$ sends $(g(a,x,y),a)$ to the trivial $\cE$ with the standard $N$-reduction at $0$ and the $N^-$-reduction $g(a,x,y)^{-1}N^{-}$ at $\infty$ and Higgs field
\begin{equation*}
\ph=\left(\Ad(g(a,x,y)^{-1})\mat{}{}{t}{}+\mat{}{-a}{1}{}+\mat{}{t^{-1}}{}{}\right)\frac{dt}{t}=\mat{xyt}{-y^{2}t-a+t^{-1}}{x^{2}t+1}{-xyt}\frac{dt}{t}.
\end{equation*}
Taking determinant, we see that the composition
\begin{equation*}
\cT_{G}\incl \cM_{G}\to \cA_{G}\cong \AA^{1}
\end{equation*}
sends $(g(a,x,y),a)$ to $-(t+(x^{2}-a)+t^{-1})(dt/t)^{2}$, i.e., sending it to $x^{2}-a\in \AA^{1}$.

In particular, the stratum $\cM_{G,\{1\}}$ corresponds to the locus $x=0$ in $\cT_G$. Then $y$ satisfies $y^2=-1$. The locus $y=i$ consists of the trivial bundle $\cO^2$ with the standard $N$-reduction at $0$ and the $N^-$-reduction $\mat{}{-i}{-i}{}N^{-}$ at $\infty$ and Higgs field $\ph=\mat{}{t-a+t^{-1}}{1}{}\frac{dt}{t}$, giving a section of the Hitchin map $f_G$. We conclude that $\cM_{G,\{1\}}$ consists of two disjoint sections of $f_G$. 
\end{exam}

\begin{exam}\label{exam: SL2hxi}
We make explicit the embedding 
\begin{equation*}
    i_\vn: \cT_{T}=T^*T\incl \cM_G.
\end{equation*} 
It sends $(h,\xi)\in T\times \frt^*$ to the trivial $G$-bundle $\cE_0$ with standard $N$-reduction at $0$ and $hN^{-}$ at $\infty$, and the Higgs field $\ph=\psi_{0}+\xi\frac{dt}{t}+\Ad^*(h^{-1})\psi_{\infty}$.

For $G=\SL_{2}$, with the choices of $\psi_0$ and $\psi_\infty$ as in \eqref{psi SL2}, $i_{\vn}$ sends $(\mat{h}{}{}{h^{-1}}, \mat{\xi}{}{}{-\xi})$ (where $h\in \Gm$ and $\xi\in \AA^1$) to the trivial bundle with Higgs field $\mat{\xi}{h^{-2}+t^{-1}}{h^2t+1}{-\xi}\frac{dt}{t}$. Taking determinant we get $-(t+(\xi^{2}+h^{2}+h^{-2})+t^{-1})(dt/t)^{2}$. Hence the composite $T^{*}T\incl \cM_{G}\to \cA_{G}\cong \AA^{1}$ sends  $(h,\xi)$ to $\xi^{2}+h^{2}+h^{-2}\in \AA^{1}$. 
\end{exam}

\sss{Geometry of $\cM_{G,J}$}\label{sss:geom MGJ} We describe each stratum $\cM_{G,J}$ of $\cM^\c_G$ in more explicit terms. By Proposition \ref{p:i_J cover}, $\cM_{G,J}$ is isomorphic to the restriction of $\cT_{L_J}$ over the closed substack $N_J^-\bs N_J^- w_J B_J/N_J\subset N_J^-\bs L_J/N_J$.

There is a canonical action of $T^*Z(L_J)$ on $\cT_{L_J}$: $(z,\z)\in T^*Z(L_J)\cong Z(L_J)\times \frz(L_J)$ acts by $(h,a)\mt (hz,a+\z)$. In particular, there is a canonical action of $T^*Z(L_J)$ on $\cM_{G,J}$.

Restricting $\b$ to $\cM_{G,J}$ we get a map
\begin{equation*}
    \cM_{G,J}\to \Bun^{w_J}_G(\Iu)\cong N_J^-\bs N_J^- w_J B_J/N_J\to w_J T
\end{equation*}
where the last map realizes $w_J T$ as the coarse moduli space of $N_J^-\bs N_J^- w_J B_J/N_J$. By Lemma \eqref{l:rel wJ}, the image of the above map is a $Z(L_J)$-torsor that we denote by $\cW_J\subset w_J T$. We denote the induced map by
\begin{equation*}
    \b_J: \cM_{G,J}\to \cW_J.
\end{equation*}
On the other hand, we have a projection $\chi_J: \cT_{L_J}\to \frc_J^*$.

\begin{lemma}\label{l:MGJ}
    The map
    \begin{equation*}
        \cM_{G,J}\xr{(\b_J, \chi_J \c i_J^{-1})}\cW_J\times \frc^*_J
    \end{equation*}
    is a $T^*Z(L_J)$-equivariant isomorphism. Here the two factors of  $T^*Z(L_J)\cong Z(L_J)\times \frz(L_J)$ acts on the two factors of $\cW_J\times \frc^*_J$ respectively by translation.
\end{lemma}
\begin{proof}

    It suffices to show that $\cM_{G,J}\to \frc^*_J$ is a $Z(L_J)$-torsor. We  identify $\frl_J$ with $\frl^*_J$ using an $L_J$-invariant symmetric blinear form on $\frl_J$, so that $\psi_{0,J}\in \op_{j\in J}\frl_{J, -\a_j}$ is a regular nilpotent element, and $\psi_{\infty,J}\in \op_{j\in J}\frl_{J, \a_j}$ is also regular nilpotent. Let $\k: \frc^*_J\to \frl^*_J\cong \frl_J$ be a Kostant slice defined by an $\sl_2$-triple $(\psi_{0,J},-2\rho_J^\vee,f)$, so that the image of $\k$ lies in $\psi_0+\frn_J^\bot=\psi_0+\frb_J$. We can rewrite \eqref{cT desc} as
    \begin{equation*}
        \cT_{L_J}=\{(N_J^-h, a)\in N_J^-\bs L_J\times \frc^*_J|\Ad(h)(\k(a))\in -\psi_{\infty,J}+\frb^-_J\}.
    \end{equation*}
    Under this identification, $\cM_{G,J}$ corresponds to the subscheme of $\cT_{L_J}$ defined by the condition that $N_J^-h\in N_J^-\bs (N_J^-w_JB_J)\cong w_JT$. Therefore 
    \begin{equation*}
        \cM_{G,J}\cong \{(\dot w_J,a)\in w_JT\times \frc^*_J|\Ad(\dot w_J)(\k(a))\in -\psi_{\infty,J}+\frb^-_J\}.
    \end{equation*}
    For each $a\in \frc^*_J$, $\k(a)\in\psi_{0,J}+\frb_J$. For $\dot w_J\in w_JT$, $\Ad(\dot w_J)(\k(a))\in \Ad(\dot w_J)\psi_{0,J}+\Ad(\dot w_J)\frb_J=\Ad(\dot w_J)\psi_{0,J}+\frb^-_J$. Now $\Ad(\dot w_J)\psi_{0,J}$ is a sum of nonzero elements from each simple root space of $\frl_J$. Therefore, $(\dot w_J, a)\in \cM_{G,J}$ if and only if $\Ad(\dot w_J)\psi_{0,J}=-\psi_{\infty,J}$. Such $\dot w_J$ clearly form a torsor for $Z(L_J)$.
\end{proof}

\begin{remark}
    The last paragraph of the proof of Lemma \ref{l:MGJ} gives the following description of $\cW_J$:
    \begin{equation*}
        \cW_J\cong \{\dot w_J\in w_JT|\Ad^*(\dot w_J)\psi_{0,J}=-\psi_{\infty, J}\}.
    \end{equation*}
    Any point $\dot w_J\in \cW_J(k)$ gives a section to the $Z(L_J)$-torsor $\cM_{G,J}\to \frc^*_J$, hence trivializes it; the same section can be viewed as a section to $\cT_{L_J}\to \frc^*_J$, hence trivializes $\cT_{L_J}$ as a $\cJ_{L_J}$-torsor. The resulting isomorphism $\dot w_J: \cJ_{L_T}\isom \cT_{L_J}$ can be seen directly from the definition of $\cT_{L_J}$ in \eqref{eq: cT_LJ}. Namely, conjugation by  $\dot w_J$ identifies $(N_J, \psi_{0,J})$ with $(N^-_J,-\psi_{\infty, J})$, hence left multiplication by $\dot w_J$ induces an isomorphism 
    $$\dot w_J: \cJ_{L_J}\cong T^*((N_J,\psi_{0,J})\bs L_J/(N_J,\psi_{0,J}))\isom T^*((N^-_J,-\psi_{\infty,J})\bs L_J/(N_J,\psi_{0,J}))=\cT_{L_J}.$$
\end{remark}

\sss{Other components of $\cM_G$}
The above discussions about $\cT_{L_J}$ and $\cM_{G,J}$ can be extended to other components of $\cM_G$ with essentially the same proofs. We will state the results only. 

For $\om\in \Om$, let $\cM^\om_G$ be the preimage of the $\om$-component of $\Bun_G(\Iu)$ under $\b$. It is stratified into the union of $\cM^\om_{G,J}$, the preimage of $\Bun^{w_J\om}_G(\Iu)$ under $\b$. Similarly define $\cM^\om_{G,\subset J}$

Let $J'=\om^{-1}(J)$. Define
\begin{equation*}
    \cT^\om_{L_J}=T^*((N^-_J, -\psi_{\infty, J})\bs L_J\om/(N_{J'}, \psi_{0,J'})).
\end{equation*}
Here $L_J\om$ makes sense (without choosing a lifting of $\om$) because $T\om$ is a well-defined $T$-coset. We have the map
\begin{equation*}
    \chi^\om_{J}: \cT^\om_{L_J}\to \psi_{0,J}+\frn_J^\bot\to \frc_J^*.
\end{equation*}

Then $\cT^\om_{L_J}$ is a left torsor under $\cJ_{L_J}$ and a right torsor under $\cJ_{L_{J'}}$. 

Restricting $\b$ to $\cM^\om_{G,J}$ we get a map
\begin{equation*}
    \cM^\om_{G,J}\to \Bun^{w_J\om}_G(\Iu)\to w_J\om T.
\end{equation*}
Denote the image of this map by $\cW^\om_J$, which is a $(Z(L_J),Z(L_{J'}))$-bitorsor by Lemma \ref{l:rel wJ}. We thus get a map $\b_J^\om: \cM^\om_{G,J}\to \cW^\om_J$. We omit the proof of the following.

\begin{prop}\label{p:MGJ om} Let $\om\in \Om$.
    \begin{enumerate}
        \item There is a canonical symplectomorphism $i^\om_J:\cT^\om_{L_J}\cong \cM^\om_{G,\subset J}$.
        \item The natural map
        \begin{equation*}
            \cM^\om_{G,J}\xr{(\b_J^\om, \chi^\om_J\c (i^\om_J)^{-1})}\cW^\om_J\times \frc^*_J
        \end{equation*}
        is an isomorphism.
        \item $\cM^\om_G$ is connected.
    \end{enumerate}
\end{prop}

\begin{cor}\label{c:MG scheme}
    The algebraic space $\cM_G$ is a scheme.
\end{cor}
\begin{proof}
    By Proposition \ref{p:MGJ om}, $\cM_G$ has an open cover by $\cT^\om_{L_J}$, which are torsors for $\cJ_{L_J}$ hence themselves schemes.
\end{proof}

\subsection{A horizontal completely integrable system}\label{subsec: horizontal,integrable}

In this subsection we will modify the map $\cM^\c_G\to \Bun^\c_G(\Iu)$ to construct another completely integrable system structure on $\cM^\c_G$. This integrable system will be used in \S\ref{s:Fukaya} to produce geometric structures such as Liouville 1-forms and Weinstein sectors on $\cM^\c_G$, which are necessary to define the Fukaya wrapped category of $\cM^\c_G$.

\sss{The toric variety $\PP$}\label{sssToricPP}
Let $\Sigma_{\wt{I}}$ be the fan in $\XX_*(T)_{\RR}$ defined by the set of closed cones $\{\s_J\}_{J\sft \wt I}$, where
\begin{align*}
\sigma_J=\sum\limits_{j\in J}\RR_{\geq 0} \alpha^\vee_j. 
\end{align*}
For each cone $\sigma_J$, let $\sigma_J^\vee=\{\xi\in \xch(T)_\RR|\j{\xi, \a_j^\vee}\ge0, \forall j\in J\}\}$ be the dual cone in $\xch(T)_{\RR}$, and $\s_J^\bot=\{\xi\in \xch(T)_\RR|\j{\xi, \a_j^\vee}=0, \forall j\in J\}$.

Then $\Sig_{\wt I}$ defines a $T$-toric variety $\PP:=\PP_{\Sig_{\wt I}}$. It is the union of affine open subschemes $U_J=\Spec \CC[\xch(T)\cap\s_J^\vee]$. It is also stratified into $T$-orbits indexed by $J\sft \wt I$:
\begin{equation*}
    \PP=\bigsqcup_{J\sft \wt I} \PP_J
\end{equation*}
where $\PP_J=\Spec k[\xch(T)\cap \s_J^\bot]$ as  a closed subscheme of $U_J=\Spec k[\xch(T)\cap \s_J^\vee]$. Note that $\PP_J$ can be identified with $(L_J)_{\ab}=L_J/L^{\der}_J$, as a quotient torus of $T$. 

When $G$ is semisimple, $\PP$ is projective.  When $G$ is simply-connected, $\PP$ is isomorphic to the weighted projective space of dimension $r$ with weights $(n_i)_{i\in \wt I}$, where $n_0=1$ and $n_i$ is the coefficient of $\a_i$ in the highest root $\th$, for $i\in I$.

\sss{The map $\BB\to \PP$}
Let $\BB\subset \Bun^\c_G(\Iu)$ be the open substack that is the union 
\begin{equation*}
    \BB=\bigcup_{w\le w_J \textup{ for some } J\sft \wt I}\Bun^{w}_G(\Iu).
\end{equation*}
By \eqref{im bJ}, we have
\begin{equation*}
    \BB\cong \colim_{J\sft \wt I}N_J^-\bs L_J/N_J,
\end{equation*}
where the transition maps are open embeddings for $J\subset J'$. 

There is a canonical map 
\begin{equation}\label{cJ}
    c_J: N_J^-\bs L_J/N_J\to U_J=\Spec k[\xch(T)\cap\s_J^\vee].
\end{equation}
Indeed, $\xch(T)\cap\s_J^\vee$ are exactly dominant weights of $T$ with respect to the Borel subgroup $B_J$ of $L_J$. For each $\l\in \xch(T)\cap \s_J^\vee$, define the matrix coefficient
\begin{equation*}
    m_\l: L_J\ni h\mt \j{v^*_\l, hv_\l}
\end{equation*}
where $v_\l\in V_\l$ is a highest weight vector in the irreducible $L_J$-module $V_\l$ with highest weight $\l$, and $v_\l^*\in V_\l^*$ is the lowest weight vector (with weight $-\l$) such that $\j{v_\l^*, v_\l}=1$. Then $m_\l$ descends to a regular function on $N_J^-\bs L_J/N_J$, and restricts to the character $\l$ on the open subset $T\subset N_J^-\bs L_J/N_J$. The map $c_J$ is characterized by $c_J^*(\l)=m_\l$ for $\l\in \xch(T)\cap \s_J^\vee$. It is clear that the maps $c_J$ are compatible with inclusions of $J$'s, therefore they glue to give a map
\begin{equation*}
    c: \BB\to \PP.
\end{equation*}
By Lemma \ref{l:rel wJ}, the image of $\b: \cM^\c_G\to \Bun^\c_G(\Iu)$ is contained in $\BB$, therefore we get a map 
\begin{equation*}
    b: \cM^\c_G\xr{\b} \BB\xr{c} \PP.
\end{equation*}
In particular,the restriction of $b$ to $\cM_{G,\vn}\cong T^*T$ is the projection $T^*T\to T$. Therefore, $b$ is also a completely integrable system.

\sss{Fibers of $b$}\label{sss:fiber b}
The restriction of $b$ to $\cM_{G,\subset J}\cong \cT_{L_J}$ is 
\begin{equation*}
    b_{\subset J}: \cM_{G,\subset J}\cong \cT_{L_J}\to N^-_J\bs L_J/N_J\xr{c_J}\Spec k[\xch(T)\cap\s_J^\vee].
\end{equation*}
In fact, $\cM_{G,\subset J}$ is the preimage of the open subset $\Spec k[\xch(T)\cap\s_J^\vee]\subset \PP$ under $b$. The reduced structure of $b^{-1}(\PP_J)$ is $\cM_{G,J}$. The restriction of $b$ to the stratum $\cM_{G,J}$ is
\begin{equation*}
    b_{J}: \cM_{G,J}\xr{\b_J} \cW_J\xr{\mu_J} \PP_J\cong (L_J)_{\ab},
\end{equation*}
where $\mu_J$ above is the composition
    \begin{equation}\label{mu J}
        \mu_J:  \cW_J\subset w_JT\subset L_J\surj (L_J)_\ab\cong \PP_J.
    \end{equation}
    Since $\cW_J$ is a $Z(L_J)$-torsor, which also acts on $(L_J)_{\ab}$ by translation compatibly, we conclude that $\mu_J: \cW_J\to (L_J)_\ab$ is a torsor for
    \begin{equation*}
        \ker(Z(L_J)\to (L_J)_{\ab})=Z(L_J^\der).
    \end{equation*}
    By Lemma \ref{l:MGJ}, we conclude that the fibers of $b_J$ are isomorphic to the product of $\frc^*_J$ with a torsor for the finite group $Z(L_J^\der)$.

\begin{prop}\label{p:nu finite flat}
The map 
    \begin{equation}\label{fb}
    \nu=(f_G, b): \cM^\c_G\to \cA_G\times \PP
\end{equation}
is finite flat of degree $|W|$.
\end{prop}
\begin{proof}
    The map $\nu$ is proper since $f_G$ is proper. From the description of $b_{\subset J}$ we see that $b_{\subset J}$ is affine since $\cT_{L_J}$ is. Therefore $b$ is affine, and the same is true for $\nu$. We conclude that $\nu$ is proper and affine, hence finite.

    On the other hand, since $\cM_G$ is smooth of the same dimension as $\cA_G\times \PP$, the finite map $\nu$ is flat.

    It remains to show that $\nu$ has degree $|W|$. It suffices to compute the degree of the restriction of $\nu$ to $\cA_G\times \PP_\vn\cong \cA_G\times T$, which is $\nu_\vn: T^*T\to \cA_G\times T$. We may restrict further over $1\in T$ and compute the degree of the flat map
    \begin{equation*}
        \nu_{\vn,1}: \frt^*\to \cA_G.
    \end{equation*}
    By Example \ref{exam: SL2hxi}, $\nu_{\vn,1}$ sends $\xi\in \frt^*$ to the point in $\cA_G$ given by the $r$-tuple $f_i(\psi_0+\xi\frac{dt}{t}+\psi_\infty)$, $1\le i\le r$. 
    
    We may assume $G$ is almost simple. Now enlarge $\frt^*$ to $\frt^*\times V_0^*\times V_\infty^*$ by allowing $\psi_0$ and $\psi_\infty$ to vary in $V_0^*$ and $V_\infty^*$, and enlarge $\cA_G$ to $\frc^*\times V_0^*\times V_\infty^*$. The maps $\nu_{\vn,1}$ for various choices of $(\psi_0,\psi_\infty)$ combine to give a map 
    \begin{equation*}
        \mu: \frt^*\times V_0^*\times V_\infty^*\to \frc^*\times V_0^*\times V_\infty^*
    \end{equation*}
    that sends $(\xi,\psi_0,\psi_\infty)$ to $(b_1,b_2,\cdots, b_r, \psi_0,\psi_\infty)$, where $b_i$ are determined by (compare \eqref{fi ph})
    \begin{equation*}
        f_i(\psi_0+\xi\frac{dt}{t}+\psi_\infty)=\begin{cases}
            b_i\left(\frac{dt}{t}\right)^{d_i}, & 1\le i\le r-1;\\
            (a_0(\psi_0)t^{-1}+b_r+a_\infty(\psi_\infty)t)\left(\frac{dt}{t}\right)^{h}, & i=r.
        \end{cases}
    \end{equation*}
    Then $\mu$ is equivariant with respect to the $\Gm$-actions that is scaling on $\frt^*$, $V_0^*$ and $V_\infty^*$, and the usual weighted action on $\frc^*$. Both the source and the target of $\mu$ are smooth of the same dimension, hence $\mu$ is flat and quasi-finite. Since the $\Gm$-actions on both the source and the target are contracting to a single point, $\mu$ is indeed finite flat. Restricting $\mu$ over $(0,0)\in V_0^*\times V_\infty^*$, it becomes the obvious map $\frt^*\to \frc^*$ which has degree $|W|$, therefore $\mu$ has degree $|W|$ and the same is true for $\nu_{\vn,1}$ and $\nu$.
\end{proof}

\begin{remark}
    From the degree calculation of $\nu$, we can see that the stratum $\cM_{G,J}$ is usually not the scheme-theoretic fiber of $b_J$ over $\PP_J\cong (L_J)_{\ab}$. Indeed, one can show, using the same argument of Proposition \ref{p:nu finite flat}, that the map $\nu_J: \cM_{G,J}\cong \cW_J\times \frc_J^*\to (L_J)_{\ab}\times \cA_G$ has degree $|Z(L_J^\der)|\times |W/W_J|$. Therefore, the scheme-theoretic preimage $b^{-1}(\PP_J)$ has multiplicity
    \begin{equation*}
        \frac{|W_J|}{|Z(L_J^\der)|}.
    \end{equation*}
\end{remark}

%%%%%%%%%%%%%%%%%%%%%%%%%%%%%%%%%%%%%%%%%%%%%%%%%

\section{Mirror symmetry for $\cM_{G}$: statement and motivation}

In this section we formulate the main result of the paper, a homological mirror symmetry between $\cM_G$ and the group-theoretic regular centralizer for the Langlands dual $\dG$. We explain how our result fits into the ramified geometric Langlands program.

\subsection{Group-theoretic regular 
centralizer}\label{ss:JG}

In this section we make general discussions about the regular centralizer group scheme for a connected reductive group $H$ over $\CC$. At times we will use a choice of a maximal torus $T\subset H$, a pair of opposite Borel subgroups $B,B^-$ containing $T$ with unipotent radical $N, N^-$. 

\sss{Stacky Steinberg base}\label{sss:stacky Steinberg}

By the {\em Steinberg base} of $H$ we mean the GIT quotient $H\sslash H$. If $T\subset H$ is a maximal torus with Weyl group $W$, then inclusion $T\incl H$ induces an isomorphism $T\sslash W\isom H\sslash H$.

We now introduce a stacky version of $H\sslash H$. Let $\nu: H_1\to H$ be any central isogeny such that $H_1^\der$ is simply-connected. 
The kernel $\ker(\nu)$ is a finite central subgroup of $H_1$ that acts on $H_1\sslash H_1$ by multiplication. Define $S_H$ to be the Deligne-Mumford stack
\begin{equation*}
    \cS_H:=(H_1\sslash H_1)/\ker(\nu).
\end{equation*}

\begin{lemma}\label{l:SH can}
    The stack $\cS_H$ is independent of the choice of the central isogeny $\nu: H_1\to H$ (where $H_1^\der$ is simply-connected) up to a canonical isomorphism.
\end{lemma}
\begin{proof}
    Any two such $H_1$ and $H_2$ can be both covered by a group of the form $H^\sc\times A$ where $A$ is a torus. Therefore it suffices to compare the constructions in the case $H_1=H^\sc\times A$ and $H_1$ and $H_2$ are related by a central isogeny $\pi: H_1\to H_2$. It then suffices to show that
    \begin{equation*}
        s_\pi: H_1\sslash H_1\cong H^\sc\sslash H^\sc\times A\to H_2\sslash H_2
    \end{equation*}
    is a $\G:=\ker(\pi)$-torsor. For this we may rename $H_2$ by $H$ hence assuming $H^\der$ is simply-connected. 

    Since $H^\der$ is simply-connected, it has fundamental representations with characters $\{\chi_i\}_{i\in I}$. These representations can be extended to $H$; we choose extensions of each $\chi_i$ and denote its character by $\wt\chi_i$ (a class function on $H$). Then we have an isomorphism
    \begin{equation}\label{SH sc}
        H\sslash H\isom \AA^I\times H^{\ab}
    \end{equation}
    where the first projection is given by $(\wt\chi_i)_{i\in I}$. Precomposing the above isomorphism with $s_\pi$ we get
    \begin{equation}\label{s'pi}
        s'_\pi: H^\sc\sslash H^\sc\times A\to \AA^I\times H^{\ab}
    \end{equation}
    sending $(h,a)$ to $(\{\chi_i(h)\wt\chi_i(a)\}_{i\in I}, \ov\pi(a))$. Here $\ov\pi: A\to H^{\ab}$ is the isogeny of tori whose kernel is identified with $\ker(\pi)=\G$. Note that $(\chi_i)_{i\in I}$ gives an isomorphism $H^\sc\sslash H^\sc\isom \AA^I$, therefore the map $s'_\pi$ is a $\ker(\ov\pi)=\G$-torsor. This shows that $s_\pi$ is also a $\G$-torsor.
\end{proof}

\begin{exam}
    For $H=\PGL_n$, we may take $H_1=\SL_n$, and get that 
    \begin{equation*}
        \cS_{\PGL_n}\cong \AA^{n-1}/\mu_n
    \end{equation*}
    as a quotient stack, where $\z\in \mu_n$ acts on the $n-1$ coordinates of $\AA^{n-1}$ by multiplication by $\z,\z^2,\cdots, \z^{n-1}$.
\end{exam}

By definition there is a canonical morphism $s_H: H\to \cS_H$ that factors through
\begin{equation}\label{canon to St base}
    \ov s_H:\frac{H}{H^{\ad}}\to \cS_H.
\end{equation}

Let $\wt H$ be the Grothendieck alteration of $H$: it is the moduli space of pairs $(h,B)$ where $h\in H$ and $B$ is a Borel subgroup of $H$ that contains $h$. The $\pi: \wt H\to H$ be the natural map. Let $\bT$ be the abstract Cartan of $H$, we have a map $\wt s_H:\wt H\to T$ sending $(h,B)$ to the image of $h$ under the canonical map $B\to \bT$.

\begin{lemma}\label{l:Cart wt H reg}
    The following diagram is commutative
    \begin{equation}\label{wt H comm}
        \xymatrix{\wt{H} \ar[d]^{\pi}\ar[r]^{\wt s_H} & \bT\ar[d]\\
        H\ar[r]^{s_H} & \cS_H}
    \end{equation}
    Restricting to the regular locus we get a Cartesian diagram \footnote{It would not be Cartesian if we replace $\cS_{H}$ by the naive base $H\sslash H$.}
    \begin{equation}\label{Hreg Cart}
        \xymatrix{\wt{H^\reg} \ar[d]^{\pi}\ar[r]^{\wt s_H} & \bT\ar[d]\\
        H^{\reg}\ar[r]^{s_H} & \cS_H}
    \end{equation}
\end{lemma}
\begin{proof}
    Let $\nu: H_1=H^\sc\times (ZH)^\c\to H$ be the central isogeny. It suffices to prove the statements for some $H_1$ since the diagrams for $H$ are obtained from the corresponding diagrams of $H_1$ by quotienting by $\ker(\nu)$ on each term. 

    The case $H_1=H^\sc\times (ZH)^\c$ then easily reduces to the case of $H^\sc$. Therefore we reduce to the case where $H$ itself is simply-connected.  In this case, $\cS_H=H\sslash H$, and \eqref{wt H comm} is clearly commutative. Now consider the map induced by the diagram \eqref{Hreg Cart}
    \begin{equation*}
        \wt\pi: \wt{H^\reg}\to H^\reg\times_{\cS_H}\bT.
    \end{equation*}
    By \cite[Theorem 1.2, Theorem 1.5]{St}, the map $s_H|_{H^\reg}: H^\reg\to \cS_H$ is smooth, hence $H^\reg\times_{\cS_H}\bT$ is also smooth, and in particular normal. Since $\pi|_{\wt H^\reg}$ and $\bT\to \cS_H$ are finite, $\wt\pi$ is also finite. Moreover, $\wt\pi$ is birational since it is an isomorphism when restricted to the regular semisimple locus. Therefore $\wt\pi$ has to be an isomorphism. 
\end{proof}

\sss{Steinberg cross-section}\label{sss:St section}

Recall a Coxeter element in the Weyl group $W$ is a product $w=s_{i_1}\cdots s_{i_r}$, where $s_{i_1}, s_{i_2},\cdots, s_{i_r}$ is any ordering of the set of simple reflections in $W$.

Let $w$ be a Coxeter element in $W$ and  $\dot w\in N_H(T)$ be a lifting of $w$.
Let
\begin{equation*}
    S_{\dot w}=N\dot w\cap \dot w N^{-}
\end{equation*} 
be the {\em Steinberg cross-section} defined in \cite[8.9]{St}. Then $S_{\dot w}$ is an affine space of dimension $r$  consisting of regular elements of $H$. When $H$ is semisimple, $S_{\dot w}$ meets every regular conjugacy classes finitely many times.

When $H$ is simply-connected, the map $S_{\dot w}\to \cS_H=H\sslash H$ is an isomorphism. 

For semisimple $H$, the map $S_{\dot w}\to \cS_H$ is a $\pi_1(H)$-torsor.

For a reductive $H$, the Steinberg cross-section $S_{\dot w}$ lies in a single fiber of $H\to H^\ab$. To cover all regular conjugacy classes of $H$, we consider the map
\begin{equation}\label{St section gen}
    \mu: S_{\dot w}\times (ZH)^\c\to \cS_H
\end{equation}
given by multiplication followed by projection to $\cS_H$. Since $S_{\dot w}\times (ZH)^\c$ maps isomorphically to $\cS_{H_1}$ for $H_1=H^\sc\times(ZH)^\c$, we conclude that $\mu$ is a $\ker(H_1\to H)$-torsor.

%%%%%%%%%%

\sss{Group-theoretic regular centralizer: the case where $H^\der$ is simply-connected}\label{sss:J sc}

For a connected reductive group $H$ over $\CC$, let $H^\reg\subset H$ be the open subset of regular elements. Let $I^{H}_{H^\reg}=\{(x,y\in H\times H^\reg|xyx^{-1}=y\}$ be the centralizer group scheme over $H^\reg$. Then $I^{H}_{H^\reg}$ is a group scheme over $H^\reg$. 

\begin{lemma}\label{l:group J sc}
    If $H^\der$ is simply-connected, then $I^{H}_{H^\reg}$ is a smooth and commutative group scheme over $H^\reg$. The group scheme 
    $I^{H}_{H^\reg}$ has a canonical descent $J^H_H$ over $\cS_H$. In particular, $J^H_H$ is a smooth and commutative group scheme over $\cS_H$.
\end{lemma}
\begin{proof}
    We first show that $I^H_{H^\reg}$ is smooth over $H^\reg$. By \cite[Theorem 1.2, Theorem 1.5]{St} and \eqref{SH sc}, the map $s_H|_{H^\reg}: H^\reg\to \cS_H$ is smooth, and each $k$-fiber of $s_H|_{H^\reg}$ is a single (regular) conjugacy class in $H$. This implies that the map $a: H\times H^{\reg}\to H^\reg\times_{\cS_H}H^\reg$ given by $a(y,x)=(yxy^{-1},x)$ is smooth. Base change $a$ along the diagonal map $\D_H: H^\reg\incl H^\reg\times_{S_H}H^\reg$ we obtain the projection $I^{H}_{H^\reg}\to H^\reg$, which is therefore smooth.
    
    We next show the commutativity of $I^{H}_{H^\reg}$. Over the regular semisimple locus $H^\rs\subset H^\reg$, $I^H_{H^\reg}|_{H^\rs}$ is a torus over $H^\rs$, hence commutative. Since $H^\rs$ is dense in $H^\reg$ and $I^H_{H^\reg}$ is smooth over $H^\reg$, $I^H_{H^\reg}|_{H^\rs}$ is dense in $I^H_{H^\reg}$. This implies that $I^H_{H^\reg}$ is commutative.

    The descent statement is proved by the same argument of \cite[Lemme 2.1.1]{NgoFL}.
\end{proof}

If $H$ is simply-connected, $S_{\dot w}$ is a Steinberg cross-section for $H$ as in \S\ref{sss:St section}, then after identifying $\cS_H$ with $S_{\dot w}$, $J^H_H$ is the centralizer group scheme of $S_{\dot w}$ under $H$.

More generally, still assuming $H^\der$ is simply-connected, for any reductive group $H_\flat$ that is a central extension of $H^\ad$, $H_\flat$ acts on $H$ by conjugation (via its map to $H^\ad$), hence $I^{H_\flat}_{H^\reg}$ is defined. The obvious analog of Lemma \ref{l:group J sc} holds, and we have a group scheme
\begin{equation*}
    J^{H_\flat}_H\to \cS_H=H\sslash H.
\end{equation*}

\sss{Group-theoretic regular centralizer: general case}\label{sss:reg cent gen}
For a general connected reductive group $H$ over $\CC$, let $\nu:H_1\to H$ be a central isogeny with $H_1^\der$ simply-connected, then we define
\begin{equation*}
    J^H_H:=J^H_{H_1}/\ker(\nu)
\end{equation*}
as a group scheme over $\cS_H=(H_1\sslash H_1)/\ker(\nu)$.  Here the $\ker(\nu)$-action on $J^H_{H_1}$ is the descent of its action on $I^{H}_{H_1}$, where $z\in \ker(\nu)$ sends  $(x,y)\in H\times H_1^\reg$ to $(x,yz)$. 

The same argument as Lemma \ref{l:SH can} shows:
\begin{lemma}
    The group scheme $J^H_H$ over $\cS_H$ is canonically independent of the choice of the central isogeny $\nu: H_1\to H$ with $H_1^\der$ simply-connected. 
\end{lemma}

More generally, if $H_\flat$ is any central extension of $H^\ad$, then the group scheme
\begin{equation*}
    J^{H_\flat}_H:=J^{H_\flat}_{H_1}/\ker(\nu) \to \cS_H
\end{equation*}
is defined, and is smooth and commutative.

\subsection{Statement of main result}

The main result of this paper is the following.

\begin{theorem}\label{th:mirror}
    Let $G$ be a reductive group. Let $\cM^{\c}_{G}$ be the neutral component of $\cM_G$. Then for any compatible system of Kostant sections (cf. Definition \ref{defncompatibleKostantsec}), which always exists, there is an equivalence of dg categories
\begin{equation}\label{main equiv}
    \cW(\cM^{\c}_{G})\cong \Coh(J^{G^{\vee, \ad}}_{G^\vee}).
\end{equation}
Here $\cW(\cM^{\c}_{G})$ denotes the wrapped Fukaya category of $\cM^{\c}_{G}$, to be made precise in Section \ref{s:Fukaya}.

Moreover, if $G$ of adjoint type, then the equivalence \eqref{main equiv} sends the unique Kostant section $S_I$ in $\cT_G$ to the structure sheaf of $J_{G^\vee}^{G^{\vee, \ad}}$. 
\end{theorem}

See Proposition \ref{propMirrorconnectedZ} below about the images of Kostant sections in $\cM_G^\circ$ under the mirror equivalence for general reductive $G$.

\sss{Changing compatible systems of Kostant sections}\label{sss:change comp sys Kos sec}
The equivalence \eqref{main equiv} depends on the choice of a compatible system of Kostant sections, one for each pseudo-Levi subgroup $L_J$ of $G$. It can be calculated  that the isomorphism classes of such compatible systems form a torsor under $\cohog{1}{\triangle, \pi\frZ}\cong \coker(\xcoch(T)\to \ZZ^{\wt I})$ (see Lemma \ref{l:cohoZ}; the map $\xcoch(T)\to \ZZ^{\wt I}$ is given by pairing $\l\in \xcoch(T)$ with affine simple roots). Correspondingly, there should be an action of $\coker(\xcoch(T)\to \ZZ^{\wt I})$ on $\Coh(J^{G^{\vee, \ad}}_{G^\vee})$. We describe now such an action in the case $G$ is almost simple though we do not check that it corresponds to changing compatible systems of Kostant sections on the Fukaya side.

For $G$ almost simple, we have a canonical decomposition
\begin{equation}\label{decomp H1 piZ}
    \coker(\xcoch(T)\to \ZZ^{\wt I})\cong \ZZ^I/\xcoch(T)\op \ZZ
\end{equation}
where the copy of $\ZZ$ is the image of $\ZZ\a_0$. Note that $\ZZ^I/\xcoch(T)$ can be identified with $\xch(\pi_1(\dG))$. By definition, $J^{G^{\vee, \ad}}_{\dG}=J^{G^{\vee, \ad}}_{G^{\vee,\sc}}/\pi_1(\dG)$. Therefore $\chi\in \xch(\pi_1(\dG))$ acts on $\Coh(J^{G^{\vee, \ad}}_{\dG})\cong \Coh^{\pi_1(\dG)}(J^{G^{\vee, \ad}}_{G^{\vee,\sc}})$ by tensoring coherent sheaves with the $1$-dimensional representation $\chi$ of $\pi_1(\dG)$. 

On the other hand, the action of the other direct summand $\ZZ$ in \eqref{decomp H1 piZ} on $Coh(J^{G^{\vee, \ad}}_{\dG})$ should be induced by the ``Dehn twist" automorphism $\d: J^{G^{\vee, \ad}}_{\dG}\to J^{G^{\vee, \ad}}_{\dG}$ defined as follows. We may assume $\dG$ is simply-connected and the general case is obtained by passing to the quotient by $\pi_1(\dG)$. For $\dG$ simply-connected, we choose a Steinberg section $S_{\dot w}$ and identify $J^{G^{\vee, \ad}}_{\dG}$ with the space of pairs $(g,x)$ where $g\in S_{\dot w}$ and $x\in C_{G^{\vee,\ad}}(g)$. Then $\d(g,x)=(g,\ov{g}x)$, with $\ov{g}$ the image of $g$ in $G^{\vee,\ad}$. The name ``Dehn twist" will be justified in Remark \ref{r:Dehn}.

\subsection{Motivation from geometric Langlands}\label{ss:GL}

As we mentioned in \S\ref{ss:auto}, the definition of $\cM_G$ was motivated by automorphic representations. In fact the equivalence \eqref{main equiv} can be thought of as an instance of categorical geometric Langlands correspondence with wild ramifications. We will now make this more precise by relating the both sides of \eqref{main equiv} to certain automorphic and spectral categories respectively.

\sss{Fukaya category as an automorphic category} 

In the proof of Proposition \ref{p:count auto forms}, we argued that the space $\cC(\psi_0,\psi_\infty)$ of functions on $\Bun_G(\bI^{++}_0, \bI^{++}_\infty)$ that are eigen under $V_0\times V_\infty$ with character $\psi_q\psi_0\times\psi_q\psi_\infty$ is the space of automorphic forms on $G(F)\bs G(\AA_F)$ (where $F=\FF_q(t)$) satisfying specific local conditions. At the categorical level, working with $\Qlbar$-sheaves over a base field $k$ of characteristic $p$, this naturally leads to considering the category
\begin{equation*}
    \cD(\psi_0,\psi_\infty):=D_{(V_0\times V_\infty, \AS_{\psi_0}\bt\AS_{\psi_\infty})}(\Bun_G(\bI^{++}_0, \bI^{++}_\infty), \Qlbar).
\end{equation*}
Here, $\AS_{\psi_0}$ is the pullback of a fixed Artin-Schreier local system on $\Ga$ (which amounts to choosing a nontrivial character $\FF_p\to \Qlbar^\times$) via the map $\psi_0: V_0\to \Ga$. Similarly for $\AS_{\psi_\infty}$. Then $\AS_{\psi_0}\bt\AS_{\psi_\infty}$ is a character sheaf on $V_0\times V_\infty$, and $\cD(\psi_0,\psi_\infty)$ above consists of complexes of sheaves on $\Bun_G(\bI^{++}_0, \bI^{++}_\infty)$ that are equivariant under $V_0\times V_\infty$ with respect to the character sheaf $\AS_{\psi_0}\bt\AS_{\psi_\infty}$.

By the general philosophy of microlocal sheaf theory, when working with the base field $k=\CC$ instead, the category of sheaves on $\Bun_G(\bI^{++}_0, \bI^{++}_\infty)$ may be identified with the wrapped Fukaya category of its cotangent bundle. The equivariance condition under $(V_0\times V_\infty, \AS_{\psi_0}\bt\AS_{\psi_\infty})$ now corresponds to replacing $T^*\Bun_G(\bI^{++}_0, \bI^{++}_\infty)$ by its Hamiltonian reduction by $V_0\times V_\infty$ at the moment value $(\psi_0,\psi_\infty)$, which recovers $\cM_G$ by Remark \ref{r:MG as Ham red}. This leads to considering $\cW(\cM_G)$ as a category of automorphic nature.

\sss{A wild character variety}

Having related the left side of \eqref{main equiv} to a category of automorphic sheaves, we now relate the space $J^{\dGad}_{\dG}$ that appears on the right side of \eqref{main equiv} to a moduli space of $\dG$-local systems on $\PP^1-\{0,\infty\}$ with wild ramifications around $0$ and $\infty$ matching the level structures imposed on the automorphic side.

Around $0$, the level structure $(\bI^+_0,\psi_0)$ is a special case of the kind considered in \cite{BBAMY} coming from a homogeneous element in the loop Lie algebra. In \cite[\S4]{BBAMY}, a local Betti moduli space $\cC_{\dG,0}$ classifying $\dG$-Stokes data matching the level structure  $(\bI^+_0,\psi_0)$ was constructed, following Minh-Tam Trinh \cite{Trinh}. Let $\SS^1_0$ be the tangent circle of $\PP^1(\CC)$ (as a real manifold) at $0$. Let $\Sig_0$ be the set of Stokes sectors (as open arcs in $\SS^1_0$) for $\psi_0$, and let $R_0\subset \SS^1_0$ be the finite set of Stokes rays.  Each Stokes ray $\rho\in R_0$ determine a simple reflection $s_\rho\in I$, and the map $\rho\mt s_\rho$ is a bijection $R_0\isom I$. Then $\cC_{\dG,0}$ is the moduli stack of $(\cE; (\cF_\s)_{\s\in \Sig_0})$, where $\cE$ is a $\dG$-torsor on $\SS^1_0$ and each $\cF_\s$ is a $\dB$-reduction of $\cE$ such that the relative position between $\cF_\s$ and $\cF_{\s'}$ (both are viewed as $\dB$-reductions of the same $\dG$-torsor $\cE|_{\rho}$ that can be trivialized) for two sectors $\s,\s'$ separated by the ray $\rho\in R_0$ is $s_\rho$. We have a canonical map 
\begin{equation*}
    \om_0: \cC_{\dG,0}\to \frac{\dG}{\dG}
\end{equation*}
recording the $\dG$-torsor $\cE$ on the circle $\SS^1_0$.

Following \cite[\S4.3.3]{BBAMY}, upon choosing a Stokes sector $\s_0\in \Sig_0$, $\psi_0$ determines a Coxeter element $c_0=\prod_{\rho\in R_0}s_\rho$ (product over the Stokes rays in counterclosewise order,  starting with the one next to $\s_0$ in the counterclockwise direction). For $w\in W$, let $\wt\cC_{\dG}(w)$ be the moduli space of pairs $(g,\b)$ where $g\in \dG$ and $\b$ is a Borel subgroup of $\dG$ such that the relative position of the pair $(\b,\Ad(g)\b)$ is $w$. With the choice of $\s_0$, we have an isomorphism
\begin{equation*}
    i_{\s_0}: \cC_{\dG,0}\cong\wt\cC_{\dG}(c_0)/\dG
\end{equation*} 
where $\dG$ acts on the pair $(g,\b)$ by simultaneous conjugation. If we fix a choice of a Borel subgroup $\dB$ as $\b$, we get an isomorphism $\wt\cC_{\dG}(c_0)/\dG\cong (\dB c_0 \dB)/\dB$ (conjugation action). Changing the sector $\s_0$ will change $c_0$ by a cyclic shift of words.

Similarly, at $\infty$, $\psi_\infty$ together with the choice of a Stokes sector $\s_\infty$ determines another Coxeter element $c_\infty$. The local Betti moduli space $\cC_{\dG,\infty}$ at $\infty$ is isomorphic to $\wt\cC_{\dG}(c_\infty)/\dG$, which is in turn isomorphic to $(\dB c_\infty \dB)/\dB$.

The global moduli space $\cC_{\dG}$ that serves as the Betti version of the spectral counterpart of $\cM_{G}$ is then be glued from the two local Betti moduli spaces at $0$ and $\infty$ by imposing that their topological monodromies are inverse to each other. In other words $\cC_{\dG}$ is defined by the Cartesian diagram
\begin{equation*}
    \xymatrix{\cC_{\dG}\ar[r]\ar[d] & \cC_{\dG,0}\ar[d]^{\om_0}\\
    \cC_{\dG,\infty}\ar[r]^{\inv\c\om_\infty} & \frac{\dG}{\dG}}
\end{equation*}
Here $\om_0$ and $\om_\infty$ denote the maps recording topological monodromy around $0$ and $\infty$, and for the bottom horizontal arrow we compose $\om_\infty$ with the inversion map on $\dG$.

To write $\cC_{\dG}$ more explicitly, we choose an auxiliary datum. A {\em longitude} for $(\PP^1; 0,\infty)$ is a smooth curve connecting $0$ and $\infty$ on $\PP^1(\CC)$ that lands in the interior of a Stokes sector for $\psi_0$ near $0$, and lands in the interior of a Stokes sector for $\psi_\infty$ near $0$. A longitude $\ell$ then determines a Stokes sector $\s_0$ at $0$ and $\s_\infty$ at $\infty$, hence Coxeter elements $c_0$ and $c_\infty$. Using $i_{\s_0}$ and its analog at $\infty$, we get an isomorphism 
\begin{equation*}
    i_{\s_0,\s_\infty}: \cC_{\dG}\cong \frac{\wt\cC_{\dG}(c_0)}{\dG}\times_{\frac{\dG}{\dG}} \frac{\wt\cC_{\dG}(c_\infty^{-1})}{\dG}.
\end{equation*}
The maps $\wt\cC_{\dG}(c_0)\to \dG$ and $\wt\cC_{\dG}(c_\infty^{-1})\to \dG$ are both given by $(g,\b)\mt g$. Here we have used the isomorphism $\wt\cC_{\dG}(c_\infty)\cong \wt\cC_{\dG}(c^{-1}_\infty)$ given by $(g,\b)\mt (g^{-1}, \Ad(g)\b)$ to rewrite the factor at $\infty$. 

Let
\begin{equation*}
    \wt\cC_{\dG}(c_0,c^{-1}_\infty):=\wt\cC_{\dG}(c_0)\times_{\dG}\wt\cC_{\dG}(c_\infty^{-1}).
\end{equation*}
This is a generalization of the Steinberg variety of triples: it classifies triples $(g,\b_0, \b_\infty)$ where $g\in \dG$, $\b_0,\b_\infty$ are Borel subgroups of $\dG$, such that the relative position of $(\b_0, \Ad(g)\b_0)$ is $c_0$ and the relative position of $(\b_\infty, \Ad(g)\b_\infty)$ is $c_\infty^{-1}$ for $i=2$. The choice of the longitude $\ell$ allows us to identify the horizontal sections of a $\dG$-connection on the sector $\s_0$ and the sector $\s_\infty$ by parallel transport. Therefore we get an isomorphism
\begin{equation*}
    i_\ell: \cC_{\dG}\cong\wt\cC_{\dG}(c_0,c^{-1}_\infty)/\dG
\end{equation*}
that depends only on the starting sector $\s_0$, the ending sector $\s_\infty$ and the homotopy class of the longitude $\ell$.

The microlocal version of the ramified geometric Langlands correspondence predicts at least a full embedding 
\begin{equation}\label{rGLC}
    \cW(\cM_G)\stackrel{?}{\incl} \Coh(\cC_{\dG}).
\end{equation}

\sss{Relation with $J^{\dGad}_{\dG}$}
We next explain the relation between $\cC_{\dG}$ and the regular centralizer $J_{\dG}$. Note that $\cC_{\dG,0}, \cC_{\dG,\infty}$ and $\cC_{\dG}$ have actions by $\BB(ZG^\vee)$, i.e., every point of these stacks have automorphisms containing $Z\dG$. Dividing by the $\BB(ZG^\vee)$-actions, we obtain stacks $\cC'_{\dG,0}, \cC'_{\dG,\infty}$ and $\cC'_{\dG}$ with a canonical isomorphism 
\begin{equation*}
    \cC'_{\dG}\cong \cC'_{\dG,0}\times_{\frac{\dG}{\dGad}}\cC'_{\dG,\infty}.
\end{equation*}
Upon choosing a longitude $\ell$ connecting $0$ and $\infty$, we get a canonical isomorphism
\begin{equation}\label{i' C'}
    i'_\ell: \cC'_{\dG}\cong \wt\cC_{\dG}(c_0,c^{-1}_\infty)/\dGad.
\end{equation}

Compatible with the expected equivalence \eqref{rGLC}, we also expect an equivalence
\begin{equation*}
    \cW(\cM^\c_G)\stackrel{?}{\incl} \Coh(\cC'_{\dG}).
\end{equation*}
We now explain that the right side above is equivalence to the right side of the mirror theorem \ref{th:mirror}.

\begin{lemma}\label{l:BwB S}
    \begin{enumerate}
        \item Let $c\in W$ be a Coxeter element, then the forgetful map $\om: \wt\cC_{\dG}(c)\to \dG$ lands in the open subset of regular elements. 
        \item The composition
    \begin{equation*}
        \frac{\wt\cC_{\dG}(c)}{\dGad}\xr{\om} \frac{\dG}{\dGad}\xr{\eqref{canon to St base}} \cS_{\dG}
    \end{equation*}
    is an isomorphism of stacks.
        \item The composition 
        \begin{equation}\label{CG to SG}
            \cC'_{\dG,0}\xr{\om_0}\frac{\dG}{\dGad}\xr{\eqref{canon to St base}} \cS_{\dG}
        \end{equation}
        is an isomorphism. The same is true when $\cC_{\dG,0}$ is replaced by $\cC_{\dG,\infty}$.
    \end{enumerate}
\end{lemma}
\begin{proof}
    (1) (2) Write $H=\dG$ and $B_H=\dB$. Since $\frac{\wt\cC_{\dG}(c)}{\dGad}\cong \frac{B_H c B_H}{B_H^{\ad}}$, it suffices to show that the canonical map $\frac{B_H c B_H}{B_H^{\ad}}\to \cS_H$ is an isomorphism. 
    
    Let $\nu: H_1=H^\sc\times (ZH)^\c\to H$ be the central isogeny, and $B_{H_1}=\nu^{-1}(B_H)$. We have a commutative diagram
    \begin{equation*}
        \xymatrix{\frac{B_{H_1}cB_{H_1}}{B^\ad_H}\ar[d]\ar[r] & \cS_{H_1}\ar[d]\\
        \frac{B_{H}cB_{H}}{B^\ad_H}\ar[r] & \cS_{H}}
    \end{equation*}
    Since both vertical  maps are $\ker(\nu)$-torsors, the above diagram is Cartesian. Thus it suffices to show the statements for $H_1$ instead of $H$.
    
    Let $S_{\dot c}=N_H \dot c\cap \dot c N^-_H$ denote the Steinberg section of $H^\sc$ (where $N_H$ is the unipotent radical of $B_H$, and $\dot c$ lifts $c$). Then the inclusion $S_{\dot c}\subset N_H\dot c N_H$ induces an isomorphism 
    \begin{equation*}
        S_{\dot c}\isom \frac{N_H \dot c N_H}{N_H}.
    \end{equation*}
    Thus we have
    \begin{equation}\label{Sw ZH}
        \frac{B_{H_1}cB_{H_1}}{B^\ad_H}\cong \frac{N_{H}\dot cN_{H}\times (ZH)^\c}{N_H}\cong S_{\dot c}\times (ZH)^\c.
    \end{equation}
    In particular, elements in $B_{H_1}cB_{H_1}$ are conjugate to elements in $S_{\dot c}\times (ZH)^\c$, which are all regular. Moreover, the projection to $\cS_{H_1}$ induces an isomorphism
    \begin{equation*}
        \frac{B_{H_1}cB_{H_1}}{B^\ad_H}\cong S_{\dot c}\times (ZH)^\c\isom \cS_{H^\sc}\times (ZH)^\c=\cS_{H_1}
    \end{equation*}
    proving the second statement.

    (3) To check \eqref{CG to SG} is an isomorphism, we may choose a Stokes sector to identify $\cC'_{\dG,0}$ with $\wt\cG_{\dG}(c_0)/\dGad$, and the statement follows from (2).
\end{proof}

\begin{cor}
\begin{enumerate}
    \item There is a canonical isomorphism over $\cS_G$
    \begin{equation}\label{J to CC0}
        J^{\dGad}_{\dG}\isom \cC'_{\dG,0}\times_{\frac{\dG}{\dGad}}\cC'_{\dG,0}
    \end{equation}
    such that the unit section of $J^{\dGad}_{\dG}$ corresponds to the diagonal $\cC'_{\dG,0}$ on the right side, under the isomorphism $\cC'_{\dG,0}\isom \cS_{\dG}$ in Lemma \ref{l:BwB S}.
    \item The stack $\cC'_{\dG}$ is canonically a $J^{\dGad}_{\dG}$-bitorsor over $\cS_{\dG}$. 
\end{enumerate}
\end{cor}
\begin{proof}
(1) The map $\pi^\reg: G^{\vee,\reg}/\dGad\to \cS_{\dG}$ is a $J^{\dGad}_{\dG}$-gerb. By Lemma \ref{l:BwB S}, the map $\om^\reg_0: \cC'_{\dG,0}\to G^{\vee,\reg}/\dGad$ can be identified with a section to $\pi^\reg$, hence $\om^\reg_0$ is a $J^{\dGad}_{\dG}$-torsor. The isomorphism \eqref{J to CC0} follows.

(2) The two maps $\om^\reg_0: \cC'_{\dG,0}\to G^{\vee,\reg}/\dGad$ and $\om^\reg_\infty: \cC'_{\dG,\infty}\to G^{\vee,\reg}/\dGad$ are both $J^{\dGad}_{\dG}$-torsors. Therefore $\cC'_{\dG,0}\times_{\frac{\dG}{\dGad}}\cC'_{\dG,\infty}$ is a left $\cC'_{\dG,0}\times_{\frac{\dG}{\dGad}}\cC'_{\dG,0}\cong J^{\dGad}_{\dG}$-torsor and a right $\cC'_{\dG,\infty}\times_{\frac{\dG}{\dGad}}\cC'_{\dG,\infty}\cong J^{\dGad}_{\dG}$-torsor.
\end{proof}

\begin{cor}\label{c:C' J} There exists a longitude $\ell$ such that the Coxeter elements $c_0$ and $c_\infty$ determined by its sectors at $0$ and $\infty$ satisfy $c_0=c_\infty^{-1}$. Any such longitude determines a section of $\cC'_{\dG}\to \cS_{\dG}$, hence a trivialization of $\cC'_{\dG}$ as a $J^{\dGad}_{\dG}$-bitorsor. 
\end{cor}
\begin{proof}
    Choose any Stokes sectors $\s_0$ and $s_\infty$ inside $\SS^1_0$ and $\SS^1_\infty$, which determine Coxeter elements $c_0$ and $c_\infty$. Both $c_0$ and $c_\infty^{-1}$ are minimal length elements in the Coxeter conjugacy class of $W$, hence they differ by a cyclic shift by \cite[Corollary 4.4]{HeNie}. In other words, upon changing $\s_\infty$, we may arrange that $c_0=c_\infty^{-1}$. Take any longitude $\ell$ that gives such a pair $(c_0,c_\infty)$, $i'_\ell$ gives an isomorphism
    \begin{equation*}
        \cC'_{\dG}\cong \wt\cC_{\dG}(c_0,c_0)/\dGad\cong \cC'_{\dG,0}\times_{\frac{\dG}{\dGad}}\cC'_{\dG,0}.
    \end{equation*}
    The diagonal copy of  $\cC'_{\dG,0}$ gives the desired section to $\cC'_{\dG}\to \cS_{\dG}\cong \cC'_{\dG,0}$.
\end{proof}

\begin{remark}[Dehn twist]\label{r:Dehn} Let $\ell_1$ and $\ell_2$ be two longitudes determining the same sectors $(\s_0,\s_\infty)$, hence the same pair of Coxeter elements $(c_0,c_\infty)$. Assume $\ell_1$ and $\ell_2$ differ by a Dehn twist of the cylinder $\PP^1(\CC)-\{0,\infty\}$. The isomorphisms $i'_{\ell_1}$ and $i'_{\ell_2}$ in \eqref{i' C'} then differ by the automorphism of $\wt\cC_{\dG}(c_0,c_\infty)$ that sends $(g,\b_1,\b_2)$ to $(g,\b_1, \Ad(g)\b_2)$. In particular, suppose $c_\infty^{-1}=c_0$, then the two trivializations of $\cC'_{\dG}$ as a left $J_{\dG}^{\dGad}$-torsor differ by the automorphism $\d$ on $J_{\dG}^{\dGad}$ introduced in \S\ref{sss:change comp sys Kos sec}. This justifies calling $\d$ the ``Dehn twist".

\end{remark}

%%%%%%%%%%%%%%%%%%%%%%%%%%%%%%%%%%%%%%%%%%%%%%%%%%%
\section{Point-counting, cohomology and P=W for $\cM_G$}\label{s:counting}
%%%%%%%%%%%%%%%%%%%%%%%%%%%%%%%%%%%%%%%%%%%%%%%%%%%

When the base field is a finite field $\FF_q$, we give an explicit formula for the number of $\FF_q$-points on $\cM_G$ in \S\ref{ss:MG Fq} and relate that to the counting of automorphic forms with specific local conditions in \S\ref{ss:auto}. 

Over any algebraically closed base field $k$, we compute the cohomology of $\cM^\c_G$ and express it in terms of root-theoretic data in \S\ref{ss:coho MG}. The result shows strong similarity with point-counting on the group regular centralizer $J^{G^\ad}_{G}$, which was done by Lusztig. To explain the similarity, in \S\ref{ss:P=W}, we speculate the relationship between $\cM^\c_G$ and $J^{G^\ad}_{G}$ by formulating a conjectural non-abelian Hodge correspondence between them, and a P=W conjecture between their cohomology groups.

This section is not used in the rest of the paper.

\subsection{Points over $\FF_{q}$}\label{ss:MG Fq} 
We work with the finite base field $k=\FF_q$, and assume $G$ is almost simple and split.
We assume $q$ is large enough so that
\begin{equation}\label{large q}
    \mbox{For each $J\sft \wt I$, $\cW_J(\FF_q)\ne \vn$, and $Z(L_J)$ has $\FF_{q}$-points on every connected component.}
\end{equation}
Here $\cW_J$ is the $Z(L_J)$-torsor introduced in \S\ref{sss:geom MGJ}.

Let $\{n_{i}\}_{i\in \wt I}\in \NN^{\wt I}$ be the Dynkin labeling (i.e., then for $i\in I$, $n_{i}$ is the coefficient of $\a_{i}$ in the highest root $\th$, and $n_{0}=1$). For $d\in \NN$, let $\wt I_{d}\subset \wt I$ be the set of $i\in \wt I$ such that $d|n_{i}$, and let
\begin{equation*}
    r_d:=|\wt I_d|.
\end{equation*}
For example, $r_1=|\wt I|=r+1$. For $m\in \NN$, let $[m]_{q}=(q^{m}-1)/(q-1)$.

\begin{prop}\label{p:MG Fq} Assume $G$ is simple adjoint and that \eqref{large q} holds. 
\begin{enumerate}
    \item We have
\begin{equation}\label{MG Fq}
|\cM^{\c}_{G}(\FF_{q})|=q^{r}\sum_{d\in \NN, \wt I_{d}\ne\vn}[r_d]_{q}\ph(d).
\end{equation}
Here $\ph(-)$ is the Euler function.
    \item Let $\om \in \Om$, and assume $q$ satisfies the following variant of \eqref{large q}:
\begin{equation}\label{large q om}
    \mbox{For each $J\sft \wt I$, $\cW^\om_J(\FF_q)\ne \vn$.}
\end{equation}
Then we also have 
\begin{equation}\label{MG om Fq}
|\cM^{\om}_{G}(\FF_{q})|=q^{r}\sum_{d\in \NN, \wt I_{d}\ne\vn}[r_d]_{q}\ph(d).
\end{equation}
\end{enumerate}
\end{prop}
\begin{proof}
We give the proof of (1); the proof (2) is the same. For any $J\sft \wt I$, by Lemma \ref{l:MGJ} we have $|\cM^{\c}_{G,J}(\FF_{q})|=|\cW_J(\FF_q)|\cdot |\frc_J^*(\FF_q)|$ since $\frc_J^*$ is an affine space of dimension $r$. By \eqref{large q}, $|\cW_J(\FF_q)|=|Z(L_J)(\FF_q)|=(q-1)^{r-|J|}|\pi_{0}(Z(L_{J}))|$.  Moreover, since $\xch(ZL_J)=\xch(T)/\Span\{\a_j|j\in J\}$, and $G$ is adjoint so $\xch(T)=\Span\{\a_i|i\in I\}$, we see that the torsion part of $\xch(ZL_J)$ has cardinality equal to the gcd of $\{n_{i}\}_{i\in \wt I-J}$ (when $J\subset I$, $\xch(ZL_J)$ is torsion-free, and $\gcd\{n_{i}\}_{i\in \wt I-J}=1$ since $n_0=1$). Since $\pi_0(ZL_{J})$ is dual to $\xch(ZL_J)_{\tors}$, we see that $|\pi_{0}(ZL_{J})|=\gcd\{n_{i}\}_{i\notin J}$. Therefore
\begin{equation*}
|\cM^{\c}_{G,J}(\FF_{q})|=q^{r}(q-1)^{r-|J|}\gcd\{n_{i}\}_{i\in \wt I-J}.
\end{equation*}
Sum over all $J\sft \wt I$, or equivalently sum over non-empty subsets $K=\wt I -J\subset \wt I$, we get
\begin{equation*}
|\cM^{\c}_{G}(\FF_{q})|=q^{r}\sum_{K\subset \wt I, K\ne\vn}(q-1)^{|K|-1}g_{K}\end{equation*}
where $g_{K}:=\gcd\{n_{i}\}_{i\in K}$. The right side can be written as
\begin{eqnarray*}
&&q^{r}\sum_{d\in \NN}d\sum_{\vn\ne K\subset \wt I, g_{K}=d}(q-1)^{|K|-1}\\
&=& q^{r}\sum_{d\in \NN}\left(\sum_{\vn\ne K\subset \wt I, d|g_{K}}(q-1)^{|K|-1}\right)
\ph(d).
\end{eqnarray*}
Now the term $\sum_{\vn\ne K\subset \wt I, d|g_{K}}(q-1)^{|K|-1}$ is the same as
\begin{equation*}
\sum_{\vn\ne K\subset \wt I_{d}}(q-1)^{|K-1|}=\sum_{j=1}^{r_d}(q-1)^{j-1}\binom{r_d}{j}=\begin{cases}[r_d]_{q},  & \wt I_{d}\ne\vn; \\ 0, & \mbox{otherwise.}\end{cases}
\end{equation*}
The formula for $|\cM^\c_G(\FF_q)|$ follows.
\end{proof}

In the next subsection we will explain the point-counting formula for $\cM^\c_G$ by computing its cohomology.

\begin{remark} The same argument for \eqref{MG Fq} can be used to compute $|\cM_{G}(\FF_{q})|$ for $G$ not necessarily adjoint. For example, if $G=\SL_{n}$ and $q\equiv1\mod n$,  we get
\begin{equation*}
|\cM_{\SL_{n}}(\FF_{q})|=q^{n-1}\sum_{d|n}d\ph(d)[n/d]_{q}.
\end{equation*}
The summation over $d|n$ can be interpreted as over endoscopic groups of $\SL_n$, consistent with the endoscopic decomposition phenomenon of $\SL_n$-Hitchin fibrations.
\end{remark}

\subsection{Cohomology of $\cM^\c_G$}\label{ss:coho MG} In this subsection we work with any algebraically closed base field $k$. We assume $G$ is semisimple. We give an explicit formula for the $\Qlbar$-cohomology of $\cM^\c_G$ in this subsection.

Recall the map $b:\cM_G^\c\to \PP$, and let
\begin{equation}
    \cF=\bR b_!\Qlbar[2r](r).
\end{equation}
By the discussion in \S\ref{sss:fiber b}, the fibers of $b$ are finite disjoint unions of affine spaces of dimension $r$. Therefore $\cF$ is an ordinary sheaf on $\PP$ (concentrated in degree zero), constructible with respect to the stratification $\PP=\cup_{J\sft \wt I}\PP_J$. The restriction of $\cF_J=\cF|_{\PP_J}$ is a local system canonically isomorphic to $\mu_{J!}\Qlbar$, where $\mu_J: \cW_J\to \PP_J\cong (L_J)_{\ab}$ is defined in \eqref{mu J}. In particular, $\cF_J$ is non-canonically isomorphic to the direct image of constant sheaf along the finite projection $ZL_J\to (L_J)_{\ab}$, from which we see that $\cF_J$ is a semisimple local system.  

Let $\cF^0_J\subset \cF_J$ be the trivial direct summand (maximal trivial sub local system). 

\begin{lemma}
    There is a unique subsheaf $\cF^0\subset \cF$ whose restriction to $\PP_J$ is $\cF^0_J$.
\end{lemma}
\begin{proof}
Write $\cF_J=\cL_J\op \cF^0_J$, where $\cL_J$ does not contain any trivial subquotient.  Let $i_J:\PP_J\incl \PP$ be the embedding. We define $\cF^0$ to be the kernel of the canonical map
\begin{equation}
    \cF\to \bigoplus_{J\sft \wt I}\bR^0 i_{J*}\cF_J\to \bigoplus_{J\sft \wt I}\bR^0i_{J*}\cL_J.
\end{equation}
Taking restriction to $\PP_J$ we get
\begin{equation}\label{ker FJ}
    \cF^0|_{\PP_J}= \ker(\cF_J\to \bigoplus_{J\supset J'}(\bR^0 i_{J'*}\cL_{J'})|_{\PP_{J}})
\end{equation}
By Lemma \ref{l:toric nontriv loc sys}, the restrictions of $\bR^0 i_{J'*}\cL_{J'}$ to each boundary stratum does not contain trivial local systems as subquotients, therefore $\cF^0_J\subset \cF^0|_{\PP_J}$. On the other hand, the $J'=J$ term of the right side of \eqref{ker FJ} is the projection $\cF_J\surj \cL_J$, which implies $\cF^0|_{\PP_J}=\cF^0$.
\end{proof}

\sss{Digression: sheaves on toric varieties} We need the following general observation for toric varieties. Let $Y$ be a $T$-toric variety, with open dense $T$-orbit $j: Y_0\incl Y$. Let $\cL$ be a tame rank one local system on $Y_0$, corresponding to a character $\chi_\cL: \pi_1(Y_0)^{\tame}\to \Qlbar^\times$. For any $T$-orbit $Y_\a\subset Y$, there is a canonical map of fundamental groups $\r_\a: \pi_1(Y_0)\to \pi_1(Y_\a)$ (choosing base points identify $Y_0$ and $Y_\a$ as quotients of $T$). 
\begin{lemma}\label{l:toric nontriv loc sys}
    Let $\cL$ be a tame rank one local system on $Y_0$ as above. Then any irreducible local system appearing as a subquotient of $(\bR^i j_*\cL)|_{Y_\a}$ for some $i$ must be tame, and the corresponding character $\chi: \pi_1(Y_\a)^\tame\to \Qlbar^\times$ satisfies $\chi\c\r_\a=\chi_\cL$.

    In particular, if $\cL$ is nontrivial, then for any $T$-orbit $Y_\a$ and any $i\in\ZZ$, $(\bR^i j_*\cL)|_{Y_\a}$ does not contain the trivial local system as a subquotient.
\end{lemma}
\begin{proof}
    It is standard to reduce to the case where $\cL$ has finite monodromy. Also we may replace $T$ by a quotient so that $Y_0\cong T$. Now $\chi_\cL: \pi_1(T)^\tame\to \Qlbar^\times$ has finite image. Let $p_0: \wt T\to T$ be an isogeny from another torus such that $p_0^*\cL$ is trivial, i.e., $\chi|_{\pi_1(\wt T)^\tame}$ is trivial. Let $\wt Y$ be the $\wt T$-toric variety constructed using the same fans as $Y$. We have a finite map $p: \wt Y\to Y$ restricting to $p_0$ over $Y_0\cong T$. Let $p_\a: \wt Y_\a\to Y_\a$ be the restriction of $p$ over $Y_\a$. 
    
    Let $\G_0=\ker(p_0)$. Note that $\chi_\cL$ factors through $\G_0\cong \pi_1(Y_0)^\tame/\pi_1(\wt Y_0)^\tame$. Also, $\G_0$ acts on $\wt Y$ and $p$ is $\G_0$-invariant.  Let $\wt j: \wt T\incl \wt Y$ be the open embedding. Now $\G_0$ acts on $\cE:=\bR p_*\bR \wt j_*\Qlbar\cong \bR j_* p_{0*}\Qlbar$, which contains $\bR j_*\cL$ as the eigen summand for character $\chi_\cL\in \Hom(\G_0,\Qlbar^\times)$. Restricting $\cE$ to $Y_\a$ we get $\cE|_{Y_\a}\cong \bR p_{\a *}((\bR\wt j_*\Qlbar)|_{Y_\a})$, which is a successive extension of $\bR p_{\a *}\Qlbar$ and shifts (because the cohomology sheaves of $(\bR\wt j_*\Qlbar)|_{Y_\a}$ are constant). This implies $\cE|_{Y_\a}$ is tame, hence the same is true for $(\bR j_*\cL)|_{Y_\a}$. The local system $(\bR^i j_*\cL)|_{Y_\a}$ corresponds to a monodromy representation $r:\pi_1(Y_\a)^\tame\to \GL(V)$ (where $V$ is the stalk of $\bR j_*\cL$ at some $y_\a\in Y_\a$) that factors through the quotient $\G_\a:=\pi_1(Y_\a)^\tame/\pi_1(\wt Y_\a)^\tame$. We have a canonical map $\ov\r_\a: \G_0\cong \pi_1(Y_0)^\tame/\pi_1(\wt Y_0)^\tame\to \G_\a=\pi_1(Y_\a)^\tame/\pi_1(\wt Y_\a)^\tame$ compatible with $\r_\a: \pi_1(Y_0)^\tame\to \pi_1(Y_\a)^\tame$. Now the action of $\G_0$ on the stalk $(\bR^i j_*\cL)|_{y_\a}$ is via the composition $\G_0\to \G_\a\xr{r} \GL(V)$. On the other hand, the $\G_0$-action on $(\bR^i j_*\cL)|_{Y_\a}$ is by the character $\chi_\cL$, which implies that the composition $\G_0\xr{\ov\r_\a} \G_\a\xr{r} \GL(V)$ is $\chi_\cL\cdot \id_V$, implying the lemma.
\end{proof}

Note that $\cF_J\cong \mu_{J!}\Qlbar$ decomposes according to $\pi_0(\cW_J)$, which is a torsor under $\pi_0(ZL_J)$, and each summand has a rank one trivial sub local system. Therefore $\pi_0(ZL_J)$ acts on $\cF^0_J$, so that its stalks are free rank one modules over $\Qlbar[\pi_0(ZL_J)]$. For each $\chi\in \pi_0(ZL_J)^*=\Hom(\pi_0(ZL_J), \Qlbar^\times)$, let $\cF^0_{J,\chi}$ be the eigen sub local system of $\cF^0_J$ under the $\pi_0(ZL_J)$-action with eigenvalue $\chi$. We have a decomposition
\begin{equation}
    \cF^0_J\cong \bigoplus_{\chi\in \pi_0(ZL_J)^*}\cF^0_{J,\chi}.
\end{equation}
Let $i_J: \PP_J\incl \PP$ be the inclusion, and let
\begin{equation}
    \ov{\cF^0_{J,\chi}}=\bR^0 i_{J*}\cF^0_{J,\chi}.
\end{equation}
Since $\cF^0_{J,\chi}$ is a trivial local system on $\PP_J$, $\ov{\cF^0_{J,\chi}}$ is a trivial local system on $\ov\PP_J$.

Let $\frS_G$ be the set of pairs $(J,\chi)$ where $J\sft \wt I$ and $\chi\in \pi_0(ZL_J)^*=\Hom(\pi_0(ZL_J), \Qlbar^\times)$. For $J\subset J'$ we have an inclusion 
$L_{J}\subset L_{J'}$ as a Levi subgroup, hence an inclusion of the centers $ZL_{J'}\subset ZL_{J}$. It thus induces a map on the character groups
\begin{equation}
    \r^{J}_{J'}: \pi_0(ZL_{J})^*\to \pi_0(ZL_{J'})^*.
\end{equation}

We define a partial order on $\frS_G$ given by $(J,\chi)\le (J',\chi')$ if $J\subset J'$ and $\chi'=\r^{J}_{J'}(\chi)$. Let $\frS^{\min}_G$ be the set of minimal elements in $\frS_G$.

\begin{lemma}\label{l:unique min char}
\begin{enumerate}
    \item For $J\subset J'\sft \wt I$, the map $\r^{J}_{J'}$ is injective.
    \item For $J_1, J_2\sft \wt I$, we have a Cartesian diagram
    \begin{equation*}
        \xymatrix{\pi_0(ZL_{J_1\cap J_2})^*\ar[r]^{\r^{J_1\cap J_2}_{J_1}}\ar[d]_{\r^{J_1\cap J_2}_{J_2}} &  \pi_0(ZL_{J_1})^* \ar[d]^{\r^{J_1}_{J_1\cup J_2}}\\
        \pi_0(ZL_{J_2})^*\ar[r]^{\r^{J_2}_{J_1\cup J_2}} &  \pi_0(ZL_{J_1\cup J_2})^* }
    \end{equation*}
    \item For each $(J',\chi')\in \frS_G$, there is a unique $(J,\chi)\in \frS^{\min}_G$ such that $(J,\chi)\le (J',\chi')$.
\end{enumerate}
\end{lemma}
\begin{proof}
(1) Up to a Tate twist, we have $\pi_0(ZL_{J})^*\cong (\xch(T)/\ZZ\Phi_J)_\tors$ (where $\ZZ\Phi_J$ is the root lattice for $L_J$). The map $\r^{J}_{J'}$ for $J\subset J'$ is induced by the projection $p^{J}_{J'}: \xch(T)/\ZZ\Phi_J\to \xch(T)/\ZZ\Phi_{J'}$. Since $\ker(p^{J}_{J'})$ is torsion free, it induces an injective map on the torsion parts.

(2) We have an exact sequence
\begin{equation*}
    0\to \xch(T)/\ZZ\Phi_{J_1\cap J_2}\xr{(p^{J_1\cap J_2}_{J_1},p^{J_1\cap J_2}_{J_2} )}\xch(T)/\ZZ\Phi_{J_1}\op \xch(T)/\ZZ\Phi_{J_2}\xr{p^{J_1}_{J_1\cup J_2}-p^{J_2}_{J_1\cup J_2}} \xch(T)/\ZZ\Phi_{J_1\cup J_2}\to 0 
\end{equation*}
Taking torsion parts we get a left exact sequence
\begin{equation*}
    0\to \pi_0(ZL_{J_1\cap J_2})^*\xr{(\r^{J_1\cap J_2}_{J_1},\r^{J_1\cap J_2}_{J_2} )}\pi_0(ZL_{J_1})^*\op\pi_0(ZL_{J_2})^*\xr{\r^{J_1}_{J_1\cup J_2}-\r^{J_2}_{J_1\cup J_2}}\pi_0(ZL_{J_1\cup J_2})^*.
\end{equation*}

(3) If $(J_1,\chi_1)\le (J',\chi')$ and $(J_2,\chi_2)\le (J',\chi')$, then by (2) there is a unique $\chi\in \pi_0(Z_{J_1\cap J_2})^*$ such that $(J_1\cap J_2,\chi)\le (J_1,\chi_1)$  and $(J_1\cap J_2,\chi)\le (J_2,\chi_2)$. From this we get the uniqueness of the minimal $(J,\chi)\le (J',\chi')$.
\end{proof}

\begin{prop}\label{p:decomp F0}
    There is a canonical decomposition
    \begin{equation*}
        \cF^0\cong \bigoplus_{(J,\chi)\in \frS^{\min}_G} \ov{\cF^0_{J,\chi}}.
    \end{equation*}
    where each $\ov{\cF^0_{J,\chi}}$ is a rank one trivial local system on $\ov{\PP_J}$.
\end{prop}
\begin{proof}
    For each $(J,\chi)\in \frS^{\min}_G$, the adjunction map $\cF^0\to i_{J*}\cF^0_J$ followed by the projection to $i_{J*}\cF^0_{J,\chi}$ gives a map $\cF^0\to i_{J*}\cF^0_{J,\chi}$, which factors through $\bR^0i_{J*}\cF^0_{J,\chi}=\ov{\cF^0_{J,\chi}}$. Therefore we get a canonical map
    \begin{equation}\label{F0 maps to decomp}
        \cF^0\to \bigoplus_{(J,\chi)\in \frS^{\min}_G} \ov{\cF^0_{J,\chi}}.
    \end{equation}
    To show this is an isomorphism, it suffices to check it after restricting to each $\PP_{J'}$. For $(J,\chi)\in\frS^{\min}_G$ and $J'\supset J$, the restriction $i_{J'}^*\ov{\cF^0_{J,\chi}}$ is canonically isomorphic to $\cF^0_{J',\chi'}$ where $\chi'=\chi|_{\pi_0(ZL_{J'})}$, and the canonical map $\cF^0_{J'}\to i_{J'}^*\ov{\cF^0_{J,\chi}}\cong \cF^0_{J',\chi}$ obtained by restricting \eqref{F0 maps to decomp} to $\PP_{J'}$ and projecting to the $(J,\chi)$-summand is the projection to the $\chi'$-summand of $\cF^0_{J'}$. Therefore, to show \eqref{F0 maps to decomp} is an isomorphism over $\PP_{J'}$, it suffices to show that the sum of the projections
    \begin{equation}
        \cF^0_{J'}\to \bigoplus_{(J,\chi)\in \frS^{\min}_G} \cF^0_{J',\r^J_{J'}(\chi)}
    \end{equation}
    is an isomorphism, i.e. the map
    \begin{equation}
        \coprod_{(J,\chi)\in \frS^{\min}_G}\pi_0(ZL_J)^*\xr{\coprod_{J,\chi} \r^{J}_{J'}}\pi_0(ZL_{J'})^*
    \end{equation}
    is a bijection. In other words, we need to show that for each $(J',\chi')\in \frS_G$, there is a unique $(J,\chi)\in \frS^{\min}_G$ such that $(J,\chi)\le (J',\chi')$. This is precisely Lemma \ref{l:unique min char}(3). 
\end{proof}

\begin{cor} Recall that $G$ is semisimple. 
The total cohomology $\cohog{*}{\cM^\c_G}$ (as a graded vector space) has a canonical decomposition
\begin{equation}\label{decomp coho MG}
    \cohog{*}{\cM^\c_G}\cong \bigoplus_{(J,\chi)\in \frS^{\min}_G}\cohog{*}{\ov\PP_J, \ov{\cF^0_{J,\chi}}}[-2|J|](-|J|),
\end{equation}
where each summand is isomorphic to $\cohog{*}{\ov\PP_J}[-2|J|](-|J|)$ (so that the top degree is $2r$) after trivializing $\ov{\cF^0_{J,\chi}}$. In particular, $\cohog{*}{\cM^\c_G}$ is concentrated in even degrees and pure with Poincar\'e polynomial
\begin{equation*}
    \sum_{i}\dim \cohog{i}{\cM^\c_G}t^i=\sum_{(J,\chi)\in \frS^{\min}_G} t^{2|J|}[r+1-|J|]_{t^2}.
\end{equation*}
\end{cor}
\begin{proof}
    We first observe that the inclusion $\cF^0\incl \cF$ induces an isomorphism on cohomology
    \begin{equation}\label{F0 F same coho}
        \cohog{*}{\PP, \cF^0}\isom \cohog{*}{\PP, \cF}
    \end{equation}
    Indeed, by excision we reduce to the claim that the inclusion $\cF^0_J\incl \cF_J$ induces an isomorphism $\cohoc{*}{\PP_J, \cF^0_J}\isom\cohoc{*}{\PP_J, \cF_J}$, which is true since nontrivial tame local systems on the torus $\PP_J$ have zero (compactly supported or not) cohomology.

    The isomorphism \eqref{F0 F same coho} implies that $\cohog{*}{\PP, \cF^0}\isom \cohog{*}{\PP, \cF}\cong \cohoc{*}{\cM^\c_G}[2r](r)$. The decomposition in Proposition \ref{p:decomp F0} after taking 
    global sections gives the decomposition of $\cohoc{*}{\cM^\c_G}[2r](r)$. The desired decomposition of $\cohog{*}{\cM^\c_G}$ is obtained by Poincar\'e duality.
\end{proof}

When $G$ is simple adjoint, $\frS^{\min}_G$ consists of pairs $(\wt I-\wt I_d, \chi)$ where $\wt I_d\ne\vn$ and $\chi$ is a primitive character of the cyclic group $\pi_0(ZL_{\wt I-\wt I_d})$ of order $d$. From this we get:

\begin{cor} Assume $G$ is simple adjoint. Then we have
\begin{equation}\label{decomp coho MG adj}
    \cohog{*}{\cM^\c_G}\cong \bigoplus_{d\in \NN, \wt I_d\ne \vn} \cohog{*}{\ov\PP_{\wt I-\wt I_d}}^{\op\ph(d)}[-2(r+1-r_d)](-r-1+r_d)
\end{equation}
and
\begin{equation*}
    \sum_{i}\dim\cohog{i}{\cM^\c_G}t^i\cong \sum_{d\in \NN, \wt I_d\ne \vn} t^{2(r+1-r_d)}[r_d]_{t^2}\ph(d).
\end{equation*}
\end{cor}

\sss{Perverse filtration}
The proper map $f^\c_G: \cM^\c_G\to \cA_G$ gives a perverse filtration on the direct image complex $\bR f^\c_{G!}\Qlbar$, hence on the cohomology of $\cM^\c_G$. More precisely, let ${}^p\t_{\le i}\bR f^\c_{G*}\QQ$ be the perverse truncations of the direct image complex $\bR f^\c_{G*}\QQ$ of the Hitchin fibration $f^\c_G:\cM^\c_G\to \cA_G$. Recall $r$ is the rank of $G$. Define the perverse filtration on $\cohog{*}{\cM^\c_G,\QQ}$ by
\begin{equation*}
    P_{i}\cohog{*}{\cM^\c_G,\QQ}:=\Im\left(\cohog{*}{\cA_G, {}^p\t_{\le i+r}\bR f^\c_*\QQ}\to \cohog{*}{\cA_G, \bR f^\c_{*}\QQ}\right)\subset \cohog{*}{\cM^\c_G,\QQ}.
\end{equation*}
Then the associated graded $\Gr^P_{i}\cohog{*}{\cM^\c_G,\QQ}$ is zero unless $0\le i\le 2r$.

We make a precise prediction on the perverse filtration on $\cohog{*}{\cM^\c_G}$.

\begin{conj}\label{c:coho MG} Assume $G$ is semisimple. The perverse filtration on $\cohog{*}{\cM^\c_G}$ is compatible with the decomposition \eqref{decomp coho MG}. Moreover, for each $(J,\chi)\in \frS^{\min}_G$, we have
\begin{equation}
    P_i\left(\cohog{*}{\ov\PP_J, \ov{\cF^0_{J,\chi}}}[-2|J|](-|J|)\right)=\t_{\le i+|J|}\left(\cohog{*}{\ov\PP_J, \ov{\cF^0_{J,\chi}}}[-2|J|](-|J|)\right), \quad i\in \ZZ.
\end{equation}
Here $\t_{\le n}$ denotes the truncation in the usual t-structure of $D(\Vect)$.

In particular, for $0\le i\le 2r$,
\begin{equation}
    \Gr^P_i\cohog{*}{\cM^\c_G}\cong \bigoplus_{(J,\chi)\in \frS^{\min}_G, |J|\le\min\{i,2r-i\}}\cohog{i-|J|}{\ov\PP_J, \ov{\cF^0_{J,\chi}}}(-|J|).
\end{equation}
\end{conj}

\begin{remark}
    If Conjecture \ref{c:coho MG} holds and $G$ is simple adjoint, then for $0\le j\le r$ and $j\le i\le 2j$, we have
    \begin{equation}\label{dim GrP}
        \dim \Gr^P_i\cohog{2j}{\cM^\c_G}=\sum_{d\in \NN, \wt I_d\ne \vn, 2j-i=r-r_d+1}\ph(d).
    \end{equation}
    Otherwise $\Gr^P_i\cohog{2j}{\cM^\c_G}=0$.
 \end{remark}

\begin{exam}\label{ex:E8 MG}
    When $G=E_{8}$, we have
\begin{eqnarray*}
r_1=9, r_{2}=5, r_3=3, r_4=2, r_5=1, r_6=1.
\end{eqnarray*}
The Betti numbers of each summand of \eqref{decomp coho MG adj} are distributed as follows:
\begin{equation}
    \begin{array}{c|ccccccccc}
         d\bs 2j & 0 & 2 & 4 & 6 & 8 & 10 & 12 & 14 & 16 \\
         \hline
         1 & 1 & 1 & 1 & 1 & 1 & 1 & 1 & 1 & 1 \\
         2 &   &   &   &   & 1 & 1 & 1 & 1 & 1 \\
         3 &   &   &   &   &   &   & 2 & 2 & 2 \\
         4 &   &   &   &   &   &   &   & 2 & 2 \\
         5 &   &   &   &   &   &   &   &   & 4 \\
         6 &   &   &   &   &   &   &   &   & 2
    \end{array}
\end{equation}
Here the column labeled by $2j$ stands for $\cohog{2j}{\cM^\c_G}$, and the row labeled by $d$ is the contribution of the $d$-summand of \eqref{decomp coho MG adj}.

On the other hand, assuming Conjecture \ref{c:coho MG}, the dimensions of the perverse graded pieces of $\cohog{*}{\cM^\c_G}$ are distributed as follows:
\begin{equation}\label{perv E8}
    \begin{array}{c|ccccccccccccccccc}
         d\bs i & 0 & 1 & 2 & 3 & 4 & 5 & 6 & 7 & 8 & 9 & 10 & 11 & 12 & 13 & 14 & 15 & 16 \\
         \hline
         1 & 1 && 1 && 1 && 1 && 1 && 1 && 1 && 1 && 1 \\
         2 &   &   &   &   & 1 && 1 && 1 && 1 && 1 &   &   &   &   \\
         3 &   &   &   &   &   &   & 2 && 2 && 2 &   &   &   &   &   &  \\
         4 &   &   &   &   &   &   &   & 2 && 2 &   &   &   &   &   &   &  \\
         5 &   &   &   &   &   &   &   &   & 4 &   &   &   &   &   &   &   &   \\
         6 &   &   &   &   &   &   &   &   & 2 &   &   &   &   &   &   &   &   
    \end{array}
\end{equation}
Here the the column labeled by $i$ stands for $\Gr^P_i\cohog{*}{\cM^\c_G}$.
\end{exam}

\subsection{Point-counting and cohomology of $J^{G^{\ad}}_{G^{\sc}}$}\label{ss:JG Fq} 
Following the notations of \S\ref{sss:J sc},  let $J^{G^{\ad}}_{G^{\sc}}$ be the centralizer group scheme of the Steinberg section in the simply-connected form $G^\sc$ while the centralizer is taken in the adjoint group $G^{\ad}$.

In \cite[Table in p.158]{L-Coxeter}, Lusztig computed the cardinality of $J^{G^{\ad}}_{G^{\sc}}(\FF_{q})$. His result can be reformulated as:
\begin{equation}\label{JG Fq}
|J^{G^{\ad}}_{G^{\sc}}(\FF_{q})|=\sum_{d\in \NN}q^{r+1-r_d}[r_d]_{q^{2}}\ph(d).
\end{equation}

\begin{exam}  When $G=E_8$, using the values of $r_d$ in Example \ref{ex:E8 MG}, the right side of \eqref{JG Fq} is
\begin{eqnarray*}
&&[9]_{q^{2}}+q^{4}[5]_{q^{2}}+q^{6}[3]_{q^{2}}\cdot 2+q^{7}[2]_{q^{2}}\cdot 2+q^{8}[1]_{q^{2}}\cdot 4+q^{8}[1]_{q^{2}}\cdot 2\\
&=&1+q^{2}+2q^{4}+4q^{6}+2q^{7}+10q^{8}+2q^{9}+4q^{10}+2q^{12}+q^{14}+q^{16}
\end{eqnarray*}
which is consistent with Lusztig's calculation in \cite[Table in p.158]{L-Coxeter}. The coefficients of $q$ above matches the column sums of Table \eqref{perv E8}.  
\end{exam}

\begin{remark}
    In \cite{L-Coxeter}, Lusztig observes that $J^{G^{\ad}}_{G^{\sc}}$ has the same point-counting as $X_c\times X_c/G(\FF_q)$, where $X_c$ is the Deligne-Lusztig variety for a Coxeter element $c\in W$. The summation in \eqref{JG Fq} is consistent with the finer information given in \cite[Table (7.3)]{L-Coxeter} concerning the Frobenius eigenvalues on $\cohoc{*}{X_c}$. For example, the Frobenius eigenvalues on $\cohoc{r}{X_{f},\Qlbar}$ are exactly given by
\begin{equation*}
\z q^{(r+1-r_d)/2}, \mbox{ where $\z\in \mu_{d}$ is primitive and $d\in \NN$ such that $\wt I_{d}\ne\vn$}.
\end{equation*}
\end{remark}

By considering a model of $J^{G^{\ad}}_{G^{\sc}}$ over $\ZZ$, point-counting over finite fields determines the Euler characteristics of the weight filtration on the compactly supported cohomology of the complex fiber of $J^{G^{\ad}}_{G^{\sc}}$. Passing to $\cohog{*}{J^{G^{\ad}}_{G^{\sc}}}$ by duality, we get the following consequence.

\begin{cor}\label{c:wt fil J}
    Consider $J^{G^{\ad}}_{G^{\sc}}$ over $\CC$. Let $W_i\cohog{*}{J^{G^{\ad}}_{G^{\sc}}, \QQ}$ be the weight filtration on the cohomology of $J^{G^{\ad}}_{G^{\sc}}$. Then $\Gr^W_i\cohog{*}{J^{G^{\ad}}_{G^{\sc}},\QQ}$ has zero Euler characteristic for $i$ odd, and for $0\le i\le 2r$
    \begin{equation}\label{chi W HJ}
        \chi(\Gr^W_{2i}\cohog{*}{J^{G^{\ad}}_{G^{\sc}},\QQ})=\sum_{d\in \NN, |i-r|\le r_d-1}\ph(d).
    \end{equation}
\end{cor}

Note that $q^{r+1-r_d}[r_d]_{q^{2}}$ is palindromic about $r$ as a polynomial of $q$, hence the right side of \eqref{JG Fq} as a polynomial of $q$ is palindromic and unimodal. 
Moreover, the formula \eqref{JG Fq} differs from \eqref{MG Fq} only in reassigning exponents of $q$ in each summand. 
We now give a conjectural explanation of these phenomena.

\subsection{Non-abelian Hodge package for $\cM_G$ and P=W conjecture}\label{ss:P=W}

As explained in \S\ref{ss:GL}, $\cM_G$ can be viewed as a Dolbeault moduli space for $G$-bundles on $\PP^1$ with specific ramifications at $0$
and $\infty$, while $J_G$ can be viewed as a Betti moduli space for $G$-local systems with the same ramification conditions. Following the  general paradigm of Simpson correspondence \cite{Simpson}, we propose

\begin{conj}\label{c:HK} Let $G$ be semisimple over $\CC$ and $G^\sc$ be its simply-connected cover.
Then there exists a hyperK\"ahler structure on $\cM^\c_G$ such that $J^{G}_{G^\sc}$ as a complex manifold is obtained from $\cM^\c_G$ by  hyperK\"ahler rotation. In particular, there is a homeomorphism $\cM^\c_G\simeq J^{G}_{G^\sc}$ as $C^\infty$-manifolds.
\end{conj}

Motivated by the P=W conjecture \cite{dCHM} for the usual Hitchin moduli space and character variety, we propose a P=W conjecture for $\cM_G$. 

\begin{conj}[P=W]\label{c:P=W} Let $G$ be semisimple over $\CC$.
    Under the homeomorphism $\cM^\c_G\simeq J^{G}_{G^\sc}$ in Conjecture \ref{c:HK}, the induced isomorphism on cohomology
    \begin{equation*}
        \cohog{*}{\cM^\c_G,\QQ}\isom \cohog{*}{J^G_{G^\sc}, \QQ}
    \end{equation*}
    sends $P_{i}\cohog{*}{\cM^\c_G,\QQ}$ isomorphically to $W_{2i}\cohog{*}{J^G_{G^\sc},\QQ}$, for $i=0,1,\cdots, 2r$.
\end{conj}

Based on the P=W conjecture and Conjecture \ref{c:coho MG}, we make more precise predictions on the cohomology of $J^G_{G^\sc}$.

\begin{conj}\label{c:coho J}
    Let $G$ be semisimple over $\CC$. Then 
    \begin{enumerate}
        \item $\cohog{i}{J^G_{G^\sc},\QQ}=0$ for odd $i$.
        \item For $0\le j\le r$ and $j\le i\le 2j$, we have
        \begin{equation}
            \dim \Gr^W_{2i}\cohog{2j}{J^G_{G^\sc},\QQ}=\sum_{d\in \NN, \wt I_d\ne \vn, 2j-i=r-r_d+1}\ph(d).
        \end{equation}
        Otherwise $\Gr^W_{2i}\cohog{2j}{J^G_{G^\sc},\QQ}=0$.
        In particular, $\dim \cohog{2}{J^G_{G^\sc},\QQ}=1$ and it is pure of weight $4$. \footnote{This follows from the observation that $r_1=r+1$ and $r_d<r$ for $d>1$.}
        \item Let $\y\in \cohog{2}{J^G_{G^\sc},\QQ}$ be a nonzero element (which has weight $4$ by part (2)). Then the cup product action of $\y$ on $\cohog{*}{J^G_{G^\sc},\QQ}$ satisfies the hard Lefschetz property with respect to the weight filtration: for $0\le i\le r$, cupping with $\y^i$ induces an isomorphism
        \begin{equation}
            \cup \y^i: \Gr^W_{2r-2i}\cohog{*}{J^G_{G^\sc},\QQ}\isom \Gr^W_{2r+2i}\cohog{*+2i}{J^G_{G^\sc},\QQ}.
        \end{equation}
    \end{enumerate}
\end{conj}
The hard Lefschetz property follows from Conjecture \ref{c:P=W} and the usual hard Lefschetz for the perverse filtration on $\bR f_{G*}\QQ$.

%%%%%%%%%%%%%%%%%%%%%%%%%%%%%%%%%%%%%%%%%%%%%%%%%%%%%
\subsection{Counting automorphic forms}\label{ss:auto}
The definition of the space $\cM_{G}$ was originally motivated by considerations of a specific kind of automorphic forms for the function field  $F=\FF_{q}(t)$ of $\PP^{1}_{\FF_{q}}$, which we now explain. We work over the base field $k=\FF_q$ and assume $G$ is split semisimple.

Fix a nontrivial additive character $\psi_{q}: \FF_{q}\to \Qlbar^{\times}$. For a closed point $x\in |\PP^{1}|$, let $F_{v}, \cO_{v}$ be the local field at $v$ and its valuation ring.

\sss{Simple supercuspidals} Recall the notion of {\em simple supercuspidal representations} of the $p$-adic group $G(F_{0})$ \cite{GrossReeder}. It is an irreducible admissible representation $\pi_{0}$ of $G(F_{0})$ (over $\Qlbar$) such that $\pi_{0}$ contains an eigenvector for $\bI^{+}_{0}$ on which it acts through a generic character $\psi_{q}\c\psi_0: \bI^{+}_{0}\xr{\psi_0} \FF_{q}\xr{\psi_{q}}\Qlbar^{\times}$. Genericity means $\psi_{0}$ is nontrivial on each affine simple root space. For a simple supercuspidal $\pi_{0}$ of $G(F_{0})$, the eigenspace $\pi_{0}^{(\bI^{+}_{0}, \psi_q\psi_{0})}$ is one-dimensional. 

\sss{Automorphic representations with local conditions} Fix generic characters $\psi_{0}$ and $\psi_{\infty}$ of $\bI^{+}_{0}$ and $\bI^{+}_{\infty}$ respectively. We consider the set $\Pi(\psi_{0},\psi_{\infty})$ of automorphic representations  $\pi$ of $G(\AA_{F})$ such that
\begin{enumerate}
\item $\pi_{v}$ is unramified for $v\ne 0,\infty$;
\item $\pi_{0}$ is a simple supercuspidal representation of $G(F_{0})$ with respect to $(\bI_{0}^{+}, \psi_{0})$;
\item $\pi_{\infty}$ is a simple supercuspidal representation of $G(F_{\infty})$ with respect to $(\bI_{\infty}^{+}, \psi_{\infty})$.
\end{enumerate}
Any such $\pi$ is necessarily cuspidal hence appears discretely in the automorphic spectrum of $G(\AA_{F})$. Let $m(\pi)$ be the multiplicity of $\pi$ in $C_c(G(F)\bs G(\AA_{F}), \Qlbar)$. 

The next result is not needed for the rest of the paper; it serves as a motivation for counting the number of points on $\cM_G$.

\begin{prop}\label{p:count auto forms}
We have
\begin{equation}\label{MGFq sum mult}
\sum_{\pi\in \Pi(\psi_{0}, \psi_{\infty})}m(\pi)=q^{-r}|\cM_{G}(\FF_{q})|.
\end{equation}
In particular, if $G$ is simple adjoint, and $q$ satisfies \eqref{large q} and \eqref{large q om} for every $\om \in\Om=\pi_1(G)$, then
\begin{equation}\label{sum mult formula}
\sum_{\pi\in \Pi(\psi_{0}, \psi_{\infty})}m(\pi)=|\pi_{1}(G)|\sum_{d\in \NN, \wt I_{d}\ne\vn}[r_d]_{q}\ph(d).
\end{equation}
\end{prop}
\begin{proof}
    Let $\cC(\psi_0,\psi_\infty)$ be the space of $\Qlbar$-valued functions on $\Bun_G(\bI^{++}_0, \bI^{++}_\infty)(\FF_q)$ that are eigen under the actions of $V_0\times V_\infty=\bI_0^+/\bI^{++}_0\times \bI_\infty^+/\bI^{++}_\infty$ with eigencharacter $\psi_q\psi_0\times \psi_q\psi_\infty$. 
    
    For each $\pi\in \Pi(\psi_{0}, \psi_{\infty})$ realizes as a sub-$G(\AA_F)$-module of $C_c(G(F)\bs G(\AA_F), \Qlbar)$, the following line
    \begin{equation*}
        \ell_\pi=\pi_0{(\bI_0^+, \psi_q\psi_0)}\ot\left(\bigotimes_{v\neq 0,\infty}\pi_v^{G(\cO_v)}\right)\ot  \pi_0{(\bI_\infty^+, \psi_q\psi_\infty)}
    \end{equation*}
    lies in $\cC(\psi_0,\psi_\infty)$. Indeed,
    \begin{equation*}
        \cC(\psi_0,\psi_\infty)=\bigoplus_{\pi\in \Pi(\psi_{0}, \psi_{\infty})}\ell_\pi\ot \Hom_{G(\AA_F)}(\pi, C_c(G(F)\bs G(\AA_F))).
    \end{equation*}
    In particular, 
    \begin{equation}\label{dimC sum mult}
        \dim\cC(\psi_0,\psi_\infty)=\sum_{\pi\in \Pi(\psi_{0}, \psi_{\infty})}m(\pi).
    \end{equation}

    To compute $\dim\cC(\psi_0,\psi_\infty)$, observe that every function $f\in \cC(\psi_0,\psi_\infty)$ has a well-defined support in $\Bun_G(\Iu)(\FF_q)$. A point $\xi\in \Bun_G(\Iu)(\FF_q)$ is called {\em $(\psi_0,\psi_\infty)$-relevant} if there exists a nonzero function $f\in\cC(\psi_0,\psi_\infty)$ whose support is the singleton $\{\xi\}$. Thus
    \begin{equation}\label{dimC rel points}
        \dim\cC(\psi_0,\psi_\infty)=|\mbox{$(\psi_0,\psi_\infty)$-relevant points (up to isom) in $\Bun_G(\Iu)(\FF_q)$}|.
    \end{equation}
    Note that the right side is counting isomorphism classes, not weighted by the reciprocal of automorphisms.

    For $\xi\in \Bun_G(\Iu)(\FF_q)$, let $\Aut(\xi)$ be its automorphism group, viewed as an algebraic group over $\FF_q$. Evaluation at $0$ (upon choosing trivializations of $\xi$ at $0$) gives a homomorphism $\ev_0:\Aut(\xi)\to \bI_0^+$, whose composition with the projection to $V_0$ is independent of the trivialization, giving a canonical homomorphism
    \begin{equation*}
        \ov\ev_0: \Aut(\xi)\to V_0.
    \end{equation*}
    Similarly we have $\ov\ev_\infty: \Aut(\xi)\to V_\infty$. It is easy to check from the definitions that $\xi$ is $(\psi_0,\psi_\infty)$-relevant if and only if the pullback of the character $\psi_q\psi_0\bt\psi_q\psi_\infty$ along the homomorphism
    \begin{equation*}
        \Aut(\xi)(\FF_q)\xr{(\ov\ev_0,\ov\ev_\infty)} V_0(\FF_q)\op V_\infty (\FF_q)
    \end{equation*}
    is trivial. This is the equivalent to that the composition
    \begin{equation}\label{Aut psi}     \Aut(\xi)\xr{(\ov\ev_0,\ov\ev_\infty)}V_0\op V_\infty\xr{\psi_0+\psi_\infty}\Ga
    \end{equation}
    is trivial. Note that $\Aut(\xi)$ is a unipotent group since it injects into  $\bI_0^+$ by evaluating at $0$. Therefore $\uAut(\xi)$ is unipotent. Therefore \eqref{Aut psi} being trivial is equivalent to the triviality of its derivative (where $\aut(\xi)=\Lie\Aut(\xi)$)
    \begin{equation*}
        \a: \aut(\xi)\xr{(\ov\frev_0,\ov\frev_\infty)}V_0\op V_\infty\xr{\psi_0+\psi_\infty}\FF_q.
    \end{equation*}
    Dualizing, we see that $\xi$ is $(\psi_0,\psi_\infty)$-relevant if and only if the image of $(\psi_0,\psi^*_\infty)$ under the map
    \begin{equation*}
        \a^*: V_0^*\op V_\infty^*\to \aut(\xi)
    \end{equation*}
    is zero.

    We next claim that 
    \begin{equation}\label{rel iff im b}
        \mbox{$\xi$ is $(\psi_0,\psi_\infty)$-relevant $\iff$ $\xi$ is in the image of $\b: \cM_G(\FF_q)\to \Bun_G(\Iu)(\FF_q)$.}
    \end{equation} 
    Moreover, for each $(\psi_0,\psi_\infty)$-relevant $\xi$, we have
    \begin{equation}\label{qr}
        |b^{-1}(\xi)(\FF_q)|\cdot |\Aut(\xi)(\FF_q)|^{-1}=q^r.
    \end{equation}
    Now \eqref{qr} implies that the number of $(\psi_0,\psi_\infty)$-relevant points of $\Bun_G(\Iu)(\FF_q)$ is equal to  $q^{-r}|\cM_G(\FF_q)|$. Combined with \eqref{dimC rel points} and \eqref{dimC sum mult}, we get \eqref{MGFq sum mult}.

    It remains to prove \eqref{rel iff im b} and \eqref{qr}.  Consider the $V_0\times V_\infty$-torsor $\pi: \Bun_G(\bI^{++}_0, \bI^{++}_\infty)\to \Bun_G(\Iu)$. Let $\wt\xi\in \Bun_G(\bI^{++}_0, \bI^{++}_\infty)(\FF_q)$ be any lift of $\xi$. The long exact sequence attached to the map of cotangent complexes induced by $\pi$ reads
    \begin{equation}\label{long ex cot}
        0\to T^*_\xi\Bun_G(\Iu)\to T^*_{\wt\xi}\Bun_G(\bI^{++}_0, \bI^{++}_\infty)\xr{\mu} V_0^*\op V_\infty^*\xr{\a^*}\aut(\xi)
    \end{equation}
    Here $T^*$ means classical cotangent spaces, and $\mu$ is the moment map for the $V_0\times V_\infty$-action. The fiber $b^{-1}(\xi)$ is $\mu^{-1}(\psi_0, \psi_\infty)$. By the long exact sequence \eqref{long ex cot}, $\mu^{-1}(\psi_0, \psi_\infty)$ is non-empty if and only if $\a^*(\psi_0,\psi_\infty)$, which we have shown to be equivalent to $\xi$ being  $(\psi_0,\psi_\infty)$-relevant. This proves \eqref{rel iff im b}. When $\xi$ is $(\psi_0,\psi_\infty)$-relevant, $b^{-1}(\xi)$ is a torsor for $T^*_\xi\Bun_G(\Iu)$. Using that $\Aut(\xi)$ unipotent, we have
    \begin{equation*}
        |b^{-1}(\xi)(\FF_q)|\cdot |\Aut(\xi)(\FF_q)|^{-1}=|T^*_\xi\Bun_G(\Iu)|\cdot |\aut(\xi)|^{-1}.
    \end{equation*}
    Now the right side is $q^{\chi(\LL_\xi)}$ where $\LL_\xi$ is the fiber of the cotangent complex of $\Bun_G(\Iu)$ at $\xi$. Since $ \Bun_G(\Iu)$ is smooth of dimension $r$, we get $\chi(\LL_\xi)=r$, and \eqref{qr} follows.

    Finally, \eqref{sum mult formula} follows from Proposition \ref{p:MG Fq}.
\end{proof}

%%%%%%%%%%%%%%%%%%%%%%%%%%%%%
%%%%%%%%%%%%%%%%%%%%%%%%%%%%%%
\section{Fukaya category of $\cM_G$}\label{s:Fukaya}

Let $\omega=\omega_{\cM_G^\circ}$ be the \emph{real part} of the canonical holomorphic symplectic form $\omega_{\alg}$ on $\cM_G^\circ$. 
Our first goal in this section is to construct suitable \emph{real} Liouville 1-forms $\vartheta$ on $\cM_G^\circ$ so that $(\cM_G^\circ, \vartheta)$ is a Liouville manifold (with $d\vartheta=\omega_{\cM_G^\circ}$). We show that all such constructed Liouville manifold structures are canonically equivalent (up to a contractible space of deformations). This enables us to set the wrapped Fukaya category $\cW(\cM_G^\circ)$ of $\cM_G^\circ$ in a canonical way. The construction of $\vartheta$ occupies \S\ref{subsec: Liouville}-\S\ref{subsecLiouvillefromZflat}, and the main tool is by ``modifying" the canonical Liouville 1-forms on $\cT_{L_J}, J\subset_{ft}\wt{I}$, by constructing (various versions of) ``admissible" vector fields $\sfZ_\flat$ on the moment polytope $\trg$ of $\PP$. We further show that for a nice subspace of $\vartheta$, $\cM_G^\circ$ admits a Weinstein manifold structure (cf. \S\ref{subsecWeinMfldMG}). 

Our second goal is to construct convenient Weinstein sectorial coverings $\{\cT_{L_J}^{\ovl{\cU}_J}\}_{J\subset_{ft}\wt{I}}$ of $\cM_G^\circ$ (cf. \S\ref{subsec: sectorial,covering}). Using results from Appendices \ref{subsecSectdeformJG} and \ref{subsecSectInclu}, we give a canonical (homotopy) colimit diagram to calculate $\cW(\cM_G^\circ)$ , whose terms are the (partially) wrapped Fukaya categories of $\ovl{\cT}_{L_J}, J\subset_{ft}\wt{I}$ (cf. Corollary \ref{cor: cWcMG,twocolimits}). This colimit diagram is crucially used in the proof of Mirror Theorem \ref{th:mirror} (cf. \S\ref{secProofMirror}).

\subsection{Grading data on $\cM_G^\circ$}

For any holomorphic symplectic manifold $(M,\omega_{M;\CC})$ of complex dimension $n$, let $\omega_M=\Real\omega_{M;\CC}$. Then we have $c_1(M, \omega_M)\in H^2(M,\ZZ)$ is trivial since it is classified by  
\begin{align*}
M\lrar \BB\Sp\lrar \BB U\lrar \BB^2 \ZZ,
\end{align*}
and $H^2(\BB\Sp, \ZZ)=0$ (cf. \cite{CKNS} for more details). Here $\Sp$ is the stable compact symplectic group (the stabilization of the compact form of $\Sp(2n,\CC)$).
There exists a canonical trivialization (up to homotopy) induced from the canonical trivialization of $\BB\Sp\lrar \BB^2 \ZZ$, since $H^1(\BB\Sp, \ZZ)=0$. Let $\JJ$ be any almost complex structure on $M$ that is compatible with $\omega_M$, which forms a contractible space. Let 
\begin{align}\label{eqkappaTMJJ}
\kappa=(\Lambda^{n}(TM, \JJ))^\vee. 
\end{align}
The above canonical trivialization gives a canonical trivialization $\alpha_M$ of $\kappa^{\otimes 2}$ (up to homotopy) that is used to assign a squared phase map to any Lagrangian submanifold $\sfL$ in $M$, i.e. $\alpha_M^{\sfL}: \sfL\lrar S^1$. For any holomorphic Lagrangian $\sfL_{\hol}$, the squared phase map is always nullhomotopic since it is classified by 
\begin{align*}
\sfL_{\hol}\lrar \Sp/U\lrar U/O\lrar \BB\ZZ, 
\end{align*}
and $H^1(\Sp/U, \ZZ)=0$. Each trivialization amounts to a lifting of $\alpha_M^{\sfL_{\hol}}$ to $\sfL_{\hol}\to \RR$, and the space of liftings forms a torsor of $H^0(\sfL_{\hol}, \ZZ)$. In fact, we can choose $\alpha_M$ so that $\alpha_M^{\sfL_{\hol}}$ has constant phase $1\in S^1$ for any holomorphic Lagrangian $\sfL_{\hol}$, and then 
 the grading on $\sfL_{\hol}$ is always an integer. This is discussed in \S\ref{appconstantgrading} based on \emph{compatible (almost) quaternionic structures}. 

\begin{remark}
For $G$ semisimple, by Lemma \ref{lemmapi1cMG}, we have both $H^1(\cM_G^\circ, \ZZ)=0$ and $H^1(\cT_G, \ZZ)=0$. Therefore, the canonical grading datum on $\cM_G^\circ$ (resp. $\cT_G$) is the unique grading datum on $\cM_G^\circ$ (resp. $\cT_G$). 
\end{remark}

\subsection{Real valued Liouville 1-forms on $\cM_G$}\label{subsec: Liouville}

For simplicity, we will mostly assume $G$ is simply connected. It will become clear that all constructions below are compatible with any central isogeny of $G$. Hence all results here are valid for $\cM_G^\circ$, for semisimple groups. 

We would like to specify for $\cM_G$ an exact symplectic structure, so that we can study (certain versions of) the Fukaya category on it. We cannot hope for a holomorphic primitive 1-form for the holomorphic symplectic form because of the following example. Our first goal in this section is to construct suitable real valued Liouville 1-forms $\vartheta$ (in a convenient way) on $\cM_G^\circ$, which admit Weinstein manifold structures.   

\begin{exam}\label{exam: omegaSL2}
Let $G=\SL_2$. 
To simplify notations, we will denote $L_{\{0\}}$ by $L^0$. 
Using coordinates from Example \ref{exam: SL2coord} and \ref{exam: SL2hxi} (and the cyclic symmetry of the affine Dynkin diagram), for the affine open cover $\cM_G=\cT_{G}\cup \cT_{T}\cup \cT_{L^0}$, we have: 
\begin{align}
\nonumber \cO(\cT_{G})&\hookrightarrow \cO(\cT_{T})\\
\nonumber (a, x, y) &\mapsto (-\xi^2-h^{-2}, h, h\xi),
\end{align}
\begin{align}
\nonumber \cO(\cT_{L^0})&\hookrightarrow \cO(\cT_{T})\\
\label{exam: axyprime}(a', x', y') &\mapsto (-\xi^2-h^{2}, h^{-1}, -h^{-1}\xi).
\end{align}
On the other hand, since $\cM_G\rightarrow \cA_G$ has projective fibers, 
we have inside $\cO(\cT_{T})$,  
$\cO(\cT_G)\cap \cO(\cT_{L^0})=\cO(\cM_G)=\CC[\xi^2+h^2+h^{-2}]$. 

Let $\vartheta_G=\frac{dy}{x}$, $\vartheta_{\vn}=\xi\frac{d h}{h}$ and $\vartheta_0=\frac{dy'}{x'}$ be the standard primitives of the restriction of $-\omega$  on $\cT_{G},\cT_T, \cT_{L^0}$ respectively. 
We have $\vartheta_G-\vartheta_{\vn}=d\xi$ and $\vartheta_0-\vartheta_{\vn}=-d\xi$. Since $H^1(\cT_G, \CC)=H^1(\cT_{L^0}, \CC)=0$ (cf. Lemma \ref{lemmapi1cMG} below) to construct a global $\vartheta$, we must find functions $f_1\in \cO(\cT_G), f_0\in \cO(\cT_{L^0})$ such that $(\vartheta_G-df_1)|_{\cT_T}=(\vartheta_0-df_0)|_{\cT_T}$. The latter is equivalent $f_1|_{\cT_T}-f_0|_{\cT_T}=2\xi+C$ for some constant $C$. Without loss of generality, we may assume $C=0$.

On the other hand, it is easy to see (using that $f_G$ is flat and proper and each Hitchin fiber has arithmetic genus $1$) that $\bR^1 f_{G,*}\cO_{\cM_G}\cong \cO_{\cA_G}\cdot \xi$, where $\xi\in \cO(\cT_T)$ gives a generating section in $\Gamma(\cA_G, \bR^1 f_{G,*}\cO_{\cM_G})\cong H^1(\cM_G, \cO)$, using the affine covering $\{\cT_G, \cT_{L^0}\}$. 
This implies that the sought-for $f_1$ and $f_0$ do \emph{not} exist. 

\end{exam}

\begin{remark}
\item[(i)] Let $G=\SL_2$. It follows from the above that $H^1(\cM_G, \cO)\cong \CC[\xi^2+h^2+h^{-2}]\cdot \xi\cong \bigoplus_{i=0}^{\infty}\CC\cdot \xi^{2i+1}$ as a quotient of $\cO(\cT_T)$, The same holds for $\cM_{\PGL_2}^\circ$. \\

\item[(ii)] Let $G=\PGL_2$. The canonical map $\eta: \Omega^2(\cM_G^\circ)/d\Omega^1(\cM_G^\circ)\rightarrow H^2(\cM_G^\circ, \CC)\cong \CC$ is zero. This is because $H_2(\cM_G^\circ, \CC)$ is generated by any Hitchin fiber (by looking at the SYZ fibration), and any algebraic 2-form vanishes on the Hitchin fibers.
On the other hand, using the spectral sequence associated to the holomorphic de Rham complex and the affine covering $\{\cT_G, \cT_{L^0}\}$ and using (i) and $H^1(\cM_G^\circ,\CC)=0$ (cf. Lemma \ref{lemmapi1cMG} below), one gets an isomorphism $H^1(\cM_G^\circ, \cO)\overset{\sim}{\to} \Omega^2(\cM_G^\circ)/d\Omega^1(\cM_G^\circ)$. Its inverse map is given by 
\begin{align*}
&\Omega^2(\cM_G^\circ)/d\Omega^1(\cM_G^\circ)\longrightarrow H^1(\cM_G^\circ, \cO)=\cO(\cT_{T})/d_{\check{\text{C}}\text{ech}}\left(\cO(\cT_{G})\oplus \cO(\cT_{L^0})\right)\cong \bigoplus_{i=0}^{\infty}\CC\cdot \xi^{2i+1}\\
&(d\varphi_0, d\varphi_1)\in \Omega^2(\cT_{L^0})\oplus \Omega^2(\cT_{G}) |d\varphi_1-d\varphi_0=0\text{ on }\cT_T \mapsto \int(\varphi_1-\varphi_0).  
\end{align*}
Note that the morphism is well defined since $\varphi_1-\varphi_0$ is an exact 1-form on $\cT_T$ (indeed, taking any 1-cycle $C$ in $\cT_T$, since it is null homologous both in $\cT_G$ and $\cT_{L^0}$, by Stokes theorem and $\eta=0$, we get $\int_C(\varphi_1-\varphi_0)=0$). By Example \ref{exam: omegaSL2}, $[\omega]_{\alg}$ is mapped to $2\xi$. It is direct to see that the same holds for $\cM_{\SL_2}$. \\

\item[(iii)] We can consider a deformation of $\cM_G, G=\SL_2$ by changing the map (\ref{exam: axyprime}) to $(a', x', y') \mapsto (-(\xi+\epsilon)^2-h^2, h^{-1}, -h^{-1}(\xi+\epsilon))$ for some $\epsilon\in \CC$, which amounts to introducing a shift in the cotangent fibers when gluing $\cT_G$ and $\cT_{L^0}$ along $\cT_T$. Denote the new space by $\cM_G^\epsilon$. We still have the holomorphic symplectic form $\omega$ on $\cM_G^\epsilon$. However, its image in $H^2(\cM_G^\epsilon, \CC)$ is not zero (assuming $\epsilon\neq 0$), i.e. it is not exact as a smooth 2-form. This means that $\cM_G^\epsilon$ is not algebraically isomorphic to $\cM_G$, although they are diffeomorphic. 
\end{remark}

\sss{Convention on Liouville 1-forms} 
From now on, by default a Liouville 1-form on an exact symplectic manifold is always a \emph{real} Liouville 1-form (unless otherwise specified) even though the symplectic manifold might be holomorphic symplectic. 
For a cotangent bundle $T^*X$, use $\vartheta_{T^*X}=-\bp d\bq$ to denote the standard Liouville 1-form on $T^*X$, with any (local) Darboux coordinates $(\bq, \bp)$.

\sss{A tropical manifold structure on the moment polytope $\triangle$}\label{sss: set-up, triangle}
Let $\triangle\subset \XX^*(T)\otimes_\ZZ \RR$ be the moment polytope with respect to an ample $T$-equivariant line bundle on  $\PP$. 
Let $\triangle_J, J\subset_{ft}\wt{I}$, be (the interior of) the face that is normal to $\alpha_j^\vee, j\in J$ (equivalently, it is the interior of the moment image of $\PP_J$). 
Composing with the moment map $\mu: \PP\rightarrow \triangle$, we get a real morphism 
\begin{align}\label{eq:pi_triangle}
\pi_\triangle: \cM_G^\circ\longrightarrow \triangle. 
\end{align}
For each $J$, the fundamental weights for the Lie algebra $\frl_J^\der$ determines a canonical map $\pi'_J: N_J^-\backslash L_J/N_J\longrightarrow \AA^J\sslash \pi_1(L_J^\ad)\times Z(L_J)/Z(L_J^\der)$ (which is the composition of $c_J$ in \eqref{cJ} and the quotient of $U_J$ by $Z(L_J^\der)$). Let $Z(L_J)_{0,\RR}$ be the real form of $Z(L_J)_0=Z(L_J)^\circ$.  
We fix a homeomorphism $\nu_J: \triangle_{\subset J}:=\bigcup_{J^\dag\subset J}\triangle_{J^\dag}\overset{\sim}{\to} \RR^J_{\geq 0}\times Z(L_J)_{0,\RR}$, so that the composition $\nu_J\circ\pi_\triangle|_{\cT_{L_J}}$ is equal to 
\begin{align}\label{eq: b_J, RR}
b_J^\RR: \cT_{L_J}\longrightarrow N_J^-\backslash L_J/N_J\overset{\pi_J'}{\longrightarrow}\AA^J\sslash \pi_1(L_J^\ad)\times Z(L_J)/Z(L_J^\der)\overline{}{\longrightarrow} \RR^J_{\geqsl 0}\times Z(L_J)_{0,\RR}, 
\end{align}
where the last map is the standard one that takes norm on each factor of $\AA^J$ (which descends to $\AA^J\sslash \pi_1(L_J^\ad)$) and projects $Z(L_J)/Z(L_J^\der)$ onto  $Z(L_J)_{0,\RR}$. Similarly, let $b_{J,\der}^{\RR}: \cT_{L_J^\der}\rightarrow \RR_{\geqslant 0}^{J}$ be the projection $b_J^\RR$ in the factor $\cT_{L_J^\der}$.
Recall that $b_{J,\der}^{\RR}$ is equivariant with respect to the canonical $\RR_+$-action on $\cT_{L_J^\der}$ (that scales the symplectic form by weight $1$) and the $\RR_+$-action on $\RR_{\geqsl 0}^J$ that scales the $j$-th factor by $-\omega_j^J(\sfh_0^J)<0$, where $\omega_j^J$ is the fundamental weight with respect to the root system defined by $J$ and $\sfh_0^J$ is the sum of all positive coroots for the same root system. Let $\sfZ_{r,J}$ denote the vector field on $\RR_{\geqsl 0}^J$ from differentiating the weighted $\RR_+$-action, which is the same as $\sfZ_{\bm}$ for $\bm=(\omega_j^J(\sfh_0^J))_{j\in J}$ as in Definition \ref{defngoodadm}. 

We further identify $Z(L_J)_{0,\RR}$ with $\frz_{J,\RR}$ through the logarithmic map and fix a linear coordinate $\bq^{J}=(q^J_1,\cdots, q^J_{r-|J|})$ on $\frz_{J,\RR}$, so that now we view $\nu_J$ as an identification of manifolds with corners
\begin{align}\label{eqnuJr}
\nu_J: \triangle_{\subset J}\overset{\sim}{\longrightarrow} \RR_{\geqsl 0}^J\times \RR^{r-|J|}. 
\end{align}

Note that for any $J_1\subset J_2$, $\nu_{J_2}\circ \nu_{J_1}^{-1}$ is of the form (\ref{eq: nu_beta,alpha,compose}) (see \cite[Lemma 3.3]{J} for an explicit formula). Then $\{\nu_J, J\subset_{ft}\wt{I}\}$ gives the charts of a tropical manifold (cf. Definition \ref{def: tropicalmfld}) with corners structure on $\triangle$.

\sss{An admissible open covering of $\triangle$ and restrictions of standard Liouville 1-forms}\label{subsubsec: adm,triang}

\begin{figure}[htbp]
\begin{tikzpicture}
\coordinate (A1) at (0, 3);
\coordinate (A2) at (-1.73*1.5, -1.5);
\coordinate (A2) at (1.73*1.5, -1.5);
\coordinate (MP2) at ({-1.73*1.5*13/42}, {-1.5*13/42+3*29/42});
\coordinate (MPP2) at ({1.73*1.5*13/42}, {-1.5*13/42+3*29/42});
\coordinate (MQ2) at ({-1.73*1.5*13/42+0.8*1.73/2}, {-1.5*13/42+3*29/42-0.8*0.5});
\coordinate (MQ2prime) at ({-1.73*1.5*13/42+0.6*1.73/2}, {-1.5*13/42+3*29/42-0.6*0.5});
\coordinate (MQQ2) at ({1.73*1.5*13/42-0.8*1.73/2}, {-1.5*13/42+3*29/42-0.8*0.5});
\coordinate (MQQ2prime) at ({1.73*1.5*13/42-0.6*1.73/2}, {-1.5*13/42+3*29/42-0.6*0.5});
\coordinate (P1) at ({-1.73*1.5*5/7}, {-1.5*5/7+3*2/7});
\coordinate (P2) at ({-1.73*1.5*2/7}, {-1.5*2/7+3*5/7});
\coordinate (Q1) at ({-1.73*1.5*5/7+0.6*1.73/2}, {-1.5*5/7+3*2/7-0.6*0.5});
\coordinate (Q2) at ({-1.73*1.5*2/7+0.6*1.73/2}, {-1.5*2/7+3*5/7-0.6*0.5});
\coordinate (PP1) at ({1.73*1.5*5/7}, {-1.5*5/7+3*2/7});
\coordinate (PP2) at ({1.73*1.5*2/7}, {-1.5*2/7+3*5/7});
\coordinate (QQ1) at ({1.73*1.5*5/7-0.6*1.73/2}, {-1.5*5/7+3*2/7-0.6*0.5});
\coordinate (QQ1prime) at ({1.73*1.5*5/7-0.3*1.73/2}, {-1.5*5/7+3*2/7-0.3*0.5});
\coordinate (QQ2) at ({1.73*1.5*2/7-0.6*1.73/2}, {-1.5*2/7+3*5/7-0.6*0.5});
\coordinate (QQ2prime) at ({1.73*1.5*2/7-0.3*1.73/2}, {-1.5*2/7+3*5/7-0.3*0.5});
\coordinate (LQQ1) at ({1.73*1.5*5/7-0.32*1.73/2}, {-1.5*5/7+3*2/7-0.32*0.5});
\coordinate (LQQ2) at ({1.73*1.5*2/7-0.32*1.73/2}, {-1.5*2/7+3*5/7-0.32*0.5});
\coordinate (PPP1) at ({-1.73*1.5*3/7}, {-1.5});
\coordinate (PPP2) at ({1.73*1.5*3/7}, {-1.5});
\coordinate (QQQ1) at ({-1.73*1.5*3/7}, {-1.5+0.6});
\coordinate (QQQ2) at ({1.73*1.5*3/7}, {-1.5+0.6});
\coordinate (L1) at ({-1.73*1.5*2/3}, {-1.5*2/3+3*1/3});
\coordinate (L2) at ({-1.73*1.5*1/3}, {-1.5*1/3+3*2/3});
\coordinate (U1) at ({-1.73*1.5*2/3+0.4*1.73/2}, {-1.5*2/3+3*1/3-0.4*0.5});
\coordinate (U2) at ({-1.73*1.5*1/3+0.4*1.73/2}, {-1.5*1/3+3*2/3-0.4*0.5});
\coordinate (LL1) at ({1.73*1.5*2/3}, {-1.5*2/3+3*1/3});
\coordinate (LL2) at ({1.73*1.5*1/3}, {-1.5*1/3+3*2/3});
\coordinate (UU1) at ({1.73*1.5*2/3-0.4*1.73/2}, {-1.5*2/3+3*1/3-0.4*0.5});
\coordinate (UU2) at ({1.73*1.5*1/3-0.4*1.73/2}, {-1.5*1/3+3*2/3-0.4*0.5});
\coordinate (LLL1) at ({-1.73*1.5*1/3}, {-1.5});
\coordinate (LLL2) at ({1.73*1.5*1/3}, {-1.5});
\coordinate (UUU1) at ({-1.73*1.5*1/3}, {-1.5+0.4});
\coordinate (UUU2) at ({1.73*1.5*1/3}, {-1.5+0.4});
\coordinate (AA1) at (0,3*0.8);
\coordinate (AA2) at ([rotate=120] AA1);
\coordinate (AA3) at ([rotate=240] AA1);
\draw[thick, dashed, cyan] (P1)--(Q1)--(Q2)--(P2);
\draw[thick, dashed, cyan] (PP1)--(QQ1)--(QQ2)--(PP2);
\draw[thick, dashed, cyan] (PPP1)--(QQQ1)--(QQQ2)--(PPP2);
\draw[thick] (1.73*1.5,-1.5) to (0,3);
\draw[thick] (-1.73*1.5, -1.5) to (0,3);
\draw[thick] (-1.73*1.5,-1.5) to  (1.73*1.5, -1.5);
\draw[thick, red!60!black, rounded corners=3pt, dashed] (MP2)--(MQ2)--(MQQ2)--(MPP2);
\draw[rotate=120, thick, red!60!black, rounded corners=3pt, dashed] ([rotate=120] MP2)--([rotate=120]  MQ2)--([rotate=120]  MQQ2)--([rotate=120]  MPP2);
\draw[rotate=240, thick, red!60!black, rounded corners=3pt, dashed] ([rotate=240] MP2)--([rotate=240]  MQ2)--([rotate=240]  MQQ2)--([rotate=240]  MPP2);
\draw[thick, dashed, orange, rounded corners=3pt] (AA1)--(AA2)--(AA3)--(AA1);
\fill[fill=cyan!60, opacity=0.5] (MP2)--(MQ2prime)--([rotate=120]  MQQ2prime)--([rotate=120] MPP2);
\fill[fill=cyan!60, opacity=0.5] (PP1)--(QQ1prime)--(QQ2prime)--(PP2);
\end{tikzpicture}
\caption{An illustration of an admissible covering, with each $\cU_J$ enclosed by a dashed curve. The filled region on the left represents $\cU_{J_1}^\flat$ for some $J_1$. The filled region on the right represents $\cU_{J_2}^\sharp$ for some $J_2$.}\label{figure: triangle_wtI}
\end{figure}
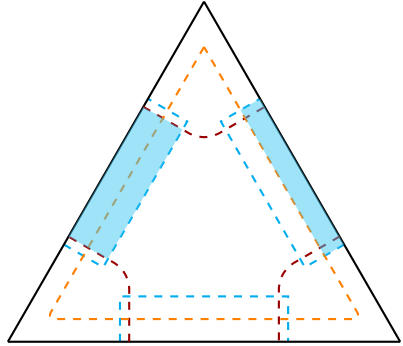

Let $\frU=\{(\cC_J, \cI_J)\}_{J\in \PfI}$ be an admissible open covering of $\triangle$ (cf. Definition \ref{def: adm cover} (i)) with the properties listed in Lemma \ref{lemma: exist admissible} (see Figure \ref{figure: triangle_wtI} for an illustration). Let $\cC\cI_\trg$ be the poset of all admissible open coverings of $\trg$. 

For each $J\in \PfI$ we fix a linear coordinate $\bq^{J}=(q^J_1,\cdots, q^J_{r-|J|})$ on $\frz_{J;\RR}$, which is regarded as a logarithmic coordinate on $Z(L_J)_{0,\RR}$. Let $\bp^{J}=(p^J_1,\cdots, p^J_{r-|J|})$ be the dual coordinate on $T^*Z(L_J)_{0,\RR}$. 
Let $\vartheta_{\cT_{L_{J}^\der}}$, $\vartheta_{T^*Z(L_{J})_{0,\cpt}}$ and $\vartheta_{T^*Z(L_{J})_{0,\RR}}$ be the standard Liouville 1-forms on $\cT_{L_{J}^\der}$, $T^*Z(L_{J})_{0,\cpt}$ and $T^*Z(L_{J})_{0,\RR}$, respectively. Let $\vartheta_{\cT_{L_J}}$ denote the standard Liouville 1-form  on $\cT_{L_J}$, i.e. it is the descent of $\vartheta_{\cT_{L_{J}^\der}}+\vartheta_{T^*Z(L_{J})_{0,\cpt}}+\vartheta_{T^*Z(L_{J})_{0,\RR}}$ through the finite quotient map $\cT_{L_{J}^\der}\times T^*Z(L_J)_0\to \cT_{L_J}$. 

Recall that for any $J_1\subsetneq J_2$, under the canonical open embedding $\cT_{L_{J_1}}\hookrightarrow \cT_{L_{J_2}}$, 
\begin{align}\label{eqvarthetacTJ2J1}
\vartheta_{\cT_{L_{J_2}}}|_{\cT_{L_{J_1}}}=\vartheta_{\cT_{L_{J_1}}}+d\left(\vartheta_{T^*Z(L_{J_1}^{J_2})_{0,\RR}}(V^Z_{J_2, J_1})\right),
\end{align}
where (1) $V^Z_{J_2, J_1}$ is the constant translation vector field on $Z(L_{J_1}^{J_2})_{0,\RR}$ defined by 
\begin{align*}
V^Z_{J_2, J_1}=\sfh_0^{J_2}-\sfh_0^{J_1}, 
\end{align*}
whose natural lifting to $T^*Z(L_{J_1}^{J_2})_{0,\RR}$ is denoted by the same notation; (2) the sum is with respect to the splitting 
\begin{align*}
\cT_{L_{J_1}^\der}\times^{Z(L_{J_1}^\der)}T^*Z(L_{J_1})_{\cpt}\times T^*Z(L_{J_1}^{J_2})_{0,\RR}\times T^*Z(L_{J_2})_{0,\RR},
\end{align*}
where $Z(L_{J_1}^{J_2}):=Z(L_{J_1})\cap L_{J_2}^\der$. Let $V^\frz_{J_2, J_1}$ denote the constant translation vector field on $\frz_{J_1;\RR}^{J_2}:=\frz_{J_1;\RR}\cap \frl_{J_2}^\der$.

\subsection{Admissible vector fields $\sfZ_\flat$ on $\trg$}\label{sss: constr vartheta}

Using Lemma \ref{lemma: exist admissible} and its proof, there exist sufficiently thin and wide admissible open coverings $\{(\cC_J, \cI_J)\}_{J\in \cP_{ft}(\wt{I})}$ and $\{(\cC_J, \cI_J')\}_{J\in \cP_{ft}(\wt{I})}$  of $\triangle$ both satisfying the conditions in that lemma, such that $\cI_J'\Supset \cI_J$ and $\cI_J'$ is sufficiently large relative to $\cI_J$.

Let $\cP_{ft}(\wt{I})$ be the poset of subsets of $\wt{I}$ of finite type. For any $k\in \ZZ$, let $\wt{I}_{>k}$ (resp. $\wt{I}_{\geqsl k}, \wt{I}_k$, etc.) be the subposet of $\cP_{ft}(\wt{I})$ consisting of $J$ with $|J|>k$ (resp. $|J|\geqsl k,\ |J|=k$, etc.). For any subset $A$ in a topological space $X$, let $\Nb(A)$ denote a sufficiently small (but not specified) open neighborhood of $A$.

\sss{Definitions of admissible vector fields $\sfZ_\flat$ on $\trg$}
Let us make the following two definitions of admissibility/rigidity about vector fields $\sfZ_\flat$ on $\trg$. These will be used for defining the desired Liouville 1-forms on $\cM_G$ later. 
In the following, by a \emph{good} admissible open covering of $\trg$, we mean a good admissible open covering with respect to $\sfZ_{r,J}, J\in \cP_{ft}(\wt{I})$, as in Definition \ref{defngoodadm}. Recall the notations from \S\ref{ssscUflatsharp}.

\begin{defn}\label{defnZsemiadm}
Let $\sfZ_\flat$ be a smooth vector field on the tropical manifold $\trg$. For any $J\in \cP_{ft}(\wt{I})$, let $\sfZ_{J,\flat}^\frz$ denote $(\nu_J)_*(\sfZ_{\flat}|_{\trg_J})$, which is a smooth vector field on $\frz_{J;\RR}$. 
Let $\frU=\{(\cC_J, \cI_J)\}_{J\in\PfI}\prec\frU'=\{(\cC'_J, \cI_J')\}_{J\in\PfI}$ be good admissible coverings, with $\cI_J'\Supset \cI_J$. For any $\delta>0$, we say $\sfZ_\flat$ is \emph{$(\delta,\frU')$-semi-admissible} if the following hold
\begin{itemize}
\item[(ZSAa)] for each $J\in \cP_{ft}(\wt{I})$, 
\begin{align*}
(\nu_J)_*(\sfZ_{\flat}|_{(\cU_J')^\sharp})=\sfZ_{r,J}\times \sfZ_{J,\flat}^{\fz},
\end{align*}
with respect to the splitting $\nu_J((\cU_J')^\sharp)\subset \nu_J(\cU_J')=\cC_J\times\cI_J'$.

\item[(ZSAb)] Let $\sfI_J=\vartheta_{T^*\frz_{J;\RR}}(\sfZ_{J,\flat}^\frz)$. Then the Hamiltonian vector field of $\sfI_J$ satisfies 
\begin{align}\label{eqXsfIJ}
|X_{\sfI_J}(\|\bp^J\|^2)|\leqsl \delta\|\bp^J\|^2,
\end{align}
where the norm is induced from the Killing form. 
\end{itemize}

We say $\sfZ_\flat$ is \emph{$\frU$-rigid} if it  satisfies condition (ZR)  which consists of two sub-conditions:  
\begin{itemize}
\item[(ZRi)] $\sfZ_{J, \flat}^{\fz}|_{\Nb(\ovl{\cI}_J)}=0$;
\item[(ZRii)]  for any $J'\supsetneq J$, on $\nu_J\left(\Nb(\ovl{\cU}_{J'}\cap \trg_{J})\right)\subset \fz_{J;\RR}$, the vector field $\sfZ^\frz_{J, \flat}$ is the (restriction of the) pullback (i.e. horizontal lifting) of a smooth vector field $\sfZ_{J, \flat}^{J', \fz}$ on $\frz_{J;\RR}^{J'}$ along the obvious projection
\begin{align}\label{eqfrzJRprojJJprime}
\frz_{J;\RR}=\frz_{J;\RR}^{J'}\times \frz_{J';\RR}\to \frz_{J;\RR}^{J'}. 
\end{align}
\end{itemize}
We say $\sfZ_\flat$ is \emph{weakly $\frU$-rigid} if it  satisfies condition (ZR') which consists of (ZRi) and 
\begin{itemize}
\item[(ZRii')] for any $J'\supsetneq J$, on $\nu_J\left(\Nb(\ovl{\cU}_{J'}\cap \trg_{J})\right)\subset \fz_{J;\RR}$, the vector field $\sfZ^\frz_{J, \flat}$ has zero component in the $\frz_{J';\RR}$-factor under the splitting of tangent spaces with respect to 
$\frz_{J;\RR}=\frz_{J;\RR}^{J'}\times \frz_{J';\RR}.$
\end{itemize}

For any $(\delta,\frU')$-semi-admissible vector field $\sfZ_\flat$ on $\trg$ that is $\frU$-rigid, we say $\sfZ_\flat$ is \emph{$(\delta, \frU', \frU)$-rigid}. We also say  $\sfZ_\flat$ is \emph{$(\delta, \frU)$-rigid} if there exists $\frU'\succ \frU$ making $\sfZ_\flat$ $(\delta, \frU', \frU)$-rigid. Similar convention applies to weakly $\frU$-rigid vector fields. 
\end{defn}

\begin{remark}\label{rmksplitting}
\begin{itemize}
\item[(i)]
Under the assumption of Definition \ref{defnZsemiadm}, observe that $\cU_{J'}=\bigcup_{J\subset J'}\cU_{J'}\cap (\cU_J')^\sharp$ and $\cU_{J'}\cap (\cU_J')^\sharp\subset \cU_{J'}\cap \cU_J'=\nu_J^{-1}(\cC_J\times \wt{\cI}_{J',J})$, for some open subset $\wt{\cI}_{J',J}\subset \frz_{J;\RR}$. Then (ZSAa) and (ZRii) for $J\subsetneq J'$ implies that on $\cU_{J'}\cap (\cU_J')^\sharp$, 
\begin{align*}
(\nu_J)_*(\sfZ_\flat|_{\cU_{J'}\cap (\cU_J')^\sharp})=\sfZ_{r,J}\times \sfZ_{J, \flat}^{J',\frz}\times 0
\end{align*}
with respect to the splitting $\RR_{\geqsl 0}^{J}\times \frz_{J;\RR}^{J'}\times \frz_{J';\RR}$. Combining with (ZRi), this implies that on the entire $\cU_{J'}$, 
\begin{align*}
(\nu_{J'})_*(\sfZ_\flat|_{\cU_{J'}})=\sfZ_{r,J'}^\dm\times 0
\end{align*}
for some vector field $\sfZ_{r,J'}^\dm$ on $\Nb(\ovl{\cC}_J)\subset \RR_{\geqsl 0}^J$, with respect to the splitting $\RR_{\geqsl 0}^{J'}\times \frz_{J';\RR}$. 

\item[(ii)] Assume $\sfZ_\flat$ is $(\delta, \frU', \frU)$-rigid. 
Using $\cU_{J'}=\bigcup_{J\subset J}\cU_{J'}\cap\cU_J$ and $\partial^\circ\cC_{J'}$ is $\cC$-adapted to $\{(\cC_{J}, \cI_J)\}_{J\subsetneq J}$, the observation in (i) implies that $\sfZ_{r,J'}^{\dm}$ is stratified on $\ovl{\cC}_{J'}$.

\item[(iii)] For $\sfZ_\flat$ satisfying (ZSAa) and (ZRii'), by a similar observation as in (i) and (ii), we see that $\sfZ_{\flat}|_{\ovl{\cU}_J}$ is stratified for all $J\in \cP_{ft}(\wt{I})$. 

\item[(iv)] An identity that we will frequently use is that for any $J\subsetneq J'$, 
\begin{align*}
(\nu_{J, J'}^{-1})_*\left((\sfZ_{r,J'}\times 0)|_{\Imagine(\nu_{J,J'})}\right)=Z_{r, J}\times V_{J', J}^{\frz}\times 0,
\end{align*}
where the left-hand-side (resp. right-hand-side) is with respect to the splitting $\RR_{\geqsl 0}^{J'}\times\frz_{J';\RR}$ (resp. $\RR_{\geqsl 0}^{J}\times\frz_{J;\RR}^{J'}\times\frz_{J';\RR}$). 

\item[(v)] For a maximal $J\in \cP_{ft}(\wt{I})$, only changing $\cC_J$ has \emph{no} effect on semi-admissibility or rigidity condition on vector fields. This is due to $(\cU_J')^\sharp$ is independent of $\cC_J$ for a maximal $J$ (when other members in the admissible covering stay the same), and in $\trg_{\subset J}$, every vector field $\sfZ_\flat$ automatically satisfies (ZRi) and (ZRii) for a maximal $J$. 

\end{itemize}

\end{remark}

\begin{defn}\label{defnZadm}
Let $\sfZ_\flat$ be a vector field on the tropical manifold $\trg$. 
Let $\frU=\{(\cC_J, \cI_J)\}_{J\subsetneq\wt{I}}$ be a good admissible covering. We say $\sfZ_{\flat}$ is \emph{strongly $(\delta,\frU)$-admissible} if the following hold
\begin{itemize}
\item[(ZAa)]  for each $J\in \cP_{ft}(\wt{I})$, 
\begin{align*}
(\nu_J)_*(\sfZ_{\flat}|_{\cU_J})=\sfZ_{r,J}\times \sfZ_{J,\flat}^{\frz},
\end{align*}
with respect to the splitting $\nu_J(\cU_J)=\cC_J\times\cI_J$, where $\sfZ_{J,\flat}^{\frz}$ is convex along $\partial \cI_J$ (cf. \S\ref{sssMWC}) and  $\sfZ_{J,\flat}^{\frz}$ has no zero outside $\cI_J$;

\item[(ZAb)] This is the same as (ZSAb) in Definition \ref{defnZsemiadm};

\item[(ZAc)] $\sfZ_{J,\flat}^\frz$ is gradient-like with respect to some exhausting Morse function $f_J^{\frz}: \frz_{J;\RR}\to [0,\infty)$, i.e. $f_J^{\frz}$ is proper and Morse, and there exists a function $\ep_J:\frz_{J;\RR}\to \RR_{>0}$ such that 
\begin{align*}
\sfZ_{J,\flat}^{\frz}f_J^{\frz}\geqsl \ep_J(|\sfZ_{J,\flat}|^2+|df_J^{\frz}|^2)
\end{align*}
with respect to some Riemannian metric on $\frz_{J;\RR}$. 
\end{itemize}

We say $\sfZ_\flat$ is \emph{$(\delta,\frU)$-admissible} if it satisfies (ZAa) and (ZAb). We also say $\sfZ_\flat$ is \emph{$\delta$-admissible} (resp. $\delta$-\emph{rigid}, \emph{strongly $\delta$-admissible}) if it is $(\delta,\frU)$-admissible (resp. $(\delta, \frU',\frU)$-rigid, 
strongly $(\delta,\frU)$-admissible) for some good admissible covering $\frU'$ and $\frU$ as above.

\end{defn}

Note that in the above definition, $\sfZ_{r;J}$  is automatically concave along $\partial^\circ\cC_J$ due to the good condition on $\frU$. 
Given a $\delta$-admissible (resp. strongly $\delta$-admissible) $\sfZ_\flat$ on $\trg$, let $\CI_{\sfZ_\flat,\delta}$ (resp. $\CI_{\sfZ_\flat,\delta}^{\sigma}$) be the collection of good admissible coverings $\frU$ with respect to which $\sfZ_\flat$ is $(\delta,\frU)$-admissible (resp. strongly $(\delta,\frU)$-admissible). Similarly, for a $(\delta,\frU)$-rigid $\sfZ_\flat$, let $\cC\cI\frU_{\sfZ_\flat,\delta}$ be the collection of $\frU'$ making $\sfZ_\flat$ $(\delta, \frU', \frU)$-rigid.

\subsection{Existence results about admissible $\sfZ_\flat$}\label{ssexisZflat}

\begin{lemma}\label{lemmasemitoadm}
\begin{itemize}
\item[(i)] Let $\frU^1\succ \frU^2$ be admissible open coverings of $\trg$. Then $\frU^2$-semi-admissible $\Longrightarrow$ $\frU^1$-semi-admissible. Hence for any  $(\delta,\frU)$-rigid $\sfZ_\flat$,    $\cC\cI\frU_{\sfZ_\flat,\delta}\subset \cC\cI_\trg$ is cofinal.  

\item[(ii)] Let $\frU^1\succ\frU^2$ be admissible open coverings of $\trg$. Then $(\delta, \frU^1)$-rigidity $\Longrightarrow$ $(\delta,\frU^2)$-rigidity.

\item[(iii)] For any $(\delta, \frU)$-semi-admissible $\sfZ_\flat$ and any $\frU'\succ \frU$, there exists a good $\frU''\succ \frU'$ such that $\sfZ_\flat$ is $(\delta, \frU'')$-admissible. 

\end{itemize}
\end{lemma}

\begin{proof}
(i) This follows from $(\cU_J^1)^\sharp\subset \bigcup_{J'\supset J}(\cU_{J'}^2)^\sharp$, for any $J\in\PfI$ (cf. Lemma \ref{lemmaAfrU12sharpflat} (i)). 

(ii) This follows from $\cU_J^2\subset \bigcup_{J^\dagg\subset J}\cU_J^1$ (cf. Lemma \ref{lemmaAfrU12sharpflat} (iii)) and Remark \ref{rmksplitting} (i).

(iii) We define $\frU''_{>k}:=\{(\cC_J'',\cI_J'')\}_{J\in \wt{I}_{>k}}$ inductively (in the decreasing order of $k$). Since the condition (ZSAb)=(ZAb) is independent of $\frU$, we only need to check (ZAa) for $\frU''_{>k}$. 
First, for any $J\in \wt{I}_{r}$, choose any sufficiently small good $\cC_J''\subset\cC_J'$ with respect to $\sfZ_{r,J}$. Since $\nu_J^{-1}(\cC_J'')\subset \cU_J^\sharp$, (ZSAa) for $\cU_J^\sharp$ implies (ZAa) for $\cU''_J$. 

Suppose we have defined $\frU''_{>k}$ such that (1) $\frU''_{>k}$ is a good admissible open covering of $\cU''_{>k}:=\bigcup_{J\in \wt{I}_{>k}}\cU''_J$---which gives a tubular neighborhood of $\trg_{>k}:=\bigcup_{J\in \wt{I}_{>k}}\trg_J$,  
with $\cC_J''\subset \cC_J'$ and $\cI_J''\supset \cI_J'$;
(2) (ZAa) holds for each $\cU''_J$. Then for any $J\in \wt{I}_k$, the intersection $\cU_{>k}''\cap \trg_J$ satisfies that 
\begin{itemize}
\item $\cU_{>k}''\cap \ovl{\trg_J}$ inherits an admissible covering of a tubular neighborhood of $\partial \trg_J$ in $\ovl{\trg_J}$;
\item $\sfZ_{J,\flat}^\frz$ is concave along $\partial(\nu_J (\cU_{>k}''\cap \trg_J))\subset \frz_{J;\RR}$. 
\end{itemize}

It is then easy to see that there exist arbitrarily large $\cI''_J\supset \frz_{J;\RR}-\ovl{\nu_J (\cU_{>k}''\cap \trg_J)}$ such that 
 \begin{itemize}
 \item $\sfZ_{J,\flat}^\frz$ is convex along $\partial \cI''_J$; 
 \item $(0, \cI''_J)$ together with the induced admissible covering of $\cU_{>k}''\cap \ovl{\trg_J}$ gives an admissible covering of $\ovl{\trg_J}$. 
\end{itemize}
Choose any good sufficiently small $\cC_J''$ such that 
$\cU_J''=\nu_{J}^{-1}(\cC_J''\times \cI_J'')\subset \bigcup_{J'\supset J}\cU_{J'}^\sharp$. Then 
(ZSAa) for $\frU$ implies (ZAa) for $\cU''_J$, and  $\frU_{>k-1}$ is a good admissible open covering of $\cU''_{>k-1}$---which gives a tubular neighborhood of $\trg_{>k-1}$. Thus, the inductive step follows. 
\end{proof}

\begin{remark}
Similarly to the definition of $\delta$-admissibility, we can define $\delta$-semi-admissibility by requiring that $\sfZ_\flat$ is $(\delta,\frU)$-semi-admissible with respect to some good admissible $\frU$. But Lemma \ref{lemmasemitoadm} shows that  $\delta$-semi-admissible is equivalent to $\delta$-admissible. 
\end{remark}

\begin{lemma}\label{lemmarigidZ}
Let $\delta>0$. For sufficiently thin and wide good admissible open coverings $\frU=\{(\cC_J, \cI_J)\}_{J\subsetneq\wt{I}}\prec\frU'=\{(\cC_J, \cI_J')\}_{J\subsetneq\wt{I}}$, i.e. $\cC'_J=\cC_J$, with $\cI_J'\Supset \cI_J$ and $\cI_J'$ is sufficiently large relative to $\cI_J$, there exists a $(\delta,\frU',\frU)$-rigid vector field $\sfZ_\flat$ on $\trg$. 
\end{lemma}

\begin{proof}
For any subposet $A\subset \cP_{ft}(\wt{I})$ and $J\in A$, Let 
\begin{align*}
\cU_J^{\ddagger; A}=\cU_J'-\nu_J^{-1}\left(\Nb(\partial^\circ C_J)\times \Nb(\bigcup_{A\ni J'\supsetneq J}\ovl{\cI'_{J', J}})\right)-\bigcup_{A\ni J^\dagg\subsetneq J}\nu_{J^\dagg}^{-1}\left(\Nb(\partial^\circ C_{J^\dagg})\times \Nb(\ovl{\cI'_{J, J^\dagg}})\right). 
\end{align*}
Here $\cI'_{J', J}$ is as in Lemma \ref{lemma: admissible, cI contractible} for $\frU'$.
Let $\cU^{\ddagger; >k}_J=\cU^{\ddagger; \wt{I}_{>k}}_J$. 
Let $\cU_{>k}^\ddagger=\bigcup_{J\in \wt{I}_{>k}}\cU^{\ddagger; >k}_J$. Let $\cU_{>k}=\bigcup_{J\in \wt{I}_{>k}}\cU_{J}=\bigcup_{J\in \wt{I}_{>k}}\cU'_{J}$. 
By the sufficiently thin and wide condition on $\frU$ and $\frU'$, we can assume that for any $J\in \wt{I}_k$,  
\begin{align*}
\nu_J^{-1}(\cI_J)\supset \ovl{\cU_{>k}-\cU_{>k}^{\ddagger}}\cap \trg_J.
\end{align*}
Note that $\cU^{\ddagger}_{>-1}=\trg$. 
We construct a $(\delta,\frU')$-semi-admissible vector field $\sfZ_\flat$ that is $\frU$-rigid, by the following inductive steps on $-1\leqsl k\leqsl r-1$ in the decreasing order:

\begin{itemize}

\item[(i)] First, when $k=r$, for each vertex of $\triangle$ indexed by a maximal $J$, let  $\sfZ_\flat^{>r-1}|_{\cU_J'}=(\nu_J^{-1})_*\sfZ_{r,J}$, which clearly satisfies all the conditions for the subposet $J\in \wt{I}_{>r-1}$ in (ZSAa), (ZSAb) and (ZR). 

\item[(ii)] Given $0\leqsl k\leqsl r-1$, suppose we have defined $\sfZ_\flat^{>k}$ on 
$\cU^{\ddagger}_{>k}$, so that it satisfies (ZSAa), (ZSAb) and (ZR) for the subposet $\wt{I}_{>k}$, where $(\cU_J')^\sharp$ are replaced by $(\cU_J')^{\sharp, \ddagger;>k}:=(\cU_J')^{\sharp;\wt{I}_{>k}}\cap \cU^\ddagger_{>k}$ and we denote the terms corresponding to $\sfZ_{J,\flat}^{\frz}$ (resp. $\sfZ_{J,\flat}^{J', \frz}$) by $\sfZ_{J,\flat}^{\frz, >k}$ (resp. $\sfZ_{J,\flat}^{J', \frz,>k}$).

\item[(iii)] 
Now for any $J\in \wt{I}_{k}$ and any $J'\supsetneq J$, the restricted vector field $(\sfZ_\flat^{>k})_{J', J}:=\sfZ_{\flat}^{>k}|_{(\cU'_{J'})^{\sharp, \ddagger;>k}\cap \cU'_J}$ satisfies that 
\begin{align*}
(\nu_J)_*(\sfZ_\flat^{>k})_{J',J}=\sfZ_{r;J}\times V_{J', J}^\frz \times\sfZ_{J',\flat}^{\frz,>k}
\end{align*}
under the splitting $\RR_{\geqsl 0}^J\times \frz_{J;\RR}^{J'}\times \frz_{J';\RR}$. Since $V_{J', J}^\frz$ is a constant vector field, (ZSAb) holds for $\sfZ_{J',\flat}^{\frz, >k}$ implies that it holds for $V_{J', J}^\frz \times \sfZ_{J',\flat}^{\frz}$ on $T^*\left(\nu_J\left((\cU'_{J'})^{\sharp, \ddagger;>k}\cap \trg_J\right)\right)$ as well. Write 
$\nu_J\left(\cU'_J\cap \cU_{J'}^{\ddagger,>k}\right)=\cC_J\times \cI_{J',J}^{\ddagger, >k}$, and 
$\nu_J\left(\cU'_J\cap \cU_{>k}^{\ddagger}\right)=\cC_J\times \cI^{\ddagger,>k}_{J}$. Then the restriction of $\sfZ_\flat^{>k}$ on $\cU'_J\cap \cU_{>k}^{\ddagger}$  satisfies the splitting condition (ZSAa) in which we denote the vector field on the factor $\cI^{\ddagger,>k}_{J}$ by $\sfZ_{J,\flat}^{\frz, >k}$. 
It is easy to verify that $\sfZ_{J,\flat}^{\frz, >k}$ satisfies (ZR) with $\ovl{\cU}_{J'}$ replaced by $\ovl{\cU_{J'}^{\ddagger,>k}}$ for any $J'\supsetneq J$. 

Choose $\sfZ_{J,\flat}^{\frz, >k-1}$ on $\frz_{J;\RR}$ to be a smooth vector field satisfying 
\begin{itemize}
\item[(1)] $\sfZ_{J,\flat}^{\frz, >k-1}|_{\Nb(\ovl{\cI_J})}=0$ and $\sfZ_{J,\flat}^{\frz, >k-1}|_{\Nb(\frz_{J;\RR}-\cI_J')}=\sfZ_{J,\flat}^{\frz, >k}|_{\Nb(\frz_{J;\RR}-\cI_J')}$;

\item[(2)] (ZRii) holds for any $J'\supsetneq J$ with $\ovl{\cU}_{J'}$ replaced by $\ovl{\cU_{J'}^{\ddagger,>k}}$;

\item[(3)] (ZSAb) holds for $\sfZ_{J,\flat}^{\frz, >k-1}$. 
\end{itemize}
\end{itemize}
Such $\sfZ_{J,\flat}^{\frz, >k-1}$ exists since 
\begin{itemize}
\item $\cI_J'$ is sufficiently large relative to $\cI_J$ (so then (1) and (3) can be easily fulfilled);

\item to satisfy (2), we choose $\sfZ_{J,\flat}^{\frz, >k-1}$ on $\cI^{\ddagger,>k}_{J',J}$ inductively in the increasing order of $|J'|$---this is achievable since for $J''\supsetneq J'\supsetneq J$, $\sfZ_{J,\flat}^{\frz, >k-1}|_{\cI^{\ddagger,>k}_{J'',J}}$ has ``more freedom" than $\sfZ_{J,\flat}^{\frz, >k-1}|_{\cI^{\ddagger,>k}_{J',J}}$. 
\end{itemize}

Lastly, define 
\begin{align*}
\sfZ_{\flat}^{>k-1}=\begin{cases}
\sfZ_\flat^{>k},&\text{ on }\cU_{>k}^\ddagger-\bigcup_{J\in \wt{I}_k}\cU'_J;\\
(\nu_J^{-1})_*(\sfZ_{r;J}\times \sfZ_{J,\flat}^{\frz, >k-1}), &\text{ on }\cU_J^{\ddagger, >k-1} \text{ for }J\in \wt{I}_{k}.
\end{cases}
\end{align*}
Then $\sfZ_{\flat}^{>k-1}$ on $\cU_{>k-1}^{\ddagger}$ satisfies (ZSAa), (ZSAb) and (ZR) for the subposet $\wt{I}_{>k-1}$. 
 Note that the two formulas above on the right-hand-side do \emph{not} agree on (a non-empty subset of) $\nu_{J}^{-1}\left(\partial^\circ C_{J}\times \cI^{\ddagger,>k}_{J}\right), J\in \wt{I}_{k}$, which is the reason why we put the region $\cU^{\ddagger, >k-1}_J$ instead of $\cU_J'$ on the second line. This finishes the inductive step. 
\end{proof}

The assumption in Lemma \ref{lemmarigidZ} requires that $\cI_J$ is not ``too large" for $\cC_J$. This is needed in  
applying the induction argument in a convenient way. However, the following corollary (as a direct consequence of the lemma) says we do \emph{not} need this assumption on $\frU$ to admit a $(\delta, \frU)$-rigid vector field.

\begin{cor}\label{corcofinalrigid}
For any admissible open covering $\frU^1$ of $\trg$ and $\delta>0$, there exist a $(\delta,\frU^1)$-rigid vector field $\sfZ_\flat$ on $\trg$. 
\end{cor}

\begin{proof}
For any $\frU^1$, there exists $\frU'\succ\frU\succ\frU^1$ such that $\frU'$ and $\frU$ satisfy the assumption in Lemma \ref{lemmarigidZ}. The statement then follows from Lemma \ref{lemmasemitoadm} (ii). 
\end{proof}

\begin{lemma}\label{lemmaadmZ}
\begin{itemize}
\item[(i)]
Let $\delta>0$. For any good admissible open covering $\frU^1=\{(\cC_J^1, \cI^1_J)\}_{J\in\PfI}$, there exists a good $\frU\succ \frU^1$ and a strongly $(\delta,\frU)$-admissible vector field $\sfZ_\flat$ on $\trg$. 

\item[(ii)] If $\sfZ_\flat$ is $(\delta,\frU)$-admissible (resp. strongly $(\delta,\frU)$-admissible), then for any admissible $\frU^1\succ \frU$, there exists a good admissible $\frU'\succ \frU^1$ such that $\sfZ_\flat$ is $(\delta,\frU')$-admissible (resp. strongly $(\delta,\frU')$-admissible).
\end{itemize}
\end{lemma}

\begin{proof}
(i) Like in the proof of Lemma \ref{lemmasemitoadm} and Lemma \ref{lemmarigidZ}, we define $\frU_{>k}=\{(\cC_J, \cI_J)\}_{J\in \wt{I}_{>k}}$ and the vector field $\sfZ_\flat^{>k}$ on $\cU_{>k}$ inductively. Since in each inductive step, we will actually need $\sfZ_\flat^{>k}$ on a neighborhood of $\ovl{\cU}_{>k}$, we also choose $\cC'_J\Supset \cC_J$ such that $\frU'_{>k}:=\{(\cC'_J,\cI_J)\}_{J\in\wt{I}_{>k}}$ is a good admissible open covering $\cU_{>k}'$---a tubular neighborhood of $\trg_{>k}$ such that $\frU'_{>k}\succ \frU_{>k}^1$, and define $\sfZ_\flat^{>k}$ on $\cU'_{>k}$. 

First, for $J\in \wt{I}_r$, choose any good $\cC_J\Subset \cC_J'\subset \cC_J^1$ and define $\sfZ_\flat^{>r-1}=\sfZ_{r,J}|_{\cU'_J}$. Suppose we have defined $\frU_{>k}\succ \frU_{>k}'$ and $\sfZ_\flat^{>k}$ on $\cU'_{>k}$ such that $\frU'_{>k}\succ \frU^1_{>k}$, and  (ZAa)--(ZAc) holds for $\sfZ_\flat^{>k}$ with $J\in\wt{I}_{>k}$ and the splitting condition in (ZAa) also holds on $\cU'_J$.  

Similarly to the inductive step in the proof of Lemma \ref{lemmasemitoadm}, for any $J\in \wt{I}_k$, we have $\sfZ_{J,\flat}^{\frz,>k}:=(\nu_{J})_*(\sfZ_{\flat}^{>k}|_{\cU'_{>k}\cap \trg_J})$ satisfies that it is concave along $\partial(\nu_J(\cU_{>k}\cap \trg_J))$. Hence we  
can choose arbitrarily large $\cI_J$ satisfying 
\begin{itemize}
 \item $\sfZ_{J,\flat}^{\frz,>k}$ is convex along $\partial \cI_J$; 
 \item the pair $(0, \cI_J)$ together with the induced admissible covering of $\cU_{>k}\cap \ovl{\trg_J}$ gives an admissible covering of $\ovl{\trg_J}$. 
\end{itemize}
Observe that $\sfZ_{J,\flat}^{\frz,>k}$ has no zero outside $\cI_J$ and for any compact subset $K\supset \cI_J$ in $\frz_{J;\RR}$, there exists $t>0$ such that $\varphi^t_{\sfZ_{J,\flat}^{\frz,>k}}(K-\cI_J)\cap K=\vn$. Using the same observation about the condition (ZSAb) as in (iii) of the proof of Lemma \ref{lemmarigidZ}, we see that for sufficiently thin and wide $\frU'_{>k}$, there exists an extension of $\sfZ_{J,\flat}^{\frz,>k}$ to $\sfZ_{J,\flat}^{\frz}$ on $\frz_{J;\RR}$ and a function $f_J^\frz$ such that (ZAb) and (ZAc) hold. Now choose sufficiently small $\cC_J'\Supset \cC_J$ and sufficiently large $\cI_J$ for $J\in \wt{I}_k$ such that $\frU^1_{>k-1}\prec\frU'_{>k-1}$ is a good admissible open covering $\cU_{>k-1}'$, and 
define
\begin{align*}
&\sfZ_{\flat}^{>k-1}=\begin{cases}(\nu_J^{-1})_*(\sfZ_{r,J}\times \sfZ_{J,\flat}^{\frz}), &\text{ on }\cU_J'\text{ for }J\in \wt{I}_k;\\
\sfZ_{\flat}^{>k}, &\text{ on }\cU_{>k}'. 
\end{cases}
\end{align*}
This is clearly well defined on $\cU'_{>k-1}$ and satisfies (ZAa) as well. This finishes the inductive step. 

(ii) This is obvious from the definition of (strongly) admissibility of $\sfZ_\flat$ and using (iii) in the proof of Lemma \ref{lemmarigidZ} about the splitting condition. 
\end{proof}

\begin{cor}\label{corsfZflatspecial}
Choose any point $\sfc_J\in \trg_J$ for each $J\in \cP_{ft}(\wt{I})$ and any $\delta>0$. There exists a strongly $\delta$-admissible $\sfZ_\flat$ on $\trg$ such that the zero points of $\sfZ_{\flat}$ are exactly $\{\sfc_J, J\in \PfI\}$. We will refer to any such $\sfZ_\flat$ as a \emph{special} strongly $\delta$-admissible vector field on $\trg$. 
\end{cor}
\begin{proof}
We can choose $\sfZ_\flat$ so that in each inductive step in the proof of Lemma \ref{lemmaadmZ} (i), $\sfZ_{J,\flat}^\frz$  has exactly one zero point, which happens at $\frc_S$ and which (necessarily) gives the minimum of $f_J^\frz$. This is achievable as long as $\frU$ is sufficiently thin and wide. 
\end{proof}

\begin{cor}\label{corZflatCI}
Let $\delta>0$. 
\begin{itemize}
\item[(i)] 
Let $\frU, \frU'$ be as in Definition \ref{defnZsemiadm}. Assume $\cC'_J=\cC_J$ and $\cI'_J$ is sufficiently large relative to $\cI_J$, then the space of $(\delta, \frU', \frU)$-rigid vector fields is nonempty and convex, hence contractible. The same holds for the subspace of $\frU$-rigid and weakly $\frU$-rigid vector fields. 

\item[(ii)] The space of $\delta$-admissible (resp. $(\delta, \frU)$-rigid) vector fields is nonempty and convex. 

\item[(iii)] Given a $\delta$-admissible (resp. strongly $\delta$-admissible) $\sfZ_\flat$, the inclusion of the subposet $\CI_{\sfZ_\flat, \delta}$ (resp. $\CI^\sigma_{\sfZ_\flat, \delta}$) into $\CI_\trg$ is cofinal. 
\end{itemize}
\end{cor}

\begin{proof}
These are direct consequences of Lemma \ref{lemmarigidZ}, Corollary \ref{corcofinalrigid} and Lemma \ref{lemmaadmZ} (and the proof of Lemma \ref{lemmaadmZ} for part (ii) on the convexity of $\delta$-admissible vector fields). Alternatively, for part (ii) on the convexity of $\delta$-admissible vector fields, one can use Lemma \ref{lemmasemitoadm} (i) to show the convexity of $\delta$-semi-admissible vector fields and then use the equivalence between $\delta$-semi-admissibility and $\delta$-admissiblility. 
\end{proof}

\subsection{The Liouville 1-form for a semi-admissible $\sfZ_\flat$}\label{subsecLiouvillefromZflat}

\begin{defnlemma}\label{defthetaZflat}
Let $\sfZ_\flat$ be a $(\delta, \frU')$-semi-admissible vector field on $\trg$. Then $\sfZ_\flat$ determines a unique Liouville 1-form $\vartheta$ on $\cM_G$ satisfying 
\begin{align}\label{eqdefnvarthetaJ}
\vartheta|_{\pi_\trg^{-1}((\cU_{J}')^{\sharp})}=\vartheta_{\cT_{L_J^\der}}+\vartheta_{T^*Z(L_J)_{0,\cpt}}+\left(\vartheta_{T^*\cI_J'}+
d\left(\vartheta_{T^*\cI_J'}(\sfZ_{J,\flat}^{\frz})\right)\right),\forall J\in \cP_{ft}(\wt{I}), 
\end{align}
with respect to the natural splitting, in which we canonically identify $\cI_J$ as an open subset of $Z(L_J)_{0;\RR}$. This is independent of the choice of $\frU'$ with respect to which $\sfZ_\flat$ is $(\delta, \frU')$-semi-admissible. 
We say $\vartheta$ is the \emph{canonical Liouville 1-form} for $\sfZ_\flat$.
\end{defnlemma}

\begin{proof}
Let $\vartheta_J$ denote the formula on the right-hand-side of \eqref{eqdefnvarthetaJ}.
Let $(\cU')^\sharp_{\supset J}=\bigcup_{J'\supset J}(\cU_{J'}')^\sharp$. Then $\vartheta_J$ is well defined on $\pi_\trg^{-1}((\cU')^\sharp_{\supset J})\cap \cT_{L_J}$. 
To show $\vartheta$ is well defined, it suffices to show that for any $J_1\subsetneq J_2$, 
\begin{align}\label{eqvarthetaJ2J1cUcT}
&\vartheta_{J_2}|_{\pi_\trg^{-1}((\cU_{J_2}')^{\sharp})\cap \cT_{L_{J_1}}}=\vartheta_{J_1}|_{\pi_\trg^{-1}((\cU_{J_2}')^{\sharp})\cap \cT_{L_{J_1}}}. 
\end{align}
But this follows from \eqref{eqvarthetacTJ2J1} and Remark \ref{rmksplitting} (iv). 

Let $\vartheta^{\frU'}$ denote the $\vartheta$ defined using $\frU'$. 
To show that $\vartheta^{\frU'}$ is independent of the choice of $\frU'$, 
it suffices to show that for any  $\frU''\succ\frU'$, $\vartheta^{\frU''}=\vartheta^{\frU'}$ (e.g. by Lemma \ref{lemmasemitoadm} and that $\CI_{\sfZ_\flat, \delta}$ is filtered). 
Using $(\cU_J'')^{\sharp}\subset (\cU')^\sharp_{\supset J}$ and \eqref{eqvarthetaJ2J1cUcT}, we have 
\begin{align*}
\vartheta^{\frU''}|_{\pi_\trg^{-1}((\cU_J'')^\sharp)}=\vartheta^{\frU'}|_{\pi_\trg^{-1}((\cU_J'')^\sharp)}. 
\end{align*}
Thus $\vartheta^{\frU''}=\vartheta^{\frU'}$.
\end{proof}

\begin{remark}\label{remarkvarthetaadmiss}
\begin{itemize}
\item[(i)] It is clear that if $\sfZ_\flat$ is $(\delta, \frU)$-admissible, then the canonical Liouville 1-form and its Liouville vector field are given by
\begin{align}
\label{eqdefnvarthetaJadm}&\vartheta|_{\pi_\trg^{-1}(\cU_{J})}=\vartheta_{\cT_{L_J^\der}}+\vartheta_{T^*Z(L_J)_{0,\cpt}}+\left(\vartheta_{T^*\cI_J}+
d\left(\vartheta_{T^*\cI_J}(\sfZ_{J,\flat}^{\frz})\right)\right),\forall J\in \cP_{ft}(\wt{I}), \\
\label{eqdefnsfZJadm}&\sfZ|_{\pi_\trg^{-1}(\cU_{J})}=\sfZ_{\cT_{L_J^\der}}+\sfZ_{T^*Z(L_J)_{0, \cpt}}+\left(\sfZ_{T^*\cI_J}-X_{\sfI_J}\right). 
\end{align}

\item[(ii)] We say a Liouville 1-form $\vartheta$ on $\cM_G$ is \emph{$(\delta,\frU')$-semi-admissible} (resp. \emph{$(\delta,\frU')$-admissible}, \emph{$\delta$-admissible}, etc.) if it is the canonical Liouville 1-form for a $(\delta,\frU')$-semi-admissible (resp. $(\delta,\frU')$-admissible, $\delta$-admissible, etc.) vector field $\sfZ_\flat$ on $\trg$. By the formula \eqref{eqdefnvarthetaJ}, it is clear that $\sfZ_{\flat}$ is uniquely determined by such a $\vartheta$. Thus, we say $\sfZ_{\flat}$ is the \emph{canonical vector field}  on $\trg$ for $\vartheta$. In this case, denote $\CI_{\sfZ_\flat, \delta}$ (resp. $\CI_{\sfZ_\flat, \delta}^\sigma$) also by $\CI_{\vartheta, \delta}$ (resp. $\CI_{\vartheta, \delta}^\sigma$). 

\end{itemize}
\end{remark}

\subsection{Weinstein manifold structure(s) on $(\cM_G, \vartheta)$}\label{subsecWeinMfldMG}

Let $0<\delta< 1$. Let $\sfZ_\flat$ be a strongly $\delta$-admissible vector field on $\trg$. Let $\vartheta$ be the canonical Liouville 1-form for $\sfZ_\flat$. 

\sss{An inductive construction of (generalized) Weinstein subdomains of $(\cM_G, \vartheta)$}\label{sssindWeinsteindomain}
Choose any $\frU\in \CI_{\sfZ_\flat, \delta}$. We define for each $0\leqsl k\leqsl r$ a Liouville subdomain $\frD_k$ (whose boundary might have corners) in $\frF$ satisfying:
\begin{itemize}
\item there exists a Morse-Bott function making $\frD_0$ a Weinstein domain, and $\frD_{k+1}$ is obtained from attaching Weinstein handles to $\frD_k$;

\item the Liouville completion of $\frD_{r}$ is $\cM_G$, thus $\cM_G$ has a (generalized) Weinstein handle decomposition.  

\end{itemize}
The construction is very similar to the construction given in \cite{J} of Liouville subdomains in the symplectic hyperpersurface $\frF$ in $\cH$ (cf. \S\ref{subsec: Hsm}). We give the details of the construction in the current setting for the convenience of the reader. 
 
Choose $R_0\gg R_1\gg\cdots\gg R_r\gg 1$. For any $J\in \wt{I}_k$, let 
\begin{align*}
&\frD^{\flat, \ovl{\cI}}_{J, R_k}:=\{(\bq^J, \bp^J)\in \ovl{\cI}_J\times \frz_{J;\RR}^*: \|\bp^J\|\leqsl R_k\}\subset T^*\ovl{\cI}_J,\\
&\frD^{Z,\cpt}_{J,R_k}:=\{(z_J^\cpt, \xi_{J,\perp}^\cpt)\in Z(L_J)_{\cpt}\times \frz_{J; \cpt}^*: \|\xi_{J,\perp}^{\cpt}\|\leqsl R_k\}\subset T^*Z(L_J)_{\cpt},\\
&\frD_{\ep, \ovl{\cC}_J}^J:=\frD_{\ep, \ovl{\cC}_J}^{L_J^\der}\subset \cT_{L_J^\der} \text{ as in \eqref{eqfrDepC}}. 
\end{align*}
Let 
\begin{align*}
&\frD^{Z}_{J,R_k}:=\frD_{J, R_k}^{Z, \cpt}\times \frD^{\flat, \ovl{\cI}}_{J, R_k},\\
&\frD_{\ep, R_k}^J=\frD_{\ep, \ovl{\cC}_J}^{J}\times^{Z(L_J)_{\cpt}}\frD^Z_{J, R_k}. 
\end{align*}
For any $J$, consider the proper function 
\begin{align}\label{eqwtchiJnorm}
\|\wt{\chi}_J\|^2: \pi_\trg^{-1}(\ovl{\cU}_J)=\cT^{\ovl{\cC}_J}_{L_J^\der}\times^{Z(L_J^\der)}T^*Z(L_J)_\cpt\times T^*\ovl{\cI}_J&\longrightarrow [0,\infty)\\
\nonumber (g_J, \xi^+_J, \xi_J^-; z_J^\cpt, \xi_{J,\perp}^\cpt; \bq^J, \bp^J)&\mapsto \|\chi_{L_J^\der}(\xi_J^+)\|^2+\|\xi_{J,\perp}^\cpt\|^2+\|\bp^J\|^2,
\end{align}
where $\xi_J^+\in \psi_{0,J}+\frb_J$ and $\xi_J^-\in -\psi_{\infty,J}+\frb_J^-$. 

First, using \eqref{eqdefnsfZJadm} and (ZAa)--(ZAc), it is easy to see that $\frD^{\flat, \ovl{\cI}}_{J,R_k}$ is a Weinstein domain, thus $\frD^{Z}_{J,R_k}$ is a (generalized) Weinstein domain (whose boundary has corners). Indeed, let 
\begin{align*}
f_{J,\flat}(\bq^J, \bp^J):=f_J^\frz(\bq^J)+\|\bp^J\|^2,
\end{align*}
which is an exhausting Morse function on $T^*\frz_{J;\RR}\cong T^*Z(L_J)_{0,\RR}$. For sufficiently small $\delta>0$, the Liouville vector field $\sfZ_{T^*\frz_{J;\RR}}':=\sfZ_{T^*\frz_{J;\RR}}-X_{\sfI_J}$ (corresponding to $\vartheta'_{T^*\frz_{J;\RR}}:=\vartheta_{T^*\frz_{J;\RR}}+
d\left(\vartheta_{T^*\frz_{J;\RR}}(\sfZ_{J,\flat}^{\frz})\right)$) is gradient-like with respect to $f_{J,\flat}$. Since $\sfZ_{T^*\frz_{J;\RR}}'$ is convex along $\partial \frD^{\flat, \ovl{\cI}}_{J,R_k}$ and has no zero outside of $\frD^{\flat, \ovl{\cI}}_{J,R_k}$, $\frD_{J,R_k}^{\flat, \ovl{\cI}}$ is a Weinstein domain whose completion is $(T^*\frz_{J;\RR},\vartheta'_{T^*\frz_{J;\RR}})$. 

Define inductively
\begin{align*}
&\frD_0=\frD_{\vn, R_0}^Z,\\
&\frD_{k+1}=\frD_k\cup\bigcup_{J\in \wt{I}_{k+1}}\frD_{R_{k+1}, R_{k+1}}^J. 
\end{align*}

Let $(\sfZ_{J,\flat}^\frz)_0=\{\bq^J_m\}_m$ be the set of zeros of $\sfZ_{J,\flat}^\frz$. Let $d_{J,m}$ be the dimension of the stable manifold of $\bq^J_m$ in $\frz_{J;\RR}$.  

\begin{lemma}\label{lemmafrDkWeinstein}
\begin{itemize}
\item[(i)] If $\sfZ_\flat$ is $(\delta, \frU)$-admissible, then for $R_0\gg R_1\gg\cdots \gg R_r\gg 1$, each $\frD_k$ has convex boundary with respect $\sfZ$ \eqref{eqdefnsfZJadm}.

\item[(ii)] If $\sfZ_\flat$ is $(\delta, \frU)$-admissible, then the completion of $\frD_k$ in $\cM_G$ gives $\cM_{G}^{\leqsl k}:=\pi_\trg^{-1}(\trg_{\leqsl k})$, where $\trg_{\leqsl k}=\bigcup_{J\in\wt{I}_{\leqsl k}}\trg_J$. 
\item[(iii)] If $\sfZ_\flat$ is strongly $(\delta, \frU)$-admissible, then each $\frD_k$ is a Weinstein subdomain. 
\end{itemize}
\end{lemma}

\begin{proof}
(i) Let $J\in \wt{I}_k$. Since $\frD_{\ep, \ovl{\cC}_J}^{J}$ has convex boundary along $\partial^{\out}\frD^J_{\ep, \ovl{\cC}_J}:=\frD_{\ep, \ovl{\cC}_J}^J\cap \|\chi_{L_J^\der}\|^{-1}(\ep)$ and concave boundary along $\partial^{\inin}\frD_{\ep, \ovl{\cC}_J}^{J}:=\frD^{J}_{\ep, \ovl{\cC}_J}\cap \pi_\trg^{-1}(\partial^\circ\cC_J)$ and $\frD_{J, R_k}^{Z}$ has convex boundary, it suffices to show that 
\begin{align*}
&\partial^{\inin}\frD^J_{R_k, R_k}=\partial^{\inin}\frD^J_{R_k, \ovl{\cC}_J}\times \frD_{J, R_k}^Z\subset \frD_{k-1}^\circ.
\end{align*}
 By assumption, $\partial^{\inin}\frD^J_{R_k, R_k}$ is a compact region whose projection to $\trg$ is $\nu_J^{-1}(\partial^\circ\cC_J\times \ovl{\cI}_J)$. Since $\nu_J^{-1}(\partial^\circ\cC_J\times \ovl{\cI}_J)\subset \cU_{\subsetneqq J}=\bigcup_{J^\dagg\in \wt{I}_{\subsetneqq J}}\cU_{J^\dagg}$, and for each $J^\dagg\subsetneq J$, $\|\wt{\chi}_{J^\dagg}\|^2$ \eqref{eqwtchiJnorm} is an exhausting function, the assertion is immediate. Note that the argument only requires $\sfZ_\flat$ to be $(\delta, \frU)$-admissible. 

(ii) We just need to check 
\begin{itemize}
\item[(1)] There is no zero of $\sfZ$ in $\cM_G^{\leqsl k}-\frD_k$;
\item[(2)] For any point $y\in \cM_G^{\leqsl k}$, the backward flow of $\sfZ$ will eventually take it inside $\frD_k$. 
\end{itemize}
By the assumption (ZAa) on $\sfZ_\flat$, we know that the only zeros of $\sfZ_\flat|_{\trg_{\leqsl k}}$ lie in $\cU_{\leqsl k}$. Moreover, for each $J\in \wt{I}_{\leqsl k}$,  $\sfZ_\flat|_{\ovl{\cU}_J}$ has zeros exactly at $\nu_J^{-1}(0\times (\sfZ_{J,\flat}^\frz)_0)$. Thus the zero locus of $\sfZ|_{\pi_\trg^{-1}(\ovl{\cU}_J)}$ (under the splitting \eqref{eqwtchiJnorm}) is exactly
\begin{align*}
\Gamma_J^\psi\times^{Z(L_J^\der)}Z(L_J)_\cpt\times (\sfZ_{J,\flat}^\frz)_0\cong Z(L_J)_{\cpt}\times (\sfZ_{J,\flat}^\frz)_0,
\end{align*}
where $\Gamma_J^\psi$ is the set of Kostant sections of $\cT_{L_J^\der}$ (which is a torsor of $Z(L_J^\der)$), and the factors $Z(L_J)_\cpt$ and $(\sfZ_{J,\flat}^\frz)_0$ represent points in the corresponding zero-sections in the cotangent bundles. 
This shows that $\frD_{\ep, R_k}^J$ contains all the zeros of $\sfZ|_{\pi_\trg^{-1}(\ovl{\cU}_J)}$, and (1) follows.

By the dynamics of $\sfZ_\flat$, we know for any $y_\flat\in \trg_{\leqsl k}$, $\varphi_{\sfZ_\flat}^{-t}(y)$ is contained in $\nu_J^{-1}(0\times \cI_J)$ for a unique $J\in \wt{I}_{\leqsl k}$ and for all $t\gg 1$. 
If $\sfZ_\flat$ is strongly $(\delta, \frU)$-admissible, then  $y_\flat$ is on the unstable manifold of a unique zero of $\sfZ_{\flat}$. 
Since the forward flow of $\sfZ|_{\pi_\trg^{-1}(\ovl{\cU}_J)}$ scales $\|\chi_{L_J^\der}\|$ and $\|\xi_{J,\perp}^\cpt\|$ by some positive weights and $\sfZ(\|\bp^J\|^2)\geqsl (2-\delta)\|\bp^J\|^2$, (2) follows as well.

(iii) We prove by induction on $k$. This is clear for $k=0$ since we assume $\sfZ$ to be strongly $(\delta, \frU)$-admissible. Suppose we have proved the case for $\frD_k$. Let 
\begin{align*}
\frD_{k+1}^\ep=\frD_k\cup \bigcup_{J\in \wt{I}_{k+1}}\frD_{\ep, R_{k+1}}^J
\end{align*}
for $\ep>0$ sufficiently small. By the same proof of (i) and (ii), we see that $\frD_{k+1}^\ep$ has convex boundary whose completion is $\cM_G^{\leqsl k+1}$.

Note that $\frD_{J,R_{k+1}}^{\flat, \ovl{\cI}}$ is a Weinstein domain that has a Weinstein handle decomposition indexed by $(\sfZ_{J,\flat}^{\frz})_0$. Let $\frD_{\bq^J_m}$ be a Weinstein handle associated with $\bq^J_m\in(\sfZ_{J,\flat}^{\frz})_0$. Then by Lemma \ref{lemmafrDepChandle} for the group $L_J^\der$, $\frD_{k+1}^\ep$ is obtained from $\frD_k$ by attaching the product handles 
\begin{align*}
\frD_{\ep, \ovl{\cC}_J, \bq^J_m}^J:=\frD_{\ep, \ovl{\cC}_J}^J\times^{Z(L_J^\der)}\frD_{J, R_{k+1}}^{Z, \cpt}\times \frD_{\bq^J_m},\quad J\in \wt{I}_{k+1}, 
\end{align*}
whose connected components 
\begin{align}\label{eqfrDJgpsi}
\frD^{J, [g_\psi]}_{\ep, \ovl{\cC}_J, \bq^J_m}, \quad [g_\psi]\in \Gamma_J^\psi\times^{Z(L_J^\der)}\pi_0(Z(L_J)), 
\end{align}
are Morse-Bott Weinstein handles of index $2|J|+r-|J|+d_{J, m}=r+|J|+d_{J, m}$. 

Lastly, since $\frD_{k+1}^\ep\subset \frD_{k+1}$ have the same completion under $\sfZ$, we conclude that $\frD_{k+1}$ is a Weinstein domain as well. The inductive step then follows. 
\end{proof}

\begin{theorem}\label{thmMGWeinstein}
\begin{itemize}
\item[(i)] For any $\delta$-admissible Liouville 1-form $\vartheta$ on $\cM_G$, 
$(\cM_G, \vartheta)$ is a Liouville manifold. For any $\frU\in \CI_{\vartheta, \delta}$, the corresponding $\frD_r$ gives a Liouville domain (whose boundary has corners) whose completion is $\cM_G$. 

\item[(ii)]
If $\vartheta$ is strongly $\delta$-admissible, then there exists an exhausting Morse-Bott function $f$ on $\cM_G$, with respect to which, $(\cM_G, \vartheta)$ is a generalized Weinstein manifold. For any $\frU\in \CI^{\sigma}_{\vartheta, \delta}$, the generalized Weinstein handles are given by \eqref{eqfrDJgpsi} for $J\in \cP_{ft}(\wt{I})$, for any fixed $\ep>0$ sufficiently small, whose cores are diffeomorphic to $D^{2|J|}\times Z(L_J)_{0, \cpt}\times D^{d_{J, m}}$.  
\end{itemize}
\end{theorem}
\begin{proof}
This is a direct consequence of Lemma \ref{lemmafrDkWeinstein}. 
\end{proof}

\begin{remark}\label{remarkMGthetadeform}
The space of $\delta$-admissible Liouville 1-forms $\vartheta$ is nonempty and convex. By the construction of $\frD_r$, it is clear that one can choose a single $\frD_r$ such that for every $\delta$-admissible $\vartheta$, it gives a Liouville domain whose completion is $\cM_G$.
By Gray stability, for any $\delta$-admissible family $\vartheta_t, 0\leqsl t\leqsl 1$, the deformation is conjugate to a compactly supported deformation of Liouville manifolds. 
Thus, the Liouville manifold structures $(\cM_G, \vartheta)$  are all deformation equivalent and they form a contractible space.

\end{remark}

\begin{remark}\label{remSpecialZflatHandle}
If we use a special $\sfZ_\flat$ as in Corollary \ref{corsfZflatspecial}, whose canonical Liouville 1-form $\vartheta$ is referred as special, then for each $J\in \PfI$, the generalized Weinstein handles of $\cM_G$ are indexed by $\Gamma_J^\psi\times^{Z(L_J^\der)}\pi_0(Z(L_J)), J\in \PfI$, and each of them has core diffeomorphic to $D^{2|J|}\times Z(L_J)_{0,\cpt}$. 
\end{remark}

%%%%%%%%%%%%%%%%%%%%%%%%%%%%%%%%%%%%%%%%%%%%%%%%%%%%%%%%%%%%%%%%%%%%%%%%%%%%%%%%%%%%%%%%%

\subsection{Weinstein sectorial coverings of $\cM_G^\circ$}\label{subsec: sectorial,covering}

Let $\delta>0$ be small. Choose any admissible open covering $\frU$ of $\trg$ and any $(\delta, \frU)$-rigid Liouville 1-form $\vartheta$ on $\cM_G^\circ$. Recall the notations from \S\ref{sssdigposets} and \S\ref{sssAppFukSecIncl}. 

\begin{prop}
The collection $\{\cT_{L_J}^{\ovl{\cU}_J}=\pi_\triangle^{-1} (\overline{\cU}_J)\}_{J\in \PfI}$ forms a sectorial covering of $(\cM_G^\circ, \vartheta)$. Moreover, for $\vec{J}=(J^1\supsetneq \cdots\supsetneq J^k)\in \sfC(\PfI)$, 
$(\cT_{L_{J_k}}^{\ovl{\cU}_{\vec{J}}},\vartheta_{\vec{J}}:=\vartheta|_{\cT_{L_{J_k}}^{\ovl{\cU}_{\vec{J}}}})$ is equivalent to the Weinstein sector $(\ovl{\cT}_{L_{J_k}}, \vartheta_{\ovl{\cT}_{L_{J_k}}})$ up to (a contractible space of) deformations of Liouville sectors.  
\end{prop}

\begin{proof}
Applying Lemma \ref{lemmarigidproductsector} to the reductive group $L_{J_k}$ and $\vartheta_{J_k}\in \frs_{\delta, \rigid}^{\frU}(\cT_{L_{J_k}}^{\ovl{\cU'_{J_k}}})$, we get the product sector structure on $\cT_{L_{J_k}}^{\ovl{\cU'_{J_k}}}$, where $\cU'_{J_k}=\nu_{J_k}^{-1}(\cC_{J_k}\times \Nb(\ovl{\cI}_{J_k}))$. 

By Lemma \ref{lemma: exist admissible}, $\ovl{\cI}_{\vec{J}}$ is cut out by transverse intersecting hypersurfaces in $\Nb(\ovl{\cI}_{J_k})$, which correspond to $\cC$-adapted hypersurfaces in $\nu_{J_k}^{-1}\left(\Nb(\ovl{\cC}_{J_k})\times \Nb(\ovl{\cI}_{J_k})\right)$. 
As $\cT_{L_{J_k}}^{\ovl{\cU}_{\vec{J}}}$ is cut out further from $\cT_{L_{J_k}^{\der}}\times^{Z(L_{J_k}^\der)}T^*Z(L_{J_k})_{\cpt}\times T^*\ovl{\cI}_{\vec{J}}$ in the factor $\cT_{L_{J_k}^{\der}}$ by $(b_{\RR, \der}^{J_k})^{-1}(\partial^\circ \cC_{J_k})$, 
the first part of the proposition follows.

The second part follows from the contractibility of $\frs_\delta(\cT_{L_{J_k}}^{\ovl{\cU}_{\vec{J}}})$ (cf. \S\ref{sssSpcSectStr}) and the first row of deformations in the diagram \S\ref{sssSummdeform} for the reductive group $L_{J_k}$. 
\end{proof}

\begin{prop}\label{propWeinsteinSectorial}
The sectorial covering $\{\cT_{L_J}^{\ovl{\cU}_J}=\pi_\triangle^{-1} (\overline{\cU}_J)\}_{J\in \PfI}$ is a Weinstein sectorial covering in the sense of \ref{defnWeinsteinsectorcovering}.
\end{prop}

\begin{proof}
The verification of the conditions in \ref{defnWeinsteinsectorcovering} follows from a combination of the same argument as for Example \ref{exam: Qi, WeinsteinSectorial} (when we look at the covering $\{\ovl{\cU}_J\}_{J\in \PfI}$ of $\trg$) and similar arguments in the proof that  each $\cT_{L_J}^{\ovl{\cU}_J}$ is a product Weinstein sector up to deformations. We omit the details. 
\end{proof}

\subsection{$\cW(\cM_G^\circ)$ as a homotopy colimit}\label{subsec: cWM_G,colimit}

Let $\Pr^L$ (resp. $\Pr^R$, $\Pr^{L,R}$) be the $\infty$-category of presentable $\infty$-categories with functors that are left (resp. right, both left and right) adjoints. Let $\Pr^L_{\st}$ (resp. $\Pr^R_{\st}$, $\Pr^{L,R}_{\st}$) be the corresponding full subcategory of presentable stable $\infty$-categories. 
Let $\CatEx$ be the $\infty$-category of small stable $\infty$-categories which are idempotent complete. We will use the same assumptions as in \S\ref{subsec: sectorial,covering}. Using the Weinstein sectorial covering $\{\cT_{L_J}^{\ovl{\cU}_J}, J\in \PfI\}$ (up to deformations) of $\cM_G^\circ$ (Proposition \ref{propWeinsteinSectorial}) and sectorial descent for wrapped Fukaya categories (Theorem \ref{thm: sectorialcover,cW}), we get the natural functor 
 \begin{align*}
 \Phi: \sfC(\PfI)^{op}&\longrightarrow \CatEx \\
 \vec{J}=(J^1\supsetneq\cdots\supsetneq J^k)&\mapsto \cW(\cT_{p(\vec{J})}^{\ovl{\cU}_{\vec{J}}}, \vartheta_{\vec{J}})
 \end{align*}
from \eqref{eqPhibulletN}, that induces an equivalence 
\begin{align*}
\cW(\cM_G^\circ, \vartheta)\simeq \colim_{\sfC(\PfI)^{op}}\Phi.
\end{align*}  
Moreover, by Lemma \ref{lemmaAppcWJicY}, we get the following. 
\begin{cor}\label{cor: cWcMG,twocolimits}
The functor $\Phi$ descends to a canonical functor (up to a contractible space of choices)
\begin{align*}
\Phi': \PfI&\longrightarrow \CatEx\\
J&\mapsto \cW(\ovl{\cT}_{L_J}, \vartheta_{\ovl{\cT}_{L_J}}),
\end{align*}
and we have 
\begin{align*}
\cW(\cM_G^\circ, \vartheta)\simeq \colim_{\PfI}\Phi'.
\end{align*}
Equivalently, we have 
\begin{align*}
\Ind\cW(\cM_G^\circ, \vartheta)\simeq \lim_{J\in \PfI^{op}}\Ind\cW(\ovl{\cT}_{L_J}, \vartheta_{\ovl{\cT}_{L_J}}).
\end{align*}
\end{cor}

\begin{proof}
Set $\Phi'=\cW^{\hat{\frU}}_{\delta, \rigid}$. The corollary follows from \eqref{eqPhiWrigU} and the remark after  Lemma \ref{lemmaAppcWJicY}. 
\end{proof}

%%%%%%%%%%%%%%%%%%%%%%%%%%%%%%%%%%%%%%%%%%%%%%%%%%%%%%%%%%%%%%%
\section{Quasi-coherent sheaves on certain Weil restrictions}\label{s:Weil res}
%%%%%%%%%%%%%%%%%%%%%%%%%%%%%%%%%%%%%%%%%%%%%%%%%%%%%%%%%%%%%%%

In this section we prove a general statement (Theorem \ref{th:ff}) describing quasi-coherent sheaves on certain schemes involving Weil restrictions. It will be applied to describe quasi-coherent sheaves on the regular centralizer group scheme for a reductive group.

\subsection{Admissibility} Let $k$ be an algebraically closed field of characteristic zero.  Let $W$ be a finite group. Let $\Rep_{k}(W)$ be the category of finite-dimensional $k$-representations of $W$.

Fix a unital commutative algebra object $C\in \textup{CommAlg}(\Rep_{k}(W))$ such that as a $W$-module, $C$ is isomorphic to the regular representation $k[W]$. In particular, $k\subset C$ is identified with the $W$-invariants in $C$.

Let $V\in\Rep_{k}(W)$. Consider the map
\begin{align*}
\wt\d_{j}: \Lambda^j(C\otimes V)&\longrightarrow C\otimes\Lambda^jV\\
(c_1\otimes v_1)\wedge\cdots\wedge(c_j\otimes v_j)&\mapsto (c_1\cdots c_j)\otimes v_1\wedge\cdots\wedge v_j.
\end{align*}
Let 
\begin{equation}\label{dj}
\delta_j: \Lambda^j((C\otimes V)^W)\overset{\sim}{\longrightarrow}(C\otimes\Lambda^jV)^W
\end{equation}
be the restriction of $\wt\d_{j}$.

\begin{defn}\label{def:adm rep}
The representation $V$ is called {\em $C$-admissible} if $\d_{j}$ is an isomorphism for all $j\ge0$.
\end{defn}

The following is clear form the definition. 
\begin{lemma}
If $\dim V=1$, then $V$ is $C$-admissible. 
\end{lemma}

\begin{lemma}\label{l:adm triv}
Let $V_{0}$ be a vector space with trivial action of $W$, and let $V\in \Rep_{k}(W)$. Then $V$ is $C$-admissible if and only if $V\op V_{0}$ is.
\end{lemma}
\begin{proof}
It suffices to treat the case $V_0=k$, the trivial representation of $W$. In this case, $\Lambda^j((C\otimes (V\op k))^W)=\Lambda^j((C\otimes V)^W\op k)=\oplus_{i=0}^j\Lambda^j((C\otimes V)^W)$. Similarly, $(C\otimes\Lambda^j(V\op k))^W=\oplus_{i=0}^j(C\otimes\Lambda^iV)^W$. The map $\d_{j}^{V\op k}$ for $V\op k$ is the direct sum of $\d^V_{i}$ for $V$, where $i$ runs from $0$ to $j$. From this we see that $\d^V_i$ are isomorphisms for all $i\ge0$ if and only if $\d^{V\op k}_j$ are isomorphisms for all $j\ge0$.
\end{proof}

The main example of admissible representations we will encounter is the following.
\begin{lemma}\label{l:ref adm} Suppose $W$ is a finite Weyl group with reflection representation $V$. Let $C$ be the coinvariant algebra for $W$, namely
\begin{equation*}
C=\Sym(V^{*})/(\Sym(V^{*})^{W}_{+})
\end{equation*}
where $\Sym(V^{*})^{W}_{+}$ is the positive degree part of the $W$-invariants on $\Sym(V^{*})$.  Then $V$ is $C$-admissible.
\end{lemma}
\begin{proof} Let $S=V\sslash W$, an affine space of the same dimension as $V$. Let $\pi: V\to S$ be the projection. Pullback of forms $\pi^*\Om^j_S\to \Om^j_V$ induces a map of coherent sheaves
\begin{equation*}
    \th_j: \Om^j_S\to (\pi_*\Om^j_V)^W.
\end{equation*}
We claim that $\th_j$ is an isomorphism. Indeed, since $\pi$ is flat, $\pi_*\Om^j_V$ and hence $(\pi_*\Om^j_V)^W$ is also locally free over $S$. Thus it suffices to check $\th_j$ is an isomorphism away from a codimension $\ge 2$ subscheme of $S$. Let $V^\circ\subset V$ be the open subset that is the complement of intersections of two root hyperplanes, and $S^\circ:=\pi(V^\circ)$. Then $S-S^\circ$ has codimension $2$. We reduce to checking that $\th_j|_s$, the restriction of $\th_j$ to $\Spec \cO_{S,s}$, is an isomorphism for $s=\pi(v)\in S^{\circ}$. When $\Stab_W(v)$ is trivial, $\pi$ is \'etale over $s$ hence $\th_j|_s$ is an isomorphism. Otherwise, $\Stab_W(v)$ is generated by a reflection $r\in W$, and we can identify $(\pi_*\Om^j_V)^W\ot \cO_{S,s}$ with $(\Om^j_V\ot \cO_{V,v})^{r}$ by restriction. Thus we reduce to showing $\Om^j_S\ot \cO_{S,s}\to (\Om^j_V\ot \cO_{V,v})^r$ is an isomorphism. Decompose $
V=V^r\op V'$, where $V'=V^{r=-1}$ is one-dimensional, we easily reduce to showing that $\th_j$ is an isomorphism for the Weyl group $\ZZ/2\ZZ=\j{r}$ and its reflection representation $V'$. In that case, the only nontrivial statement is for $j=1$. Let $x$ be a linear coordinate on $V'$ so that $x^2$ is the coordinate for $S'=V'\sslash \j{r}$, then $\Om_{S'}\to (\pi_*\Om_{V'})^r$ sends $d(x^2)$ to $2xdx$, which freely generates $(\pi_*\Om_{V'})^r$ as a $k[x^2]$-module. This proves that $\th_j$ is an isomorphism in general.

Now consider the restriction of the isomorphism $\th_j$ to the fiber at the closed point $0\in S$, i.e., applying $(-)\ot_{\cO_S}k(0)=(-)\ot_{\Sym(V^*)^W}\Sym(V^*)^W/(\Sym(V^*)^W_+)$ on both sides of $\th_j$. Note that $(\pi_*\Om^j_V)^W\ot_{\cO_{S}}k(0)\cong (\Sym(V^*)\ot \wedge^jV^*)^W\ot_{\cO_{S}}k(0)\cong (C\ot \wedge^jV^*)^W$. We obtain that
\begin{equation}\label{ovthj}
    \ov\th_{j}: \wedge^i(T^*_0S)\isom (\pi_*\Om^j_V)^W\ot_{\cO_{S}}k(0)\cong (C\ot \wedge^j(V^*))^W.
\end{equation}
For $j=1$ we obtain $T^*_0S\cong (C\ot V^*)^W$, under which \eqref{ovthj} is the same as \eqref{dj} for $(C,V^*)$. We conclude that $V\cong V^*$ is $C$-admissible.
\end{proof}

\subsection{Statement of main result}\label{ss:setup ff}

\sss{Setup}\label{sss:ff setup} Let $k$ be a field of characteristic $0$. Let $W$ be a finite group. Let $Q$ be a finite type $k$-scheme, and let $S=Q\sslash W$. Consider the quotient map
\begin{equation*}
\pi: Q\to S.
\end{equation*}
Then $\cO_{Q}$ is a coherent sheaf of $\cO_{S}[W]$-modules. We assume that $\cO_{Q}$ is locally free over $\cO_{S}[W]$. In particular, $\pi$ is finite flat.

Let $X$ be another $k$-scheme of finite type with a $W$-action and let $\xi: X\to Q$ be a quasi-projective $W$-equivariant morphism. We can form the Weil restriction
\begin{equation*}
\cR_{X}=R_{Q/S}(X).
\end{equation*}
Since $\pi: Q\to S$ is finite and locally free, $X\to Q$ is quasi-projective, $\cR_{X}$ is represented by a $S$-scheme that is separated and of finite type by \cite[\S7.6, Theorem 4, Proposition 5]{BLR}. 

The actions of $W$ on $Q$ and $X$ induce an action of $W$ on $\cR_{X}$. Define the $S$-scheme
\begin{equation*}
J_{X}=(\cR_{X})^{W}.
\end{equation*}

Since Weil restriction $R_{Q/S}$ is right adjoint to base change $(-)\times_{S}Q$, we have a co-unit map
\begin{equation*}
\cR_{X}\times_{S}Q\to X
\end{equation*}
that restricts to a map
\begin{equation*}
\vk: J^{(1)}_{X}:=J_{X}\times_{S}Q\to X.
\end{equation*}
Note that $\vk$ is $W$-equivariant with respect to the $W$-action on the second factor $Q$ and its action on $X$.

Consider the diagram
\begin{equation*}
\xymatrix{ J_{X} & \ar[l]_-{\pr_{1}} J_{X}\times_{S}(Q/W) \ar[r]^-{\vk} & X/W
}
\end{equation*}
We have the functor
\begin{equation}\label{prk}
\Phi_{X}:=\vk_{*}\pr_{1}^{*}: \QCoh(J_{X})\to \QCoh(X/W)=\QCoh(X)^W.
\end{equation}

We introduce some more notations. For $q\in Q(K)$ (where $K$ is a field), let $W_{q}\subset W$ be its stabilizer. Since the finite $K$-scheme $\pi^{-1}(\pi(q))$ is a disjoint union of thick points indexed by the $W$-orbits of $q$. Let $C_{q}$ be the local ring of $\pi^{-1}(\pi(q))$ at $q$. Then $W_{q}$ acts on $C_{q}$ by algebra automorphisms. By our assumption, $\cO_{Q}$ is locally free as an $\cO_{S}[W]$-module, which implies that $C_{q}\cong K[W_{q}]$ as a $W_{q}$-module.  

Let $X_q$ be the base change of $X$ along $q: \Spec K\to Q$. Then $W_q$ acts on $X_q$. We have natural maps
\begin{equation*}
    X_q^{W_q}\xr{\io_q} X_q\xr{i_q} X.
\end{equation*}

The main result of this section is:

\begin{theorem}\label{th:ff} In the above setup, assume further:
\begin{itemize}
\item The map $\xi: X\to Q$ is smooth.
\item For each field $K$ containing $k$ and $q\in Q(K)$, and $x\in X_q^{W_{q}}(K)$, the $W_{x}=W_{q}$-module $T_{x}(X_{q})$ (tangent space of $X_{q}$ at $x$) is $C_{q}$-admissible (see Definition \ref{def:adm rep}). 
\end{itemize}
Then:
\begin{enumerate}
    \item The functor $\Phi_X$ in \eqref{prk} is fully faithful.
    \item The essential image of $\Phi_X$ consists of $\cG\in \QCoh(X)^W$ such that for any field $K$ containing $k$ and $q\in Q(K)$, $i_q^*\cG\in \QCoh(X_q)^{W_q}$ lies in the essential image of $\io_{q*}: \QCoh(X_q^{W_q})\to \QCoh(X_q)^{W_q}$ that sends $\cH\in \QCoh(X_q^{W_q})$ to $\io_{q*}\cH$ with the trivial $W_q$-action.
\end{enumerate}
\end{theorem}

The rest of the section is devoted to the proof of Theorem \ref{th:ff}. The proof of full faithfulness will be given in \S\ref{ss:ff}, after some reductions. The essential image statement will be proved in \S\ref{ss:im}.

\subsection{First reductions}\label{ss:ff reduction}

We are in the setting of \S\ref{ss:setup ff}.
\begin{lemma}\label{l:Phi left}
The functor $\Phi_{X}$ admits a left adjoint
\begin{equation*}
\Phi^{\ell}_{X}=\pr_{1*}(\vk^{*}(-)\ot \pr_{2}^{*}\om_{\pi}): \QCoh(X/W)\to \QCoh(J_{X}).
\end{equation*}
Here $\pr_{2}: J_{X}\times_{S}(Q/W) \to Q/W$ is the second projection, and $\om_{\pi}\in \QCoh(Q/W)$ is the relative dualizing sheaf for $\pi: Q\to S$, which is a $W$-equivariant invertible sheaf on $Q$.

In particular, $\Phi_{X}$ is fully faithful if and only if the co-unit map
\begin{equation*}
c_{X}: \Phi^{\ell}_{X}\c\Phi_{X}\to \id_{\QCoh(J_{X})}
\end{equation*}
is an isomorphism.
\end{lemma}
\begin{proof}
Since $\pr_{1}$ is the base change of the finite flat map $\ov\pi: Q/W\to S$, we have a canonical isomorphism $\pr_{1}^{!}(-)\cong \pr^{*}_{1}(-)\ot \pr_{2}^{*}\om_{\pi}$. Therefore we may rewrite $\Phi_{X}$ as $\vk_{*}(\pr_{1}^{!}(-)\ot \pr_{2}^{*}\om^{-1}_{\pi})$. The formula for the left adjoint of $\Phi_{X}$ then follows.
\end{proof}

The isomorphism $c_{X}$ is induced by an explicit kernel sheaf $\cK$ on $J_{X}\times J_{X}$ as follows. Consider the {\em derived} fiber product of $J^{(1)}_{X}$ with itself over $X$
\begin{equation*}
J^{(2)}_{X}=(J_{X}\times_{S}Q)\times^{\bR}_{X}(J_{X}\times_{S}Q)=J^{(1)}_{X}\times_{X} J^{(1)}_{X}
\end{equation*}
equipped with a diagonal action of $W$ (which acts also on $X$). We have two projections $\oll{p}, \orr{p}: J^{(2)}_{X}/W\to J_{X}$. Let $p=(\oll{p}, \orr{p}): J^{(2)}_{X}/W\to J_{X}\times_{S} J_{X}$. Let $\oll{\frq}, \orr{\frq}: J^{(2)}/W\to Q/W$ be the projections to the first and second factors of $Q/W$. Then $c_{X}$ is induced by the kernel sheaf
\begin{equation*}
\cK=p_{*}(\orr{\frq}^{*}\om_{\pi})\in \QCoh(J_{X}\times_{S} J_{X}).
\end{equation*}

The co-unit map $c_{X}$ corresponds to a map of coherent complexes
\begin{equation}\label{kernel to diag}
\d_{X}: \cK\to \D_{J*}\cO_{J_{X}}
\end{equation}
where $\D_{J}: J_{X}\incl J_{X}\times_{S} J_{X}$ is the diagonal map. 
 
\begin{lemma}\label{l:check point}
If $c_{X}(\cF)$ is an isomorphism for any $\cF$ of the form $\cF=\ph_{*}\cO_{\Spec K}$, where $K$ is a field containing $k$ and $\ph: \Spec K\to  J_{X}$ is any morphism over $k$, then $c_{X}$ is an isomorphism.

\end{lemma}
\begin{proof}
The map $c_{X}$ is an isomorphism if and only if the map $\d_{X}$ in \eqref{kernel to diag} is an isomorphism in $\QCoh(J_{X}\times_{S} J_{X})$, if and only if $\cC=\Cone(\d_{X})\cong 0$.   By \cite[Lemma 10]{Arinkin}, it suffices to show that $\psi^{*}\cC\cong 0$ for all field-valued point $\psi: \Spec K\to J_{X}\times_{S} J_{X}$. A fortiori, it suffices to show that for any field valued point $\ph: \Spec K\incl J_{X}$, with induced map $\ph'=\ph\times \id: \Spec K\times_{S}J_{X}\to J_{X}\times_{S} J_{X}$, $\ph'^{*}\cC\cong 0$. 

By base change and projection formula, the functor $\Phi_{X}^{\ell}\c\Phi_{X}\c\ph_{*}: \QCoh(\Spec K)\to \QCoh(J_{X})$ is given by the kernel sheaf $\ph'^{*}\cK$. The map $c_{X}\c \ph_{*}:\Phi_{X}^{\ell}\c\Phi_{X}\c\ph_{*}\to \ph_{*}$ is induced from the map $\ph'^{*}\d_{X}: \ph'^{*}\cK\to \ph'^{*}\D_{J*}\cO_{J_{X}}$. Therefore the $\ph'^{*}\cC\cong 0$ if and only if $c_{X}\c \ph_{*}$ is an isomorphism of functors, if and only if $c_{X}(\ph_{*}\cF)$ is an isomorphism for all $\cF\in \QCoh(\Spec K)$, if and only if $c_{X}(\ph_{*}\cO_{\Spec K})$ is an isomorphism.
\end{proof}

\sss{Base change}\label{sss:bc}
We consider the behavior of $c_{X}$ under base change. Let $\b: S'\to S$ be any morphism from a finite type $k$-scheme $S'$. Let $Q'=Q\times_{S}S'$, $X'=X\times_{Q}Q'$, and form $R_{X'}$ and $J_{X'}$ using $(X', Q')$ in place of $(X,Q)$. We have a commutative diagram where both squares are Cartesian
\begin{equation}\label{bc Phi}
\xymatrix{J_{X'} \ar[d]^{\b_{J}} & \ar[l]_-{\pr'_{1}} J_{X'}\times_{S'}(Q'/W) \ar[d]^{\wt\b_{J}}\ar[r]^-{\vk'} & X'/W\ar[d]^{\b_{X/W}}\\
J_{X} & \ar[l]_-{\pr_{1}} J_{X}\times_{S}(Q/W) \ar[r]^-{\vk} & X/W
}
\end{equation}
Let $\Phi_{X'}=\vk'_{*}\c \pr'^{*}_{1}: \QCoh(J_{X'})\to \QCoh(X'/W)$, with left adjoint $\Phi^{\ell}_{X'}$ and co-unit map $c_{X'}: \Phi^{\ell}_{X'}\c\Phi_{X'}\to \id_{\QCoh(J_{X'})}$. 

\begin{lemma}\label{l:cX pull push}
Base change isomorphisms induce canonical isomorphisms of functors
\begin{eqnarray*}
\Phi_{X}\c \b_{J*}\cong \b_{X/W*}\c\Phi_{X'}, \quad  \b_{J*}\c\Phi^{\ell}_{X'}\cong \Phi^{\ell}_{X}\c\b_{X/W*},\\
\Phi_{X'}\c\b_{X/W}^{*}\cong \b_{J}^{*}\c\Phi_{X}, \quad \b_{J}^{*}\c\Phi^{\ell}_{X}\cong \Phi^{\ell}_{X'}\c\b_{X/W}^{*}.
\end{eqnarray*}
In particular, if $c_{X'}(\cF')$ is an isomorphism for some $\cF'\in \QCoh(J_{X'})$, then $c_{X}(\b_{J*}\cF')$ is also an isomorphism.
\end{lemma}
\begin{proof}
Since $\pr_{1}$ is flat by assumption, we have the base change isomorphism $\pr_{1}^{*}\b_{J*}\cong \wt\b_{J*}\pr'^{*}_{1}$. This implies $\Phi_{X}\c \b_{J*}\cong \b_{X/W*}\c\Phi_{X'}$. Passing to left adjoints we get $\b_{J}^{*}\c\Phi^{\ell}_{X}\cong \Phi^{\ell}_{X'}\c\b_{X/W}^{*}$.

On the other hand, since $J_{X}$ and $X$ are both flat over $S$, the Cartesian square on the right side of \eqref{bc Phi} is also derived Cartesian, hence base change for quasi-coherent sheaves holds: we have  canonical isomorphisms $\vk^{*}\c\b_{X/W*}\cong \wt\b_{J*}\c\vk'^{*}$ and $\b^{*}_{X/W}\c\vk_{*}\cong \vk'_{*}\c\wt\b^*_{J}$, which induce isomorphisms $\b_{J*}\c\Phi^{\ell}_{X'}\cong \Phi^{\ell}_{X}\c\b_{X/W*}$ (in view of the formula for $\Phi^{\ell}_{X}$ given in Lemma \ref{l:Phi left}), and $\Phi_{X'}\c\b_{X/W}^{*}\cong \b_{J}^{*}\c\Phi_{X}$. 
\end{proof} 

\subsection{Proof of full faithfulness}\label{ss:ff}

\sss{}\label{sss:Wq} Combining Lemma \ref{l:check point} and Lemma \ref{l:cX pull push}, to check that $c_{X}$ is an isomorphism, it suffices to check the same statement after making a base change along any morphisms $s: S'=\Spec K\to S$, where $K$ is a field, as in \S\ref{sss:bc}, and check that $c_{X}(\ph_{*}\cO_{\Spec K})$ is an isomorphism for any $\ph:\Spec K\to J_{X}$ over $s$. Let $q\in Q(K)$ be the image of $\ph\in J_{X}(K)$. 

After the base change, and renaming $K$ be $k$, we reduce to the situation where $S=\Spec k$, and $Q=\Spec C$, where $C\cong k[W]$ as a representation of $W$, and equipped with a $k$-point $q\in Q(k)$. We decompose the underlying space of $Q$ into singleton schemes (possibly non-reduced) $Q=\coprod_{w\in W/W_{q}}Q_{w(q)}$. Correspondingly we have a decomposition $X=\coprod_{w\in W/W_{q}}X_{w(q)}$, where $X_{w(q)}=\xi^{-1}(Q_{w(q)})$. We have $\cR_{X}=\prod_{w\in W/W_{q}}\cR_{Q_{w(q)}/k}(X_{w(q)})$, and $J_{X}=(\cR_{X})^{W}\cong (\cR_{Q_{q}/k}(X_{q}))^{W_{q}}$. We may thus reduce further to the case where $W$ is replaced with $W_{q}$, and $Q=\Spec C_{q}$ is a singleton scheme with the unique $k$-point $q$.

\sss{New setup}\label{sss:Q=C} Below we will write $Q=\Spec C$, with $C$ a local $k$-algebra with residue field $k$ and with a $W$-action such that $C\cong k[W]$ as a $W$-module. We write $J_{X}^{(1)}=J_{X}\times_{k}\Spec C$ as $J_{X,C}$. By Lemma \ref{l:check point}, to show $c_{X}$ is an isomorphism in this situation, it suffices to show that for any $\ph\in J_{X}(k)$
\begin{equation}\label{c ph}
c_{\ph}:=c_{X}(\ph_{*}\cO_{\Spec k}): \Phi^{\ell}\c\Phi(\ph_{*}\cO_{\Spec k})\to \ph_{*}\cO_{\Spec k}.
\end{equation}
is an isomorphism in $\QCoh(J_{X})$. 

Consider the following diagram where the left square is Cartesian
\begin{equation}\label{diagram: S,J_X,0,C}
\begin{tikzcd}[column sep=5em]
\Spec k\ar[d, "\varphi"]&(\Spec C)/W\ar[d, "\wt{\varphi}"]\ar[l]\ar[r, "\wt{x}"]&X/W\ar[d, equal]\\
J_{X}\ar[ur, phantom, very near end, "\urcorner"]&J_{X,C}/W\ar[l, "\pr_1"']\ar[r, "\varkappa"]&X/W. 
\end{tikzcd}
\end{equation}

Let $X_{0}$ be the fiber of $\xi: X\to Q=\Spec C$ over $\Spec k\incl \Spec C$. Then $X_{0}$ is a smooth $k$-scheme of finite type.  Let $x\in X_{0}(k)$ be unique closed point in the image of $\wt x$. Note that $x\in X_{0}^{W}$.

By the reduction steps we have done, to prove Theorem \ref{th:ff}, it remains to prove the following statement.

\begin{prop} In the situation of \S\ref{sss:Q=C}, suppose $T_{x}X_{0}$ (as a $W$-representation) is $C$-admissible, then the map $c_{\ph}$ in \eqref{c ph} is an isomorphism. 
\end{prop}
\begin{proof}We may replace $X$ by an affine neighborhood of $x$. Hence we may assume $X$ is affine, therefore $J_{X}, J_{X,C}$ are also affine.  We thus identify quasi-coherent sheaves on $J_{X}$ with $\cO(J_{X})$-modules, and similarly for $J_{X,C}$.

Let $\om_{C/k}=\Hom(C,k)$ viewed as an object in $\QCoh(\Spec C)^{W}$.

Unravelling the definitions of $\Phi$, we have $\Phi(\ph_{*}\cO_{\Spec k})\cong \vk_{*}\pr_{1}^{*}\ph_{*}\cO_{\Spec k}\cong \vk_{*}\wt\ph_{*}\cO_{\Spec C}\cong \wt x_{*}\cO_{\Spec C}$ (using flat base change). By  Lemma \ref{l:Phi left}, 
\begin{equation*}
\Phi^{\ell}\c\Phi(\ph_{*}\cO_{\Spec k})=(\vk^{*}\wt x_{*}\cO_{\Spec C}\ot \om_{C/k})^{W}.
\end{equation*}
Here we understand both sides as $\cO(J_{X})$-modules;  the right side before taking $W$-invariants is a $W$-equivariant $\cO(J_{X})\ot C=\cO(J_{X,C})$-module. The map $c_{\ph}$ is induced from the co-unit map
\begin{equation*}
\vk^{*}\wt x_{*}\cO_{\Spec C}=\vk^{*}\vk_{*}\wt\ph_{*}\cO_{\Spec C}\to \wt\ph_{*}\cO_{\Spec C}
\end{equation*}
by tensoring with $\om_{C/k}$ and taking $W$-invariants.

Now let $V=T_{x}X_{0}$. Let $I_{\wt x}\subset \cO(X)$ be the ideal of the closed subscheme given by the image of $\wt x$. Let $\fm_{x}\subset \cO(X_{0})$ be the maximal ideal at $x$. We have $I_{\wt x}\ot_{C}k\isom \fm_{x}$. We thus have a $W$-equivariant surjection
\begin{equation*}
I_{\wt x}\surj \fm_{x}\surj \fm_{x}/\fm_{x}^{2}=V^{*}.
\end{equation*}
Since $W$ is finite and $\ch(k)=0$, we can choose a $W$-equivariant section of the above surjection, $s: V^{*}\to I_{\wt x}$. By shrinking $X$ we may assume $s(V^{*})$ generate $I_{\wt x}$.  We can form the Koszul resolution of $\wt x_{*}\cO_{\Spec C}$
\begin{equation*}
\Kos(\wt x): \wedge^{n}V^{*}\ot\cO(X)\to \cdots \to \wedge^{2}V^{*}\ot \cO(X)\to V^{*}\ot \cO(X)\to 0
\end{equation*}
in degrees $-n,\cdots, 0$, where $n=\dim V$. The differentials are defined by
\begin{equation*}
\th\ot f\mt \sum_{i=1}^{n}(\th\wedge v_{i})\ot s(v^{i})f, \quad  \th\in \wedge^{*}(V^{*}), f\in \cO(X)
\end{equation*}
where $\{v_{i}\}$ is a basis for $V$ and $\{v^{i}\}$ the dual basis for $V^{*}$.

Let $M=(V\ot C)^{W}$, viewed as either a $k$-vector space or an affine $k$-scheme. By \cite[\S7.6, Proposition 5(h)]{BLR}, $J_{X}$ is smooth over $k$. A calculation using definitions show that there is a canonical isomorphism to $M\cong T_{\ph}J_{X}$. The map $s_{X}: X\to V\times_{k}\Spec C$ (induced by $s$) induces a map
\begin{equation*}
\cR_{X}=R_{C/k}(X)\to R_{C/k}(V\times_{k}\Spec C)=V\ot C
\end{equation*}
where $V\ot C$ is viewed as an affine space.  Taking $W$-invariants we get a map
\begin{equation*}
s_{J}: J_{X}\to M
\end{equation*}
that sends $\ph\in J_{X}(k)$ to $0\in M(k)$,  and induces an isomorphism $T_{\ph}J_{X}\isom T_{0}M=M$. We also have the Koszul  resolution of $\ph_{*}\cO_{\Spec k}$
\begin{equation*}
\Kos(\ph): \wedge^{n}M^{*}\ot\cO(J_{X})\to \cdots \to \wedge^{2}M^{*}\ot \cO(J_{X})\to M^{*}\ot \cO(J_{X})\to 0
\end{equation*}
in degrees $-n,\cdots, 0$, where we note that $\dim M=\dim V=n$ since $C\cong k[W]$ as a $W$-module. We have a map of complexes
\begin{equation}\label{map Kos}
(\vk^{*}\Kos(\wt x)\ot \om_{C/k})_{W}\to \Kos(\ph)
\end{equation}
given in degree $-i$ by the map
\begin{equation*}
(\wedge^{i}V^{*}\ot \cO(J_{X,C}))_{W}=(\wedge^{i}V^{*}\ot \om_{C/k})_{W}\ot \cO(J_{X})\xr{\d^{*}_{i}}\wedge^{i}M^{*}\ot \cO(J_{X})
\end{equation*}
where $\d^{*}_{i}$ is map dual to the map $\d_{i}$ in \eqref{dj}. By our assumption, $\d_{i}$ hence $\d^{*}_{i}$ are isomorphisms, we conclude that \eqref{map Kos} is an isomorphism of complexes. It thus induces a  quasi-isomorphism in $\QCoh(J_{X})$
\begin{equation*}
c'_{\ph}: (\vk^{*}\wt x_{*}\cO_{\Spec C}\ot \om_{C/k})_{W}\to \ph_{*}\cO_{\Spec k}
\end{equation*}
One can check that the $c'_{\ph}$ is essentially $c_{\ph}$, upon identifying $W$-invariants with coinvariants. This shows that $c_{\ph}$ is an isomorphism.
\end{proof}

The proof of the full faithfulness in Theorem \ref{th:ff} is now complete.

\subsection{Characterization of essential image of $\Phi_X$}\label{ss:im}

Let us denote the full subcategory of $\QCoh(X)^W$ described in part (2) of Theorem \ref{th:ff} by $\QCoh(X)^{W,\da}$. 

\begin{lemma}\label{l:im in}
    The essential image of $\Phi_X$ is contained in $\QCoh(X)^{W,\da}$.
\end{lemma}
\begin{proof}
Let $q\in Q(K)$ with image $s\in S(K)$. Let $J_{X,s}$ be the fiber of $J_{X}\to S$ over $s$, which can be identified with the fiber $J^{(1)}_{X,q}$ of $J_{X}^{(1)}\to Q$ over $q$. We have a commutative diagram
\begin{equation*}
\xymatrix{ J_{X,s}\ar[d]^{i_{J,s}} &\ar[l]_-{\pr_{1,q}}^-{\sim} J_{X,q}^{(1)}/W_{q} \ar[r]^-{\vk_{q}}\ar[d]^{\wt i_{J, q}}& X_{q}/W_{q}\ar[d]^{i_{q}}
\\
J_{X} & \ar[l]_-{\pr_{1}} J_{X}^{(1)}/W\ar[r]^-{\vk} & X/W}
\end{equation*}
Note that the right square is derived Cartesian because both $J_{X}^{(1)}$ and $X$ are flat over $Q$. Let $\cF\in \QCoh(J_{X})$. By base change, we have an isomorphism
\begin{equation}\label{iq PhiF}
i_{q}^{*}\Phi_{X}(\cF)=i_{q}^{*}\vk_{*}\pi_{1}^{*}\cF\cong \vk_{q*}\wt i_{J,q}\pr_{1}^{*}\cF^{*}=\vk_{q*}\pr_{1,q}^{*}i_{J,s}^{*}\cF.
\end{equation}
Since $\vk_{q}$ is $W_{q}$-equivariant, it factors as $\vk_{q}: J_{X,q}^{(1)}/W_{q}\xr{\wt\vk_{q}} X_{q}^{W_{q}}/W_{q}\xr{\io_{q}} X_{q}/W_{q}$, therefore the right side of \eqref{iq PhiF} is of the form $\io_{q*}\cH$ for $\cH=\wt\vk_{q*}\pr_{1,q}^{*}i_{J,s}^{*}\cF$. The $W_{q}$-equivariant structure on $\cH$ is trivial because we may identify $J_{X,q}^{(1)}$ with $J_{X, s}$, hence $W_{q}$ acts trivially on $J_{X,q}^{(1)}$. 
\end{proof}

\sss{}
To finish the proof of part (2) of Theorem \ref{th:ff}, it remains to show that every object in $\QCoh(X)^{W,\da}$ is in the essential image of $\Phi_X$. 

Let $\cG\in \QCoh(X)^{W,\da}$, and consider the unit map $u: \cG\to \Phi_X\c\Phi^\ell_X(\cG)$. Then $\Cone(u)\in \QCoh(X)^{W,\da}$ and satisfies $\Phi^\ell_X(\Cone(u))\cong 0$ because $\Phi^\ell_X\c\Phi_X\cong \id$. If we can conclude that $\Cone(u)\cong0$, then $\cG\cong \Phi_X\c\Phi^\ell_X(\cG)$, i.e. ,$\cG$ is in the essential image of $\Phi_X$. Therefore Theorem \ref{th:ff}(2) is equivalent to:
\begin{equation}\label{ess im equiv}
\mbox{Suppose $\cG\in \QCoh(X)^{W,\da}$ is such that $\Phi^\ell_X(\cG)\cong 0$, then $\cG\cong 0$.}
\end{equation}

\sss{}  We do some reductions similar to  \S\ref{sss:Wq}. To show $\cG\cong 0$, it suffices to show that for any field valued point $x: \Spec K\to X$, the derived restriction $x^*\cG\cong 0$. Let $q=\xi(x)\in Q(K)$ be the image of $x$. We apply the discussion in \S\ref{sss:bc} to the base change $q^\flat: \Spec K\to Q^\flat$, where $q^\flat=\pi(q)$. Using that the formation of $\Phi^\ell_X$ commutes with pullbacks along such a base change (Lemma \ref{l:cX pull push}), and the condition in Theorem \ref{th:ff}(2) is also stable under pullback along such a base change,  we may assume $Q^\flat=\Spec K$. A reduction similarly to \S\ref{sss:Wq} further reduces the situation to the case $Q=\Spec C_q$ is local artinian with the unique $K=k$-point $q$ fixed by $W=W_q$. The point $x$ becomes a $k$-point of $X$. 
    
    Below we denote $C_q$ by $C$. Let $X_0=\xi^{-1}(\Spec k)$ be closed fiber of $X$, and let $i_0: X_0\incl X$ be the inclusion. Now $x\in X_{0}(k)$. In this situation, $\QCoh(X)^{W,\da}$ consists of $\cG\in \QCoh(X)^{W}$ such that $i_0^*\cG$ is in the essential image of $\io_*: \QCoh(X_0^W)\to \QCoh(X_0)^W$ (direct image under $\io: X_0^W\incl X_0$ with trivial $W$-action). If $x\notin X_{0}^{W}(k)$, the fiber $x^{*}\cG$ is automatically zero. Therefore we will assume $x\in X_{0}^{W}(k)$. We need to show
\begin{equation*}
\mbox{Assume $Q=\Spec C$ and $x\in X_{0}^{W}(k)$. If $\Phi^\ell_X(\cG)\cong 0$, then $x^{*}\cG\cong 0$.}
\end{equation*}
For this statement we may replace $X$ by any affine open neighborhood of $x$. In particular, we may assume $X$ is affine, in which case $J_X$ and $J_X^{(1)}=J_{X}\times \Spec C$ are also affine. 

Later in the argument we will compare the situation of $X$ with the situation of $V_{C}=V\times_{\Spec k}\Spec C$ (where $V$ is an affine space over $k$), related by an \'etale map $X\to  V_{C}$. 

\sss{Change of $X$} Recall that $C$ is a local artinian $k$-algebra with $W$-action such that $C^{W}=k$. Suppose $\y: Y\to \Spec C$ is another $W$-equivariant morphism from a $k$-scheme of finite type, and $f: X\to Y$ is a $W$-equivariant map over $\Spec C$. We have a commutative diagram
\begin{equation}\label{JXY}
\xymatrix{ J_{X}\ar[d]^{f_{J}} &\ar[l]_-{\pr_{1,X}}  J^{(1)}_{X}/W\ar[r]^-{\vk_X}\ar[d]^{\wt f_{J}} & X/W\ar[d]^{f}\\
 J_{Y} & J^{(1)}_{Y}/W\ar[l]_-{\pr_{1,Y}}\ar[r]^-{\vk_Y} & Y/W}
\end{equation}
The left square is derived Cartesian but the right one may not be. However, the situation is much better when $f$ is \'etale. 

\begin{lemma}\label{l:JXY}
In the above situation, assume both $\xi: X\to \Spec C$ and $\y: Y\to \Spec C$ are flat. 
\begin{enumerate}
\item $f_{*}$ sends $\QCoh(X)^{W,\da}$ to $\QCoh(Y)^{W,\da}$.

\item If $f$ is \'etale, so is $f_{J}$. Moreover, the following commutative diagram is Cartesian
\begin{equation}\label{JXY Cart}
\xymatrix{ J_{X}\ar[d]^{f_{J}}\ar[r]^{\vk_{X,0}} & X_{0}^{W}\ar[d]^{f_0^W}\\
J_{Y}\ar[r]^{\vk_{Y,0}} & Y_{0}^{W}
}
\end{equation}
Here, $\vk_{X,0}$ is the restriction of $\vk_X$ over $\Spec k\incl \Spec C$, and similarly for $\vk_{Y,0}$. 

\item If $f$ is \'etale, 
then we have a natural isomorphism of functors
\begin{equation}\label{Phi left bc}
\Phi^{\ell}_{Y}\c f_{*}\cong f_{J*}\c\Phi^{\ell}_{X}: \QCoh(X)^{W,\da}\to \QCoh(J_{Y}).
\end{equation}
\end{enumerate}
\end{lemma}
\begin{proof}
(1) Let $\cG\in \QCoh(X)^{W,\da}$. Let $q\in Q(K)$. We have a commutative diagram
\begin{equation*}
\xymatrix{X_{q}^{W_{q}}\ar[d]^{\ph_{q}}\ar[r]^-{\io_{q}} & X_{q}\ar[d]^{f_{q}}\ar[r]^-{i_{q}} & X\ar[d]^{f}\\
Y_{q}^{W_{q}}\ar[r]^-{\io'_{q}} & Y_{q}\ar[r]^-{i'_{q}} & Y}
\end{equation*}
The square on the right is derived Cartesian because both $X$ and $Y$ are flat over $Q$. Therefore base change for quasi-coherent sheaves $i_{q}'^{*}f_{*}\cong f_{q*}i^{*}_{q}$ holds. 

By assumption, $i_{q}^{*}\cG=\io_{q*}\cH$ for some $\cH\in \QCoh(X_{q}^{W_{a}})$ with the trivial $W_{q}$-action. Therefore by base change,
\begin{equation*}
i_{q}'^{*}f_{*}\cG\cong f_{q*}i^{*}_{q}\cG\cong f_{q*}\io_{q*}\cH=\io'_{q*}ph_{q*}\cH,
\end{equation*}
all with the trivial $W$-action. This shows that $f_{*}\cG\in  \QCoh(Y)^{W,\da}$.

(2) From the definition of Weil restriction one easily checks that if $f$ is formally smooth (i.e., it satisfies the existence part of the lifting property from $\ov A=A/I$-points to $A$-points for any $k$-algebra $A$ and square zero ideal $I$), then the same is true for the induced map on Weil restrictions $f_{\cR}: \cR_{X}\to \cR_{Y}$. Moreover, the fiber of $f_{\cR}$ at a geometric point $y:\Spec K\to \cR_{Y}$ (with image $s\in S(K)$) is $R_{C/K}(X_{\wt y})$, where $\Spec C$ is the fiber of $Q$ over $s$, $\wt y: \Spec C\to Y$ corresponds to $y$, and $X_{\wt y}=X\times_{Y,\wt y}\Spec C$. Since $f$ is \'etale, $X_{\wt y}$ is a disjoint union of finitely many copies of $\Spec C$ (since $C$ is local artinian with an algebraically closed residue field),   we conclude that $f_{\cR}^{(-1)}(y)\cong R_{C/K}(X_{\wt y})$ is also a disjoint union of finitely many copies of $\Spec K$. Combined with the formal smoothness, we conclude that $f_{\cR}$ is \'etale.  In particular, $f_{\cR}$ is formally \'etale, i.e., it satisfies the existence and uniqueness of the lifting property from $\ov A=A/I$-points to $A$-points for any $k$-algebra $A$ and square zero ideal $I$. Passing to $W$-fixed points, we see that $f_{J}$ is also formally \'etale. Since $f_{J}$ is of finite presentation, it is \'etale.

Now $f_0: X_0\to Y_0$ is \'etale implies that $f_0^W:X_0^W\to Y_0^W$ is \'etale. Since both $f_J$ and $f_0^W$ are \'etale, to check that the commutative diagram \eqref{JXY Cart} is Cartesian, it suffices to check at the level of geometric points. Let $y: \Spec K\to J_Y$ be a geometric point of $J_Y$ with image $y_0: \Spec K\to Y_0^W$. By definition, $y$ can be viewed as a $W$-equivariant map $\wt y: \Spec (K\ot_k C)\to Y$ whose restriction to $\Spec K\incl \Spec (K\ot_k C)$ is $y_0$. The $K$-points of the fiber $f_J^{-1}(y)$ are the $W$-equivariant liftings $\wt x: \Spec (K\ot_k C)\to X$ of $\wt y$. The map
\begin{equation}\label{f fiber}
    f_J^{-1}(\wt y)(K)\to f_0^{W,-1}(y_0)(K)
\end{equation}
sends such a lifting $\wt x$ to its restriction to the reduced structure of $\Spec K\incl \Spec (K\ot_k C)$. We need to show that \eqref{f fiber} is a bijection. Since $f$ is \'etale, $\wt x$ is determined by its restriction to the reduced structure. This shows that \eqref{f fiber} is injective. Conversely, given $x_0: \Spec K\to X_0^W$ lifting to $y_0$, it extends uniquely as a map $\wt x: \Spec(K\ot_kC)\to X$ lifting $\wt y$; the uniqueness of such an extension forces it to be $W$-equivariant, i.e., $\wt x\in f_J^{-1}(\wt y)(K)$. This shows that \eqref{f fiber} is surjective, hence a bijection.

(3) Let $J'$ be the fiber product $J'=J_Y^{(1)}\times_{Y}X$, then we have a natural $W$-equviariant map over $Y$
\begin{equation*}
    j: J_X^{(1)}\to J'.
\end{equation*}
Snice both $J_X^{(1)}$ and $J'$ are \'etale over $Y$, so is $j$. Passing to the fiber over $\Spec k\incl \Spec C$, $j$ becomes the map
\begin{equation*}
    j_0: J_X^{(1)}\ot_C k=J_X\to J'\ot_Ck\cong J_Y\times_YX=J_Y\times_{Y_0}X_0=:J'_0.
\end{equation*}
By (2), $j_0$ is the inclusion
\begin{equation}\label{j0 XW}
    J_X\cong J_Y\times_{Y_0^W}X_0^W\incl J_Y\times_{Y^W_0}(X_0\times_{Y_0}Y_0^W)=J_Y\times_{Y_0}X_0.
\end{equation}
In particular, $j_0$ is a closed embedding. Therefore $j_0$ is both open and closed. Hence $j$ is both open and closed (being \'etale). We can summarize the situation by a diagram where all maps are $W$-equivariant
\begin{equation*}
    \xymatrix{J_X^{(1)}\ar@{^(->}[r]^j\ar[dr]_{\wt f_J}\ar@/^2pc/[rr]^{\vk_X} & J'\ar[d]^{f'}\ar[r]^{\vk'} & X\ar[d]^f\\
    & J_Y^{(1)}\ar[r]^{\vk_{Y}} & Y}
\end{equation*}

For $\cG\in \QCoh(X)^{W}$, we have   isomorphisms by base change
\begin{equation}\label{ffbc}
    \vk_Y^*f_*\cG\cong f'_*\vk'^*\cG, \quad \wt f_{J,*}\vk^*\cG \cong f'_*j_*j^*\vk'^*\cG.
\end{equation}
Now suppose $\cG\in \QCoh(X)^{W,\da}$. We claim that the unit map
\begin{equation}\label{jj}
    \vk'^*\cG\to j_*j^*\vk'^*\cG
\end{equation}
is an isomorphism. 

\begin{claim} $\cG$ is set-theoretically supported on $X_0^W$, i.e., $\cG|_U\cong 0$ where $U=X-X_0^W$.
\end{claim}
\begin{proof}[Proof of Claim]
    Indeed, letting $i_{U,0}: U_0=U\times_{\Spec C}\Spec k\incl U$ be the inclusion, we have $i^*_{U,0}(\cG|_U)\cong (i_0^*\cG)|_{U_0}$, and the latter is zero because by assumption $i_0^*\cG$ is supported on $X_0^W$. Now we observe that $i_{U,0}^*: \QCoh(U)\to \QCoh(U_0)$ is conservative. Indeed, it suffices to treat the case where $U=\Spec A$ is affine, in which case $U_0$ is defined by a nilpotent ideal $I\subset A$. We have a finite filtration of $A$ by $A/I$-modules, namely $I^n/I^{n+1}$ for $n=0,1,\cdots, N$ (where $N$ is such that $I^{N+1}=0$). Now suppose $M\in \QCoh(U)\cong D(A\lmod)$ is such that $M\ot^{\bL}_A(A/I)\cong 0$, then $M$ is a successive extension of $M\ot^{\bL}_A(I^n/I^{n+1})\cong M\ot^{\bL}_A(A/I)\ot^{\bL}_A/I(I^n/I^{n+1})\cong 0$, hence $M\cong0$. 
\end{proof}

Recall that $j$ is both open and closed. To show \eqref{jj} is an isomorphism, it suffices to show that $\vk'^*\cG|_{J'-J^{(1)}_X}\cong 0$. By \eqref{j0 XW}, the underlying space of $J'-J^{(1)}_X$ does not meet $X_0^W$. Therefore \eqref{jj} is an isomorphism. By \eqref{ffbc}, we conclude with a natural isomorphism
\begin{equation*}
    \wt f_{J,*}\vk_X^*\cong \vk_Y^*f_*: \QCoh(X)^{W,\da}\to \QCoh(J_Y^{(1)}).
\end{equation*}
Applying $\pr_{1,Y*}$ to both sides gives the desired isomorphism \eqref{Phi left bc}.
\end{proof}

\sss{Reduction to vector space case}
Let $V=T_{x}X_{0}$, a $k$-vector space which by assumption is a $C$-admissible $W$-module. Let $\fm_{x}$ (resp. $\wt\fm_{x}$) be the maximal ideal of $\cO_{X_{0},x}$ (resp. $\cO_{X, x}$), then we have a $W$-equivariant surjection $\wt\fm_{x}\surj \fm_{x}\surj \fm_{x}/\fm_{x}^{2}=V^{*}$. Choose a $W$-equvariant section $s: V^{*}\to \wt\fm_{x}$ to this surjection. By shrinking $X$ we may assume $s(V^{*})$ consists of regular functions on $X$, therefore giving a map $\s: X\to V$ that maps $x$ to $0$. Let $f=(\s, \xi): X\to V_{C}:=V\times_{k}\Spec C$. By construction, $\s|_{X_{0}}$ is \'etale at $x$, hence $f$ is also \'etale at $x$. By shrinking $X$ further (but still affine and contains $x$), we may assume $f$ is \'etale. 
	
	Note that $X=V_{C}$ satisfies the assumptions of Theorem \ref{th:ff}. Indeed, we only need to check the second condition for $x\in V^{W}(k)$, in which case $T_{x}(X_0)=V$ and it is $C$-admissible by assumption.
	
	We continute with the proof of \eqref{ess im equiv}. Since $\cG\in \QCoh(X)^{W,\da}$, we have $f_{X*}\cG\in \QCoh(V_{C})^{W,\da}$ by Lemma \ref{l:JXY}(1). Since
    $\Phi_{X}^{\ell}(\cG)\cong 0$, by Lemma \ref{l:JXY}(3), we have
\begin{equation*}
\Phi_{V_{C}}^{\ell}(f_{*}\cG)\cong f_{J*}\Phi^{\ell}_{X}(\cG)\cong 0.
\end{equation*}
If \eqref{ess im equiv} holds for $V_{C}$, then $f_{*}\cG\cong 0$. Since $f$ is a map between affine schemes, $f_{*}\cG\cong 0$ implies $\cG\cong 0$. Therefore the proof of Theorem \ref{th:ff}(2) reduces to the special case $X=V_{C}$, which we prove below.

\begin{lemma} Let $V$ be a $C$-admissible finite-dimensional $k$-representation of $W$, viewed as an affine space over $k$. Consider the situation $X=V_{C}=V\times\Spec C$. Then
\begin{enumerate}
\item $J_{V_{C}}$ is the affine space given by the $k$-vector space $M=(V\ot C)^{W}$.
\item Suppose $\cG\in \QCoh(V_{C})^{W,\da, \hs}$ is in the heart of the usual $t$-structure on $\QCoh(V_{C})^{W,\da}$, then $N:=\G(V_C,\cG)^{W}$ has a canonical structure of a $\Sym(M^{*})$-module, making it an object $N\in \QCoh(J_{V_{C}})$. There is a canonical isomorphism $\Phi_{V_{C}}(N)\cong \cG$.
\item The essential image of $\Phi_{V_{C}}$ is $\QCoh(V_C)^{W,\da}$.
\end{enumerate}
\end{lemma}
\begin{proof}
(1) It is a basic exercise that $R_{C/k}(V_{C})$ is the affine space attached to the $k$-vector space $V\ot C$. Taking $W$-invariants we get the desired description of $J_{V_{C}}$.

(2) Let $L=\Gamma(V_C,\cG)$ as a $\Sym(V^*)\ot C$-module. The condition $\cG\in \QCoh(V_C)^{W,\da}$ implies that $\Tor^C_i(L,k)$ has trivial $W$-action for all $i$. We claim that in this case the canonical $C$-linear map
\begin{equation}\label{aL}
    a_L: L^W\ot C\to L
\end{equation}
is an isomorphism. 

Indeed, taking reduction mod the maximal ideal $\fm_C$ of $C$, $a_L$ becomes the canonical map $L^W\to L\ot_C k$, which is surjective because the $W$-action on $L\ot_C k$ is trivial, and taking $W$-invariants is an exact functor. Since $\fm_C$ is nilpotent, we conclude that $a_L$ is also surjective. 

Now let $L_1=\ker(a_L)$, a $C$-submodule of $L^W\ot C$. Taking $W$-invariants for the short exact sequence
\begin{equation}\label{L1L}
    0\to L_1\to L^W\ot C\to L\to 0,
\end{equation}
we get $L_1^W=0$. On the other hand, applying $\ot^{\bL}_C k$ to \eqref{L1L}, we get a short exact sequence
\begin{equation*}
    0\to \Tor^C_1(L,k)\to L_1\ot_C k\to L^W \to L\ot_C k\to 0.
\end{equation*}
Since the $W$-actions on both $\Tor^C_1(L,k)$ and $L^W$ are trivial, the same is true for $L_1\ot_C k$. Repeating the argument from the previous paragraph, we see that $L_1$ is generated by $L_1^W$ as a $C$-module. Since $L_1^W=0$, we have $L_1=0$ and $a_L$ is an isomorphism.

Now we give $N=L^W$ a module structure over $\Sym(M^*)$. The commuting actions of $\Sym(V^*)$ and $C$ on $L\cong N\ot C$ gives a $W$-equivariant $k$-algebra homomorphism
\begin{equation*}
    \Sym(V^*)\to \End_C(N\ot C)\cong \End_k(N)\ot C.
\end{equation*}
This in particular induces a $W$-equivariant $k$-linear map
\begin{equation*}
    V^*\to \End_k(N)\ot C, \quad \mbox{or} \quad  V^*\ot C^*\to \End_k(N).
\end{equation*}
By $W$-equivariance, the above map factors through the $W$-coinvariants, giving a $k$-linear map
\begin{equation}\label{M*action}
    M^*=(V^*\ot C^*)_W\to \End_k(N).
\end{equation}
It remains to argue that \eqref{M*action} gives {\em commuting} actions of elements of $M^*$ on $N$. For a $k$-vector space $U$, we use $\Alt^i(U)\subset U^{\ot i}$ to denote the subspace of $U^{\ot i}$ on which $S_i$ acts through the sign character. The pairwise commutators of the action of $M^*$ on $N$ can be packaged into a $k$-linear map
\begin{equation*}
    c: \Alt^2(M^*)\to \End_k(N).
\end{equation*}
Note the canonical map
\begin{equation}\label{delta 2}
    \d^*_2: (\Alt^2(V^*)\ot C^*)_W\to \Alt^2((V^*\ot C^*)_W)=\Alt^2(M^*)
\end{equation}
that is dual to the map $\d_j$ of \eqref{dj}. 
Precomposing $c$ with \eqref{delta 2} we get
\begin{equation*}
    c'=c\c\d^*_2: (\Alt^2(V^*)\ot C^*)_W\to \End_k(N).
\end{equation*}
Precomposing $c'$ with the projection from $\Alt^2(V^*)\ot C^*$, we obtain a map
\begin{equation*}
    c'': \Alt^2(V^*)\to \End_k(N)\ot C=\End_C(N\ot C)
\end{equation*}
that is the pairwise commutator of the action of $V^*$ on $N\ot C$. By assumption, the action of $V^*$ extends to $\Sym(V^*)$, hence $c''=0$, which implies $c'=0$. Since $\d^*_2$ is an isomorphism by assumption, we conclude that $c=0$, i.e., the map \eqref{M*action} extends to an action of $\Sym(M^*)$.

Finally, a direct calculation shows that $\Phi_{V_C}(N)\cong N\ot C$, equipped with the $W$-action on $C$, the obvious $C$-action and the action of $\Sym(V^*)$ induced by
\begin{equation*}
    V^*\xr{\id_{V^*}\ot \coev_C} V^*\ot C^*\ot C\to M^*\ot C\to \End_C(N\ot C).
\end{equation*}
The map $a_L:N\ot C\isom L$ in \eqref{aL} is easily seen to be a $W$-equivariant $\Sym(V^*)\ot C$-module map. Since $a_L$ is an isomorphism, we conclude that $\Phi_{V_C}(N)\cong N\ot C\cong L$ as objects in $\QCoh(V_C)^{W,\da}$.

(3) From (2), we know that $\QCoh(V_{C})^{W,\da, \hs}$ is in the essential image of $\Phi_{V_{C}}$. We follow a similar argument as in the proof of \cite[Proposition 4.23 (6)]{Gan}. Since the $t$-structure on $\QCoh(V_{C})^{W}$ is both left and right-complete, and $\Phi_{V_{C}}$ is a fully faithful continuous $t$-exact functor, for any $\cG\in  \QCoh(V_{C})^{W}$, writing $\cG\cong \lim_m\colim_n\tau^{\geq m} \tau^{\leq n}\cG$, we see that if all cohomology sheaves $\cH^j(\cG)$ are in the essential image of $\Phi_{V_{C}}$, so are all the finite extensions $\tau^{\geq m} \tau^{\leq n}\cG$ (using $\Phi_{V_{C}}$ is fully faithful), and hence $\cG$ is in the essential image of $\Phi_{V_{C}}$ (using $\Phi_{V_{C}}$ preserves both limits and colimits). 

\end{proof}

%%%%%%%%%%%%%%%%%%%%%%%%%%%%%
\section{Quasi-coherent sheaves on regular centralizers}\label{s:QCoh J} %%%%%%%%%%%%%%%%%%%%%%%%%%%%%%%%%

\subsection{Regular centralizer and Weil restriction}\label{ss:J to R}

Let $G$ be a connected reductive group over $\CC$ with maximal torus $T$ and Weyl group $W$.

Let $\nu: G\to G_\flat$ be a central quotient. We have the regular centralizer group scheme $J^{G_\flat}_G$ over $\cS_G$ as defined in \S\ref{sss:reg cent gen}.

\sss{The canonical map $j$} Let $T$ be the abstract Cartan of $G$, and let $T_\flat$ be that of $G_\flat$. Let $W$ be the abstract Weyl group that acts on both $T$ and $T_\flat$. Then we have a canonical map
\begin{equation}\label{J to R}
    j^{G_\flat}_G: J^{G_\flat}_G\to \left(R_{T/\cS_G}(T_\flat\times T)\right)^W=:J_{T_\flat, T}
\end{equation}
where the action of $W$ on $R_{T/\cS_G}(T_\flat\times T)$ is induced from its action on both $T$ and $T_\flat$. To define the map \eqref{J to R}, it suffices to give a $W$-equivariant map
\begin{equation}\label{can j}
    (J^{G_\flat}_G)\times_{\cS_G}T\to T_\flat.
\end{equation}
Here the $W$-action on the left side is only on the $T$-factor. 

To construct \eqref{can j}, we first consider the case $G$ is simply-connected. Let $S=S_{\dot w}\subset G$ be a Steinberg cross-section as in \S\ref{sss:St section}, so that the projection $S\to \cS_G$ is an isomorphism. We identify $J^{G_\flat}_G$ with the universal centralizer group scheme of $S$ under the conjugation action of $G_\flat$. Recall the Cartesian diagram
\begin{equation*}
        \xymatrix{\wt{G^{\reg}} \ar[d]^{\pi}\ar[r] & T\ar[d]\\
        G^{\reg}\ar[r]^{s_G} & \cS_G}
    \end{equation*}
Restricting it over $S\subset G^\reg$, we get an isomorphism
\begin{equation*}
    \wt S:=S\times_{G^\reg}\wt {G^{\reg}}\isom T.
\end{equation*}
The above map is defined as follows: $\wt S$ classifies pairs $(g,B)$ where $g\in S$ and $B$ is a Borel subgroup of $G$ containing $g$. The above map sends $(s,B)$ to the image of $s$ under the canonical map $B\surj T$. Using this isomorphism we may identify $J^{G_\flat}_G\times_{\cS_G}T$ with
\begin{equation*}
    \wt J:=J^{G_\flat}_G\times_{S}\wt S,
\end{equation*}
which classifies triples $(g,y,B)$ where $g\in \cS$, $y\in G_\flat$ commuting with $g$, and $B$ a Borel subgroup of $G$ containing $g$. Since $g$ is regular, $C_G(g)\subset B$, hence $C_{G_\flat}(g)\subset B_\flat=\nu(B)$. The map \eqref{can j} sends $(g,y,B)$ to the image of $y\in B_\flat$ under the canonical projection $B_\flat\surj T_\flat$. It is easy to see that this assignment is $W$-equivariant when restricted to the open dense locus where $g$ is regular semisimple, hence it is $W$-equivariant everywhere. This completes the construction of \eqref{can j} for simply-connected groups.

Next consider the case where $G=G^\sc\times A$ where $A$ is a torus. In this case, we can repeat the above construction, with $S$ replaced by $S_{\dot w}\times A$.

Finally, for general reductive $G$, consider the central isogeny $\nu: G_1:=G^\sc\times (ZG)^\c\to G$. The previous paragraph gives a $W$-equivariant map
\begin{equation}\label{can j wt}
    (J_{G_1}^{G_\flat})\times_{\cS_{G_1}}\wt T\to T_\flat
\end{equation}
which is also $\ker(\nu)$-invariant for the diagonal action of $\ker(\nu)$ on the left side. Therefore \eqref{can j wt} descends to the desired map \eqref{can j}. This completes the construction of $j^{G_\flat}_G$ in \eqref{J to R}.

\sss{Variant of $j$}
Let $H\subset G$ be a Levi subgroup. Note that $H^\der$ is also simply-connected. Let $H_\flat=H/\ker(\nu)$. We have a generalization of the map $j^{G_\flat}_G$ to a canonical map of group schemes over $\cS_H$
\begin{equation*}
    j^{G_\flat, H_\flat}_{G,H}: J^{G_\flat}_G\to R_{\cS_H/\cS_G}(J^{H_\flat}_H)^{W(G,H)}.
\end{equation*}
Here $W(G,H)=N_G(H)/H$. The construction of $j^{G_\flat, H_\flat}_{G,H}$ depends on the choice of a parabolic subgroup $P$ that contains $H$ as a Levi factor. Indeed, let $\wt G_P$ be the moduli space of pairs $(g,Q)$ where $g\in G$ and $Q$ is a parabolic subgroup of $G$ containing $g$ and conjugate to $P$. Then have an isomorphism $\cS_H\cong S\times_{G^\reg}\wt G^\reg_P$, i.e., the moduli of pairs $(g,Q)$ where $g\in S$ and $Q\ni g$ is a parabolic subgroup in the class $P$. We thus have a map $\wt j: J^{G_\flat}_G\times_{\cS_G}\cS_H\cong J^{G_\flat}_G\times_{G^{\reg}}\wt G^{\reg}_P\to J^{H_\flat}_H$ sending $(g,y,Q)$ (where $g\in S$, $y\in C_{G_\flat}(g)$ and $Q\ni g$ is a parabolic in class $P$) to $(\ov g, \ov y)\in J^{H_\flat}_H$: here we choose a surjection $\pi: Q\surj H$ realizing $H$ as the Levi quotient of $Q$ (so $\pi$ is well-defined up to $H$-conjugation), $\ov g=\pi(g)\in H^\reg$, and $\ov y=\pi_\flat(y)\in H_\flat$ (noting that $y$ has to lie in $Q_\flat=Q/\ker(\nu)$, and $\pi_\flat: Q_\flat\to H_\flat$ is induced from $\pi$), so that the pair $(\ov g,\ov y)$ is well-defined up to simultaneous conjugation by $H$, hence defining a point in $J^{H_\flat}_H$.

The following commutative diagram can be verified from the construction
\begin{equation}\label{jHG}
    \xymatrix{J^{G_\flat}_G\ar[rrr]^{j^{G_\flat}_G} \ar[d]^{j^{G_\flat, H_\flat}_{G,H}} && & R_{T/\cS_G}(T_\flat\times T)^W\ar@{^(->}[d]\\
     R_{\cS_H/\cS_G}(J^{H_\flat}_H)^{W(G,H)}\ar[rr]^-{R_{\cS_H/\cS_G}(j^{H_\flat}_H)} && R_{\cS_H/\cS_G}(R_{T/\cS_H}(T_\flat\times T)^{W_H})^{W(G,H)}\ar@{=}[r] & R_{T/\cS_G}(T_\flat\times T)^{N_W(W_H)}}
\end{equation}
 
\begin{prop}\label{p:J to R open} Assume $G^\der$ is simply-connected. For any central quotient $\nu: G\to G_\flat$, the map $j^{G_\flat}_G$ as in \eqref{J to R} is an open embedding whose image is characterized as follows: for any $k$-scheme $S'$, an $S'$-point of $J_{T_\flat, T}=\left(R_{T/\cS_G}(T_\flat\times T)\right)^W$ given by a $W$-equivariant map $S'\times_{\cS_G}T\to T_\flat$ is in the image of $j^{G_\flat}_G$ if and only if for any (abstract) root $\a$ of $(G,T)$ with corresponding reflection $s_\a\in W$, the composition
\begin{equation}\label{S'}
    S'\times_{\cS_G}T^{s_\a}\to T_\flat^{s_\a}\xr{\a}\mu_2
\end{equation}
is the constant map to the identity element in $\mu_2$.

In particular, if $\a^\vee$ is primitive in $\xcoch(T_\flat)$ (for example if $G_\flat^\der$ is simply-connected), the map $j^{G_\flat}_G$ is an isomorphism.
\end{prop}
\begin{proof}
    The analogous statement for the regular centralizer group scheme of the Lie algebra is proved in \cite{DG}. We will adapt the proof given in \cite[Prop.2.4.7]{NgoFL} to the group situation.

    To simplify notations, we denote $J^{G_\flat}_G$ by $J$ and $\left(R_{T/\cS_G}(T\times T)\right)^W$ by $R$. Let $R'\subset R$ be the sub-functor defined by the condition specified in the statement of the proposition. We claim that $R'$ is open in  $R$. Indeed, taking $S'=R$ in \eqref{S'}, we obtain a map $\e_\a: R\times_{\cS_G} T^{s_\a}\to \mu_2$. Let $\wt R^-_\a\subset R\times_{\cS_G} T^{s_\a}$ be the preimage of $-1$ under the map $\e_\a$, and let $R^-_\a\subset R$ be the image of $\wt R^-_\a$. Since $T^{s_\a}\to \cS_G$ is finite, $R^-_\a\subset R$ is closed. Finally, $R'$ is the complement of $\cup_{\a\in \Phi}R^-_\a$ ($\a$ runs over the set $\Phi$ of roots of $G$), hence open. 

    We claim that the image of $j^{G_\flat}_G$ lands in $R'$. Indeed, by the definition of $J^{G_\flat}_G$ as the quotient of $J^G_G/\ker(\nu)$, it suffices to check the statement for $G_\flat=G$, in which case $R'=R$ and the statement is obvious.
    
    We now show that the induced map $j': J\to R'$ is an isomorphism. We use $R'^{G_\flat}_G$ and $R'_G$ to emphasize the dependence of $R'$ on $G_\flat$. It is easy to see that $R'^{G_\flat}_G\cong R'_G/\ker(\nu)$. Since $J^{G_\flat}_G\cong J^G_G/\ker(\nu)$, it suffices to show that $j'$ is an isomorphism when $G=G_\flat$, which we will assume from now on. Our goal is thus to show that $j=j^G_G$ is an isomorphism. 

    By Lemma \ref{l:group J sc}, $J$ is 
    smooth over $\cS_G$. Also $R$ is smooth over $\cS_G$ since Weil restriction preserves formal smoothness, and taking $W$-invariants preserves smoothness. Therefore, it suffices to show $j$ is an isomorphism over an open subset $U\subset \cS_G$ whose complement has codimension at least $2$. Let $T^{\vs}\subset T$ be the union of intersections $T^{s_\a}\cap T^{s_\b}$, where $(\a,\b)$ runs over pairs of roots such that $\a\ne \pm \b$. Then $T^\vs$ has codimension $2$ in $T$. Let $S^\vs_G\subset \cS_G$ be the image of $T^\vs$, then $S^\vs_G$ has codimension $2$, and we may take $U=\cS_G-\cS_G^\vs$. 

    To check $j|_U$ is an isomorphism, it suffices to check that $j_a: J_a\to R_a$ is an isomorphism for each geometric point $a\in U$. Indeed, both $J$ and $R$ are flat over $\cS_G$, so by the fibral criterion of flatness,  the fact that $j_a$ is an isomorphism for all geometric points $a\in U$ implies that $j|_U$ is flat, \'etale and hence an isomorphism. 

    Now let $a\in U(k)$ be a geometric point. If $a$ does not lie in the discriminant divisor of $\cS_G$, then $j_a$ is clearly an isomorphism. Now if $a$ does lie on the discriminant divisor of $\cS_G$, we may assume $a$ is the image of $\wt a\in T^{s_\a}(k)$ for some root $\a$. Moreover, since $\wt a\notin T^\vs$, such $\a$ is unique up to sign. Choose an embedding $T\incl G$ as a maximal torus, and identify $\wt a$ as a semisimple element of $G$. Let $H=C_G(\wt a)$. Since $G^\der$ is simply-connected, $H$ is connected reductive, and is in fact a Levi subgroup of $G$ whose roots with respect to $T$ are $\pm\a$. Moreover $H^\der$ is also simply-connected (so isomorphic to $\SL_2$). Let $a_H\in \cS_H(k)$ be the image of $\wt a$. Consider the projections
    \begin{equation*}
        \xymatrix{T\ar[r]^{\chi^T_H}\ar@/_2ex/[rr]_{\chi^T_{G}} & \cS_H \ar[r]^{\chi^H_G} & \cS_G}
    \end{equation*}
    By our assumption that $\wt a\notin T^\vs$,  $\chi^H_G$ is \'etale over $a$, $\chi^{H,-1}_G(a)$ is a discrete scheme that can be identified with the cosets $W/W_H=W/\j{s_\a}$, with $a_H$ being the unit coset. Let $\chi^{T,-1}_G(a)=\Spec C_a$ and $\chi^{T,-1}_H(a)=\Spec C_{a_H}$. Then
    \begin{equation}\label{isom C}
        \Spec C_a=W\times^{\j{s_\a}}\Spec C_{a_H}.
    \end{equation}
    Taking the fibers of every space in the diagram \eqref{jHG} we get a commutative diagram
    \begin{equation*}
        \xymatrix{J_a\ar[r]^-{j_a}\ar[d]^{j^{H}_{G,a}|_{a_H}} & (R_{C_a/k}(T_\flat\times\Spec C_a))^W\ar[d]^{b}\\
        J_{H,a_H}\ar[r]^-{j_{H,a_H}} & (R_{C_{a_H}/k}(T_\flat\times \Spec C_{a_H}))^{s_\a}}
    \end{equation*}
    The isomorphism $b$ is induced from the isomorphism \eqref{isom C}. Thus it suffices to check that $j_{H,a_H}$ is an isomorphism. This further reduces to the case $H=\SL_2$ and $a\in \cS_H$ is the image of a regular unipotent element, for which one checks that $j_{H,a_H}$ is an isomorphism explicitly.
\end{proof}

\subsection{Descent for equivariant quasi-coherent sheaves on $T$}

We keep the notations from \S\ref{ss:J to R}. We shall consider the problem of descent of quasi-coherent sheaves along the canonical map
\begin{equation*}
    p_G: T/W\to \cS_G
\end{equation*}
as well as an affine analog of this problem.

\begin{lemma}\label{l:WJ des}
The pullback $p_G^{*}: \QCoh(\cS_G)\to \QCoh(T)^{W}$ is fully faithful. An object $\cF\in \QCoh(T)^{W}$ lies in the essential image of $p_G^{*}$ if and only if for any reflection $r_\a\in W$ with respect to a root $\a$, letting $\io_{\a}$ be the inclusion of $\ker(\a)\incl T$, the following holds:
\begin{equation}\label{ref action G}
    \mbox{$r_\a$ acts trivially on each cohomology sheaf of the derived restriction $\io_{\a}^{*}\cF$.}
\end{equation}   
\end{lemma}
\begin{proof} 

Fully faithfulness:  since $\cO_{\cS_G}$ generates $\QCoh(\cS_G)$, it suffices to check that $p^*_G$ induces a quasi-isomorphism 
\begin{equation}\label{End isom}
    \bR\End_{\cS_G}(\cO_{\cS_G})\to \bR\End_{T}(\cO_T)^W.
\end{equation}
Both sides above are concentrated in degree zero, and the above map can be identified with the pullback of regular functions
\begin{equation*}
    \cO(\cS_G)\to \cO(T)^W.
\end{equation*}
Recall $\cS_G=(T_1\sslash W)/\Gamma$ for a central isogeny $G_1\to G$ with kernel $\Gamma$ (with $G_1^\der$ simply-connected). Therefore $\cO(\cS_G)\cong \cO(T_1)^{W\times \Gamma}=\cO(T)^W$. This proves that \eqref{End isom} is an isomorphism and hence $p^*_G$ is fully faithful.

Now we characterize the essential image of $p_{G}^{*}$. first assume $G^\der$ is simply-connected. In this case we have $\cS_G=T\sslash W$. This case is a variant of \cite[Theorem 1.1]{Lon} (see also \cite[Prop 2.15(5)]{Gan}). First, the argument of \cite[Lemma 1.6]{Gan} shows that $\cF$ is in the essential image of $p_{G}^{*}$ iff for all $x\in T(\CC)$, the action of the stabilizer $W_{x}$ on $H^{*}(i_{x}^{*}\cF)$ is trivial, where $i_x:\{x\}\incl T$. Since $G^\der$ is simply-connected, $W_{x}$ is generated by reflections. Therefore the above condition is equivalent to that for any $x\in A(\CC)$ and any reflection $r\in W_{x}$, $r$ acts on $H^{*}(i_{x}^{*}\cF)$ trivially.
Clearly this is equivalent to that for any reflection $r_\a\in W$, $r_\a$ acts on each cohomology sheaf of $i_{r_\a}^{*}\cF$ trivially, where $i_{r_\a}: T^{r_\a}\incl T$ is the inclusion of $r_\a$-fixed points. Since each coroot $\a^\vee$ is a primitve element in $\xcoch(T)$, we have $\ker(\a)=T^{r_\a}$, hence the above condition is further equivalent to \eqref{ref action G}.

Now consider the general case. Choose a central isogeny $\nu:G_1\to G$ with kernel $\Gamma$ such that $G_1^\der$ is simply-connected. Let $T_1\subset G_1$ be the preimage of $T$. Now $\Gamma$ acts on both $\QCoh(\cS_{G_1})$ and $\QCoh(\wt T)^W$, such that
\begin{equation*}
    p_{G_1}^*: \QCoh(\cS_{G_1})\to \QCoh(\wt T)^W
\end{equation*}
is $\Gamma$-equivariant. Taking $\Gamma$-equivariant objects on both sides above we recover $p^*_G$. Since $p_{G_1}^*$ is fully faithful, we conclude that $\Im(p^*_G)$ is the category of $\Gamma$-equivariant objects in $\Im(p_{G_1}^*)$. Therefore, $\cF\in \QCoh(T)^W$ lies in the essentially image of $p_G^*$ if and only if $\nu^*\cF\in \QCoh(\wt T)^W$ lies in the essential image of $p_{G_1}^*$. By the previous paragraph, the latter happens if and only if for each root $\a$,
\begin{equation}\label{ref action wt G}
    \mbox{$r_\a$ acts trivially on the cohomology sheaves of $\wt\io_{\a}^*\nu^*\cF$ for each root $\a$,}
\end{equation}
where $\wt\io_\a: \wt{\ker(\a)}:=\ker(\a: \wt T\to \Gm)\incl \wt T$. 

We claim that \eqref{ref action wt G} is equivalent to \eqref{ref action G}. Indeed, let $\nu_\a: \wt{\ker(\a)}\to \ker(\a)$ be the restriction of $\nu$. Then $\nu_\a$ is a $\j{r_\a}$-equivariant $\Gamma$-torsor. If \eqref{ref action G} holds, then $r_\a$ acts trivially on cohomology sheaves of $\nu^*_\a\io_\a^*\cF\cong \wt\io_{\a}^*\nu^*\cF$, which implies \eqref{ref action wt G}. Conversely, if \eqref{ref action wt G} holds, then $r_\a$ acts trivially on cohomology sheaves of $\nu_{\a*}\wt\io_{\a}^*\nu^*\cF\cong \nu_{\a*}\nu^*_\a\io_\a^*\cF$, which contains $\io_\a^*\cF$ as a summand. Therefore \eqref{ref action G} holds. This concludes the characterization of the image of $p^*_G$.
\end{proof}

Next we move on to an affine analog of the above lemma.

\sss{$\QCoh(T_{\flat}\times T)^W$ as a limit}

Let $T_\flat$ be another torus with $W$-action with a given $W$-equivariant surjective map $T_\flat\to T_\ad$. Consider the category
\begin{equation*}
    \QCoh(T_\flat\times T)^W
\end{equation*}
where $W$ acts diagonally. 

Denote
\begin{equation*}
    \wt W_\flat:=\xch(T_\flat)\rtimes W.
\end{equation*}
This contains the affine Weyl group
\begin{equation*}
    \Wa:=\xch(T_{\ad})\rtimes W=(\mbox{root lattice of $G$})\rtimes W
\end{equation*}
as a normal subgroup.

Since $T_{\flat}=\Spec \CC[\xch(T_{\flat})]$, we have
\begin{equation*}
\QCoh(T_\flat\times T)\cong \QCoh(T)^{\xch(T_{\flat})}
\end{equation*}
with respect to the trivial action of $\xch(T_{\flat})$ on $T$. Adding $W$-equivariance we get
\begin{equation}\label{Wa equiv}
\QCoh(T_\flat\times T)^W\cong \QCoh(T)^{\wt W_\flat}
\end{equation}
where the action of $\tilW_\flat$ factors through the usual action of $W$ via $\tilW_\flat\surj W$.

For each $J\sft \wt I$, restricting the equivariance to the subgroup $W_{J}\subset \Wa\subset \wt W_\flat$ to get a functor
\begin{equation}\label{eq: r_J}
r_{J}: \QCoh(T_\flat\times T)^W\cong \QCoh(T)^{\tilW_\flat}\to \QCoh(T)^{W_{J}}.
\end{equation}
These functors induce a canonical symmetric monoidal functor
\begin{equation}\label{dT to lim}
\QCoh(T_\flat\times T)^{W,\otimes_c}\to \lim_{J\sft \wt I}\QCoh(T)^{W_{J},\otimes}.
\end{equation}
where $\otimes_c$ is the convolution monoidal structure for the group stack $(T_\flat\times T)/W$ over $T/W$.

\begin{lemma}\label{l:lim Wa}
    When $G_\flat=G_\ad$ (hence $T_\flat=T_\ad$ and $\tilW_\flat=\Wa$), the functor \eqref{dT to lim} is an equivalence.
    Forgetting the (symmetric) monoidal structure, we get a canonical equivalence between $\QCoh((T_\flat\times T)/W)$ and $ \lim_{J\sft \wt I}\QCoh(T)^{W_{J}}$. 
   
\end{lemma}

\begin{proof}
Let $\CAlg(\Pr^L)$ denote the $\infty$-category of commutative algebra objects in $\Pr^L$. Since $\CAlg(\Pr^L)\to \Pr^L$ preserves limits, it suffices to check that the functor is an equivalence in $\Pr^L$. 

By \cite[Corollary 4.2]{PL} we have
\begin{equation}\label{eq: BBWJ,BBWa}
\colim_{J\sft \wt I}\BB W_{J}\isom\BB \Wa
\end{equation}
in the category of $\infty$-groupoids. Let $\widehat{\Cat}_\infty$ (resp. $\Gpd_\infty$) be the $\infty$-category of (not necessarily small) $\infty$-categories (resp. $\infty$-groupoids), and let $\PfI$ be the poset of finite type subsets of $\wt{I}$. Let $X\to \PfI$ be the left fibration classifying $\PfI\to \Gpd_\infty, J\mapsto \BB W_J$. 
Then \eqref{eq: BBWJ,BBWa} is equivalent to the natural isomorphism $ |X|\overset{\sim}{\to}\BB W_\aff$. 
Let $F: \BB\Wa\to \widehat{\Cat}_\infty$ be the functor that corresponds to the $\Wa$-action on $\QCoh(T)$. Taking $\lim_{\BB \Wa} F$ is equivalent to $\colim_{\BB\Wa} F^{op}$, where $F^{op}: \BB \Wa\cong \BB\Wa^{op}\to \wh{\Cat}_\infty^{op}$. Let $F_X^{op}: X\to |X|\ovs{\sim}{\to}\BB W_\aff\ovs{F^{op}}{\to} \wh{\Cat}_\infty^{op}$ be the natural functor. Let $F_J^{op}$ be the restriction of $F^{op}_X$ to each fiber $X_J\simeq \BB W_J$ of $X$ over $J\in\PfI$. 
By functoriality, $F_X^{op}$ is the one naturally induced from $F^{op}$ through the inclusions $W_J\subset \Wa$. 

Then 
\begin{align}\label{eq: colimBWJ,Fop}
\lim_{J\sft \wt I}\QCoh(T)^{W_{J}}\simeq \underset{J\sft \wt I}{\colim}\ \underset{\BB W_J}{\colim}F_J^{op}\simeq \underset{X}{\colim}F_X^{op}\simeq \underset{\BB \Wa}{\colim}F^{op}\simeq \QCoh(T)^{\Wa}, 
\end{align}
in which we use \cite[Proposition 4.3.3.10]{Lu1} to get the first equivalence and \cite[Corollary 4.1.2.6]{Lu1} to get that $X\to |X|$ is cofinal. Note that one can also take colimit of $F$ and $F_J$ in \eqref{eq: colimBWJ,Fop} instead of their opposite versions, however, one needs to replace the target of $F$ and $F_J$ by the $\infty$-category of presentable (or presentable stable) $\infty$-categories $\Pr^L$ since the inclusion $\Pr^L\to \widehat{\Cat}_\infty$ does not preserve colimits (but it preserves limits). 
\end{proof}

\sss{Descent under affine Weyl group} 

For each $J\sft \wt I$, we have
the functor
\begin{equation*}
    p_J^*=p_{L_J}^*: \QCoh(\cS_{L_J})\to \QCoh(T)^{W_J}
\end{equation*}
which is a full embedding by Lemma \ref{l:WJ des}. Passing to limits we get a fully faithful embedding
\begin{equation}\label{lim pJ}
    \lim p_J^*: \lim_{J\sft \wt I}\QCoh(\cS_{L_J})\incl \lim_{J\sft \wt I}\QCoh(T)^{W_J}\cong \QCoh(T_\ad\times T)^W.
\end{equation}

\begin{defn}\label{defnQCohWddaggT}
    Let $\QCoh(T_\flat\times T)^{W,\dda}\subset \QCoh(T_\flat\times T)^{W}$ be the full subcategory consisting of those $\cF$ such that for any $J\sft \wt I$, $r_J(\cF)\in \QCoh(T)^{W_J}$ lands in the essential image of $p_J^*$.
\end{defn}

To state the next result, we need some notation. For each root $\a$, let 
\begin{equation*}
    \io_\a: T_\flat\times \ker(\a)\incl T_\flat\times T
\end{equation*}
be the inclusion.  Let $\ker(\a)_\flat=\ker(\a: T_\flat\to \Gm)$. Let 
\begin{equation*}
    \io'_\a: \ker(\a)_\flat\times \ker(\a)\incl T_\flat\times \ker(\a)
\end{equation*}
be the inclusion. We have a functor
\begin{equation*}
    \io'_{\a*}: \QCoh(\ker(\a)_\flat\times\ker(\a))\to \QCoh(T_\flat\times \ker(\a))^{\j{r_\a}}
\end{equation*}
sending $\cG\in \QCoh(\ker(\a)_\flat\times\ker(\a))$ to $\io'_{\a*}$ with the trivial $\j{r_\a}$-equivariant structure. 

\begin{prop}\label{p:aff des}
An object $\cF\in \QCoh(T_\flat\times T)^{W}$ lies in $\QCoh(T_\flat\times T)^{W,\dda}$ if and only if the following condition is satisfied for all roots  $\a$:
\begin{equation}\label{res a supp}
\mbox{$\io_\a^{*}\cF$ is in the essential image of $\io'_{\a,*}$.}
\end{equation}
\end{prop}
\begin{proof} We first show that $\cF\in \QCoh(T_\flat\times T)^{W,\dda}$ if and only if the following a priori weaker version of \eqref{res a supp} holds:
    \begin{equation}\label{coho res a supp}
        \mbox{Each cohomology sheaf of $\io_\a^{*}\cF$ is in the essential image of $\io'_{\a,*}$.}
    \end{equation}

    Let $\cG\in \QCoh(T)^{\wt W_\flat}$ be the object corresponding to $\cF$ under the equivalence \eqref{Wa equiv}. Let $\Phi_{\aff,re}$ be the set of affine real roots of the loop group $LG$. For $J\sft \wt I$, the roots of $L_J$ (with respect to $T$) are viewed as a subset $\Phi_J\subset \Phi_{\aff,re}$. Each $\wt\a\in \Phi_{\aff,re}$ corresponds to a reflection $r_{\wt\a}\in W_J\subset \Wa$. We have the projection $\Phi_{\aff,re}\surj \Phi(G,T)$ by taking the linear part of an affine root.

    By Lemma \ref{l:WJ des}, $r_J(\cF)$ lies in the image of $p_J^*$ if and only if for any root $\wt\a\in \Phi_J$ with projection $\a\in \Phi(G,T)$, the action of $r_{\wt\a}$ on the cohomology sheaves of $\io_{\a}^*\cG$ is trivial. Since any affine real root can be transformed to a root in $\Phi_J$ for some $J\sft \wt I$, we see that $\cF\in \QCoh(T_\flat\times T)^{W,\dda}$ if and only if
    \begin{equation}\label{aff ref triv action}
        \mbox{For any $\wt\a\in \Phi_{\aff,re}$, $r_{\wt\a}$ acts on the cohomology sheaf of $\io_\a^{*}\cG$ trivially.}
    \end{equation}
    Here $\a$ is the linear part of $\wt\a$. 

    For a fixed finite root $\a\in \Phi(G,T)=\Phi_I$, letting $\wt\a$ vary over all affine roots $\wt\a$ with linear part $\a$, the corresponding reflections $r_{\wt\a}$ are of the form $t_{n\a}\c r_\a $ for all $n\in\ZZ$, where $t_{n\a}\in \Wa$ is the translation by $n\a\in \xch(T_\ad)$. Thus \eqref{aff ref triv action} for all $\wt\a$ with linear part $\a$ is equivalent to that each cohomology sheaf of $\io_\a^*\cG$ has trivial action by $r_\a$ and $t_{\a}$, which is equivalent to saying that each cohomology sheaf of $\io_\a^*\cF\in \QCoh(T_\flat\times \ker(\a))^{\j{r_\a}}$ is scheme-theoretically supported on $\ker(\a)_\flat\times \ker(\a)$ with the trivial $r_\a$-action. This is exactly the condition \eqref{coho res a supp}. This shows that $\cF\in \QCoh(\QCoh(T_\flat\times T)^{W,\dda}$ if and only if \eqref{coho res a supp} holds.

    Finally we show \eqref{coho res a supp} implies \eqref{res a supp}, hence they are equivalent. For this it suffices to show that $\io'_{\a*}$ is exact and fully faithful. To simplify notation, let $Y=T_\flat\times \ker(\a)$ and $\io=\io'_\a: Z:=\ker(\a)_\flat\times \ker(\a)\incl Y$. Exactness is clear since $\io$ is an embedding. For full faithfulness it suffices to check it for the generating object $\cO_Z$ that the map
    \begin{equation}\label{End OZ}
        \bR\End_{Z}(\cO_{Z})\to \bR\End_{Y}(\io_*\cO_Z)^{r_\a}
    \end{equation}
    is a quasi-isomorphism. Note the right side is $\bR\Hom(\io^*\io_*\cO_Z,\cO_Z)^{r_\a}$. Koszul resolution calculation gives 
    \begin{equation*}
        \cH^{-i}(\io^*\io_{*}\cO_Z)=\un\Tor_i^{\cO_{Y}}(\cO_Z,\cO_Z)=\begin{cases}
            \cO_Z, & i=0;\\
            \cN^\vee_{Z/Y}, & i=1;\\
            0, & i\ne 0,1.
        \end{cases}
    \end{equation*}
    Since $r_\a$ acts on the conormal bundle $\cN^\vee_{Z/Y}$ by $-1$, the counit map $\io^*\io_*\cO_Z\to \cO_Z$ induces an isomorphism $\bR\Hom_Z(\cO_Z,\cO_Z)^{r_\a}\isom\bR\Hom_Z(\io^*\io_*\cO_Z,\cO_Z)^{r_\a}$, which implies that \eqref{End OZ} is an isomorphism.
\end{proof}

\subsection{Quasi-coherent sheaves on regular centralizers}\label{ss:QCoh J}

We keep the notations from \S\ref{ss:J to R}. Recall the open embedding $\frj=j^{G_\flat}_G$ from \eqref{J to R} which induces a functor
\begin{equation*}
    \frj_*: \QCoh(J^{G_\flat}_G)\to \QCoh((R_{T/\cS_G}(T_\flat\times T))^W)=\QCoh(J_{T_\flat, T}).
\end{equation*}
On the other hand, we apply the setup in \S\ref{sss:ff setup} to $Q=T$, $S=\cS_G$,  and $\xi: X:=T_{\flat}\times T\to T=Q$ with the diagonal $W$-action. Then $J_X$ becomes $J_{T_\flat, T}$. The functor $\Phi_X$ in \eqref{prk} becomes
\begin{equation*}
    \Phi_{T_\flat,T}: \QCoh(J_{T_\flat, T})\to \QCoh(T_\flat\times T)^W.
\end{equation*}
Consider the composite functor
\begin{equation}\label{jR gen}
\Psi: \QCoh(J^{G_\flat}_G)\xr{\frj_*}\QCoh(J_{T_\flat, T})\xr{\Phi_{T_\flat, T}}\QCoh(T_\flat\times T)^W.
\end{equation}
The functor $\Psi$ \eqref{jR gen} is symmetric monoidal with respect to $\QCoh(J_{G^\vee}^{G^{\vee, \ad}})^{\otimes_c}$ and $\QCoh((T^{\vee, \ad}\times T^\vee)/W)^{\otimes_c}$ 
(e.g. one can use \cite[Chapter 5, 5.3.1]{GR} and \cite[Corollary 2.14]{JCorr}). 

The goal of this section is prove:

\begin{theorem}\label{th:main QCoh}
    The functor $\Psi$ above is fully faithful with essential image $\QCoh(T_\flat\times T)^{W,\dda}$.
\end{theorem}

The rest of the subsection is devoted to the proof of this theorem. We first consider a special case.

\sss{Simply-connected derived group case} We assume $(G_\flat)^\der$ is simply-connected. Note that $G^\der\to G_\flat^\der$ is a central isogeny, therefore $G^\der$ is also simply-connected. In this case, $\frj$ is an isomorphism by Proposition \ref{p:J to R open}. Therefore, in view of Proposition \ref{p:aff des},  Theorem \ref{th:main QCoh} in this case is equivalent to:

\begin{prop}\label{p:Im Psi}
    Suppose $(G_\flat)^\der$ is simply-connected. Then $\Phi_{T_\flat,T}$ is fully faithful with essential image consisting of $\cF\in \QCoh(T_\flat\times T)^W$ satisfying \eqref{res a supp} for each root $\a\in \Phi(G,T)$.
\end{prop}
\begin{proof} We shall apply Theorem \ref{th:ff} to $X=T_\flat\times T\to T=Q$ and $Q^\flat=\cS_G=T\sslash W$. 
We first check that the assumptions of Theorem \ref{th:ff} are satisfied. 

First we claim that the map $\pi:T\to T\sslash W$ satisfies our assumption: i.e., $\cO(T)$ is a locally free $\cO(T)^{W}[W]$-module of rank one. For this it suffices to check that for any geometric point $\t\in (T\sslash W)(K)$, the fiber $\pi^{-1}(\t)=\Spec C_{\t}$ is isomorphic to $K[W]$ as a $W$-module. Let $t\in T(K)$ be in the fiber of $\t$, with stabilizer $W_{t}$. Then the localization of $C_{\t}$ at $t$ is $C_{t}\cong \cO_{T, t}/((\cO_{T,t}^{W_{t}})_{+})$, where $\cO_{T, t}$ is the local ring of $T\times\Spec K$ at $t$, so that $\cO_{T,t}^{W_{t}}$ is the local ring of $(T\sslash W)_{K}$ at $\t$, and $(\cO_{T,t}^{W_{t}})_{+}$ is its maximal ideal. It is easy to see that $C_{t}$ is isomorphic to the coinvariant algebra $K[\frt]/(K[\frt]^{W_{t}}_{+})$ of $W_{t}$ with respect to the representation $\frt=\Lie T$. Since $G^{\der}$ is simply-connected, $W_{t}$ is a reflection subgroup of $W$, and $\frt$ is the reflection representation of $W_{t}$, up to  
trivial modules of $W_{t}$. This implies that $C_{t}\cong K[W_{t}]$ as $W_{t}$-modules. Therefore $\cO(T)$ is a locally free $\cO(T)^{W}[W]$-module of rank one. 

We next check the admissibility condition for the relative tangent spaces of $\xi: X=T_\flat\times T\to T$.  Let $K$ be a field containing $k$ and $t\in T(K)$. We use the notations $W_{t}$ and $C_{t}$ introduced in the previous paragraph. Let $x=(t', t)\in X^{W_{t}}(K)$ (so that $t'\in T_{\flat}^{W_{t}}(K)$). We need to check that $T_{\xi, x}\cong \frt_\flat\ot K$ (where $\frt_\flat=\Lie T_\flat$) is a $C_{t}$-admissible $W_{t}$-module. By assumption, as $W$-modules, $\frt_{\flat}$ and $\frt$ are isomorphic up to trivial $W$-modules. Therefore $\frt_{\flat}\ot K$ is isomorphic to the reflection representation of $W$ up to a trivial $W$-module, hence as a $W_{t}$-module, $T_{\xi, x}= \frt_\flat\ot K$ is also isomorphic to the reflection representation of $W_t$ up to a trivial module. By Lemma \ref{l:adm triv} and Lemma \ref{l:ref adm}, $\frt_{\flat}\ot K$ is a $C_{t}$-admissible $W_{t}$-module. The assumptions of Theorem \ref{th:ff} are checked. In particular, $\Phi_{T_1, T_2}$ is fully faithful.

Now we show that the essential image of $\Phi_{T_\flat, T}$ characterized in Theorem \ref{ss:ff}(2) is the same as the one characterized by \eqref{res a supp}. By Theorem \ref{ss:ff}(2) the image of $\Phi_{T_\flat, T}$ consists of $\cF$ such that for any geometric point $t\in T(K)$, the derived restriction $i_t^*\cF\in \QCoh(T_{\flat,K})^{W_t}$ satisfies
\begin{equation}\label{itF image}
    \mbox{$i_t^*\cF$ lies in the essential image of $\io_{t*}: \QCoh(T_{\flat,K}^{W_t})\to \QCoh(T_{\flat,K})^{W_t}$.}
\end{equation}
We claim that \eqref{itF image} is equivalent to
\begin{equation}\label{coho itF image}
    \mbox{Each cohomology sheaf of $i_t^*\cF$ lies in the essential image of $\io_{t*}$.}
\end{equation}
Indeed, the argument for the equivalence between \eqref{coho itF image} and \eqref{itF image} is the same as the argument of the equivalence between \eqref{coho res a supp} and \eqref{res a supp}. The point is that the $W_t$-action on $\wedge^i(N^\vee_{T^{W_t}_\flat/T_\flat})$ for $i>0$ does not have trivial factors. Here the normal bundle $N^\vee_{T^{W_t}_\flat/T_\flat}$ is the trivial bundle whose fiber is identified with the reflection representation of $W_t$.

We further claim that \eqref{coho itF image} is equivalence to:
\begin{equation}\label{Wt act triv on global sections}
    \mbox{$W_t$ acts trivially on $\bR\Gamma(T_{\flat, K}, i_t^*\cF)$.}
\end{equation}
Indeed,  \eqref{coho itF image} obviously implies \eqref{Wt act triv on global sections} by taking global sections. Conversely, assume \eqref{Wt act triv on global sections} holds for $t\in T(K)$. Since $T_\flat$ is affine, each cohomology sheaf $\cH^i(i_t^*\cF)$ satisfies that its global sections $M^i:=\Gamma(T_{\flat, K}, \cH^i(i_t^*\cF))$ has trivial $W_t$-action. Now $M^i$ is a module over the associative algebra $\cO(T_{\flat,K})\#W_t$. A simple calculation shows that if $W_t$ acts trivially on $M^i$, then $f-w(f)$ acts on $M^i$ by zero for any $f\in \cO(T_{\flat,K})$, therefore the $\cO(T_{\flat,K})$-action on  $M^i$  factors through  $\cO(T_{\flat,K}^{W_t})$, and we conclude that \eqref{coho itF image} holds.

The same argument for the equivalence between \eqref{itF image} and \eqref{Wt act triv on global sections} shows that \eqref{res a supp} is equivalent to:
\begin{equation}\label{ra act triv on global sections}
    \mbox{$r_\a$ acts trivially on $\bR\Gamma(T_{\flat}\times \ker(\a), \io_\a^*\cF)$.}
\end{equation}
Here the $r_\a$-action is induced from the $r_\a$-action on both $T_{\flat, K}$ and the $r_\a$-equivariant structure on $\io_\a^*\cF$.

It remains to show that \eqref{ra act triv on global sections} for all roots $\a$ is equivalent to \eqref{Wt act triv on global sections} for all geometric points $t\in T(K)$.

Suppose \eqref{res a supp} holds for all roots $\a$. Let $t\in T(K)$. Since $G^\der$ is simply-connected, for any geometric point $t\in T$, $W_t$ is generated by reflections by \cite[Ch. II, Theorem 4.2]{SpSt}. Therefore, to see \eqref{Wt act triv on global sections}, it suffices to check that for each $r_\a\in W_t$, $r_\a$ acts trivially on $\bR\Gamma(i_t^*\cF)$. Now $r_\a\in W_t$ implies $t\in T^{r_\a}=\ker(\a)$ (since $\a^\vee$ is primitive), and \eqref{res a supp} implies that $\bR\Gamma(\io_\a^*\cF)$ has trivial $r_\a$-action. Now we have an $r_\a$-equivariant quasi-isomorphism of $\cO(T_{\flat,K})$-modules
\begin{equation}\label{global section tensor}
    \bR\Gamma(i_t^*\cF)\cong \bR\Gamma(\io_\a^*\cF)\ot^{\bL}_{\cO(\ker(\a))}K(t)
\end{equation}
where $K(t)=K$ is to emphasize it is the residue field of $t$. We conclude that the action of $r_\a$ on $\bR\Gamma(i_t^*\cF)$ is also trivial. This being true for all reflections $r_\a\in W_t$, we conclude that \eqref{Wt act triv on global sections} holds.

Conversely, suppose \eqref{Wt act triv on global sections} holds. We may assume $\cF$ is cohomologically bounded from above. Let $\a\in \Phi(G,T)$. The involution $r_\a$ on $M=\bR\Gamma(\io_\a^*\cF)$ decomposes it into eigenspaces $M^+\op M^-$, where $r_\a$ acts trivially on $M^+$ and by $-1$ on $M^-$. View $M^-$ as a quasi-coherent sheaf on $T_{\flat}$. Taking $-1$ eigenspaces under $r_\a$ on both sides of \eqref{global section tensor}, we see that $M^{-}\ot^{\bL}_{\cO(\ker(\a))}K_t=(\bR\Gamma(i_t^*\cF))^-=0$ for any geometric point $t$ of $T$. By \cite[Lemma 10]{Arinkin}, this implies that $M^-=0$. This finishes the proof.

\end{proof}

The proof of Theorem \ref{th:main QCoh} will be reduced to the case $G^\der_\flat$ is simply-connected by the following two steps.

\sss{Reduction: changing $G$}\label{sss:red Gder sc}

Now $G$ is a connected reductive group over $\CC$ and $\nu: G\to G_\flat$ is a central quotient.

Let $\pi: G_1\to G$ be a central isogeny  and let $T_1=\pi^{-1}(T)$. We claim that if Theorem \ref{th:main QCoh} holds for $(G_1, G_\flat)$, then it holds for $(G,G_\flat)$. This can thus be applied to the case where $G_1^\der$ is simply-connected. 

Let $\wt\Psi:\QCoh(J_{G_1}^{G_\flat})\to \QCoh(T_\flat\times \wt T)^W$ be the analog of $\Psi$ for $(G_1,G_\flat)$. Let $\pi_J: J_{G_1}^{G_\flat}\to J_G^{G_\flat}$ and $\Pi_T=\id_{T_\flat}\times \pi_T: T_\flat\times \wt T\to T_\flat\times T$ be the projections induced by $\pi$, which are both $\Gamma$-torsors. In the following commutative diagram, all squares are Cartesian and all vertical maps are $\Gamma$-torsors
\begin{equation*}
    \xymatrix{J_{G_1}^{G_\flat}\ar[r]^-{\wt\frj}\ar[d]^{\pi_J} & J_{T_\flat, \wt T}\ar[d] & \ar[l]_-{\pr_1} J_{T_\flat, \wt T}\times_{\cS_{G_1}}(\wt T/W) \ar[r]^-{\wt\vk}\ar[d] & (T_\flat\times \wt T)/W\ar[d]^{\Pi_{T}}\\
    J_{G}^{G_\flat}\ar[r]^-{\frj} & J_{T_\flat, T} & \ar[l]_-{\pr_1} J_{T_\flat, \wt T}\times_{\cS_{G}}(T/W) \ar[r]^-{\vk} & (T_\flat\times T)/W
    }
\end{equation*}
The base change isomorphisms give two commutative diagrams, one for $(\pi_{J*}, \Pi_{T*})$ and one for $(\pi_J^*, \Pi_T^*)$
\begin{equation}\label{wt Psi and Psi}
\xymatrix{\QCoh(J_{G_1}^{G_\flat})\ar[r]^-{\wt\Psi}\ar@<1ex>[d]^{\pi_{J*}} & \QCoh(T_\flat\times \wt T)^W\ar@<1ex>[d]^{\Pi_{T*}}\\
\QCoh(J_{G}^{G_\flat})\ar[r]^-{\Psi} \ar@<1ex>[u]^{\pi_{J}^*} & \QCoh(T_\flat\times T)^W \ar@<1ex>[u]^{\Pi_{T}^*}}
\end{equation}
The functor $\wt\Psi$ is $\Gamma$-equivariant. Taking $\Gamma$-invariants we recover $\Psi$. If $\wt\Psi$ is fully faithful, so is $\Psi$. Moreover, $\QCoh(T_\flat\times T)^{W,\dda}$ coincides with $(\QCoh(T_\flat\times T)^{W,\dda})^\Gamma$. By assumption, the essential image of $\wt\Psi$ is $\QCoh(T_\flat\times T)^{W,\dda}$, therefore, to show that the image of $\Psi$ is $\QCoh(T_\flat\times T)^{W,\dda}$, it remains to show that $\cF\in \QCoh(T_\flat\times T)^W$ lies in the image of $\Psi$ if and only if $\Pi_T^*\cF\in \QCoh(T_\flat\times \wt T)^W$ lies in the image of $\wt\Psi$, where $\pi_T:\wt T\to T$ is the restriction of $\pi$. The commutativity of \eqref{wt Psi and Psi} for pullbacks shows that if $\cF\in \Im(\Psi)$, then $\Pi_T^*\cF\in \Im(\wt\Psi)$. Conversely, if $\Pi_T^*\cF\in \Im(\wt\Psi)$, hence $\Pi_T^*\cF\cong \wt\Psi(\cG)$ for some $\cG\in \QCoh(J_{G_1}^{G_\flat})$, then by \eqref{wt Psi and Psi}, $\Pi_{T*}\Pi_T^*\cF\cong \Pi_{T*}\wt\Psi(\cG)\cong \Psi(\pi_{J*}\cG)$ is in the image of $\Psi$. Now $\cF$ is a summand of $\Pi_{T*}\Pi_T^*\cF$ (where $\xch(\Gamma)$ acts trivially), and since $\Psi$ is fully faithful and $\QCoh(J_G^{G_\flat})$ is Karoubian, $\Im(\Psi)$ is closed under taking direct summands, we conclude that $\cF\in \Im(\Psi)$.

\sss{Reduction: changing $G_\flat$}\label{sss:red Gflat} 
Suppose $\nu: G\to G_\flat$ factors as $G\xr{\nu_\sh} G_\sh\xr{\mu} G_\flat$, where $\mu$ is a central isogeny. Then the projections $J_{G}^{G_\sh}\to J_{G}^{G_\flat}$ and $T_\sh\times T\to T_\flat\times T$ are both $\ker(\mu)$-torsors. The same argument as in \S\ref{sss:red Gder sc} shows that if Theorem \ref{th:main QCoh} holds for $(G,G_\sh)$, then it also holds for $(G,G_\flat)$. 

\sss{Finish of the proof of Theorem \ref{th:main QCoh}} After \S\ref{sss:red Gder sc}, we may assume $G^\der$ is simply-connected. We may even assume that $G=G^\der\times A$ where $A$ is a torus, with $G^\der$ simply-connected. Let $\ov A=\nu(A)\subset G_\flat$ be the image of $A$ and take $G_\sh=G^\der\times \ov A$. We have a factorization $\nu: G\to G_\sh\xr{\mu} G_\flat$ with $\mu$ a central isogeny. Applying \S\ref{sss:red Gflat} we may replace $G_\flat$ by $G_\sh$, which has a simply-connected derived group. In this case Theorem \ref{th:main QCoh} is equivalent to Proposition \ref{p:Im Psi}. We are done. \qed

Applying Theorem \ref{th:main QCoh} to $G_\flat=G_\ad$, and combining with Lemma \ref{l:lim Wa}, we obtain:
\begin{cor}\label{c:QCoh J as lim}
    There is a canonical symmetric monoidal equivalence
    \begin{equation*}
        \QCoh(J_G^{G_\ad})^{\otimes_c}\simeq \lim_{J\sft \wt I}\QCoh(\cS_{L_J})^{\otimes}.
    \end{equation*}
Forgetting the monoidal structure, this gives a canonical equivalence
\begin{align*}
 \QCoh(J_G^{G_\ad})\simeq \lim_{J\sft \wt I}\QCoh(\cS_{L_J}).
\end{align*}
\end{cor}

%%%%%%%%%%%%%%%%%%%%%%%%%%%%%
\section{Proof of the Mirror Theorem \ref{th:mirror}}\label{secProofMirror}%%%%%%%%%%%
%%%%%%%%%%%%%%%%%%%%%%%%%%%%%%

In this section, we give the proof of the Mirror Theorem \ref{th:mirror}. 
Using results from \S\ref{s:Fukaya}--\S\ref{s:QCoh J}, the remaining step is to construct an equivalence between the canonical limit diagram for calculating $\Ind\cW(\cM_G^\circ)$ (\S\ref{subsec: cWM_G,colimit}) and the canonical limit diagram for $\QCoh(J_{G^\vee}^{G^{\vee, \ad}})$ in Corollary \ref{c:QCoh J as lim}. As stated in the Mirror Theorem, such an equivalence depends on a choice of a compatible system of Kostant sections (cf. Definition \ref{defncompatibleKostantsec} below).

We begin with a review of HMS for Lie algebra regular centralizers in \S\ref{ss: reveiw of J}. 
Building on this, in \S\ref{sspi1cTHMG}--\S\ref{subsecCanonicalME}, for a reductive group $H$, we use monadicity on wrapped Fukaya categories (for a finite central quotient $H\to \Hf$ with $Z(\Hf)$ connected) and the canonical $\pi Z(H)$-action on $\cW(\cT_H)$ to establish 
the canonical  $\pi Z(H)$-equivariant mirror equivalence $\Phi_H$ \S\ref{subsecCanonicalME}. Then in \S\ref{ssPhiHrestriction}, we establish the compatibility of $\Phi_H$ with respect to restriction to a Levi subgroup $L\subset H$ containing $T$, which is encoded in the diagram \eqref{eqWPHIndWQCoh}. For $\cM_G^\c$, we assemble these diagrams into the natural commutative diagram \eqref{eqfrMPfrIWfrQfrP} of coCartesian fibrations over $\PfI^{op}$. We also show the existence of a compatible system of Kostant sections $\sigma$ in Lemma \ref{defncompatibleKostantsec}. 
Lastly, for any $\sigma$, we get an equivalence between the respective canonical limit diagrams for $\Ind\cW(\cM_G^\circ)$ and $\QCoh(J_{G^\vee}^{G^{\vee, \ad}})$ over $\PfI^{op}$. Taking limit categories on both sides, we get a desired mirror equivalence $\Phi_\sigma$. This finishes the proof of the Mirror Theoerem in \S\ref{ssfinishingpfMT}.

\subsection{Mirror symmetry for Lie algebra regular centralizers (review of \cite{J})}\label{ss: reveiw of J}

Let $H$ be any reductive group with a maximal torus $T$, a complete set of simple roots $I_H$ and Weyl group $W$.
To simplify notations, we will use either $\cW(\ovl{\cT}_H)$ or $\cW(\cT_H)$ to denote the wrapped Fukaya category of the Liouville sector $\ovl{\cT}_H$ (which is also a generalized Weinstein sector) defined in \S\ref{sssOpenWeinJG} and \S\ref{sssYJcJLJ}.

We recall the main result from \cite{J} for any reductive group $H$. For any finite set $S$, let $\cP(S)$ (resp. $\cP^\dagg(S)$) denote the power set of $S$ (resp.  $\cP(S)-\{S\}$). Recall $\cP_{ft}(\wt{I})=\{J\subset_{ft}\wt{I}\}\subsetneq \cP(\wt{I})$.

\sss{Main results from \cite{J}}

We state the HMS results and parabolic induction patterns from \cite{J}, for reductive groups with connected centers. We give details in later sections on establishing HMS for general reductive groups and parabolic induction patterns based on this. Although the latter has been established for $\cJ_{H}$ in \cite{J}, it is slightly different from the current setting where $\cT_{H}$ is a torsor of $\cJ_{H}$ (cf. Remark \ref{remarkPhiHcJH} for the relation between HMS for them), and this difference becomes more significant when we consider gluing wrapped Fukaya categories of them for calculating $\cW(\cM_G^\c)$ (cf. Corollary \ref{cor: cWcMG,twocolimits}). 

In the following, to simplify notations, we will assume all Kostant sections have been cylindricalized when they are considered objects in the relevant Fukaya categories. 

\begin{theorem}[HMS for groups with connected center]\cite[Theorem 7.7]{J}\label{thmJconnectedcenter}
Let $H$ be any reductive group with connected center. Fix any Kostant section $S_{H}$ in $\cT_H$, we get 
a natural equivalence
\begin{align*}
\Upsilon_{S_{H}}: \cW(\cT_H)&\overset{\sim}{\longrightarrow} \Coh(T^\vee\sslash W), \\
S_{H}&\mapsto \cO_{T^\vee\sslash W},\\
\sfL&\mapsto \cO(T^\vee\sslash W)=\End(S_H)^{op}\car\Hom_{\cW(\cT_H)}(S_H, \sfL).
\end{align*}
We will call $\Upsilon_{S_{H}}$ the mirror equivalence \emph{relative to $S_H$}. 
\end{theorem}
\begin{proof}[Sketch of proof]
Here is a sketch of the proof based on Proposition \ref{propcPIcWcY}. The case for $H^\flat=H^\ad\times H_\ab$ follows directly from  Proposition \ref{propcPIcWcY} using the K\"unneth formula for product sectors (cf. \cite{GPS2}). The case for $H$ follows from monadicity results between (ind-completion of) wrapped Fukaya categories (cf. \S\ref{subsecMonadicity} for more details on the same line of arguments). 
\end{proof}

For any standard Levi subgroup $L\subset H$, let $\res^H_L: \IndW(\cT_H)\to \IndW(\cT_L)$ (which actually preserve compact objects by the following theorem). 

\begin{theorem}[Parabolic Induction pattern for groups with connected centers]\cite[Proposition 7.8]{J}\label{thmJparabolicconnectedcenter}
Let $H$ be any reductive group with connected center (then all standard Levi subgroups $L$ have connected centers as well). For any Kostant section $S_H$, $\res^H_L(S_H)$ is isomorphic to any Kostant section $S_L$ in $\cT_L$, and the morphism 
\begin{align*}
\cO(T^\vee\sslash W)=\End_{\cW(\cT_H)}(S_H)\lrar \End(\res^H_L(S_H))=\End_{\cW(\cT_L)}(S_L)=\cO(T^\vee\sslash W_L) 
\end{align*}
is given by the natural inclusion $\cO(T^\vee\sslash W)\hookrightarrow \cO(T^\vee\sslash W_L)$ of invariant functions. 

\end{theorem}

\subsection{Fundamental groups of $\cT_H$ and $\cM_G^\circ$}\label{sspi1cTHMG}
 We will need the following lemma on fundamental groups. 

\begin{lemma}\label{lemmapi1cMG}
\begin{itemize}
\item[(i)]
Let $H$ be any reductive group with a maximal torus $T$ and a complete set of simple roots $I_H$. There is a canonical isomorphism $\pi_1(\cT_H)\cong \pi_1(H)$. 

\item[(ii)]
In the case of a semisimple group $G$, for any $J\sft\wt{I}$, we have the natural commutative diagram
\begin{equation}\label{eqLemmapi1cMG}
\begin{tikzcd}
\pi_1(\cT_{L_\vn})\ar[d, "{\wr}"]\ar[r]&\pi_1(\cT_{L_J})\ar[r]\ar[d, "{\wr}"]&\pi_1(\cM_G^\circ)\ar[d,"{\wr}"]\\
\pi_1(T)\ar[r]&\pi_1(L_J)\ar[r]&\pi_1(G),
\end{tikzcd}
\end{equation}
which is compatible with any inclusions $J\subset J'$. 
\end{itemize}
\end{lemma}
\begin{proof}
Using the Weinstein handle attachment description of $\cM_{G}^\circ$ (cf. Remark \ref{remSpecialZflatHandle}) and similarly for $\cT_H$, one starts from $T^*T\overset{h.e.}{\simeq}T_{\cpt}$ and then attach handles to it to obtain $\cM_G^\circ$ and $\cT_H$, respectively. The $2$-handles for $\cM_{G}^\circ$ (resp. $\cT_H$) can be chosen as $D^2_{j,[\zeta]}:=\chi_{\{j\}}^{-1}([0])^\circ\times \{(\zeta,0)\}\subset \cT_{L_{\{j\}}^\der}\times^{Z(L_{\{j\}}^\der)}T^*Z(L_{\{j\}})\cong\cT_{L_{\{j\}}}$, for $j\in \wt{I}$ (resp. $j\in I_H$) and $[\zeta]\in \pi_0 Z(L_{\{j\}})$, where $\chi_{\{j\}}^{-1}([0])^\circ$ is any connected component of $\chi_{\{j\}}^{-1}([0])$. 
The attaching map of $\partial D^2_{j,[\zeta]}$ boils down to a calculation for groups of semisimple rank $1$, and it is easy to see that it 
gives $\alpha_j^\vee\in \XX_*(T)=\pi_1(T)$. Hence the canonical morphism $\pi_1(T)\cong\pi_1(T^*T)\to\pi_1(\cT_G)\to \pi_1(\cM_{G}^\circ)$ identifies both $\pi_1(\cT_G)$ and $\pi_1(\cM_{G}^\circ)$ with $\XX_*(T)/\XX_*(T_{\sc})$ (and the same holds for $\cT_H$ where $\XX_*(T_{\sc})$ is replaced by the coroot lattice). Using the canonical  morphism $\pi_1(T)\to \pi_1(G)$ to identify $\pi_1(G)$ with $\XX_*(T)/\XX_*(T_{\sc})$, we get the canonical isomorphism $\pi_1(\cT_G)\cong\pi_1(\cM_{G}^\circ)\cong \pi_1(G)$ (and the same for $\cT_H$). 

Regarding the diagram \eqref{eqLemmapi1cMG}, the isomorphism $\pi_1(\cT_{L_J})\overset{\sim}{\to}\pi_1(L_J)$ follows from the case $H=L_J$, where the attaching 2-cells are indexed by $j\in J$. Hence we have the commutative diagram that is compatible with  any inclusion $J\subset J'$. 
\end{proof}

\subsection{Monadicity on wrapped Fukaya categories}\label{subsecMonadicity}

\sss{Description of $\cW(\cT_{H})$ from monadicity}
Let $H^\flat$ be any finite central quotient of $H$ such that $\pi_0 Z(H^\flat)=1$.
Let $\pi_{\cT,H}: \cT_H\to \cT_{H^\flat}$ and $\nu_\flat: T\to T^\flat$ be the obvious qoutient maps. 
Then we have the natural adjunction between the wrapped Fukaya categories:  
\begin{equation}\label{eqFJLFJRflat}
\begin{tikzcd}
\cW(\cT_{H})\ar[r, shift right=0.2em, "F_*"']&\cW(\cT_{H^\flat})\ar[l, shift right=0.2em, "F^*"'],
\end{tikzcd}
\end{equation}
where $F^*$ is both the left and the right adjoint of $F_*$. 
Let  $\frT_\flat:=F_*\circ F^*$ be the monad acting on $\Ind\cW(\cT_{H^\flat})$. Since $F_*$ is clearly conservative, we have the canonical equivalence
\begin{align}\label{eqFstarBBL}
F_*: \Ind\cW(\cT_{H})\overset{\sim}{\longrightarrow} \Mod_{\frT_{H, \flat}}(\Ind\cW(\cT_{H^\flat})), 
\end{align}
by Barr-Beck-Lurie theorem.

%%%%%%%%%%%%%%%%%%%%%%%%%%%%%%%%%%%%%%%%%%%%%%%%%

\sss{Functorial monadicity}\label{subsecMonadmuSh}

Let $\pi(T)$ be the ordinary fundamental groupoid of $T$ (which is equivalent to the $\infty$-groupoid associated with the space $T$).  
The adjunction \eqref{eqFJLFJRflat} is naturally compatible with restriction and co-restriction functors among $\IndW$ for the sectorial coverings $\{\cT_{L_J}^{\ovl{\cU}_J}, J\in \cP(I_H)\}$ and $\{\cT_{L^\flat_J}^{\ovl{\cU}_J}, J\in \cP(I_H)\}$ of $\cT_H^{\ovl{\cV}}$ and $\cT_{\Hf}^{\ovl{\cV}}$, respectively.  
In particular, we have the natural diagram
\begin{equation}\label{diagrammuShcovering}
\begin{tikzcd}[column sep=6em]
\Fun\left(\pi (T^\flat)^{op}, \Mod(\CC)\right)\simeq\IndW(\cT_{L_\vn^\flat})\ar[r, "\cores"', shift right=0.2em]\ar[d, "F_\vn^*"',shift right=0.2em]&\IndW(\cT_{H^\flat})\ar[l, "\res"',shift right=0.2em]\ar[d, "F^*"',shift right=0.2em ]\\
\Fun\left(\pi (T)^{op}, \Mod(\CC)\right)\simeq\IndW(\cT_{L_\vn}) \ar[u, "F_{\vn;*}"', shift right=0.2em]\ar[r, "\cores"', shift right=0.2em]&\IndW(\cT_{H})\ar[l,"\res"',shift right=0.2em ]\ar[u, "F_*"', shift right=0.2em]
\end{tikzcd},
\end{equation}
in which the functor $\res$ (resp. $\cores$) respects the vertical adunctions. 
The identification on the left column above is to be consistent with \S\ref{sssDefpiZHaction}, and this is consistent with the standard identification 
\begin{align*}
\cW(T^*T)^{op}\simeq \Loc^w(T)\simeq \Perf(\CC[\pi_1(T)]),
\end{align*}
where $\Loc^w(T)$ is the full subcategory of $\Loc(T)$ consisting of compact objects. Although we have the standard isomorphism $\pi (T)^{op}\simeq \pi (T)$, we will keep the identification in \eqref{diagrammuShcovering}, since this will make things consistent in \S\ref{subsecIWQstd} (cf. \S \ref{sssconventionoppositealg} for more explanations). 

Now we use notations from \S\ref{sssRegCovSke}. 
Since the regular covering $\cT_{H}\to \cT_{H^\flat}$ classified by \eqref{eqKJflattoPpi1cT}
restricts to a  regular covering $\cT_{L_J}\to \cT_{L_J^\flat}$ classified by 
\begin{align}
&\psi_{L_J, \pi_1}: \pi_1(\cT_{L_J^\flat}, z_\vn^\flat)\cong \pi_1(L_J^\flat)\twoheadrightarrow K=\ker(L_J\to L_J^\flat)\car\pi_{\Lambda, H}^{-1}(y_\vn^{J,\flat}), 
\end{align}
by functoriality, to understand the induced monadicity between $\Ind\cW(\cT_{L_J})$ and $\Ind\cW(\cT_{L^\flat_J})$ for any $L_J$, it suffices to understand the case for $H$.

\sss{Description of $\frT_\flat=F_*F^*$}\label{sssfrTflatF}
 
 Now consider the local system $\cL_K=(\pi_{\cT, H})_*\CC_{\cT_H}\in \Loc(\cT_{H^\flat})^{\bd, \hs}$. This is classified by 
\begin{align}\label{eqpi1HflattoK}
\pi_1(\cT_{H^\flat}, y_\vn^\flat)\cong \pi_1(H^\flat)\twoheadrightarrow K\car \cO(\pi_{\Lambda, H}^{-1}(\yvnf))=\cL_{K}|_{\yvnf}.
\end{align}
Using the obvious algebra structure on $\cO(\pi_{\Lambda, H}^{-1}(\yvnf))$ (as a representation of $K$), $\cL_K$ is naturally a commutative algebra object in $\Loc(\cT_{H^\flat})^{\bd, \hs;\otimes}$. Note if we use $\yvnf=\zvnf$, then we have a canonical identification $\pi_{\Lambda, H}^{-1}(\zvnf)=K$, but this is not needed anywhere.

\begin{prop}\label{propfrTflatcLKbd}
The monad $\frT_\flat$ is canonically identified with the tensor action of $\cL_K\in \Alg(\Loc(\cT_{H^\flat})^{\bd, \hs; \otimes})$ (cf. \S\ref{appendsubseclocalsystemaction}). 
\end{prop}

\begin{proof}
One can use the same proof as \cite[Lemma 7.4]{J} to concretely realize $F_*$ and $F^*$ between wrapped Fukaya categories, where $F^*$ is easy to realize but to describe $F_*$ one needs to restrict to a certain full subcategory $\cW(\cT_H)^s$ that generates $\cW(\cT_{H})$, e.g. the full subcategory of (connected and) simply connected Lagrangians or contractible Lagrangians of $\cW(\cT_{H})$. Then it is direct to check on $\cW(\cT_H)^s$ that $\frT_\flat$ identifies with the tensor action of $\cL_K$. 
\end{proof}

\begin{remark}
\begin{itemize}
\item[(i)]
The functor $F_*$ and $F^*$ are respectively represented by the Lagrangian correspondences
\begin{align*}
&\cG_{*}:=\Graph(\pi_{\cT, H})\subset \cT_{H}^-\times \cT_{H^\flat}\\
&\cG^*:=\Graph(\pi_{\cT, H})^t\subset \cT_{H^\flat}^-\times \cT_{H}
\end{align*}
Then $\frT_\flat=F_*F^*$ is given by the algebraic convolution of $\cG^*$ followed by $\cG_*$. It is direct to check the geometric convolution $\cG^*\circ \cG_*$ is well defined and gives $\bL_{\pi}=(\Delta_{\cT_{H^\flat}}, b, \cL_{K})$. Here we use \S\ref{sssGeometricComposition} to have a ``nice" geometric composition. Assuming the geometric composition agrees with the algebraic composition, then by \S\ref{ssstensorLoc}, $\bL_{\pi}$ gives the tensor action of $\cL_K$. Since we are not sure about the status of composition of Lagrangian correspondences for stopped Liouville manifolds, we only state this argument informally.  

\item[(ii)] If we identify wrapped Fukaya categories with microlocal sheaves, then $\frT_\flat=(\pi_{\Lambda, H})_*\pi_{\Lambda, H}^*$ is clearly given by the tensor action of $\cL_K$ on microlocal sheaves. 
\end{itemize}
\end{remark}

Using standard representation theory of finite abelian groups, we have the canonical ring isomorphism between the group algebra of $K^*$ and $\cO(K)$ (which can be viewed as a version of Fourier transform): 
\begin{align*}
\CC[K^*]&\overset{\sim}{\longrightarrow}\cO(K)\\
e_{\chi}&\mapsto \chi,\quad \chi\in K^*\\
\frac{1}{|K^*|}\sum_{\chi\in K^*}\kappa_j(\chi)^{-1}e_\chi&\mapsto \delta_{\kappa_j},\quad \kappa_j\in K
\end{align*}
in which the second (resp. third) line gives the decomposition of both sides as representations of $K$ (resp. $K^*$). 
It follows that the local system $\cL_K$, classified by \eqref{eqpi1HflattoK}, decomposes into invertible local systems 
\begin{align*}
\cL_K=\bigoplus_{\chi\in K^*}\cL_\chi,
\end{align*}
where each $\cL_\chi$ is classified by 
\begin{align}\label{eqcLchiclass}
\pi_1(\cT_{H^\flat}, \yvnf)\cong \pi_1(H^\flat)\twoheadrightarrow K\circlearrowright \CC\cdot \chi\subset \cL_K|_{\yvnf}. 
\end{align}
The algebra structure on $\cL_K$ as an object in  $\Loc(\cT_{H^\flat})^{\bd, \hs}$ is equivalent to the data of a local system-valued character of $K^*$: 
\begin{align}\label{eqphicLK}
\phi_{\cL_K}: K^*&\lrar \Loc(\cT_{H^\flat})^{\bd, \hs;\otimes}\\
\nonumber \chi&\mapsto \cL_\chi,
\end{align}
where the structure maps corresponding to multiplication in $K^*$ is the canonical isomorphism $\cL_{\chi_1}\otimes \cL_{\chi_2}\overset{\sim}{\to}\cL_{\chi_1\chi_2}$. 

For any $\infty$-category $\cC$, let $\cC^{\Gpd}$ be the $\infty$-groupoid associated with $\cC$, i.e. it is the $1$-full subcategory of $\cC$ whose $1$-morphisms consist of invertible morphisms in $\cC$.

\begin{prop}\label{lemmaKstarmuSh}
We have the canonical $K^*$-action on $\cW(\cT_{H^\flat})$, given by 
\begin{align}
K^*&\longrightarrow \Aut(\cW(\cT_{H^\flat}))\\
\nonumber\chi&\mapsto\quad \otimes \cL_{\chi},
\end{align}
which corresponds to the monad $\frT_\flat$ (by taking the categorical version of group algebra object in $\End(\cW(\cT_{H^\flat}))$). 

The automorphism group of this $K^*$-action is canonically identified with $K$, through the automorphism of $\cL_\chi$ by multiplication with the scalar $\chi(\kappa)$ for any $\kappa\in K$. Thus the canonical $K^*$-action gives a canonical full embedding as $\infty$-groupoids, $\BB K\hrar \Fun\left(K^*,  \End(\cW(\cT_{H^\flat}))^{\otimes}\right)^{\Gpd}$.

\end{prop}

\begin{proof}
This directly follows from Proposition \ref{propfrTflatcLKbd} and the discussions above. 
\end{proof}

Combining \eqref{eqFstarBBL} and Proposition \ref{lemmaKstarmuSh}, we directly get 

\begin{cor}\label{corWTHWTHfKstar}
The functor $F_*: \Ind\WTH\to \Ind\WTHf$ canonically identifies 
\begin{align*}
\Ind\WTH\simeq \Mod_{\frT_\flat}(\Ind\WTHf)\simeq  \Ind\WTHf^{K^*}.
\end{align*}
\end{cor}

%%%%%%%%%%%%%%%%%%

\subsection{Framings of the $K^*$-action on $\cW(\cT_{H^\flat})$}

To deal with the automorphism group $K$ of the canonical $K^*$-action, we make the following definitions.

\sss{The groupoid $\Xi(K^*)$ (resp. $\wt{\Xi}(K^*)$) of line characters (resp. framed line characters) of $K^*$}
\begin{defn}\label{defnlinecharKstar}
\begin{itemize}
\item[(i)] A \emph{line-valued character of $K^*$} (or line character for short) is a character of $K^*$ valued in $\Pic(\CC)$\footnote{Here $\Pic(\CC)$ is the usual Picard groupoid consisting of $\CC$-lines in homological degree $0$.}
\begin{align*}
\ell: K^*\lrar \Pic(\CC)=(\CCline)^{\Gpd, \otimes}.
\end{align*}
Concretely, $\ell$ is a collection of one-dimensional  $\CC$-modules $\ell_\chi, \chi\in K^*$, with a compatible collection of isomorphisms, called the structure maps, 
\begin{align}\label{eqdefellchi12}
\ell_{\chi_1}\otimes_\CC\ell_{\chi_2}\ovs{\sim}{\lrar} \ell_{\chi_1\chi_2},
\end{align}
Let $\ell^{\triv}$ be the trivial line character. 

\item[(ii)] A \emph{framing $\eta$} of a line character $\ell$ of $K^*$ is a trivialization of $\ell$, i.e. an isomorphism $\ell$ to $\ell^\triv$. Concretely, $\eta$ consists of $\eta_\chi\in \ell_\chi-\{0\}$ for each $\chi\in K^*$, such that $\eta_{\chi_1}\otimes \eta_{\chi_2}\mapsto \eta_{\chi_1\chi_2}$ under the given structural maps. 

\end{itemize}
\end{defn}

The groupoid of line characters of $K^*$ is 
\begin{align}\label{eqXiKstar}
\Xi(K^*):=\Maps(K^*, \Pic(\CC))\simeq \BB K,
\end{align}
where $\Maps$ is taken in the ordinary category of commutative groupoids (in this case, it is equivalent to be taken in the $\infty$-category of commutative  $\infty$-groupoids). 
The groupoid of framed line characters $\wt{\Xi}(K^*)$, i.e. a line character equipped with a framing, forms the universal cover of $\Xi(K^*)\simeq \BB K$.

\sss{Definition of framings of the $K^*$-action on $\cW(\cT_{H^\flat})$}\label{sssframingKstWTHf}

Let $\pi \THf:=\pi_{\leqsl 1}\THf$ be the ordinary fundamental groupoid.  Let $\ell_{\pf, \chi}=\cL_\chi|_{\pf}$. We can rewrite $\phi_{\cL_K}$ \eqref{eqphicLK} as a functor of groupoids
\begin{align}\label{eqphiLKpiTHfXiKst}
\phi_{\cL_K}: \pi \THf&\lrar \Xi(K^*)\\
\nonumber \pf&\mapsto \ell_{\pf}=(\ell_{\pf, \chi})_{\chi\in K^*}
\end{align}

\begin{defn}\label{defnframeKst}
A \emph{framing} of the $K^*$-action on  $\cW(\cT_{H^\flat})$  (cf. Proposition \ref{lemmaKstarmuSh}) \emph{at any given $\pf\in \THf$} is a framing of $\ell_{\pf}$ as a line character of $K^*$. We call the data of  the $K^*$-action on  $\cW(\cT_{H^\flat})$ together with a framing at $\pf$ a \emph{framed} $K^*$-action at $\pf$. 
\end{defn}

\begin{remark}\label{remframingKactionmuSh}
For any $\pf\in \THf$, the category of framed $K^*$-actions at $\pf$ forms a universal cover of the space $\BB K$ in Proposition \ref{lemmaKstarmuSh}. 
The most relevant framings of the $K^*$-action on  $\cW(\cT_{H^\flat})$ are based at $\yf\in \UfS$, which will induce canonical mirror equivalences, cf.  
 \S\ref{subsecIWQstd} and especially Corollary \ref{corIWQzfeta} below.

\end{remark}

\begin{lemma}\label{lemmacanonicalgenerator}
Let $\pf\in \THf$. Any $K$-labeling $\mathbi{l}_p: \pi_{\cT, H}^{-1}(\pf)\to K$, i.e. an isomorphism as $K$-torsors,  that sends $p\in \pi_{\cT, H}^{-1}(\pf)$ to $1\in K$, determines a framing of the $K^*$-action on $\cW(\cT_{H^\flat})$ at $\pf$:  
\begin{align*}
\eta^p_\chi=\sum_{\kappa\in K}\chi(\kappa)^{-1}\delta_{\mathbi{l}_{p}^{-1}(\kappa)}\in \ell_{\pf, \chi},  
\end{align*}
where $\delta_{\mathbi{l}_p^{-1}(\kappa)}$ is the delta function at the point labeled by $\kappa$. This gives a $K$-equivariant bijection between $\pi_{\cT, H}^{-1}(\pf)$ and the set of framings of the  $K^*$-action on $\cW(\cT_{H^\flat})$ at $\pf$. 
\end{lemma}

\begin{proof}
The proof is straightforward. 
\end{proof}

\begin{prop}\label{propwtphiLK}
There is a natural Cartesian diagram of groupoids
\begin{equation*}
\begin{tikzcd}[row sep=0.2em]
\pi\cT_{H}\ar[r, "\wt{\phi}_{\cL_K}"]\ar[ddd]&\wt{\Xi}(K^*)\ar[ddd]\\
\ &\  \\
\ &\  \\
\pi\THf\ar[r, "\phi_{\cL_K}"]&\Xi(K^*)\\
\text{where}\quad \wt{\phi}_{\cL_K}: \pi\cTH\ar[r]& \wt{\Xi}(K^*)\\
p\ar[r, mapsto]& (\ell_{p^\flat}, \eta^p), 
\end{tikzcd},
\end{equation*}
and $\pf$ is the image of $p$ in $\THf$. 
\end{prop}

\begin{proof}
This follows from Lemma \ref{lemmacanonicalgenerator}. Alternatively, this follows from the standard theory of regular covering maps (with covering group $K$). 
\end{proof}

\subsection{Framed $K^*$-actions on $\Coh(\cS_{H^{\flat,\vee}})$ and relative equivariant coherent sheaves}

\sss{Framed $K^*$-actions on $\QCoh(\cS_{H^{\flat, \vee}})$}
On the mirror side, we have the canonical $K^*$-action on $\cS_{H^{\flat,\vee}}$, induced by the natural inclusion $K^*\hookrightarrow T^{\flat, \vee}$ that is dual to 
\begin{align}\label{eqqpi1TfK}
\wt{\psi}:\XX_*(T^\flat)=\pi_1(T^\flat)\twoheadrightarrow K. 
\end{align}
More explicitly, for any $\chi\in K^*$, let 
\begin{align}\label{eqwtpsichi}
\wt{\psi}_\chi:=\chi\circ\wt{\psi}: \XX_*(T^\flat)=\pi_1(\Tf)\lrar \CC^\times. 
\end{align}
Then the $K^*$-action on $\cO(T^{\flat, \vee})$ is given by 
\begin{align}
\label{eqbtauchi}\btau_\chi: \cO(T^{\flat, \vee})&\lrar \cO(T^{\flat, \vee})\\
\nonumber \sum_{\lambda\in \XX_*(T^\flat)} a_\lambda f^\lambda&\mapsto \sum_{\lambda\in \XX_*(T^\flat)} a_\lambda \wt{\psi}_\chi(\lambda)\cdot f^{\lambda},
\end{align}
where $f^\lambda$ denote the character function on $\Tfv$ corresponding to $\lambda\in \XX_*(T^\flat)$ and $a_\lambda=0$ except for a finite set of $\lambda$.

Let  $A=\cO(\cS_{H^{\flat, \vee}})= \cO(T^{\flat, \vee})^{W}\hookrightarrow \cO(T^{\flat, \vee})$. The above determines the $K^*$-action on $A$, denoted by the same notation 
\begin{align}\label{eqavarphichia}
\btau_{\chi}: A&\overset{\sim}{\lrar} A,
\end{align} 
This induces the canonical $K^*$-action on $\QCoh(\cS_{H^{\flat, \vee}})$, given by 
\begin{align*}
K^*\ni\chi\mapsto (F\mapsto (\btau_{\chi^{-1}})^*F).
\end{align*}
By identifying $\QCoh(\cS_{H^{\flat, \vee}})$ with $\Mod(A)$, the standard $K^*$-action is given by 
\begin{align}\label{eqMtoMchi}
K^*\ni\chi\mapsto& (M\mapsto M^{\chi}),
\end{align}
where $M^\chi$ represents the $\chi$-twisted $A$-module structure on $M$, on which $A$ acts by $A\overset{\btau_\chi}{\lrar} A\lrar \End_\CC(M)$.
The automorphism of this $K^*$-action is canonically identified with $K$, where for each $\kappa\in K$ we have the automorphism of $M^\chi$ given by multiplying $\chi(\kappa)$. This gives a full embedding of the $\infty$-groupoid
\begin{align}\label{eqBKFunKQCoh}
\BB K\hrar \Fun(K^*, \QCoh(\cSHfv)^{\otimes})^{\Gpd}. 
\end{align}

The following definition is the counterpart of Definition \ref{defnframeKst} on the mirror side.

\begin{defn}\label{defframedKQCoh}
A $K^*$-action on $\QCoh(\cS_{H^{\flat, \vee}})$ \emph{associated with a line character $\ell$ of $K^*$} (Definition \ref{defnlinecharKstar}) is the  $K^*$-action on $\Mod(A)$ given by 
\begin{align}\label{eqdefMchiellchi}
K^*\ni\chi\mapsto& (M\mapsto M^{\chi}\otimes_{\CC}\ell_\chi),
\end{align}
where $A$ acts on the latter $\ell_\chi$ trivially and the isomorphism on the right-hand-side is the canonical one. 
A \emph{framing} $\eta$ of the above $K^*$-action is 
a framing $\eta$ of $\ell$. 
\end{defn}

Note that the $K^*$-action \eqref{eqdefMchiellchi} is well defined--- the composition of the action of $\chi_1$ and $\chi_2$ is canonically identified with the action of $\chi_1\chi_2$: 
\begin{align}\label{eqMchi1chi2ellchi}
M\overset{\chi_1}{\mapsto }M^{\chi_1}\otimes_\CC\ell_{\chi_1}\overset{\chi_2}{\mapsto }(M^{\chi_1})^{\chi_2}\otimes_\CC\ell_{\chi_1}\otimes_\CC \ell_{\chi_2}\cong M^{\chi_1\chi_2}\otimes_\CC\ell_{\chi_1\chi_2}. 
\end{align}
The standard $K^*$-action on $\QCoh(\cS_{H^{\flat, \vee}})$ is associated with $\ell^{\triv}:=(\CC)_{\chi\in K^*}$, with the structure map \eqref{eqdefellchi12} given by the standard multiplication map on $\CC$. 

We have the following counterpart of Remark \ref{remframingKactionmuSh} for framed $K^*$-actions on $\QCoh(\cS_{H^{\flat, \vee}})$.

\begin{remark}\label{remstdframingQCoh}
For any two framed $K^*$-actions on $\QCoh(\cSHfv)$, given by the data $(\ell^{i}, \eta^{i}), i=1,2$, we have the canonical isomorphism
\begin{align*}
\ell^1_\chi&\ovs{\sim}{\lrar} \ell^2_\chi\\
\eta^1_\chi&\mapsto \eta^2_\chi,
\end{align*}
that induces the canonical isomorphism between the two resulting $K^*$-actions on $\QCoh(\cSHfv)$.
There is also the free $K$-action: 
\begin{align}\label{eqtwistframeetakappa}
K\ni\kappa: (\ell_\chi, \eta_\chi)_\chi\mapsto  (\ell_\chi, \chi(\kappa)\cdot \eta_\chi)_\chi.
\end{align} 
Thus the collection of framed $K^*$-actions on $\QCoh(\cSHfv)$ forms a contractible ($\infty$-)groupoid with a free $K$-action. The quotient by this $K$-action gives the essential image of $\BB K$ in \eqref{eqBKFunKQCoh}. 
 \end{remark}

\sss{Equivariant quasi-coherent sheaves relative to a framing}
For any $K^*$-action on $\QCoh(\cSHfv)$ associated with $\ell$, denote the corresponding $K^*$-equivariant category by $\QCoh(\cSHfv)^{K^*;\ell}$. 
A framing $\eta$ determines a unique equivalence between this $K^*$-action on $\QCoh(\cSHfv)$ and the standard $K^*$-action, through the unique isomorphism $(\ell, \eta)\overset{\sim}{\to} (\ell^\triv, 1)$. 
In particular, this determines a unique equivalence between their respective $K^*$-equivariant categories, denoted by 
\begin{align}
\Phi_{\ell\to \std}^\eta: \QCoh(\cSHfv)^{K^*;\ell}\ovs{\sim}{\lrar} \QCoh(\cSHfv)^{K^*}.
\end{align} 

We have the notion of $\kappa$-twisted equivariant structure sheaf $\cO^\kappa$ for the standard $K^*$-action, whose equivariant structure on $M=A$, up to isomorphisms, is given by
\begin{align*}
\btau_\chi^{\kappa}: A^{\chi}&\ovs{\sim}{\lrar} A\\
a^{\chi}&\mapsto \chi(\kappa)\btau_{\chi^{-1}}(a). 
\end{align*} 
where $\btau_{\chi^{-1}}$ is given in \eqref{eqavarphichia}.

\begin{defn}\label{defnkappatwetaQCoh}
Consider the $K^*$-action on $\QCoh(\cSHfv)$ associated with $\ell$. 
For any framing $\eta$,   
define the \emph{$\kappa$-twisted equivariant structure sheaf relative to the framing $\eta$}, denoted by $\cO^{\kappa;\ell}_{\eta}$,  to be $(\Phi_{\ell\to \std}^{\eta})^{-1}(\cO^{\kappa})$ (up to isomorphisms) .
\end{defn}

Note that by definition $\cO^{\kappa;\ell}_\eta$ is isomorphic to the $A$-module $M=A$ equipped with the equivariant structure maps with respect to the $K^*$-action associated with $\ell$: 
\begin{align*}
\varphi_{\chi;\eta}^{\kappa;\ell}: A^{\chi}\otimes_\CC \ell_{\chi}&\ovs{\sim}{\lrar} A\\
a^{\chi}\otimes \eta_{\chi}&\mapsto \chi(\kappa)\btau_{\chi^{-1}}(a). 
\end{align*} 
Twisting the framing $\eta$ of $\ell$ by $\kappa_1\in K$ \eqref{eqtwistframeetakappa} results in a canonical identification 
\begin{align*}
\cO_{\eta}^{\kappa;\ell}\cong \cO_{\kappa_1\cdot \eta}^{\kappa\kappa_1;\ell},\quad \forall \kappa\in K
\end{align*}
in $\QCoh(\cSHfv)^{K^*;\ell}$.

\sss{The $\infty$-groupoid $\frG^{K^*}$}\label{sssfrGKstar}

Consider the $\infty$-groupoid $\frG^{K^*}\simeq \Fun(\BB K^*, \Pr_\st^L)^{\Gpd}$ of presentable stable $\infty$-categories $\cC$ equipped with $K^*$-actions. Let  $\IWf$ denote $\Ind\cW(\cT_{H^\flat})$ with its canonical $K^*$-action, as an object in $\frG^{K^*}$.  Let $\Qf$ denote the object $\QCoh(\cS_{H^{\flat,\vee}})$ equipped with the \emph{standard} $K^*$-action. 
A framed  $K^*$-action on $\QCoh(\SHfv)$ gives an object in $\frG^{K^*}_{/\cQ_\flat}$ in the obvious way. 
We have the natural commutative diagram (but \emph{not} Cartesian in general) of $\infty$-groupoids
\begin{equation}\label{eqcQXiwt}
\begin{tikzcd}
\wt{\Xi}(K^*)\ar[r, "\wt{\cQ}^\flat_{\Xi}"]\ar[d] &\frG^{K^*}_{/\cQ_\flat}\ar[d]\\
\Xi(K^*)\ar[r, "\cQ_{\Xi}"] &\frG^{K^*},
\end{tikzcd}
\end{equation}
where $\cQ^\flat_{\Xi}$ is the natural functor that sends $\ell\in \Xi(K^*)$ to the $K^*$-action on $\QCoh(\SHfv)$ associated to $\ell$, and $\wt{\cQ}^\flat_{\Xi}$ is the natural functor that sends $(\ell, \eta)$ to the corresponding framed $K^*$-action on $\QCoh(\SHfv)$.

\subsection{Isomorphisms from $\IWf$ to $\Qf$ in $\frG^{K^*}$ for a fixed $\yf$}\label{subsecIWQstd}

\sss{The functor $\mu_{\yf}$}\label{sssBaseptKostant}

Now we use notations from \S\ref{sssPJflatPHJ}. 
For each 
$\yf\in U^{\flat,\SS}:=U_{I_H}^{\flat, \SS}$, we have the canonical functor,
\begin{align}\label{eqmuzJflatcG}
&\mu_{\yf}:=\Hom_{\IndW(\cT_{H^\flat})}(S_{\yf}, -): \Ind\WTHf\lrar \Mod(\CC).
\end{align}
This corresponds to taking microlocal stalks at $\yf$ on the microlocal sheaf side. 
Similarly, for any $J\subset I_H$ and $\yJf\in \UJfS$, we have 
\begin{align*}
&\mu_{\yJf}:=\Hom_{\IndW(\cT_{H^\flat})}(S_{\yJf}, -): \Ind\WTHf\lrar \Mod(\CC).
\end{align*}
For any $J\subsetneq I_H$, let 
\begin{align*}
&\mu_{\yJf}^J:=\Hom_{\Ind\WTLJf}(S_{\yJf}, -): \Ind\WTLJf\lrar \Mod(\CC).
\end{align*}
Then we have the natural isomorphism of functors
\begin{align}\label{eqmuzJfres}
\mu_{\yJf}^J\circ\res^{I_H}_{J}\cong \mu_{\yJf}. 
\end{align}

\begin{remark}
Let $Y_H=\{g\in w_0T|\Ad(g)\psi_{\infty,I_H}=-\psi_{0,I_H}\}$ (for $L_J$ it was denoted by $\cW_J$) is a $ZH$-torsor. 
Using the notaions from \S\ref{sssClosureKostant}, we have a canonical identification $\cX^{H}\cap \SS^H\cong Y_H$. 
It is easy to see that the natural inclusion $U_H^{\SS}:=U_{I_H}^{\SS}\hookrightarrow \cX^{H}\cap \SS^H\cong Y_H$ is a homotopy equivalence and 
the image is $Z(H)_{\cpt}$-invariant.  
Therefore, we have a natural isomorphism $\pi U_H^\SS\ovs{\sim}{\to}\pi Y_{H}$ as $\pi Z(H)_\cpt$-torsors.  
Here and after (e.g. \S\ref{subsecfrPHtoWTH} below), 
one can safely replace $U_H^\SS$ (resp. $\pi U_H^\SS$) with $Y_H$ (resp. $\pi Y_H$). 
\end{remark}

\sss{A Diagram of mirror equivalences} 
Let 
\begin{align*}
&A_{\yf}= \End(S_{\yf})^{op}, \\
&A^\vn=\End_{\cW(\cT_{L_\vn^\flat})}(S_{y_\vn^\flat})^{op}=\cO(T^{\flat, \vee})^{op}. 
\end{align*}
Since $A_{\yf}$ is commutative and concentrated in degree $0$, for any isomorphism $S_{\yfO}\overset{\sim}{\to} S_{\yfT}$, it induces the canonical identification $A_{\yfO}= A_{\yfT}$. Alternatively, one can use the $\pi\ZHf$-action on $\WTHf$ \S\ref{subsecpiZHactionW} to show that the isomorphisms in the image of \eqref{eqrho12sigmay} for $\Hf$ canonically identifies $A_{\yfO}$ and $A_{\yfT}$. Therefore, we can set 
\begin{align}\label{eqAHf}
A_{\Hf}:=\cO(\Hfv\sslash \Hfv)^{op}=A_{\yfO}= A_{\yfT}. 
\end{align}

\sss{Convention about opposite groups and algebras}\label{sssconventionoppositealg}
 In the following, for any group or algebra homomorphism $\varphi: A\to B$, let $\varphi^{op}: A^{op}\to B^{op}$ be the opposite version of $\varphi$. 

Most of the time, the groups and algebras $A$ that we encounter are commutative, so there is the trivial identification $A^{op}=A$. However, for the algebras $\Avn$ and $\AHf$ that we will focus on, we will identify them with their respective opposite versions by a nontrivial involution---this is determined by the canonical involution on the fundamental groupoid $\pi(\Tf)^{op}\ovs{\sim}{\to} \pi(\Tf)$ (cf. \S\ref{sssiotHf}). 

Most of the time, the distinction between $\Avn$ (resp. $\AHf$) and $A^{\vn, op}$ (resp. $\AHf^{op}$) is \emph{not} important and can be ignored.  The distinction is needed when we calculate the tensor action of $\cL_\chi, \chi\in K^*$ on $\WTHf$ in \S\ref{ssstensorLocUHflat}. Since $\cL_\chi|_{\cT_{\Lfvn}}$ is classified by \eqref{eqwtpsichi}, 
while $\Ind\WTHf$ is identified with the module category of $\Avn=\CC[\pi_1(\Tf)^{op}]$, we need to use the opposite version of \eqref{eqwtpsichi} to identify the tensor action of $\cL_\chi$ with the tensor action of an $\Avn$-module. 
This becomes even more important when we calculate the $\pi \ZH$-action on the mirror categories of  $\WTHf$ and $\WTH$ in \S\ref{subsecwtPhiHPop} and \S\ref{subsecpiZHQft} (especially to give the correct statements in Proposition \ref{eqbiakappa} and Corollary \ref{corpiZHQfKstar}), where the canonical involution on the fundamental groupoid $\pi\ZHf^{op}\ovs{\sim}{\to} \pi \ZHf$ must be taken into account. The latter determines that when we identify $\cL_\chi|_{\cT_{\Lfvn}}$ as an $\Avn=\CC[\pi_1(\Tf)^{op}]$-module, we need to use the canonical involution $\pi_1(\Tf)^{op}\ovs{\sim}{\to} \pi_1(\Tf)$ (rather than the trivial identity).

\sss{The (anti-)involution between $\Avn$ (resp. $\AHf$) and $\Avnt=\cO(\Tfv)$ (resp. $\AHft=\cO(\Tfv)^W$)}\label{sssiotHf}
To simplify notations, let $\Avnt=A^{\vn, op}$ and $\AHft=\AHf^{op}$. 
Let 
\begin{align}\label{eqiotvn}
\iotvn: \Avn&\ovs{\sim}{\lrar} \Avnt\\
\nonumber f^{\lambda}&\mapsto f^{-\lambda}, \quad \forall \lambda\in \XX_*(\Tf). 
\end{align}
be the (anti-)involution that is induced from the natural involution on the fundamental groupoid $\pi(\Tf)^{op}\ovs{\sim}{\to} \pi(\Tf)$. 
This induces the involution $\iotHf$ by the commutative diagram
\begin{equation*}
\begin{tikzcd}
\Avn\ar[r, "\sim"', "\iotvn"]&\Avnt\\
(\Avn)^W=\AHf\ar[r, "\sim"', "\iotHf"]\ar[u, hook]&\AHft=(\Avnt)^W\ar[u, hook]
\end{tikzcd}. 
\end{equation*}

Let 
\begin{align}
\label{eqwtpsit}& \wt{\psi}^t:=(\wt{\psi}^{op})^{-1}: \pi_1(\Tf)^{op}\ovs{\wt{\psi}^{op}}{\lrar} K^{op}\ovs{(-)^{-1}}{\lrar} K,\\
\label{eqwtpsichit}&\wt{\psi}_{\chi}^t: \pi_1(\Tf)^{op}\ovs{\wt{\psi}^t}{\lrar}K\ovs{\chi}{\lrar}\CC^\times, 
\end{align}

Using the exact sequence
\begin{align*}
1\lrar \pi_1 (T)\lrar \pi_1(\Tf) \lrar K\lrar 1,
\end{align*}
we can make the identification of groupoids  
$\pi T\cong K/\pi_1 \Tf=K/\XX_*(\Tf)$. For any $\k_1, \k_2\in K$ and any $\gaTO\in \Hom_{\pi (T)^{op}}(\k_1, \k_2)=\Hom_{\pi T}(\k_2, \k_1)$, let $\gafTO\in  \Hom_{\pi (T)^{op}}(1, 1)=\XX_*(\Tf)^{op}$ be the projection of $\gaTO$. The following lemma is direct to check. It shows that $\wt{\psi}^t$ and $\wt{\psi}_{\chi}^t$ above are the correct ``opposite version" of $\wt{\psi}$ and $\wt{\psi}_\chi$, respectively, as it intertwines with $\iotvn$. In particular, item (i) of the following says that $\wt{\psi}^t$ is the character valued in $K$ that takes $\gafTO$ to the correct difference between its source and target in $\pi (T)^{op}$. This will be needed in Proposition \ref{eqbiakappa}. 

\begin{lemma}\label{lemmawtpsit}
\begin{itemize}
\item[(i)] For any $\gaTO\in \Hom_{\pi (T)^{op}}(\k_1, \k_2)=\Hom_{\pi (T)}(\k_2, \k_1)$ as above, let $\gaOT\in \Hom_{\pi T}(\k_1, \k_2)$ be the inverse path of $\gaTO$. Then we have  
\begin{align*}
\wt{\psi}^t(\gafTO)=\k_1^{-1}\k_2=\wt{\psi}(\iotvn(\gafTO))=\wt{\psi}(\gafOT). 
\end{align*}

\item[(ii)] Under the canonical involution $\pi_1(\Tf)^{op}\ovs{\sim}{\to}\pi_1(\Tf)$, $\cL_\chi|_{T^*\Tf}$ is the $\pi_1(\Tf)^{op}$-module classified by $\wt{\psi}_{\chi}^t$. 
\end{itemize}
\end{lemma}

Let 
\begin{align}\label{eqbtauchit}
\btau_{\chi}^t: \Avn=\CC[\pi_1(\Tf)^{op}]&\lrar \Avn=\CC[\pi_1(\Tf)^{op}]\\
\nonumber \sum_{\lambda\in \XX_*(\Tf)} a_\lambda f^\lambda&\mapsto \sum_{\lambda\in \XX_*(\Tf)} a_\lambda \wt{\psi}_{\chi}^t(\lambda) f^{\lambda}
\end{align}
be the ``opposite version" of $\btau_\chi$ \eqref{eqavarphichia}. It is easy to see that 
\begin{align}
\nonumber &\wt{\psi}_{\chi}^t\circ (\iotvn)^{-1}|_{\XX_*(\Tf)}=\wt{\psi}_{\chi}|_{\XX_*(\Tf)} \\
\label{eqtauchitvn}\Longrightarrow\ &\btau_{\chi}^t=(\iotvn)^ {-1}\circ\btau_{\chi}\circ \iotvn: \Avn\ovs{\sim}{\lrar} \Avn\\
\label{eqtauchitHf}\Longrightarrow\ &\text{(the restriction of) }\btau_{\chi}^t=(\iotHf)^{-1}\circ\btau_{\chi}\circ \iotHf: \AHf\ovs{\sim}{\lrar} \AHf. 
\end{align}

Now by Theorem \ref{thmJconnectedcenter}, the mirror equivalence relative to $S_\yf$ is given by 
\begin{align}\label{eqUpzJflmuShH}
\Upsilon_{\yf}: \Ind\WTHf&\ovs{\sim}{\lrar} \Mod(A_{\Hf})\\
\nonumber\cG&\mapsto A_{\Hf}\car \mu_{\yf}(\cG).
\end{align}
We will also refer to $\Upsilon_{\yf}$ as the mirror equivalence \emph{based at $\yf$}. 
Similarly, for any $J\subset I_H$ and $\yJf\in \UJfS$, we have 
\begin{align}\label{eqUpJzJfA}
\Upsilon^J_{\yJf}: \Ind\WTLJf&\ovs{\sim}{\lrar} \Mod(A_{\LJf})\\
\nonumber \cG&\mapsto A_{\LJf}\car \mu_{\yJf}^J(\cG). 
\end{align}

\begin{warn}
Let $E_{\yJf}=\End_{\WTHf}(S_{\yJf})^{op}$, for $J\subsetneq I_H$. 
Then one can similarly define 
\begin{align*}
\Upsilon_{\yJf}: \Ind\WTHf&\ovs{\sim}{\lrar} \Mod(E_{\yJf})\\
\nonumber\cG&\mapsto E_{\yJf}\car \mu_{\yJf}(\cG).
\end{align*}
Although we have \eqref{eqmuzJfres}, in general
\begin{align*}
\Upsilon_{\yJf}^J\circ\res^{I_H}_J\not\cong \Upsilon_{\yJf},
\end{align*}
due to $E_{\yJf}\not\cong A_{\yJf}$. This is why we put a superscript $J$ on $\Upsilon^J_{\yJf}$ \eqref{eqUpJzJfA}. 
\end{warn}

Using Theorem \ref{thmJparabolicconnectedcenter}, for any isomorphism 
\begin{align}\label{eqvarrhoSyfvres}
\varrho: S_{\yfv}\ovs{\sim}{\to}\res_\vn S_{\yf} \text{ in } \cW(\cT_{\Lfvn}), 
\end{align}
where $\res_\vn=\res^{I_H}_\vn=\res^{\Hf}_{\Tf}$, 
we have the following natural commutative diagram 
\begin{equation}\label{diagUpsilonzJflatyJ} 
\begin{tikzcd}[column sep=8em, row sep=2em]
 \Ind\WTHf\ar[d, "\res_\vn:=\res^{I_H}_\vn"']\ar[r,"{\Upsilon_{\yf}}", "\sim"']&\Mod(A_{\Hf})=\Mod(\End(S_{\yf})^{op})\ar[d, "\otimes_{\AHf}\Avn"]\\
 \Ind\cW(\cT_{L_\vn^\flat})\ar[r,"{\Upsilon^\vn_{\yf}:=\Hom(\res_\vn S_{\yf}, -)}", "\sim"',  ""{name=U, below}]\ar[r, bend right=2em, "\Upsilon^\vn:=\Upsilon^\vn_{y_\vn^\flat}"', "\sim",  ""{name=D, above}]&\Mod(A^{\vn})
\arrow[from=U, to=D, Rightarrow, "\sim", "\varrho^* "']
\end{tikzcd}, 
\end{equation}
in which $\Mod(\Avn)$ is canonically identified with $\Mod(\End(\res_\vn S_{\yf})^{op})$ for $\Upsilon^\vn_\yf$ and $\Mod(\End(S_{\yfv})^{op})$ for $\Upsilon^\vn$, respectively. Note that without loss of generality, one may take $\yfv=\zfv$.

\sss{Calculation of the action of tensoring by $\cL_\chi\in \Loc(\cT_{H^\flat})^{\bd, \hs;\otimes}, \chi\in K^*$}\label{ssstensorLocUHflat}

We have included the details of the calculation of the tensor action by any rank 1 local system $\cL$ on $\cW(\cT_{\Hf})$ in \S\ref{appendssstensorLocTHf}. One of the main results is Lemma \ref{lemmaMchicLyf}. 
In the special case that $\cL=\cL_\chi$, let $\alpha^\chi$ (resp. $\alpha^\chi_\vn$) denote $\alpha^{\cL_\chi}$ (resp. $\alpha^{\cL_\chi}_\vn$). 
For any $\Avn$ (resp. $\AHf$) module $M_\vn$ (resp. $M$), let $M_\vn^\chi$ (resp. $M^\chi$) be the $\chi$-twisted $\Avn$ (resp. $\AHf$) module defined using $\btau_\chi^t$ in place of $\btau_\chi$ \eqref{eqavarphichia}. Then Lemma \ref{lemmaMchicLyf} specializes to the following.

\begin{cor}\label{lemmaAJAvnchimuzyflat}
\begin{itemize}
\item[(i)]
The $A_{\Hf}$-module (resp. $A^{\vn}$-module) structure on $\mu_{\yf}(\cG)^{\cL_\chi}\otimes_{\CC} \ell_{\yf, \chi}$ (resp. $\mu_{\yfv}(\cG)^{\cL_\chi^\vn}\otimes_{\CC}\ell_{\yfv, \chi}$) is canonically identified $\mu_{\yf}(\cG)^\chi\otimes_{\CC} \ell_{\yf, \chi}$ (resp. $\mu_{\yfv}(\cG)^\chi\otimes_{\CC}\ell_{\yfv, \chi}$).

\item[(ii)] The automorphism $\alpha^\chi$ is \emph{canonically} identified with 
\begin{align*}
M\mapsto M^\chi\otimes_\CC \ell_{\yf, \chi},
\end{align*}
which is independent of the choice of $\varrho$ \eqref{eqvarrhoSyfvres}.  
\end{itemize}
\end{cor}

\begin{cor}\label{corUpzJflatKstarequiv}
Under $\Upsilon_{\yf}$, the induced $\Kdu$-action on $\Mod(\AHf)$ is canonically isomorphic to the $\Kdu$-action associated with the collection of lines $\ell_{\yf}=(\ell_{\yf,\chi})_\chi$ (cf. Definition \ref{defframedKQCoh}). 
\end{cor}

\begin{proof}
By Corollary \ref{lemmaAJAvnchimuzyflat} (ii), the induced automorphism $\alpha^\chi$ is canonically  isomorphic to \eqref{eqdefMchiellchi} for $\ell_\chi=\ell_{\yf, \chi}$. Since  $A_{\Hf}$ is concentrated in degree $0$ and the two $\Kdu$-actions both preserve the natural $t$-structure on $\Mod(A_{\Hf})$, it suffices to check the two $\Kdu$-actions are canonically identified after passing to cohomology modules of $A_{\Hf}$ (equivalently, the two induced $\Kdu$-actions on $\Mod(\AHf)^\hs$ canonically agree). 
This is straightforward to check. It suffices to check the structure maps
$m_{\chi,\chiO}: \alpha^{\chi}\alpha^{\chiO}\ovs{\sim}{\to}\alpha^{\chi\chiO}$ 
agrees with 
\eqref{eqMchi1chi2ellchi} after passing to cohomology modules, for $\ell_\chi=\ell_{\yf, \chi}$. By definition, for each $\cG\in \Ind\WTHf$, the structure map $m_{\chi,\chiO}$ is canonically given by 
\begin{equation*}
\begin{tikzcd}
\mu_{\yf}(\cG)\ar[drr, "\alpha^{\chi\chiO}"', mapsto]\ar[r, "\alpha^{\chiO}", mapsto]& \mu_{\yf}(\cG)^{\chiO}\otimes_{\CC} \cL_{\chiO}|_{\yf}\ar[r, "\alpha^{\chi}", mapsto]& (\mu_{\yf}(\cG)^{\chiO})^{\chi}\otimes_{\CC} \cL_{\chiO}|_{\yf} \otimes_{\CC} \cL_{\chi}|_{\yf}\ar[d, "\sim"],\\
&&\mu_{\yf}(\cG)^{\chiO\chi}\otimes_{\CC} \cL_{\chiO\chi}|_{\yf},
\end{tikzcd}
\end{equation*}
where the rightmost vertical isomorphism is defined by (1) the canonical isomorphism $\cL_{\chiO}|_{\yf} \otimes_{\CC} \cL_{\chi}|_{\yf}\cong \cL_{\chiO\chi}|_{\yf}$; (2) a canonical identification $(\mu_{\yf}(\cG)^{\chiO})^{\chi}\cong \mu_{\yf}(\cG)^{\chiO\chi}$ that 
is the identity map on the underlying $\CC$-module $\mu_{\yf}(\cG)$. 
Therefore, $m_{\chi, \chiO}$ agrees with \eqref{eqMchi1chi2ellchi} after passing to cohomology modules. 
This finishes the proof. 
\end{proof}

Recall the $\infty$-groupoid $\frGK$ and $\IWf, \Qf\in \frGK$ from \S\ref{sssfrGKstar}. 
Let $\Qft=\Mod(\AHf)\in  \frGKt$ with the standard $\Kdu$-action (given induced by $\btau^t_\chi, \chi\in K^*$) that is the opposite version of the standard $K^*$-action on $\AHf^{op}=\cO(\Tfv)^W$. We view $\IWf\in \frGKt$ with the tensor action of $\cL_{\chi}, \chi\in \Kdu$. 
Let $\cQ_{\yf}\in \frGKt$ denote the object $\Mod(\AHf)$ equipped with the $\Kdu$-action associated with $\ell_{\yf}$. 

\begin{cor}\label{corIWQzfeta}
The mirror equivalence $\Upsilon_{\yf}$ based at $\yf$ gives a canonical isomorphism from $\IWf$ to $\cQ_{\yf}$ in $\frGKt$. In particular,  any framing $\eta$ of  $\ell_{\yf}$ determines a canonical isomorphism 
\begin{align*}
\IWf\ovs{\Upsilon_\yf}{\underset{\sim}{\lrar}}\cQ_{\yf}\overset{\eta}{\underset{\sim}{\lrar}}\Qft.
\end{align*}
and this gives a $K$-torsor of isomorphisms from $\IWf$ to $\Qft$ in $\frGKt$. 
\end{cor}
\begin{proof}
This directly follows from Corollary \ref{corUpzJflatKstarequiv} and Remark \ref{remstdframingQCoh}. 
\end{proof}

Combining with Lemma \ref{lemmacanonicalgenerator}, we directly get 
\begin{cor}\label{corfiberfrgammaJf}
For any $\yf\in \UfS$, let $\bil_y: \pi_{\cT, H}^{-1}(\yf)\to K$ be the labeling that sends $y\in \pi_{\cT, H}^{-1}(\yf)$ to $1\in K$. Let $\eta^{y}$ be the framing of $\ell_{\yf}$ determined by the labeling $\bil_y$ from Lemma \ref{lemmacanonicalgenerator}. 
Then the $K$-torsor $\pi_{\cT, H}^{-1}(\yf)$ determines a $K$-torsor of isomorphisms in $\frGKt$
\begin{equation*}
\begin{tikzcd}
\wUpy: \IWf\ar[r, "\wt{\Upsilon}_{\yf}", "\sim"']&\cQ_{\yf}\ar[r, "\eta^y", "\sim"']& \Qft,\quad\text{for } y\in \pi_{\cT, H}^{-1}(\yf). 
\end{tikzcd}
\end{equation*}
\end{cor}

%%%%%%%%%%%%%%%%%%%%%%%%%%%%%%%%%%%%%%
 \subsection{The $\pi \ZHc$-equivariant functor $\nu_H: \frP_H\to \WTH$}\label{subsecfrPHtoWTH}
 
  Let $\frPH$ denote the groupoid $\pi U_{H}^{\SS}$, which is canonically a $\pi Z(H)$-torsor. 
The inclusion $U_H^{\SS, \cpt}\hookrightarrow U_{H}^{\SS}$ is a homotopy equivalence and has \emph{no} effect when we replace the definition of $\frP_H$ by $\pi U_{H}^{\SS,\cpt}$. We will use the latter definition of $\frP_H$ below only for cosmetic reasons. 

 Using the presentation from \S\ref{sssconcretepiZH}, for any $y_i=z_iy, i=1,2$ in $\frP_H$, where $y\in \frP_H$ and $z_i\in \ZHc$, we identify
 \begin{align*}
\Hom_{\frP_H}(y_1, y_2)=  \Hom_{\frP_H}(z_1y, z_2y)=\{\sigma_y(\gaOT):=y\cdot \gaOT: \gaOT\in \Hom_{\pi \ZHc}(z_1, z_2)\text{ are geodesics}\}.
 \end{align*}
This presentation is clearly independent of the choice of $y\in \frP_H$.

Using Remark \ref{remsigmayHindpresentation} (and the canonical equivalence $\WTH\simeq \cW(\cT_H')$), we get a canonical $\pi \ZHc$-equivariant functor
\begin{align}
\nonumber \nu_H: \frP_H&\lrar \WTH\\
\nonumber y&\mapsto S_y\\
\label{eqrho12sigmay}\rho_{1,2}:=\sigma_{y}(\gamma_{1,2})\in \Hom_{\frP_H}(y_1, y_2)&\mapsto \rho^H_{1,2}:=\sigma^H_{y}(\pathz_{1,2})\in \Hom_{\WTH}(S_{y_1}, S_{y_2}),
\end{align}
for any splitting $y_i=z_i y$ and $\gaOT\in \Hom_{\pi \ZHc}(z_1, z_2)$. It is clear that $\nu_H$ is well defined, and it is naturally compatible with finite central quotients 
\begin{equation*}
\begin{tikzcd}
\frP_H\ar[r, "\nu_H"]\ar[d]&\WTH\ar[d, "F_*"]\\
\frP_{H^\flat}\ar[r, "\nu_{H^\flat}"]&\WTHf
\end{tikzcd}. 
\end{equation*}

\sss{Calculation of $\nu_{\Hf}$}
For any $\sfL\in \WTH$, we have the natural functor from the canonical action of $\pi \ZHc$ on $\WTH$
\begin{align*}
a_\sfL: \pi \ZHc&\lrar \WTH\\
z&\mapsto z\cdot \sfL. 
\end{align*}
For $\sfL=S_y$, denote $a_\sfL$ by $a_y$. 
The compatibility with restriction \S\ref{ssspiZHcompatible} induces the natural commutative diagram
\begin{equation*}
\begin{tikzcd}
\pi Z(H^\flat)_\cpt\ar[r, "a_{\yf}"]\ar[d]&\WTHf\ar[d, "\res_J"]\\
\pi Z(L_J^\flat)_\cpt \ar[r, "a_{\res_JS_{\yf}}"]&\WTLJf
\end{tikzcd}. 
\end{equation*}
Recall that there is a non-canonical isomorphism $\res_JS_{\yf}\cong S_{\yJf}$ in $\WTLJf$, for any $\yJf\in \UJfS$.

\begin{lemma}\label{lemmaayfinclusion}
For any $\zfi\in \ZHfc, i=1,2$, $a_{\yf}$ induces an inclusion 
\begin{align*}
a_{\yf}: \Hom_{\pi \ZHfc}(\zfO, \zfT)\hookrightarrow  H^0(\Hom_{\WTHf}(S_{\zfO\cdot \yf}, S_{\zfT\cdot \yf}))\cong \Hom_{\WTHf}(S_{\zfO\cdot \yf}, S_{\zfT\cdot \yf}).
\end{align*}
For $\zfO=\zfT=\zf$, the inclusion is naturally identified with the composite canonical inclusion
\begin{align*}
\iota^\flat: \XX_*(Z(H^\flat))=\XX^*(H^{\flat, \vee;\ab})\hookrightarrow \cO(H^{\flat, \vee;\ab})^\times=\cO(H^{\flat, \vee}\sslash H^{\flat, \vee})^\times\hookrightarrow\cO(H^{\flat, \vee}\sslash H^{\flat, \vee})=\cO(T^{\flat, \vee}\sslash W). 
\end{align*}
\end{lemma}

\begin{proof}
We only show the case for $\zfO=\zfT=\zf$; the general case easily follows from it (it also follows from a similar argument).  This follows from Theorem \ref{thmJconnectedcenter} and \S\ref{ssspiZHcompatible}. 
Indeed, we have the commutative diagram coming from $\res_\vn$: 
\begin{equation*}
\begin{tikzcd}
\XX_*(Z(H^\flat))\ar[r, hook] \ar[d, hook]&\cO(H^{\flat, \vee;\ab})^\times=\cO(H^{\flat, \vee}\sslash H^{\flat, \vee})^\times\ar[d, hook]\ar[r, hook]&\cO(T^{\flat, \vee}\sslash W)=\End_{\WTHf}(S_{\zf\cdot \yf})\ar[d, hook]\\
\End_{\pi T^\flat_\cpt}(\zf)=\XX^*(T^{\flat, \vee})\ar[r, hook]&\cO(T^{\flat,\vee})^\times\ar[r, hook] &\cO(T^{\flat,\vee})\cong  \End_{\cW(T^*T^\flat)}(\res_\vn S_{\zf\cdot \yf})
\end{tikzcd}
\end{equation*}
where the bottom row induced from $a_{\res_\vn S_{\zf\cdot \yf}}\cong a_{\yvnf}$ is the obvious one.  
\end{proof}

\begin{cor}\label{cornuHflat}
For any $\yfi\in \frP_{\Hf}, i=1,2$, $\nu_{\Hf}$ induces an inclusion 
\begin{align*}
\nu_{\Hf}: \Hom_{\frPHf}(\yfO, \yfT)\hookrightarrow  H^0(\Hom_{\WTHf}(S_{\yfO}, S_{\yfT}))\cong \Hom_{\WTHf}(S_{\yfO}, S_{\yfT}).
\end{align*}
For $\yfO=\yfT=\yf$, the inclusion is naturally identified with the composite canonical inclusion
\begin{align*}
\iota^\flat: \XX_*(Z(H^\flat))=\XX^*(H^{\flat, \vee;\ab})\hookrightarrow \cO(H^{\flat, \vee;\ab})^\times=\cO(H^{\flat, \vee}\sslash H^{\flat, \vee})^\times\hookrightarrow\cO(H^{\flat, \vee}\sslash H^{\flat, \vee})=\cO(T^{\flat, \vee}\sslash W). 
\end{align*}
\end{cor}

\begin{proof}
For any $\yf\in \frPHf$, we have the natural commutative diagram
\begin{equation*}
\begin{tikzcd}[column sep=1em]
\zf\in\ar[d, mapsto]& \pi \ZHfc\ar[d, "\sigma_{\yf}"']\ar[drrr, "a_{\yf}"]&  &\\
\zf y^\flat\in & \frPHf\ar[rrr, "\nu_H"]& & &\WTHf
\end{tikzcd},
\end{equation*}
where $\sigma_{\yf}$ is induced from the $ \pi \ZHfc$-action on $\frPHf$. The corollary then follows directly from Lemma \ref{lemmaayfinclusion}. 
\end{proof}

\begin{remark}\label{remnuLcompatiblevarrho}
We remark that $\nu^\cL$ \eqref{eqnuupcL} is functorial in $(\yHf, \yLf;\varrho)$ in the following sense. 
 Given another pair  $(y_H^{\prime, \flat}, y_L^{\prime, \flat}; \varrho')$,  $\varphi_{\Hf}: S_{\yHf}\ovs{\sim}{\to} S_{y_H^{\prime, \flat}}$ and $\varphi_{\Lf}: S_{\yLf}\ovs{\sim}{\to} S_{y_L^{\prime, \flat}}$ such that $\res_{\Lf}^{\Hf}(\varphi_{\Hf})$ and $\varphi_{\Lf}$ are conjugated by $\varrho$ and $\varrho'$, 
Let $\gamma_{\Hf}$ represents the homotopy class of paths from $\yHf$ to $\yHpf$ corresponding to $\varphi_{\Hf}$ (and similarly for $\gamma_{\Lf}$), and let $\gamma_{\Hf}: \cL_{\yHf}\to  \cL_{\yHpf}$ also represents the parallel transport along $\gamma_{\Hf}$. Then 
one directly checks that $\nu_\varrho^\cL$ and $\nu_{\varrho'}^\cL$ are related by the following commutative diagram
 \begin{equation*}
 \begin{tikzcd}[column sep=6em]
 \cL_{\yHf}\ar[d, "\gamma_{\Hf}"', "\sim"]\ar[r, "\nu_\varrho^\cL", "\sim"']&\cL_{\yLf}\ar[d, "\gamma_{\Lf}", "\sim"']\\
  \cL_{\yHpf}\ar[r, "\nu_{\varrho'}^\cL","\sim"']&\cL_{\yLpf}
 \end{tikzcd}. 
 \end{equation*}
\end{remark}

 Let 
\begin{align}\label{eqpsiZHfK}
\psi=\wt{\psi}|_{\pi_1(\ZHf)}: \pi_1(\ZHf)= \XX_*(Z(\Hf))\lrar K
\end{align}
be the homomorphism that classifies the $K$-covering of the quotient map $Z(H)\to Z(\Hf)$, and let $\psi_\chi=\chi\circ \psi=\wt{\psi}_\chi|_{\pi_1(Z(\Hf))}$. 
The above implies that for any $\gafOT\in \XX_*(\ZHf)$, if we identify $\gafOT$ with $\iota^\flat(\gafOT)$, then 
\begin{align*}
\btau_\chi(\gafOT)(\gafOT)^{-1}=\psi_\chi(\gafOT). 
\end{align*}

Similarly, let 
\begin{align}\label{eqpsiZHfKpsit}
\psi^t=\wt{\psi}^t|_{ \pi_1(\ZHf)^{op}}: \pi_1(\ZHf)^{op}= \XX_*(Z(\Hf))^{op}\lrar K. 
\end{align}
and let $\psi^t_\chi=\chi\circ\psi^t=\wt{\psi}^t_\chi|_{\XX^*(\ZHf)^{op}}$. Then for any $\gafTO\in \XX_*(\ZHf)^{op}$, 
\begin{align*}
\btau^t_\chi(\gafTO)(\gafTO)^{-1}=\psi^t_\chi(\gafTO). 
\end{align*}

\subsection{The $\pi\ZHc^{op}$-equivariant functor $\wPhiHP :\frPH^{op}\to \mathrm{Maps}_{\frG^{\Kdu}}(\IWf, \Qft)$}\label{subsecwtPhiHPop}

\sss{The functors $\QPf$ and $\wQPf$}

\begin{lemma}\label{lemmawtcQPfG}
We have a natural commutative diagram (\emph{non}-Cartesian in general) of $\infty$-groupoids
\begin{equation*}
\begin{tikzcd}[column sep=6em]
\frPH^{op}\ar[r, "\wQPf"]\ar[d]&\frGKt_{/\Qft}\ar[d]\\
\frPHf^{op}\ar[r, "\QPf"']&\frGKt
\end{tikzcd},
\end{equation*}
where $\QPf$ (resp. $\wQPf$) is the natural functor that sends $\yf$ (resp. $y$) to $\cQ_{\yf}$ (resp. $\cQ_{\yf}$ with the framing $\eta^y$ of $\ell_{\yf}$). 

\end{lemma}

\begin{proof}
The diagram comes from the natural commutative diagram 
\begin{equation}\label{diagQPf}
\begin{tikzcd}[column sep=8em]
\frPH^{op}\ar[rr, bend left=2.3em, "\wQPf"]\ar[r, "{\wt{\phi}_{\cL_{K}}|_{\frPH^{op}}}"]\ar[d]&\wt{\Xi}(\Kdu)\ar[r, "{\wt{\cQ}^{\flat,t}_\Xi}"] \ar[d]&\frG^{\Kdu}_{/\Qft}\ar[d]\\
\frPHf^{op}\ar[rr, bend right=2em, "\QPf"]\ar[r, "{\phi_{\cL_{K}}|_{\frPHf^{op}}}"]&\Xi(\Kdu)\ar[r, "{\cQ^{\flat,t}_\Xi}"] &\frG^{\Kdu}
\end{tikzcd},
\end{equation}
where the functors in the middle of the diagram are the version of \eqref{eqphiLKpiTHfXiKst} and \eqref{eqcQXiwt} for $\AHf$ with respect to $\btau^t=(\btau^t_\chi)_{\chi\in K^*}$.  
\end{proof}

\begin{lemma}\label{lemmaphiPHinvariant}
The functor $\phi_{\cL_K}|_{\frPHf}$ (resp. $\wt{\phi}_{\cL_K}|_{\frPH}$) is naturally $\pi \ZHc$-invariant  with respect to the natural  $\pi \ZHc$-action on $\frPHf$ (resp. $\pi \ZHc$-action on $\frPH$). 
\end{lemma}

\begin{proof}
This follows from that the local systems $\cL_K=(\pi_{\cT, H})_*\CC_{\cT_H}$ and $(\pi_{\cT, H})^*\cL_K\cong \cO(K)_{\cT_H}$, equipped with their natural algebra structures,  are canonically $\pi \ZHc$-equivariant. 
\end{proof}

\begin{remark}
Note that $\phi_{\cL_K}|_{\frPHf}$ is \emph{not} $\pi \ZHfc$-invariant unless the covering $Z(H)\to Z(H^\flat)$ is a trivial covering (with covering group $K$).  
\end{remark}

Using diagram \eqref{diagQPf}, Lemma \ref{lemmaphiPHinvariant} immediately implies the following. 
\begin{cor}\label{corQPfinvariant}
The functors $\QPf$ and $\wQPf$ \eqref{diagQPf} are both naturally $\pi \ZHc^{op}$-invariant.
\end{cor}

\sss{The groupoid $\WPHf$}\label{sssWPHf}

Let $\wt{\cW}_{\frP, \Hf}$ be the full subcategory of $\WTHf$ consisting of $S_{\yf}, \yf\in \UHfSc$. 
Let $\cW_{\frP, \Hf}$ be the $\infty$-groupoid $\wt{\cW}_{\frP, \Hf}^{\Gpd}$, which is equivalent to an ordinary groupoid. We have the natural $\pi\ZHc$-equivariant functor 
\begin{align}\label{eqnuHf}
\nu_{\Hf}^{\Gpd}: \frPHf\lrar \cW_{\frP, \Hf}
\end{align}
induced from $\nu_{\Hf}$. 

Corollary \ref{cornuHflat} means the embedding 
\begin{align*}
\Hom_{\frPHf}(\yfO, \yfT)\hookrightarrow  \Hom_{\WPHf}(S_{\yfO}, S_{\yfT})
\end{align*}
gives a canonical ``basis" of the latter, in the sense that  
\begin{align*}
\Hom_{\WPHf}(S_{\yfO}, S_{\yfT})=\CC^\times\times\Hom_{\frPHf}(\yfO, \yfT).
\end{align*} 
Let $q_{\frP, \Hf}: \cW_{\frP, \Hf}\lrar \frP_{\Hf}$ be the natural functor of groupoids, which takes $S_{\yHf}$ to $\yHf$ and takes any isomorphism $\lambda\cdot \gamma: S_{y_H}\ovs{\sim}{\to} S_{y_H'}, \lambda\in \CC^\times$ to $\gamma$. In other words, $q_{\frP, \Hf}$ is the quotient on the morphism sets by the $\CC^\times$-action. 

\sss{The natural functor $\wUpPf: \frPHf^{op}\to \frG^{K^*}_{\IWf/}$}

Recall $\Upyf: \Ind\WTHf\to \Mod(A_{\Hf}), \yf\in \frPHf$, from \eqref{eqUpzJflmuShH}, where $A_{\Hf}$ as in \eqref{eqAHf}. 
We can assemble these into a canonical $\pi \ZHfc^{op}$-equivariant functor 
 \begin{align*}
\UpPf:  \WPHf^{op}&\lrar \Fun_{{\Pr}_\st^L}(\Ind\WTHf, \Mod(A_{\Hf})),\\
\yf&\mapsto \Upyf,
 \end{align*}
where for any isomorphism $\varphi^\flat\in \Hom_{\WPHf}(S_{\yfT}, S_{\yfO})$, we have the invertible natural transformation 
\begin{equation}\label{diagramtauvarphiUpsilon}
\begin{tikzcd}[column sep=8em]
\Ind\WTHf \arrow[r, bend left=20, ""{name=U, below}, "{\UpyfO}"]
\arrow[r, bend right=20, ""{name=D}, "{\UpyfT}"']
& \Mod(A_{\Hf})
\arrow[Rightarrow, from=U, to=D, "\tau'_{\varphi^\flat}", "\sim"'].
\end{tikzcd},
\end{equation}
coming from pre-composing with $\varphi^\flat$: 
 \begin{align*}
 A_{\Hf}\car\mu_{\yfO}(\cG)\overset{\varphi^\flat}{\underset{\sim}{\lrar}} A_{\Hf}\car\mu_{\yfT}(\cG). 
 \end{align*}
 Here the $\pi \ZHfc^{op}$-equivariance is automatic from the $\pi \ZHfc$-action on $\Ind\WTHf$. 
  In particular, $\UpPf$ induces a canonical $\pi \ZHfc^{op}$-equivariant functor
 \begin{align*}
\UpPf': \WPHf^{op}&\lrar ({\Pr}_\st^L)_{ \Ind\WTHf/}.
 \end{align*}
 
Now pre-composing $\UpPf'$ with $\nu_{\Hf}^{\Gpd, op}$ \eqref{eqnuHf} and using (the opposite version of) the $\pi \ZHc$-invariant functor $\phi_{\cL_K}|_{\frPHf}: \frPHf\to \Xi(K^*)$ from Lemma \ref{lemmaphiPHinvariant} (note that the latter is \emph{not} $\pi \ZHfc$-invariant in general), we get the natural $\pi \ZHc$-equivariant functor 
 \begin{align}\label{eqwUpPf}
\wUpPf: \frPHf^{op}&\lrar \frG^{\Kdu}_{\IWf/},\\
\nonumber\yf&\mapsto \wUpyf,
 \end{align}
which assembles Corollary \ref{corUpzJflatKstarequiv} for $\yf\in\frPHf$.

\sss{Definition of  $\wPhiHP$}\label{sssDefwPhiHP}

Using  \eqref{eqwUpPf} and Corollary \ref{corQPfinvariant}, we have a natural $\pi\ZHc^{op}$-equivariant functor from taking the induced functor on fiber products
\begin{equation*}
\begin{tikzcd}
\wUpPf\underset{\QPf}{\times} \wQPf:\quad \frPHf^{op}\underset{\frPHf^{op}}{\times}\frPH^{op}\ar[r]&\frG^{\Kdu}_{\IWf/}\underset{\frG^{\Kdu}}{\times} \frG^{\Kdu}_{/\cQ^t_\flat}, 
\end{tikzcd}
\end{equation*}
where the fiber product on the right-hand-side are defined from the two natural evaluation functors. This defines  
the sought-for $\pi\ZHc^{op}$-equivariant functor by canonically identifying the source and target above with the respective source and target in the following: 
\begin{align}\label{wPhiPH}
\wPhiHP:\frPH^{op}\to \Maps_{\frG^{\Kdu}}(\IWf, \Qft). 
\end{align}

\sss{Definition of $\PhiHP$}\label{sssdefPhiHP}

Let $\Qf^{t, \Kdu}:=\Mod(\AHf)^{\Kdu}$ and let $\IWf^{\Kdu}:=\Ind\WTHf^{\Kdu}$ be the associated categories of $\Kdu$-invariants in $\Pr_\st^L$. By Corollary \ref{corWTHWTHfKstar}, we have the canonical equivalence $\IWf^{\Kdu}\simeq \Ind\WTH$. We also have the canonical equivalence $\Qf^{t, \Kdu}\simeq \QCoh(\SHv)$, induced from the involution $\iota_{\Hf}^t: \AHf\ovs{\sim}{\to} \AHf^{op}$. Then $\wPhiHP$ \eqref{wPhiPH} naturally induces the $\pi\ZHc^{op}$-equivariant functor 
\begin{align}\label{PhiPH}
\PhiHP:\frPH^{op}\to \Maps_{({\Pr}_\st^L)^{\Gpd}}( \Ind\WTH,  \QCoh(\SHv)). 
\end{align}

\sss{A more concrete description of $\wPhiHP$}
 By Corollary \ref{corfiberfrgammaJf}, for each object $y\in \frPH$, we have a canonical isomorphism 
\begin{align}\label{eqwUpy}
\wUpy: \IWf\ovs{\sim}{\lrar}\Qft.
\end{align} 
in $\frG^{\Kdu}$.
 Then $\wUpy$ induces a canonical equivalence
\begin{align}\label{eqwUpyKstar}
\wUpy^{\Kdu}: \IWf^{\Kdu}\ovs{\sim}{\lrar}\Qf^{t, \Kdu}
\end{align} 
in $\Pr^L_{\st}$.

For any $y_1, y_2\in \frP_H$ and any isomorphism $\rho_{2,1}\in \Hom_{\frP_H}(y_2, y_1)$ \eqref{eqrho12sigmay},  
let $\rho_{2,1}^{H^\flat}:=F_*(\rho_{2,1}^{H})=(\rho_{2,1}^\flat)^{H^\flat}$. They determine a canonical diagram in  $\frG^{\Kdu}$,
\begin{equation}\label{diagmuSQ}
\begin{tikzcd}[row sep=0.5em, column sep=15em]
&\cQ_{\yfO}\ar[dd, "\sim"', "\wt{\Psi}_{1,2}=\wt{\Psi}_{\yfO, \yfT}"]\arrow[dr, bend left=15, "\eta^{y_1}", ""{name=U, below}]&\   \\
\IWf\ar[ur, "\wUpyfO", "\sim"']\arrow[dr,  "\wUpyfT"',"\sim"]& &\Qft\\ 
&\cQ_{\yfT}\arrow[ur, bend right=15, "\eta^{y_2}"', ""{name=D, near start}]&\   
\arrow[Rightarrow, from=U, to=D, "\sim"', "{\tau_{\rho_{2,1}}}"]
\end{tikzcd}
\end{equation}
 in which
\begin{itemize}
\item[(i)] $\wt{\Psi}_{1,2}=\wt{\Psi}_{\yfO, \yfT}:=\wUpyfT\circ\wUpyfO^{-1}$ is the canonical equivalence that makes the left triangle strictly commutative (up to a contractible space of choices);
\item[(ii)]  the invertible natural transformation marked with $\tau_{\rho_{2,1}}$, 
is induced from 
\begin{align}\label{eqrho12flatparalltransport}
A_{\Hf}\car\mu_{\yfO}(\cG)\overset{\rho_{2,1}^{\Hf}}{\underset{\sim}{\longrightarrow}} A_{\Hf}\car\mu_{\yfT}(\cG),\quad(\ell_{\yfO},\ \eta^{y_1})\overset{\rho_{2,1}^\flat}{\underset{\sim}{\longleftarrow}}(\ell_{\yfT},\ \eta^{y_2}). 
\end{align}
where the first one is the natural transformation given by pre-composition with $\rho_{2,1}^{\Hf}$ 
\begin{align}\label{eqroh21Hf}
\rho_{2,1}^{\Hf}: \Hom_{\WTHf}(S_{\yfO},-)\ovs{\sim}{\Longrightarrow} \Hom_{\WTHf}(S_{\yfT},-)
\end{align}
and the second one in \eqref{eqrho12flatparalltransport} is the parallel transport of $\ell_{\yfT}$ to $\ell_{\yfO}$ along the path $\rho_{2,1}^\flat$ that identifies the framing $\eta^{y_2}$ with $\eta^{y_1}$. 

\end{itemize}  
In particular, the above determines a canonical diagram
\begin{equation}\label{diagmuSKcQfrG}
\begin{tikzcd}[column sep=8em]
\IWf \arrow[r, bend left=20, ""{name=U, below}, "{\wUpyO}"]
\arrow[r, bend right=20, ""{name=D}, "{\wUpyT}"']
& \Qft
\arrow[Rightarrow, from=U, to=D, "\tau_{\rho_{2,1}}", "\sim"'].
\end{tikzcd}
\end{equation}
Consequently, we get a canonical diagram in $(\Pr^L_\st)^\Gpd$ by taking $\Kdu$-invariants
\begin{equation}\label{diagmuSKcQKsttau}
\begin{tikzcd}[column sep=8em]
\Ind\WTH\simeq \IWf^{\Kdu}\hspace{1em} \arrow[r, bend left=15, ""{name=U, below}, "{\wUpyO^{\Kdu}}"]
\arrow[r, bend right=15, ""{name=D}, "{\wUpyT^{\Kdu}}"']
&\hspace{2em}\Qf^{t, \Kdu}\simeq \QCoh(\SHv)^\vee
\arrow[Rightarrow, from=U, to=D, "{\tau^{\Kdu}_{\rho_{2,1}}}", "\sim"']. 
\end{tikzcd}
\end{equation}

It is clear that for any morphisms $\rho_{i+1,i}$ from $y_{i+1}$ to $y_{i}$, $i=1,2$, we have 
the canonical isomorphism
\begin{align*}
\tau_{\rho_{2,1}}\circ \tau_{\rho_{3,2}}\cong \tau_{\rho_{3,1}}, 
\end{align*}
where $\rho_{3,1}$ is the concatenation $\rho_{3,2}\bullet \rho_{2,1}$,  
and they satisfy the natural (strict) associativity relations, inherited from the  (strict) associativity relations of compositions of  \eqref{eqrho12flatparalltransport}.

\begin{remark}
\begin{itemize}
\item[(i)] The above also completely defines $\wPhiHP$, but the more ``abstract" definition of $\wPhiHP$ in  \S\ref{sssDefwPhiHP} incorporates the canonical $\pi \ZHc^{op}$-equivariance automatically.

\item[(ii)] The left half of the diagram \eqref{diagmuSQ} is different from \eqref{eqwUpPf}. 
In \eqref{diagmuSQ}, the natural transformation $\tau_{\rho_{2,1}}$ combines the natural transformation $\tau'_{\rho_{2,1}^{\flat}}$ from \eqref{diagramtauvarphiUpsilon} and the natural morphism attached to $\rho_{2,1}$ through $\wQPf$. 

\item[(iii)]
In terms of microlocal sheaves, the natural transformation $\rho_{2,1}^{\Hf}$ \eqref{eqroh21Hf} is just parallel transport of microlocal stalks along $\rho_{2,1}^\flat$ (with respect to some natural ``brane structure" on $\UHfS$). 
\end{itemize}
\end{remark}

We will show that $\PhiHP$ is independent of finite central quotients below (Corollary \ref{corPhiHySy}).

\subsection{Calculation of the induced $\pi \ZH^{op}$-action on $\Qft$ and $\mathrm{QCoh}(\SHv)$}\label{subsecpiZHQft}

Now the $\pi\ZH^{op}$-equivariant functor $\wPhiHP$ \eqref{wPhiPH} gives a canonical isomorphism in $\frG^{\Kdu}$
\begin{align}\label{eqycGcLmueta}
 \frPH^{op}\times^{\pi \ZH} \IWf&\ovs{\sim}{\lrar} \Qft,
\end{align}
where $\pi \ZH$ acts on $\frPH^{op}$ by the canonical involution $\pi Z(H)\ovs{\sim}{\to}\pi Z(H)^{op}$ (as groupoids), and 
the left-hand-side is \emph{non-canonically} $\pi \ZH^{op}$-equivariantly isomorphic to $\IWf$ (again for the latter $\pi \ZH^{op}$ acts through  $\pi Z(H)\ovs{\sim}{\to}\pi Z(H)^{op}$). 
This isomorphism is explicitly given by the following diagram
 \begin{equation*}
 \begin{tikzcd}[row sep=0.5em]
 \frPH^{op}\times^{\pi \ZH} \Ind\WTHf\times \Kdu\ar[r, "\sim"]\ar[dd]&\Mod(\AHf)\times \Kdu\ar[dd, "\emph{\tsfa}"]\ar[ddl, Rightarrow, "\sim"']\\
\  &\ \\
 \frPH^{op}\times^{\pi \ZH} \Ind\WTHf\ar[r, "\sim"]&\Mod(\AHf)\\
 (y, \cG)\ar[r, mapsto]&A_{\Hf}\car \mu_{\yf}(\cG)
 \end{tikzcd}
 \end{equation*}
 in which the vertical functors are the respective action maps of $\Kdu$ and the  invertible natural transformation is encoded by the trivialization $\ell_{\yf,\chi}\ovs{\eta^y}{\underset{\sim}{\to}} \ell_{\triv,\chi}$ for each $(y, \cG, \chi)$.

The remaining (left) $\pi\ZH^{op}$-action on the left-hand-side induces a natural $\pi\ZH^{op}$-action on $\Qft$. 
Of course, we can switch from the $\pi\ZH^{op}$-action on $\Qft$ to the $\pi\ZH$-action on $\Qft$ by the canonical isomorphism $\pi\ZH\cong \pi\ZH^{op}$. We will keep using the $\pi\ZH^{op}$-action on $\Qft$ till Corollary \ref{corpiZHQfKstar}, since (1) from the identification above it is more natural to consider the $\pi\ZH^{op}$-action on $\Qft$; (2) since $\AHf\subset \CC[\pi_1 (\Tf)^{op}]$ and consequently $\AHf^\times =\CC[\pi_1 Z(\Hf)^{op}]^{\times}$, sticking to the $\pi\ZH^{op}$-action will make the computations in Proposition \ref{eqbiakappa} and Corollary \ref{corpiZHQfKstar} more consistent. We will finally get to the $\pi\ZH$-action on $\QCoh(\SHfv)^{K^*}\simeq \QCoh(\SHv)$ in Corollary \ref{corpiZHQfKstarfinal}, by taking the opposite version of the aforementioned $\pi\ZH^{op}$-actions.

We have the complete analog of Lemma \ref{lemmawtpsit}, where $T\to \Tf$ is replaced by $\ZH\to \ZHf$. 
We can calculate the $\pi\ZH^{op}$-action on $\Qft$ explicitly using the equivalent model $ \pi \ZH^{op}\cong K/(\pi_1Z(\Hf))^{op}$ acts on 
\begin{align*}
\pi_{\Lambda, H}^{-1}(\yf)/\left(\pi_1Z(\Hf)\right)^{op}\times^{K/\pi_1Z(\Hf)}\IWf\simeq \frPH^{op}\times^{\pi \ZH} \IWf,
\end{align*}
for any fixed $\yf\in \frPHf^{op}$.

\begin{prop}\label{eqbiakappa}
The natural $K/(\XX_*(\ZHf))^{op}\cong \pi \ZH^{op}$-action on $\Qft$ induced from $\wPhiHP$ \eqref{eqycGcLmueta} acts by 
\begin{itemize}
\item[(1)] for $\kappa\in K$, it does the identity on $\Mod(\AHf)$, and induces the automorphism of the standard $\Kdu$-action on $\Mod(\AHf)$ through $\kappa^{-1}\in \Aut_{\Xi(\Kdu)}(\ell_\triv)= K$, which gives $\bia_{\kappa}$ in the diagram below: 
\begin{equation}\label{diagrambiakappa}
\begin{tikzcd}
\Kdu\times \Mod(\AHf)\ar[d, equal]\ar[r, "\tsfa"]&\Mod(\AHf)\ar[d, equal]\ar[dl, Rightarrow, "\sim", "\bia_{\kappa}"']\\
\Kdu\times \Mod(\AHf)\ar[r, "\tsfa"]&\Mod(\AHf)
\end{tikzcd}. 
\end{equation}

\item[(2)] for any morphism $\gamma_{2,1}\in \Hom_{\pi Z(H)}(\kappa_2, \kappa_1)=\Hom_{\pi Z(H)^{op}}(\kappa_1, \kappa_2)$, it induces  the invertible morphism $\tau^{\bia}_{\gamma_{2,1}}: \bia_{\kappa_1}\Rightarrow \bia_{\kappa_2}$ exhibited in the following diagram 
\begin{equation}\label{diagrammSQzstdtau}
\begin{tikzcd}[column sep=10em, row sep=4em,
execute at end picture={
\foreach \Nombre in  {A,B,...,F}
  {\coordinate (\Nombre) at (\Nombre.center);}
\fill[cof,opacity=0.3] 
 (A) -- (B) -- (C) -- cycle;
\fill[greeo,opacity=0.3] 
 (D) -- (E) -- (F) -- cycle;
\fill[pur,opacity=0.3] 
(A) -- (C) -- (F) -- (D) -- cycle;
}]
|[alias=A]| \Kdu\times \Mod(\AHf)\ar[r, "\tsfa"]\ar[d, equal]&|[alias=D]|\Mod(\AHf)\ar[d, equal, color=gray]\ar[dl, Rightarrow, color=black!60, "\bia_{\kappa_1}"', "\sim"] & \\
|[alias=B]| \Kdu\times \Mod(\AHf)
\ar[r, near end, color=black!70, "\tsfa"]\ar[dr, equal]&|[alias=E]| \Mod(\AHf) \ar[dr, equal, color=black!70]& \\
&|[alias=C]|\Kdu\times  \Mod(\AHf)\ar[r, "\tsfa"']\arrow[from=uul,  equal , near start, crossing over, ""{name=U, below}]&|[alias=F]| \Mod(\AHf)\arrow[from=uul, equal, near start, crossing over,  ""{name=X, below}] 
\arrow[from=ull, to=U, Rightarrow,  "{id\times m_{\gamma_{2,1}^{\flat}}}"',near start, "\sim"]
\arrow[from=ul, to=X, Rightarrow, color=black!80,  "{m_{\gamma_{2,1}^{\flat}}}"', near start, "\sim"]
\arrow[from=D, to=C, bend right=6em, crossing over, Rightarrow, "\bia_{\kappa_2}", "\sim"', near start]
\end{tikzcd}
\end{equation}
in which 
\begin{itemize}
\item[(2a)] the filling of the right (with arrows in gray color) and the left triangular faces are the automorphisms of $id_{\Mod(\AHf)}$ given by multiplication by $\gamma_{2,1}^{\flat}$ (which we identify with $\iota^\flat(\gamma_{2,1}^{\flat})$ cf. Corollary \ref{cornuHflat});
\item[(2b)] the filling of the back (with arrows in gray color) and top rectangular faces are given by $\bia_{\kappa_1}$ and $\bia_{\kappa_2}$, respectively; the filling of the bottom rectangular face is the trivial one;
\item[(2c)] the filling of the solid is given by the canonical isomorphism
\begin{align}\label{eqpropmak2c}
\id_{\tsfa}\circ (id\times m_{\gamma_{2,1}^{\flat}})\circ\bia_{\kappa_1}\cong \bia_{\k_2}\circ m_{\gamma_{2,1}^{\flat}}\circ \id_{\tsfa}\in \Hom_{\Fun(\Kdu\times \Mod(\AHf),\Mod(\AHf))}(\tsfa, \tsfa). 
\end{align}
\end{itemize}
\end{itemize}
\end{prop}

\begin{proof}
For any $\kappa\in K$, since it fixes $\yf$ (and does trivially on $\IWf$ by definition), it is direct from \eqref{eqycGcLmueta} that it only changes the framing $\eta^y$ to $\eta^{\kappa\cdot y}$. This action is exactly $\bia_{\kappa}$ in \eqref{diagrambiakappa}. 

Using the calculation of $\nu_{\Hf}$ Corollary \ref{cornuHflat}, we directly see that $\tau^{\bia}_{\gamma_{2,1}}: \bia_{\kappa_1}\Rightarrow \bia_{\kappa_2}$ is given by the data of 
\begin{align}\label{eqMgafTOAHftimes}
&M\ovs{\gamma_{2,1}^{\flat}}{\underset{\sim}{\lrar}} M,\quad  \forall M\in \Mod(A_{\Hf}), \text{ where }\gamma_{2,1}^\flat\in A_{\Hf}^\times=\cO(H^{\flat, \vee, \ab})^{op, \times}
\end{align}
which corresponds to (2a), and 
\begin{align*}
&(\ell_{\yf}, \eta^{\k_1\cdot y})\ovs{\gafTO}{\underset{\sim}{\longleftarrow}} (\ell_{\yf}, \eta^{\k_2\cdot y})\text{ by parallel transport}, 
\end{align*}
which induces (2b) (without using $\gafTO$). 
Let us finally check (2c). 
This amounts to checking the following commutative diagram for each $\chi\in \Kdu$ and $M\in \Mod(A_{\Hf})$.
The subscript under $\bullet$ indicates the $\AHf$-module structure on $M^\chi$ that we are using (either $M^\chi$ or the original $M$ without the twisting). Here the upper right corner (resp. the lower left corner) represents the left-hand-side (resp. right-hand-side) of \eqref{eqpropmak2c} on $(\chi, M)\in \Kdu\times \Mod(\AHf)$: 
\begin{equation*}
\begin{tikzcd}[column sep=10em]
\tsfa(\chi, M)=M^\chi\otimes\ell_{\triv, \chi}\ar[r, "1\otimes \chi(\k_1^{-1})"]\ar[d, "\left(\gafTO\underset{M^\chi}{\bullet} \right)\otimes 1=\left(\btau^t_{\chi}(\gafTO)\underset{M}{\bullet} \right)\otimes 1"']&M^\chi\otimes\ell_{\triv, \chi}\ar[d, "\left(\gafTO\underset{M}{\bullet}\right)\otimes 1=\left(\btau^t_{\chi^{-1}}(\gafTO)\underset{M^\chi}{\bullet} \right)\otimes 1"]\\
M^\chi\otimes\ell_{\triv, \chi}\ar[r, "1\otimes \chi(\k_2^{-1})"]&M^\chi\otimes\ell_{\triv, \chi}. 
\end{tikzcd}
\end{equation*}
Note the essential difference between $id_\tsfa\circ(\id\times m_{\gafTO})$ and $m_{\gafTO}\circ id_\tsfa$: for a given $M\in \Mod(\AHf)$, the former sends it to the right vertical morphism in the above diagram, for each $\chi\in K^*$, while the latter  sends it to the left vertical morphism.

This follows directly from the identity (recall $\btau^t_\chi$ from \eqref{eqbtauchit})
\begin{align}
\nonumber &\chi(\k_1^{-1}) \gafTO=\btau^t_\chi(\gafTO)\chi(\k_2^{-1})\\
\label{eqbtauchitgaTOf}\Longleftrightarrow\ &\btau^t_\chi(\gafTO)(\gafTO)^{-1}=\psi^t_\chi(\gafTO)=\chi(\k_1^{-1}\k_2).
\end{align} 
\end{proof}

\begin{remark}
Now viewing $\gafOT=(\iotHf)^{-1}(\gafTO)\in \Hom_{\pi Z(H)}(\k_1, \k_2)$, then \eqref{eqbtauchitgaTOf} translates to 
\begin{align*}
\btau_{\chi}(\gafOT)(\gafOT)^{-1}=\psi_{\chi}(\gafOT)=\chi(\k_1^{-1}\k_2). 
\end{align*}
\end{remark}

\begin{cor}\label{corpiZHQfKstar}
The natural $K/(\XX_*(\ZHf))^{op}\cong \pi \ZH^{op}$-action on $\cQ_\flat^{t, \Kdu}$ induced from $\PhiHP$ \eqref{PhiPH} acts by 
\begin{itemize}
\item[(1)] for any $\k \in K$, it gives the automorphism $\tau^{\k}=\otimes_{\cO}\cO^{\k}$ of $\cQ_\flat^{t, \Kdu}$ that twists each equivariant structure by $\k\in K=\XX^*(\Kdu)$;
\item[(2)] for any morphism $\gamma_{2,1}\in \Hom_{\pi Z(H)}(\kappa_2, \kappa_1)=\Hom_{\pi Z(H)^{op}}(\kappa_1, \kappa_2)$, it gives the natural transformation $\tau^{\k_1}\Rightarrow \tau^{\k_2}$, corresponding to 
\begin{align*}
\Hom_{\Qf^{t, \Kdu}}(\cO^{\k_1}, \cO^{\k_2})=\Hom_{\Mod(\AHf)}(\cO, \cO)^{\k_1^{-1}\k_2}=\cO(\Hfv\sslash \Hfv)^{op, \k_1^{-1}\k_2}\ni  \gaTO^{\flat},
\end{align*}
where $\Hom_{\Mod(\AHf)}(\cO, \cO)^{\k_1^{-1}\k_2}$ is the component on which $\Kdu$ acts by the character $\k_1^{-1}\k_2\in K$. 
\end{itemize}
\end{cor}

\begin{proof} 
For any $\k\in K$, let $\wt{\tau}^{\k}:  \Qf^{t,\Kdu}\to \Qf^{t,\Kdu}$ be the functor induced by $\bia_{\k}$ in Proposition \ref{eqbiakappa}. For any $M\in \Qf^{t,\Kdu}$ with the structure isomorphisms $M^\chi\ovs{F^\chi}{\to}M$, $\wt{\tau}^{\k} M$ has the same underlying $A_{\Hf}$-module $M$, but equipped with the structure isomorphisms determined by the following commutative diagram
 \begin{equation*}
\begin{tikzcd}[column sep=10em]
\wt{\tau}^{\k}(M^\chi)=M^\chi\ar[r, "F^\chi"]\ar[d, "\chi(\k^{-1})"]&M=\wt{\tau}^\k M\ar[d, equal]\\
(\wt{\tau}^{\k} M)^\chi=M^\chi\ar[r, "\chi(\k)\cdot F^\chi"]&M=\wt{\tau}^\k M.
\end{tikzcd}
\end{equation*}
This shows that we have the canonical isomorphism of functors $\wt{\tau}^{\k}\cong \tau^{\k}$. 

It is direct to trace through the diagram that for any $M$ as above and $\gaTO$, the morphism $\tau^{\k_1}(M)\to  \tau^{\k_2}(M)$ induced from the natural transformation is given by 
\begin{equation*}
\begin{tikzcd}[column sep=10em]
(\tau^{\k_1} M)^\chi=M^\chi\ar[r, "\chi(\k_1)\cdot F^\chi"]\ar[d, "\gafTO\uds{M}{\bullet}=\btau_{\chi^{-1}}^t(\gaTO^\flat)\uds{M^\chi}{\bullet}"']&M\ar[d,"\gafTO"]\\
(\tau^{\k_2} M)^\chi=M^\chi\ar[r, "\chi(\k_2)\cdot F^\chi"]&M. 
\end{tikzcd}
\end{equation*}
The diagram indeed commutes since  
\begin{align*}
\chi(\k_2)\btau^t_{\chi^{-1}}(\gaTO^\flat)=\gaTO^\flat\chi(\k_1)\ \Longleftrightarrow\ \eqref{eqbtauchitgaTOf}. 
\end{align*}
\end{proof}

Lastly, we take the opposite version of Corollary \ref{corpiZHQfKstar} using the involution $\iotHf$, and we get 
\begin{cor}\label{corpiZHQfKstarfinal}
The natural $K/\XX_*(\ZHf)\cong \pi \ZH$-action on $\cQ_\flat^{K^*}\simeq \QCoh(\SHv)$ induced from $\PhiHP$ \eqref{PhiPH} acts by 
\begin{itemize}
\item[(1)] for any $\kappa\in K$, it gives the automorphism $\tau^\kappa=\otimes_{\cO}\cO^{\kappa}$ of $\QCoh(\SHfv)^{K^*}$ that twists each equivariant structure by $\kappa\in K$;
\item[(2)] for any $\gaOT\in \Hom_{\pi \ZH}(\kappa_1, \kappa_2)$, it gives the natural transformation $\tau^{\k_1}\Rightarrow \tau^{\k_2}$, corresponding to 
\begin{align*}
\Hom_{\QCoh(\SHfv)^{K^*}}(\cO^{\kappa_1}, \cO^{\kappa_2})=\Hom_{\QCoh(\SHfv)}(\cO, \cO)^{\kappa_1^{-1}\kappa_2}=\cO(\Hfv\sslash \Hfv)^{\kappa_1^{-1}\kappa_2}\ni \gaOT^{\flat},
\end{align*}
where $\Hom_{\QCoh(\SHfv)}(\cO, \cO)^{\kappa_1^{-1}\kappa_2}$ is the component on which $K^*$ acts by the character $\kappa_1^{-1}\kappa_2\in K=\XX^*(K^*)$. 
\end{itemize}
\end{cor}

We will call the $\pi Z(H)$-action on $\QCoh(\cS_{H^\vee})$ from Corollary \ref{corpiZHQfKstarfinal} as the \emph{canonical} $\pi Z(H)$-action on $\QCoh(\cS_{H^\vee})$.

\subsection{Canonical mirror equivalence}\label{subsecCanonicalME}

For any $y_H\in \frP_{H}$ and $y_L\in \frP_{L}$, let $\yHf\in \frP_{\Hf}$ and $\yLf\in \frP_{\Lf}$ denote their respective images. 
Given an isomorphism $\varrho: S_{\yLf}\ovs{\sim}{\to}\res^{\Hf}_{\Lf}S_{\yHf}$ in $\cW(\cT_{L^\flat})$, we have a similar diagram as \eqref{diagUpsilonzJflatyJ} in $\frG^{K^*}$ (cf. Remark \ref{remgeneralizetoresHfLf}), 
\begin{equation}\label{eqresHfLfpfstar}
\begin{tikzcd}[column sep=6em]
\cI\cW_{\flat, H}\ar[r, "\res_{\Lf}^{\Hf}"]\ar[d, "\wt{\Upsilon}_{\yHf}"',"\sim"]&\cI\cW_{\flat, L}\ar[d,  "\wt{\Upsilon}_{\yLf}", "\sim"']\ar[dl, Rightarrow, "\sim"']\\
\cQ_{\yHf}\ar[r, "p_\flat^*"]&\cQ_{\yLf}
\end{tikzcd}, 
\end{equation} 
where (1) the subscripts $H, L$ in the categories mean $\cI\cW_\flat$ associated with $H, L$, respectively; (2) $p_\flat^*$ is $\otimes_{\AHf} \ALf$ which respects the $K^*$-actions on both sides through a similar diagram as \eqref{diagalphachiJvn} for $\cL=\cL_\chi, \chi\in K^*$ and $A^\vn$ replaced by $\ALf$;  (3) the natural isomorphism canonically depends on $\varrho$.  

\sss{}\label{psichivarrhoeta}
Let $\psi^\chi_\varrho: \ell_{\chi, \yHf}\ovs{\sim}{\to}\ell_{\chi, \yLf}$ denote $\nu_\varrho^{\cL_\chi}$ \eqref{eqnuupcL}. 
For any $y_H, y_L$ over $\yHf, \yLf$ respectively, the framings $\eta^{\yH}$ and $\eta^{\yL}$ identify $(\psi^\chi_\varrho)_\chi$ with $(\psi_{\varrho;\eta}^\chi)_\chi\in K=\XX^*(K^*)$. 

We say a pair $(\yH, \yL)$ is \emph{$(K, \varrho)$-compatible} if $(\psi_{\varrho;\eta}^\chi)_\chi=1\in K$. For any $(K,\varrho)$-compatible pair $(\yH, \yL)$, 
diagram \eqref{eqresHfLfpfstar} induces the following diagram by taking the framings $\eta^{y_H}, \eta^{y_L}$ and taking $K^*$-invariants 
\begin{equation}\label{eqresHLRyHyL}
\begin{tikzcd}[column sep=6em]
\Ind\cW(\cT_{H})\ar[r, "\res^H_L"]\ar[d,"\wt{\Upsilon}_{y_H}^{K^*}"']&\Ind\cW(\cT_{L})\ar[d, "\wt{\Upsilon}_{y_L}^{K^*}"]\ar[dl, Rightarrow, "\sim"']\\
\QCoh(\cS_{H^\vee})\ar[r, "p^*"]&\QCoh(\cS_{L^\vee}),
\end{tikzcd}
\end{equation}
where $p: \cS_{L^\vee}\to \cS_{H^\vee}$.

\begin{cor}\label{corPhiHySy}
We have the canonical $\pi Z(H)$-equivariant equivalence
\begin{align}\label{eqcorfrPHIndWTHQCohSH}
\Phi_H: \frP_H^{op}\times^{\pi Z(H)}\Ind\cW(\cT_H)&\ovs{\sim}{\lrar}\QCoh(\SHv), \\
\nonumber (y, S_y)&\mapsto \cO_{\cS_{H^\vee}}, 
\end{align}
where $\pi Z(H)$ acts on the left-hand-side by the canonical involution $\pi Z(H)\ovs{\sim}{\to} \pi Z(H)^{op}$ and it acts on the right-hand-side by the canonical action. This is independent of the choice of finite central quotients $\Hf$. 
\end{cor}

\begin{proof}
We have established $\Phi_H$ for any given $H\to \Hf$. 
Let us first show that $(y, S_y)\mapsto \cO_{\cS_{H^\vee}}$. By construction, the $\frT_\flat$-module structure on $F_*(S_y)=S_{\yf}$ is given by the obvious projection 
\begin{align}\label{eqFlsFusFlsSy}
F_*(F^*F_*S_y)=\cO(\pi_{\cT, H}^{-1}(\yf))\otimes_\CC S_\yf\lrar F_*(S_y)=\CC\delta_y\otimes_\CC S_\yf,
\end{align}
here we have used the obvious identification $F_*(S_{y'})=\CC\delta_{y'}\otimes_\CC S_{\yf}\cong S_{\yf}$ (the latter uses $\delta_{y'}\mapsto 1$) for any $y'\in \pi_{\cT, H}^{-1}(\yf)$. The labeling $\bl_y$ (equivalent to the framing $\eta^y$) induces a canonical identification of algebras $\cO(\pi_{\cT, H}^{-1}(\yf))\cong \CC[K^*]$. 
Recall $\wt{\Upsilon}_y$ \eqref{eqwUpy}, then \eqref{eqFlsFusFlsSy} gives the natural identification $\wt{\Upsilon}_y(\frT_\flat (F_*S_y))\cong \CC[K^*]\otimes_{\CC} A_{\Hf}$.
On the other hand, using the framing $\eta^y$ and Corollary \ref{lemmaAJAvnchimuzyflat}, we have the canonical identification $\wt{\Upsilon}_y(\frT_\flat (F_*S_y))=\bigoplus_{\chi\in K^*}A_{\Hf}^\chi$. 

Now let us show that the canonical identification (through the canonical associativity $ (F_*F^*)F_*\cong F_*(F^*F_*)$ on the abelian categories by taking hearts) of the above two $A_{\Hf}$-module structures on $\wt{\Upsilon}_y(\frT_\flat (F_*S_y))$ is identified with the left vertical arrow of the following natural commutative diagram 
\begin{equation*}
\begin{tikzcd}
\bigoplus_{\chi\in K^*}A_{\Hf}^\chi\ar[r, "\prod_{\chi\in K^*}\btau_\chi^t"]\ar[d, "\sim"',"\bigoplus_{\chi\in K^*}\btau^t_\chi"]&A_{\Hf}\ar[d, "\sim"]\\
\CC[K^*]\otimes A_{\Hf}\cong \cO(\pi_{\cT, H}^{-1}(\yf))\otimes_\CC \AHf\ar[r]& \CC\delta_y\otimes_{\CC} A_{\Hf},
\end{tikzcd}
\end{equation*}
 where the right vertical identification uses the canonical identification $\CC\ovs{\sim}{\to}\CC\delta_y, 1\mapsto \delta_y$. 
 It is clear from diagram \eqref{eqresHLRyHyL} that there is an isomorphism $\wt{\varrho}: S_{z_\vn}\ovs{\sim}{\to}\res_\vn S_{\yH}$. Then by the canonical compatibility of $F_*, F^*$ with $\res_\vn$, we have a natural commutative diagram where the vertical isomorphisms are determined by $\wt{\varrho}$
\begin{equation*}
\begin{tikzcd}
(F_*F^*)F_*(\res_\vn S_{\yH}) \ar[r, "\sim"] &F_*(F^*F_*)(\res_\vn S_{\yH})\\
(F_*F^*)F_*S_{z_\vn}\ar[u, "\sim"] \ar[r, "\sim"] &F_*(F^*F_*)S_{z_\vn}\ar[u, "\sim"']
\end{tikzcd}. 
\end{equation*}
Using diagram \eqref{diagUpsilonzJflatyJ}, the sought-for canonical isomorphism is determined by the upper row under $\Upsilon_{\yf}^\vn$ in the above diagram, through the unit map of the adjunction associated to the functor $\otimes_{\AHf} \Avn: \Mod(\AHf)\to \Mod(\Avn)$ $=\Mod(\End(\res_\vn S_\yf)^{op})$ as the right adjoint. It then suffices to check the corresponding statement for the lower row under $\Upsilon_{\yfvn}$. 
 The latter is straightforward to check, and it is in fact the natural isomorphism $\tsfa_*\pr_1^*\ovs{\sim}{\to}q^*q_*:\QCoh(\Tfv)\to \QCoh(\Tfv)$ (up to switching from $A^\vn$ to $A^{\vn, op}$) associated to the Cartesian diagram
 \begin{equation*}
 \begin{tikzcd}
K^*\times \Tfv\ar[r, "\pr_1"]\ar[d, "\tsfa"']&\Tfv\ar[d,"q"]\\
\Tfv \ar[r, "q"]&T^\vee
 \end{tikzcd}
 \end{equation*}
 applied to $\cO_{\Tfv}$.

Next, we show that $\Phi_H$ is independent of the choice of finite central quotients. For any $H\to H^\flat_1\to H^\flat_2$,  let $\frT_i$ be the monad acting on $\Ind\cW(\cT_{H^\flat_i})$ from $\cT_{H}\to \cT_{H^\flat_i}$. 
By the $\pi Z(H)$-equivariant functor $\Ind \cW(\cT_{H_1^\flat})\to \Ind\cW(\cT_{H_2^\flat})$, 
the natural commutative diagram 
 \begin{equation*}
\begin{tikzcd}
\cW_{\frP, H^\flat_1}^{op}\times \Ind\cW(\cT_{H^\flat_1})\ar[d]\ar[r]&\Mod(A_{H^\flat_1})\ar[d]\\
\cW_{\frP, H^\flat_2}^{op}\times \Ind\cW(\cT_{H^\flat_2})\ar[r]&\Mod(A_{H^\flat_2})
\end{tikzcd}
\end{equation*}
induces the following natural  commutative diagram 
  \begin{equation}\label{eqPHopIWAH1AH2}
\begin{tikzcd}
\frP_H^{op}\times^{\pi Z(H)} \Ind\cW(\cT_{H_1^\flat})\ar[d]\ar[r,"\sim"]&\Mod(A_{H_1^\flat})\ar[d]\\
\frP_H^{op}\times^{\pi Z(H)} \Ind\cW(\cT_{H_2^\flat})\ar[r, "\sim"]&\Mod(A_{H_2^\flat})
\end{tikzcd}
\end{equation}
that is automatically $\pi Z(H)^{op}$-equivariant. The left vertical arrow in \eqref{eqPHopIWAH1AH2} is naturally right-lax equivariant with respect to the monad action $\frT_i$ on $\Ind\cW(\cT_{H_i^\flat}), i=1,2$, 
 so the right vertical arrow is naturally right-lax equivariant with the standard monad acting on each of $\Mod(A_{H_i^\flat}), i=1,2$. 

Let $\Isog^{\circ}_{H/}$ be the ordinary category of finite central quotients of $H$ that have connected centers, with a morphism between two objects given by the quotient map (if there is one). 
The above argument shows that we have a natural diagram $\Isog^\circ_{H/}\times \Delta^1\to \Pr_\st^L$, where the diagram indexed by $\Isog^\circ_{H/}\times \{0\}$ (resp. $\Isog^\circ_{H/}\times \{1\}$) is given by the analog of the left column (resp. right column) of \eqref{eqPHopIWAH1AH2}, and the diagram is naturally $\pi Z(H)$-equivariant and it is compatible with the respective monad actions as above. Since $\Isog^\circ_{H/}$ is a filtered category, taking the module categories of the respective monads, we get the canonical equivalence \eqref{eqcorfrPHIndWTHQCohSH} that is independent of the choice of $H^\flat$. 
\end{proof}

Corollary \ref{corPhiHySy} implies that we have a canonical identification 
\begin{align*}
\Phi_{H, y}:=\Phi_H(y, -)\cong \wt{\Upsilon}^{K^*}_{y}. 
\end{align*}

\begin{remark}\label{remarkPhiHcJH}
If we consider $\cJ_H$ instead of $\cT_H$, then $\frP_H$ is canonically trivialized by the identity Kostant section. Thus $\Phi_H$ gives a canonical equivalence from $\cW(\cJ_H)$ to $\Coh(\cS_{H^\vee})$.  
\end{remark}

Let $[y_H]\in \pi_0\frP_H$ and let $p_{H, L}: \pi_0(\frP_H)\to \pi_0(\frP_L)$ be the natural map \S\ref{ssspi0frZtorsorfrK} (in which it was denoted as $p_{J',J}$).

\begin{prop}\label{propresHLSyH}
\begin{itemize}
\item[(i)] 
For any $\varrho: S_{\yLf}\ovs{\sim}{\to}\res^{\Hf}_{\Lf} S_{\yHf}$,  $\res^H_L(S_{y_H})\cong S_{y_L}$ in $\cW(\cT_L)$ if and only if $(\psi_{\varrho;\eta}^\chi)_\chi\in K$ belongs to $K_L:=\ker(K\to \pi_0Z(L))$. 

\item[(ii)]
$\res^H_L(S_{y_H})\cong S_{y_L}$ in $\cW(\cT_L)$ if and only if  $p_{H, L}([y_H])=[y_L]$. 
\end{itemize}

\end{prop}

\begin{proof}
(i) Given any $\yH$, choose $\yL$ so that $(\yH, \yL)$ is $(K, \varrho)$-compatible. Then we have the diagram \eqref{eqresHLRyHyL}. By Corollary \ref{corPhiHySy}, $\Phi_{H, \yH}(S_{\yH})\cong \wt{\Upsilon}_{\yH}^{K^*}(S_{\yH})\cong \cO_{\cS_{H^\vee}}$ and $\Phi_{H, \yL}(S_{\yL})\cong \wt{\Upsilon}_{\yL}^{K^*}(S_{\yL})\cong \cO_{\cS_{L^\vee}}$, therefore, $\res_L^H(S_{\yH})\cong S_{\yL}$. Since for $y_L'$ over $\yLf$, $S_{y_L'}\cong S_{\yL}$ in $\cW(\cT_L)$ if and only $y_L'\in K_L\cdot y_L$, the statement follows.

(ii) Let $(\yH, \yL)$ be as in (i) under which  we have the diagram \eqref{eqresHLRyHyL}. Consider the corresponding natural commutative diagram by taking the left adjoint of $\res_L^H$ and $p^*$, which are respectively $\cores_L^H$ and $p_*$. 
Let $\pi: \cS_{L^{\flat,\vee}}\to \cS_{L^\vee}$. Then using 
\begin{align*}
&\QCoh(\cS_{L^\vee})\ni \cF\cong p^*(\cO_{\cS_{\dH}})\Longleftrightarrow\ \pi^*\cF\cong \cO_{\cS_{L^{\flat,\vee}}}\text{ and }
p_*(\cF)\text{ has a direct summand isomorphic to }\cO_{\cS_{\dH}},
\end{align*}
we see that $S_{y_L'}\cong \res^H_L(S_{y_H}), y_L'\in \frP_L$ if and only if $\cores(S_{y_L'})$ has a direct summand isomorphic to $S_{y_H}$. 
It follows from Lemma \ref{lemmapifrKJpJ}, on which \S\ref{ssspi0frZtorsorfrK} is based, Corollary \ref{corpiZHQfKstarfinal} and \cite[Proof of Proposition 5.2]{J}, that the last equivalent condition above holds if and only if $p_{H, L}([y_H])=[y_L']$. This finishes the proof. 
\end{proof}

\sss{$\nu^\cL_\varrho$ for general $H$}
Now we have the direct generalization of $\nu^\cL_\varrho$ for any $\varrho: S_{\yL}\ovs{\sim}{\to} \res^H_L S_{\yH}$. By the last part of the proof of Proposition \ref{propfvarrhocLCtimes}, we directly get the following functoriality of $\nu^\cL_\varrho$ with respect to covering maps.

\begin{lemma}\label{lemmavarrhopicL1}
Let $\pi: H\to H_1$ be a finite central quotient. For any $\yH$ and $\yL$, let $y_{H_1}$ and $y_{L_1}$ be their respective images under the projection $\cT_{H}\to \cT_{H_1}$. For any $\varrho:  S_{\yL}\ovs{\sim}{\to} S_{\yH}$, let $\varrho_1: S_{y_{L_1}}\ovs{\sim}{\to} S_{y_{H_1}}$ be the image isomorphism. Then for any rank $1$ local system $\cL_1$ on $\cT_{H_1}$, we have the commutative diagram 
\begin{equation*}
\begin{tikzcd}
(\pi^*\cL_1)_{y_H}\ar[d, equal]\ar[r, "\nu_{\varrho}^{\pi^*\cL_1}", "\sim"']&(\pi^*\cL_1)_{y_L}\ar[d, equal]\\
(\cL_1)_{y_{H_1}}\ar[r, "\nu_{\varrho_1}^{\cL_1}", "\sim"']&(\cL_1)_{y_{L_1}}. 
\end{tikzcd}
\end{equation*}
In particular, if $\pi^*\cL_1$ is a trivial local system, then $\nu_{\varrho_1}^{\cL_1}$ is the parallel transport along any homotopy class of paths from $y_{H_1}$ to $y_{L_1}$ that lifts to a homotopy class of paths from $y_H$ to $y_L$. 
\end{lemma}

\begin{cor}\label{corvarrhoresLHpstar}
For any $\varrho:S_{y_L}\ovs{\sim}{\to} \res^H_L(S_{y_H})$, let $\varrho^\flat: S_{\yLf}\ovs{\sim}{\to} \res^{\Hf}_{\Lf}(S_{\yHf})$ be the image isomorphism. 
Then $(\yH, \yL)$ is $(K, \varrho^\flat)$-compatible, i.e. 
 $(\psi_{\varrho^\flat, \eta}^\chi)_\chi=1\in K$ for the framings $\eta^{\yH}$ and $\eta^{\yL}$ (cf. \S\ref{psichivarrhoeta}). 
In this case, we get a canonical commutative diagram
\begin{equation}\label{eqimageresHLRyHyL}
\begin{tikzcd}[column sep=6em]
\Ind\cW(\cT_{H})\ar[r, "\res^H_L"]\ar[d,"\Phi_{H, \yH}"']&\Ind\cW(\cT_{L})\ar[d, "\Phi_{L,\yL}"]\ar[dl, Rightarrow, "\sim"']\\
\QCoh(\cS_{H^\vee})\ar[r, "p^*"]&\QCoh(\cS_{L^\vee}).
\end{tikzcd}
\end{equation}
\end{cor}

\begin{proof}
Since there is the canonical identification of local systems $\pi^*\cL_K\cong \cO(K)_{\cT_H}$, we have $\nu_{\varrho}^{\pi^*\cL_\chi}: \ell_{\yHf, \chi}\ovs{\sim}{\to}\ell_{\yLf, \chi}$ takes the framing $\eta^{\yH}_\chi$ to $\eta^{\yL}_{\chi}$. Then the statement follows from Lemma \ref{lemmavarrhopicL1} and diagram \eqref{eqresHLRyHyL}. 
\end{proof}

%%%%%%%%%%%%%%%%%%%%%%%

\subsection{Compatibility of $\Phi_H$ with restrictions to a Levi $L\subset H$}\label{ssPhiHrestriction}

Now we spell out the compatibility with respect to passing to Levi subgroups. Let $L\subset H$ be a Levi subgroup containing $T$. We would like to define a functor
\begin{equation}
\t^{H}_{L}: \frP_{H}\to \frP_{L}
\end{equation}
together with a natural isomorphism making the following diagram commutative
\begin{equation}
\xymatrix{\frP_H\ar[d]_{\t^{H}_{L}}\ar[r]^-{\nu_{H}} & \cW(\cT_{H})\ar[d]^{\res^{H}_{L}}\\
\frP_L\ar[r]^-{\nu_{L}} & \cW(\cT_{L})}. 
\end{equation}
Instead of constructing $\t^{H}_{L}$ directly, we make the following constructions.

Let $\iota: \cW_{\frP, L}\hookrightarrow \cW(\cT_{L})^\Gpd$ be the $\pi Z(L)$-equivariant fully faithful embedding, and 
let $\cW(\cT_{L})^\Gpd_{\dmd}$ be the union of the connected components in $\cW(\cT_{L})^\Gpd$ that are in the image of $\iota$ on $\pi_0$. 
Let $\cW_{\frP, H, L}$ be the groupoid defined by the following Cartesian square (in $\Gpd_\infty$ or $\Cat_\infty$)
\begin{equation}\label{eqWPHLWTLGpddmd}
\begin{tikzcd}
\cW_{\frP, H, L}\ar[r, "r_H", "\sim"']\ar[d]&\cW_{\frP, H}\ar[d, "\res_L^H"]\\
\cW_{\frP, L}\ar[r, "\sim"]&\cW(\cT_{L})^\Gpd_{\dmd}\ar[ul, "\lrcorner", phantom]
\end{tikzcd},
\end{equation}
which is well defined by Proposition \ref{propresHLSyH}. 
Explicitly, we can present $\cW_{\frP, H, L}$ by the following. The objects consist of triples $(y_H, y_L;\varrho: \res^H_L S_{\yH} \ovs{\sim}{\to} S_{\yL})$ and a morphism between $(y_H, y_L, \varrho)$ and $(y_H', y_L',\varrho')$ is a pair $\varphi_H: S_{\yH}\ovs{\sim}{\to} S_{y_H'}$ and $\varphi_L: S_{y_L}\ovs{\sim}{\to} S_{y_L'}$ such that $\res_{L}^H(\varphi_H)$ and $\varphi_L$ are conjugate by $\varrho$ and $\varrho'$.

We have the following natural commutative diagram, whose left column is naturally $K^*$-equivariant with respect to the $K^*$-actions
\begin{equation*}
\begin{tikzcd}
\cW_{\frP,\Hf}^{op}\times\Ind\cW(\cT_\Hf)\ar[r]\ar[d]&\Mod(\AHf)\ar[d, "\otimes_{\AHf} \ALf"]\\
\cW(\cT_{\Lf})^{\Gpd,op}_{\dmd}\times \Ind\cW(\cT_{\Lf})\ar[r]&\Mod(\ALf)\\
\cW_{\frP, \Lf}^{op}\times \Ind\cW(\cT_{\Lf})\ar[r]\ar[u, "\sim"]&\Mod(\ALf)\ar[u, equal]
\end{tikzcd}.
\end{equation*}
By taking fiber products on both sides, we get a natural commutative diagram 
\begin{equation}\label{eqcWfrPHfLffrG}
\begin{tikzcd}
\cW_{\frP,\Hf}^{op}\times\Ind\cW(\cT_{\Hf})\ar[r]&\Mod(\AHf)\ar[d, equal]\\
\cW_{\frP,\Hf,\Lf}^{op}\times \Ind\cW(\cT_{\Hf})\ar[r]\ar[d]\ar[u, "\sim"]&\Mod(\AHf)\ar[d, "\otimes_{\AHf}\ALf"]\ar[dl, Rightarrow, "\sim"]\\
\cW_{\frP, \Lf}^{op}\times \Ind\cW(\cT_{\Lf})\ar[r]&\Mod(\ALf)
\end{tikzcd},
\end{equation}
in which the top square is strictly commutative and for any $(\yHf, \yLf, \varrho^\flat)\in \cW_{\frP,\Hf,\Lf}^{op}$ and the lower natural 2-isomorphism restricts to \eqref{diagUpsilonzJflatyJ} (after replacing $L_\vn^\flat$ by $L^\flat$ in that diagram).

\begin{lemma}\label{lemmacWfrPHLphinu}
There is a natural commutative diagram (up to a contractible space of choices)
\begin{equation*}
\begin{tikzcd}[row sep=1.5em, column sep=6em]
\cW_{\frP, \Hf}^{op}\ar[r, "\phi'_{\Hf}"]&\Xi(K^*)\simeq \BB K\\
\cW_{\frP, \Hf,\Lf}^{op}\ar[r, "\phi'_{\Hf, \Lf}"]\ar[u, "\sim"]\ar[d]&\Xi(K^*)\ar[u, equal]\ar[d, equal]\ar[dl, "\nu"', Rightarrow, "\sim"]\\
\cW_{\frP, \Lf}^{op}\ar[r,  "\phi'_{\Lf}"']&\Xi(K^*)
\end{tikzcd},
\end{equation*}
where $\phi'_{\Hf}=(\phi_{\cL_K}|_{\frP_{\Hf}^{op}})\circ q_{\frP, \Hf}^{op}$ (cf. \eqref{diagQPf}; similarly for $\phi'_{\Lf}$) and $\phi'_{\Hf, \Lf}$ is the functor that makes the upper square strictly commutative. 
The diagram naturally lifts to a commutative diagram (up to a contractible space of choices)
\begin{equation}\label{eqwtphiprimeHL}
\begin{tikzcd}[row sep=1.5em, column sep=6em]
\cW_{\frP, H}^{op}\ar[r,  "\wt{\phi}'_{H}"]&\wt{\Xi}(K^*)\simeq \pt\\
\cW_{\frP, H,L}^{op}\ar[r, "\wt{\phi}'_{H, L}"]\ar[u, "\sim"]\ar[d]&\wt{\Xi}(K^*)\ar[u, equal]\ar[d, equal]\ar[dl, Rightarrow, "\wt{\nu}"',"\sim"]\\
\cW_{\frP, L}^{op}\ar[r, "\wt{\phi}'_{L}"']&\wt{\Xi}(K^*)
\end{tikzcd}, 
\end{equation}
in which $\wt{\phi}'_{H}=(\wt{\phi}_{\cL_K}|_{\frP_{H}^{op}})\circ q_{\frP, H}^{op}$ (and similarly for $\wt{\phi}'_L$).  
\end{lemma}

\begin{proof}
For any $(y_H^\flat, y_L^\flat; \varrho^\flat)$, the natural 2-isomorphism $\nu$ is given by 
\begin{align*}
\nu_{\varrho^{\flat,-1}}^{\cL_\chi}=\psi_{\varrho^{\flat,-1}}^{\chi}: \ell_{y_H^\flat, \chi}\ovs{\sim}{\lrar} \ell_{y_L^\flat, \chi},\text{ for } \chi\in K^*
\end{align*}
(cf. \eqref{eqnuupcL}). The functoriality of $\nu_{\varrho^{\flat,-1}}^{\cL_\chi}$ with respect to any morphism $(y_H^\flat, y_L^\flat; \varrho^\flat)\to (\yHpf, \yLpf; \varrho^{\prime, \flat})$ has been explained in Remark \ref{remnuLcompatiblevarrho}. This gives the first diagram.

For the second diagram as a lifting, 
let $\wt{\nu}$ be the lifting of $\nu$ that attaches to any $(y_H, y_L; \varrho)$ the isomorphism
\begin{align*}
(\ell_{\yHf}, \eta^{\yH})\ovs{\sim}{\lrar} (\ell_{\yLf}, \eta^{\yL})
\end{align*}
in $\wt{\Xi}(K^*)$. 
This is well defined by Corollary \ref{corvarrhoresLHpstar}. The rest follows easily. 
\end{proof}

Similarly to the definition of $\PhiHP$ in \S\ref{sssDefwPhiHP} and \S\ref{sssdefPhiHP}, after combining diagram \eqref{eqcWfrPHfLffrG} and \eqref{eqwtphiprimeHL} and taking $K^*$-invariants, we get a natural commutative diagram 
\begin{equation}\label{eqWPHIndWQCoh}
\begin{tikzcd}
\cW_{\frP, H}^{op}\times^{\pi Z(H)}\Ind\cW(\cT_H)\ar[r]&\QCoh(\cS_{H^\vee})\\
\cW_{\frP, H,L}^{op}\times^{\pi Z(H)}\Ind\cW(\cT_H)\ar[r]\ar[u,"\sim"]\ar[d]&\QCoh(\cS_{H^\vee})\ar[u, equal]\ar[d, "p^*"]\ar[dl, Rightarrow, "\sim"]\\
\cW_{\frP, L}^{op}\times^{\pi Z(L)}\Ind\cW(\cT_L)\ar[r]&\QCoh(\cS_{L^\vee})
\end{tikzcd},
\end{equation}
where (1) the restriction of the top row to $\frP_{H}^{op}\times^{\pi Z(H)}\Ind\cW(\cT_H)$ naturally identifies with $\Phi_H$ by construction (and similarly for the bottom row); (2) for each $(y_H, y_L, \varrho)\in \WPHL$, the lower square restricts to \eqref{eqimageresHLRyHyL} (up to replacing $\varrho$ by $\varrho^{-1}$), by Corollary \ref{corvarrhoresLHpstar}.

\sss{The $\pi\frZ$-torsor $\frP$ over $\PfI^{op}$}
Let $\frP_{H, L}$ be the groupoid consisting of triples $(y_H,y_L,[\varrho])$ where $y_H\in \frP_H$, $y_L\in \frP_L$ and $[\varrho]$ is an isomorphism  $\varrho: \res^H_L(S_{y_H}) \ovs{\sim}{\lrar} S_{y_L}$ modulo the scaling action of $\CC^\times$. This is the same as the groupoid obtained from $\cW_{\frP, H, L}$ by taking quotient of the $\CC^\times$-action on morphism sets. In particular, we get a natural commutative diagram
\begin{equation}\label{eqWPHLfrPHL}
\begin{tikzcd}
\WPH\ar[d]&\WPHL\ar[l, "\sim"]\ar[r]\ar[d]&\WPL\ar[d]\\
\frP_H&\frP_{H, L}\ar[l, "\sim"]\ar[r]&\frP_{L}
\end{tikzcd}, 
\end{equation}
in which all vertical arrows are  taking quotient of the $\CC^\times$-action on morphism sets. Note that each row naturally carries the action of the diagram $\pi\frZ_{H,L}:=(\pi Z(H)=\pi Z(H)\rightarrow \pi Z(L))$, and the vertical arrows are equivariant with respect to those actions. Since the lower row is a $\pi\frZ_{H,L}$-torsor, a $\pi\frZ_{H,L}$-equivariant section of the vertical arrows is determined (up to a contractible space of choices) by choosing for a given object $(y_H, y_L, [\varrho])\in \frP_{H, L}$ a lifting $(y_H, y_L, \varrho: \res^H_LS_{y_H}\ovs{\sim}{\to} S_{y_L})$ in $\cW_{\frP, H, L}$, which forms a torsor of $\CC^\times$. Moreover, if we restrict the left column of diagram \eqref{eqWPHIndWQCoh} to such a section, then we get an equivalence between the two columns. 

To systematically encode the associativity for triples $H\supset L\supset M$, we consider the category $\frP$ together with a fibration
\begin{equation}\label{eqfrPtoPfI}
\frP\to \cP_{ft}(\wt I)^{op}
\end{equation}
as follows. Objects in $\frP$ are pairs $(J, y)$ where $J\sft \wt I$ and $y\in \frP_{L_J}$. A morphism $(J,y)\to (J',y')$ when $J\supset J'$ is an isomorphism $\res^{J}_{J'}S_{y}\isom S_{y'}$ in $\cW(\cT_{L_{J'}})$ modulo the scaling $\CC^\times$-action. Passing to nerves we get a left fibration of $\infty$-categories. 

The functor $J\mt \pi(ZL_{J})$ via Grothendieck construction also gives a fibration, which we denote by $\pi\frZ\to \cP_{ft}(\wt I)^{op}$. Then $\frP$ is a $\pi\frZ$-torsor.

Similarly, we get a left fibration $\frWP\to \PfI^{op}$ and a $\pi \frZ$-equivariant functor 
\begin{align}\label{eqvarpifrWPfrP}
\varpi: \frWP\to \frP  \text{ over }\PfI^{op}, 
\end{align}
whose unstraightening encodes \eqref{eqWPHLfrPHL} for each $J\supset J'$ by setting $H=L_J$ and $L=L_{J'}$. Alternatively, one can assemble the diagrams \eqref{eqWPHLWTLGpddmd} for $L_J\supset L_{J'}, J\supset J'\in \PfI$, into a functor $\PfI^{op}\lrar \Corr(\Gpds)^{\hor=\mathrm{isom}}$, where the target is the 1-full subcategory of $\Corr(\Gpds)$ whose 1-morphisms consist of correspondences with invertible horizontal maps. Since there is a canonical equivalence  $\Gpds\to  \Corr(\Gpds)^{\hor=\mathrm{isom}}$, the aforementioned functor canonically factors as $\WP: \PfI^{op}\to \Gpds$, and the left fibration $\frWP\to \PfI^{op}$ is naturally identified with the Grothendieck construction of $\WP$. 

Let $\frP^{op}\to  \cP_{ft}(\wt I)^{op}$ be the left fibration where $\frP_{L_J}, \frP_{L_J, L_{J'}}$, etc. are replaced by their respective opposite categories. Let $\frWP^{op}\to  \cP_{ft}(\wt I)^{op}$ be defined analogously. 

\begin{defn}\label{defncompatibleKostantsec}
We say a section of the left fibration $\frP\to \cP_{ft}(\wt I)^{op}$ (up to isomorphisms) is a \emph{compatible system of Kostant sections}. 
\end{defn}

\begin{lemma}\label{l:section frY}
The left fibration $\frP\to \cP_{ft}(\wt I)^{op}$ has a section. Therefore, there exists a compatible system of Kostant sections. 
\end{lemma}
\begin{proof}
The obstruction to the existence of a section lies in $\cohog{1}{\triangle, \pi\frZ}$, which vanishes by Lemma \ref{l:cohoZ}.
\end{proof}

\begin{cor}\label{corfrWPsection}
The left fibration $\frWP\to \cP_{ft}(\wt I)^{op}$ has a section. Therefore, $\varpi$ \eqref{eqvarpifrWPfrP} admits a $\pi\frZ$-equivariant section $\wh{\sigma}: \frP\to \frWP$ over $\PfI^{op}$. 
\end{cor}

\begin{proof}
Choose any section $s$ of $\frP\to \cP_{ft}(\wt I)^{op}$, the obstruction to lifting $s$ to a section $\sigma$ of $\frWP\to \cP_{ft}(\wt I)^{op}$ lies in $\cohog{2}{\trg, \CC^\times}$, which is trivial. Any lifting $\sigma$ induces a $\pi\frZ$-equivariant section $\wh{\sigma}: \frP\to \frWP$ over $\PfI^{op}$. 
\end{proof}

\begin{remark}
Since $\cohog{1}{\trg, \CC^\times}$ is trivial, the lifting $\sigma$ of $s$ is unique up to isomorphisms. Thus a compatible system of Kostant sections is also equivalent to a section of $\frWP\to \cP_{ft}(\wt I)^{op}$ (up to isomorphisms). 
\end{remark}

\sss{The coCartesian fibrations $\frIW$ and $\frQ$ over $\PfI^{op}$}

Let $\frIW\to \cP_{ft}(\wt I)^{op}$ (resp. $\frQ\to \cP_{ft}(\wt I)^{op}$) be the Grothendieck construction applied to the functor $\cP_{ft}(\wt I)^{op}\to \Pr_{\st}^{L}$ given by $J\mt\IndW(\cT_{L_J})$ (resp. $J\mt\QCoh(S_{L_J^\vee})$) under restriction functors (resp. $*$-pullbacks).

\begin{lemma}\label{lemmafrMPfrQ}
Applying similar construction of $\frWP$, the diagrams \eqref{eqWPHIndWQCoh} for $H=L_J, L=L_{J'}, J\supset J'\in \PfI$, naturally assemble into the right commutative triangle of the following diagram
\begin{equation}\label{eqfrMPfrIWfrQfrP}
\begin{tikzcd}
\frP^{op}\times^{\pi\frZ}\frIW\ar[dr]&\ar[l] \frWP^{op}\times^{\pi\frZ}\frIW\ar[d]\ar[r]&\frQ\ar[dl]\\
&\cP_{ft}(\wt I)^{op}
\end{tikzcd},
\end{equation}
which is $\pi\frZ$-equivariant with respect to the natural $\pi\frZ$-actions on them. 
\end{lemma}

\begin{proof}
For the first part on the right commutative triangle, it suffices to check that the diagram $\cQ_{\PfI^{op}}$, from assembling the right column of \eqref{eqWPHIndWQCoh}, indeed agrees with the standard diagram $\QCoh(\cS_{L_J^\vee}), J\in \PfI^{op}$, under $*$-pullback. By construction, for any $k$-simplex $J_0\supset J_1\supset\cdots\supset J_k$, we have a commutative diagram 
\begin{equation*}
\begin{tikzcd}
\QCoh(\cS_{L_{J_0}^\vee})\ar[d]\ar[r]&\QCoh(\cS_{L_{J_1}^\vee})\ar[r]\ar[d]&\cdots\ar[r]&\QCoh(\cS_{L_{J_k}^\vee})\ar[d]\\
\QCoh(\cS_{L_{J_0}^{\flat, \vee}})\ar[r]&\QCoh(\cS_{L_{J_1}^{\flat, \vee}})\ar[r]&\cdots\ar[r]&\QCoh(\cS_{L_{J_k}^{\flat, \vee}})
\end{tikzcd}
\end{equation*}
 (for a choice of $L_{J_0}^\flat$), where (1) the bottom diagram is the standard one under $*$-pullbacks; (2) the vertical functors and horizontal functors on the top row are the respective $*$-pullbacks; (3) each square indexed by $\Delta^1\times (J_i\supset J_j), i<j$, is canonically identified with the standard square with $*$-pullbacks, in view of Remark \ref{remgeneralizetoresHfLf} and Corollary \ref{corvarrhoresLHpstar}. 
 
 Note that we can pick a collection of generators for each $\QCoh(\cS_{L_{J_i}^\vee})$, say $\cO^\kappa, \kappa\in K$, such that the functors on the top diagram preserve these generators, and $\Hom_{\QCoh(\cS_{L_{J_i}^\vee})}(\cO^{\kappa_1}, \cO^{\kappa_2})$ for $\kappa_i\in K$, is concentrated in degree $0$. 
Therefore, to show that the top diagram is canonically identified with the standard diagram under $*$-pullback (up to a contractible space of choices), it suffices to show the case for every 2-simplex. Now for any 2-simplex $J_0\supset J_1\supset J_2$ of $\PfI^{op}$, let $p_{j,i}: \cS_{L_{J_j}^\vee}\to \cS_{L^\vee_{J_i}}, i<j$, denote the projection. It is easy to see that the natural 2-isomorphism $p_{2,1}^*p_{1,0}^*\ovs{\sim}{\to}p_{2,0}^*$ must be the trivial one. Indeed, for each generator $\cO_{\cS_{L_{J_0}^\vee}}^\kappa\in \QCoh(\cS_{L_{J_0}^\vee})$, the natural 2-isomorphism gives an automorphism of $p_{2,0}^*\cO_{\cS_{L_{J_0}^\vee}}^\kappa=\cO_{\cS_{L_{J_2}^\vee}}^\kappa$, which must be the identity using that the functor between abelian categories $\QCoh(\cS_{L_{J_2}^\vee})^\hs\to \QCoh(\cS_{L_{J_2}^{\flat, \vee}})^\hs$ is faithful. 

For the last part on the $\pi\frZ$-equivariance with respect to the natural $\pi\frZ$-actions on them, by a similar argument as above, it suffices to check that the induced $\pi\frZ_{H, L}$-action on the right column of \eqref{eqWPHIndWQCoh} is the standard one (up to a contractible space of choices). 
Again, by a similar argument as above (and using the established Corollary \ref{corPhiHySy}), it suffices to observe (1) from Proposition \ref{eqbiakappa} that for any $(y_H, y_L, \varrho)\in \WPHL$ and $\kappa\in K$, we have a strict commutative diagram (up to a contractible space of choices)
\begin{equation*}
\begin{tikzcd}[column sep=6em]
K^*\car\Mod(\AHf)\ar[d, "\bia_\kappa"]\ar[r, "\otimes_{\AHf} \ALf"]&\Mod(\ALf)\cal K^*\ar[d, "\bia_\kappa"]\\
K^*\car\Mod(\AHf)\ar[r, "\otimes_{\AHf} \ALf"]&\Mod(\ALf)\cal K^*,
\end{tikzcd}
\end{equation*}
where the $K^*$-actions are the standard ones and $\bia_\k$ is only doing $\kappa^{-1}\in \Aut_{\Xi(K^*)}(\ell_\triv)$; (2) taking $K^*$-invariants, we get the standard diagram as desired. 
\begin{equation*}
\begin{tikzcd}[column sep=6em]
\QCoh(\SHv)\ar[d, "\otimes\cO^\kappa"]\ar[r, "p^*"]&\QCoh(\cS_{\dL})\ar[d, "\otimes\cO^\kappa"]\\
\QCoh(\SHv)\ar[r, "p^*"]& \QCoh(\cS_{\dL})
\end{tikzcd}
\end{equation*}

This finishes the proof. 
\end{proof}

\subsection{Finishing the proof of the Mirror Theorem \ref{th:mirror}}\label{ssfinishingpfMT}

Let $\s$ be a section of $\frWP^{op}\to \PfI^{op}$ that induces a section $\wh{\s}$ of $\varpi: \frWP^{op}\to\frP^{op}$ over $\PfI^{op}$, which exists by Corollary \ref{corfrWPsection}. Then we obtain from \eqref{eqfrMPfrIWfrQfrP} 
an equivalence over $\cP_{ft}(\wt I)^{op}$
\begin{equation}
\Phi_{\sigma}: \frIW \simeq \frP^{op}\times^{\pi\frZ}\frIW\ovs{\sim}{\lrar} \frQ.
\end{equation}
Taking limits we obtain
\begin{equation}
\IndW(\cM^{\c}_{G})\simeq \lim_{J\sft \wt I}\QCoh(\cS_{L^{\vee}_{J}}).
\end{equation}
Combining with Corollary \ref{c:QCoh J as lim}, we get a desired equivalence 
\begin{align*}
\IndW(\cM^{\c}_{G})\simeq \QCoh(J_{G^\vee}^{G^{\vee, \ad}}). 
\end{align*}

Choose any $G^\flat$ with connected center and let $K=\ker(G\to G^\flat)$. Then we can canonically identify $\QCoh(J_{\dG}^{G^{\vee, \ad}})$ (resp. $\QCoh(\cS_{\dG})$) with $\QCoh(J_{G^{\flat, \vee}}^{G^{\vee, \ad}})^{K^*}$ (resp. $\QCoh(\cS_{G^{\flat, \vee}})^{K^*}$). Let $\pi: J_{\dG}^{G^{\vee, \ad}}\to \cS_{\dG}$ be the natural projection. 

\begin{prop}\label{propMirrorconnectedZ} 
For any section $\s$ of $\frWP^{op}\to \PfI^{op}$,  we have a natural commutative diagram
\begin{equation*}
\begin{tikzcd}[column sep=6em]
\cW(\cM_G^\circ)\ar[r, "\lim \Phi_\sigma", "\sim"']&\Coh(J_{\dG}^{G^{\vee, \ad}})\\
\cW(\cT_G)\ar[r, "\Phi_{\sigma(I)}", "\sim"']\ar[u, "\cores"]&\Coh(\cS_{\dG})\ar[u, "\pi^*"].
\end{tikzcd}
\end{equation*}
Under $\lim\Phi_\sigma$, the Kostant section $S_{\sigma(I)}$ (which is also a Hitchin section) is sent to the structure sheaf on $J_{\dG}^{G^{\vee, \ad}}$ (up to a canonical isomorphism). For any $\kappa\in K=\XX^*(K^*)$, $S_{\kappa\cdot \sigma(I)}$ is sent to $\cO^{\kappa}$, the structure sheaf on $J_{\dG}^{G^{\vee, \ad}}$ twisted by the character $\kappa$  (up to a canonical isomorphism). In particular, the Kostant sections $S_{\kappa\cdot \sigma(I)}, \kappa\in K$, generate $\cW(\cM_G^\circ)$. 

\end{prop}

\begin{proof}
It is clear that $\Phi_\sigma$ induces the left commutative square of the following; the functor $\pi_*$ directly follows from the right commutative square, where $\Psi$ is as in Theorem \ref{th:main QCoh}
\begin{equation*}
\begin{tikzcd}[column sep=5em]
\IndW(\cM_G^\circ)\ar[r, "\lim \Phi_\sigma", "\sim"']\ar[d, "\res"]&\QCoh(J_{\dG}^{G^{\vee, \ad}})\ar[d, "\pi_*"]\ar[r, hook, "\Psi"]&\QCoh(T^\vee)^{\Wa}\cong \QCoh(T^{\vee,\ad}\times T^\vee)^W\ar[d]\\
\IndW(\cT_G)\ar[r, "\Phi_{\sigma(I)}", "\sim"']&\QCoh(\cS_{\dG})\ar[r, hook]&\QCoh(T^\vee)^{W}
\end{tikzcd}.
\end{equation*}
Then the diagram in the proposition comes from taking the left adjoint of $\res$ and $\pi_*$, respectively, and restrict to the respective subcategories of compact objects. 
Since $\Phi_{\sigma(I)}(S_{\kappa\cdot \sigma(I)})=\cO^{\kappa}_{\cS_{G^\vee}}$, the rest follows immediately. It is clear that one only needs a Kostant section from each connected component of $\frP_G$ to generate $\cW(\cM_G^\circ)$. 
\end{proof}

%%%%%%%%%%%%%%%%%%%%%%%%%%%%%%%%%%%%%%%%%%%%%%%%%%%%%%%%

%%%%%%%%%%%%%%%%%%%%%%%%%%%%%%%%%%%%%%%%%%%
%%%%%%%%%%%%%%%%%%%%%%%%%%%%%%%%%%%%%%%%%%%

%%%%%%%%%%%%%%%%%%%%%%%%%%%%%%%%%%%%%%%%%%%%%%%%%%%%%

\appendix

\section{Tropical manifolds and admissible coverings}\label{AppSecTropical}

In this section, we introduce some convenient constructions and terminologies to be used in \S\ref{subsec: Liouville}--\ref{subsec: sectorial,covering} and Appendix \ref{secAppSectors}. The key notions are tropical manifolds and admissible coverings. 

\subsection{Tropical manifolds}

Let $M$ be a (connected) real smooth manifold with corners. Then the corner structure determines a natural stratification of $M$ by $M^\circ$, connected components of codimension $i$ corners contained in $\partial M$ for $i\geqslant 1$. Denote the stratification of $M$ by $\{S_\alpha\}_{\alpha\in A}$ and call $A$ the \emph{indexed set} of $M$.  Define a partial order on $A$ given by $\alpha<\beta$ if and only if $\overline{S}_\alpha\supset S_\beta$. For each $\alpha$, let $M_\alpha$ be the star of $S_\alpha$, i.e. $M_\alpha=\bigcup_{\beta<\alpha}S_\beta$. Clearly, $\alpha<\beta\Leftrightarrow M_\alpha\overset{\text{open}}{\subset} M_\beta$. Let $d_\alpha=\dim S_\alpha$ and $\cod_\alpha=\codim S_\alpha=d-d_\alpha$. Let $A_{\gneqq\alpha}:=\{\beta\in A: \beta\gneqq \alpha\}$ and $A_{>\alpha}:=\{\beta\in A: \beta> \alpha\}$. Similarly, we define $A_{\lneqq \alpha}$ and $A_{<\alpha}$. 

\begin{defn}\label{def: tropicalmfld}
 We say $M$ is a \emph{tropical manifold} of dimension $d$ if 
 \begin{itemize}
\item[(1)]  For each $\alpha$, we have a fixed isomorphism of manifolds with corners
 \begin{align*}
 \nu_\alpha: M_\alpha\overset{\sim}{\to} \RR_{\geqslant 0}^{cod_\alpha}\times \RR^{d_\alpha}
 \end{align*} 
which in particular identifies $S_\alpha$ with $\{0\}\times \RR^{d_\alpha}$;

\item[(2)] Let $\pi_\alpha$ be the projection $M_\alpha\to S_\alpha$ determined by $\nu_\alpha$ in the obvious way. Then for $\alpha>\beta$, we have 
\begin{align*}
\pi_\alpha|_{S_\beta}\circ \pi_\beta=\pi_\alpha|_{M_\beta},
\end{align*}
i.e. the map $\nu_{\beta, \alpha}:= \nu_\alpha|_{M_\beta}\circ \nu_\beta^{-1}$ fits into a commutative diagram of manifolds with corners
\begin{equation}\label{eqnubetaalphapinu}
\begin{tikzcd}[column sep=4em]
\RR_{\geqslant 0}^{\cod_\beta}\times\RR^{d_\beta}\ar[r, hook, "\nu_{\beta, \alpha}"]\ar[d, "\pi_\beta^\nu"']& \RR_{\geqslant 0}^{\cod_\alpha}\times\RR^{d_\alpha}\ar[d, "\pi_\alpha^\nu"]\\
\RR^{d_\beta}\ar[r, "\pi_{\alpha;\beta}^\nu:=\pi_\alpha^\nu\mid_{0^{\cod_\beta}\times \RR^{d_\beta}}"]&\RR^{d_\alpha},
\end{tikzcd}
\end{equation}
where $\pi_\alpha^\nu$ and $\pi_\beta^\nu$ are the obvious projections. 
Moreover, $\nu_{\beta, \alpha}$ is the composition 
\begin{equation}\label{eq: nu_beta,alpha,compose}
\begin{tikzcd}[row sep=1em, column sep=6em]
&\RR_{\geqslant 0}^{\cod_\beta}\times\RR^{d_\beta}\ar[r, "id\times f_{\beta,\alpha}"]&\RR_{\geqslant 0}^{\cod_\beta}\times\RR^{d_\beta}\ar[r, "(id\times \exp)"] &\RR_{\geqslant 0}^{\cod_\beta}\times\RR_{>0}^{d_\beta}\\
\ar[r, " {(\mu_{\beta,\alpha},\ \pr)} "]& \RR_{\geqslant 0}^{\cod_\alpha}\times\RR_{>0}^{d_\alpha}\ar[r, "(P_{\beta,\alpha}\times \log)"] &\RR_{\geqslant 0}^{\cod_\alpha}\times\RR^{d_\alpha},
\end{tikzcd}
\end{equation}
where $f_{\beta,\alpha}$ is an affine automorphism of $\RR^{d_\beta}$; $\exp: \RR^{d_\beta}\to  \RR^{d_\beta}_{>0}$  is the exponential map on each factor (which has inverse map $\log$); $\mu_{\beta,\alpha}: \RR_{\geqslant 0}^{\cod_\beta}\times\RR_{>0}^{d_\beta}\to \RR_{\geqslant 0}^{\cod_\alpha}$ is of the form 
\begin{align*}
(r_1^\beta,\cdots, r^\beta_{\cod_\beta}; u_1,\cdots, u_{d_\beta})\mapsto (a_1\cdot r_1^\beta\prod_i u_i^{m_{1;i}}, \cdots, a_{\cod_\beta}\cdot r_{\cod_\beta}^\beta\prod_i u_i^{m_{(\cod_\beta); i}}, u_1,\cdots, u_{d_\beta-d_\alpha}), 
\end{align*}
for some $m_{j;i}\in \QQ$ (more generally, one can assume $m_{j;i}\in \RR$) and $a_i\in \RR_{>0}$; 
$\pr: \RR^{d_\beta}_{>0}\to  \RR^{d_\alpha}_{>0}$ is the projection to the last $d_\alpha$ factors; $P_{\beta,\alpha}:\RR_{\geqslant 0}^{\cod_\alpha}\overset{\sim}{\to} \RR_{\geqslant 0}^{\cod_\alpha}$ is the automorphism given by a permutation on the $\cod_\alpha$ many coordinates. Note that for each fixed $x^\nu=(u_{d_\beta-d_\alpha+1},\cdots, u_{d_\beta})\in \RR^{d_\alpha}_{>0}$, $\mu_{\beta,\alpha}^{x_\nu}:=\mu_{\beta,\alpha}|_{\RR_{\geqslant 0}^{\cod_\beta}\times\RR_{>0}^{d_\beta-d_\alpha}\times \{x^\nu\}}$ gives a partial compactification $\RR_{\geqslant 0}^{\cod_\beta}\times\RR_{>0}^{d_\beta-d_\alpha}\hookrightarrow\RR_{\geqslant 0}^{\cod_\beta}\times\RR_{\geqsl 0}^{d_\beta-d_\alpha}=\RR_{\geqslant 0}^{\cod_\alpha}$ that is the identity inclusion on the last $d_\beta-d_\alpha$ coordinates. 
\end{itemize}
\end{defn}

\begin{remark}\label{rmk: normal,face,tropical}
\begin{itemize}
\item[(o)] Note that we do \emph{not} require $\pi_{\alpha;\beta}^\nu$ in \eqref{eqnubetaalphapinu} to be the projection to the last $d_\alpha$ factors (equivalently, $f_{\beta,\alpha}$ does \emph{not} need to factor as $f_{\beta,\alpha}^1\times id_{\RR^{d_\alpha}}$, where $f_{\beta,\alpha}^1$ is an affine automorphism on the first $d_\beta-d_\alpha$ factors), otherwise it may put too rigid conditions on tropical manifolds.\\

\item[(i)] For any $x\in S_\alpha$, let $x^\nu=\nu_\alpha(x)\in \RR^{d_\alpha}$. By fixing an affine identification $\varphi: (\pi_{\alpha;\beta}^\nu)^{-1}(x^\nu)\overset{\sim}{\to} \RR^{d_\beta-d_\alpha}$ for each $\beta\in A_{<\alpha}$,  the fiber $\pi_\alpha^{-1}(x)$ naturally inherits a tropical manifold structure of dimension $\cod_\alpha$ with indexed set $A_{<\alpha}$: the isomorphisms $\nu_\beta^{x}: M_\beta^x:=\pi_\alpha^{-1}(x)\cap M_\beta\overset{\sim}{\to} \RR_{\geqslant 0}^{\cod_\beta}\times \RR^{d_\beta-d_\alpha}$ is given by the restriction $\nu_\beta|_{M^x_\beta}:  M_\beta^x\overset{\sim}{\to}\RR_{\geqslant 0}^{\cod_\beta}\times (\pi_{\alpha;\beta}^\nu)^{-1}(x^\nu)$ composed with $\varphi$; the restriction of $\nu_{\beta,\gamma}$, for $\beta<\gamma\in A_{<\alpha}$, gives $\nu_\gamma^x|_{M_\beta^x}\circ(\nu_\beta^x)^{-1}$, hence it satisfies the condition \eqref{eq: nu_beta,alpha,compose}. \\

\item[(ii)] Each stratum closure $\overline{S}_\beta=\bigcup_{\alpha>\beta}S_\alpha$ naturally inherits a tropical manifold structure of dimension $d_\beta$ with indexed set $A_{>\beta}$. Indeed, using that $S_\beta=\nu_\beta^{-1}(0^{\cod_\beta})\times \RR^{d_\beta}$, we set  
 $\nu_\alpha^{\beta}: M_{\alpha;\beta}:=M_\alpha\cap \ovl{S}_\beta\overset{\sim}{\to} P_{\beta,\alpha}(0^{\cod_\beta}\times\RR_{\geqslant 0}^{d_\beta-d_\alpha})\times \RR^{d_\alpha}$ to be the restriction of $\nu_\alpha$, and the collection $\{\nu_\alpha^\beta: \alpha\in A_{>\beta}\}$ gives the desired tropical manifold structure.\\

\item[(iii)] For any two tropical manifolds $M$ and $N$, the product $M\times N$ naturally inherits a product tropical manifold structure. 
\end{itemize}
\end{remark}

\begin{exam}
Let $N$ be a lattice of rank $n$, with dual lattice $M$. For any complete simplicial fan $\Sigma$ in $N_\RR=N\otimes_{\ZZ} \RR$ with primitive generators of one-dimensional cones denoted by $\{v_1,\cdots, v_m\}$, one gets a tropical manifold structure on its moment polytope $\trg_\Sigma$ associated to any very ample line bundle on the toric variety $X_{\Sigma}$. The indexed set $A$ is the collection of cones $\sigma$ in the fan (which is a subset of the power set of $\{1,\cdots, m\}$). Let $T_N=N\otimes_\ZZ \GG_m$. 

First,  
we have the moment map $\mu: X_\Sigma\to M_\RR$ which induces a homeomorphism between $(X_\Sigma)_{\geqslant}:=X_\Sigma\sslash (T_N)_\cpt$ and $\trg_\Sigma$ (cf. \cite[\S 4.2]{Fulton}). 
Second,  
for any cone $\sigma$ spanned by the subcollection of primitive vectors $v_{i_1}, \cdots, v_{i_k}$ (with $i_1<\cdots<i_k$), take a rational basis of vectors $u^\sigma_1,\cdots, u^\sigma_k, u^{\sigma}_{k+1}, \cdots, u^\sigma_{n}$ in $M_\QQ$, so that $\lng v_{i_j}, u^\sigma_{\ell}\rng=\delta_{j,\ell}$ and $M_{\sigma}':=\Span_{\ZZ}\{u_{\ell}^\sigma: 1\leqsl \ell\leqsl n\}\supset M$. This identifies the open affine subvariety $X_\sigma$, associated with $\sigma$, with $\AA^{k}\times^{Z_\sigma} \GG_m^{n-k}$ (here the quotient is the geometric quotient), where $Z_\sigma=(M'_\sigma/M)^*$. Now define $\nu_\sigma: (\trg_\Sigma)_\sigma\overset{\sim}{\to} \RR_{\geqslant 0}^{k}\times \RR^{n-k}$ by the following commutative diagram
\begin{equation*}
\begin{tikzcd}
X_\sigma\ar[r,"\sim"]\ar[d, "\mu|_{X_\sigma}"] &\AA^{k}\times^{Z_\sigma} \GG_m^{n-k}\ar[d, " {(|\ |, \log{|\ |}) }"]\\
(\trg_{\Sigma})_\sigma\ar[r,"\sim","\nu_\sigma"']& \RR_{\geqslant 0}^{k}\times \RR^{n-k}
\end{tikzcd}
\end{equation*}
One directly checks that the collection of $\{\nu_\sigma\}_{\sigma\in A}$ defines a tropical manifold structure on $\trg_\Sigma$. 

\end{exam}

\sss{Conventions on manifolds with convex/concave corners}\label{sssMWC}
 
Let $M$ be an $m$-dimensional manifold with convex/concave corners. For simplicity, we will assume that $M$ is a submanifold with corners in $\RR^{N}$, for some sufficiently large $N$, in the following sense. There exists a smooth submanifold $\wt{M}\supset \ovl{M}$ in $\RR^N$ without boundary such that $M$ is locally defined in $\wt{M}$ by the intersection and/or union of loci defined by finitely many inequalities $C^\text{sub-an}(\RR^N)\ni f_i\geqsl 0, 1\leqsl i\leqsl k$, where $C^\text{sub-an}$ means the set of subanalytic functions, and the collection $\{f_i=0\}_{1\leqsl i\leqsl k}$ (which are assumed to be smooth) and $\wt{M}$ have transverse intersections.  
The corner structure on $M$ gives a Whitney stratification $\cS_M=\{\sfS_j\}_j$ on $M$ (with $M^\circ$ being the open stratum). This will always hold for any manifold with corners that we encounter later. 

We say $V$ is a \emph{smooth vector field on $M$},  if it satisfies the following: 
\begin{itemize}
\item[(1)] $V$ is the restriction to $M$ of a smooth vector field on $\Nb(M)$ in $\RR^N$; 
\item[(2)] $V$ is stratified, i.e. the restriction of $V$ to each stratum of $M$ is tangent to the stratum. 
\end{itemize}
For any open subset $U$ in a manifold $M$ with corners, 
let $\partial^\circ U=\ovl{\partial U\cap M^\circ}$. Assume $\ovl{U}\subset M$ admits a Whitney stratification $\cS_{\ovl{U}}$ satisfying that (1) $\partial^\circ U\cap \sfS_j$, for each $j$, is the union of the strata in a subset $\cS'_{\ovl{U}}\subset \cS_{\ovl{U}}$; (2) the collection of strata in $\cS_{\ovl{U}}-\cS'_{\ovl{U}}$ consists of those nonempty $\sfS_j\cap U$ (possibly not connected; note that this requires that $U=\ovl{U}-\partial^\circ U$). 
Then for a smooth vector field $V$ defined in a neighborhood of $\partial^\circ U$, we say $V$ is \emph{convex} (resp. \emph{concave}) along $\partial^\circ U$ if
for any $\sfS_j$, the pairing of $V_x$ with $N_x^*(\sfS_j\cap U)_{\sfS_j}-\{0\}$ is strictly positive (resp. negative), for $x\in \sfS_j\cap \partial^\circ U$. Here $N_x^*(\sfS_j\cap U)_{\sfS_j}\subset T^*\sfS_j$ stands for the conormal cone to $\sfS_j\cap U$ as an open subset in $\sfS_j$ defined in \cite[Definition 5.3.6]{KS}.

Recall that a map $f: M\to N$ between stratified manifolds $(M,\cS_M)$ and $(N, \cT_N)$ is said to be \emph{stratified} if (1) for every stratum $T_i\in \cT_N$, $f^{-1}(T_i)$ is a union of strata in $\cS_M$; (2) if a stratum $\sfS_j\in \cS_M$ is mapped to $T_i$, then $f|_{\sfS_j}: \sfS_j\to T_i$ is submersive. 
We say a map $f: M\to Q$ between manifolds with corners is a \emph{trivial stratified fiber bundle}, if $f$ is a stratified map with respect to the stratification on $M$ and $Q$ (induced by the respective corner structures), and there exists $q\in Q$ and a stratified preserving homeomorphism $\phi$ making the following diagram commutative
\begin{equation*}
\begin{tikzcd}
f^{-1}(q)\times Q\ar[dr,"\pr_Q"']\ar[r, "\phi"]& M\ar[d,"f"]\\
&Q
\end{tikzcd}.
\end{equation*}

\subsection{Admissible open coverings}

Let $M$ be a tropical manifold of dimension $d$. For a subset $X\subset M$, let $\partial^\circ X:=\overline{(\partial X\cap M^\circ)}\subset M$. 

\begin{defn}[Admissible coverings]\label{def: adm cover}
Let $M$ be a tropical manifold of dimension $d$, with stratification $\{S_\alpha\}_{\alpha\in A}$ and structure maps $\nu_\alpha$. 
Let $\cV$ be an open subset of $M$ with a Whitney stratification on its boundary. We say a collection $\frU=\{(\cC_\alpha, \cI_\alpha)\}_{\alpha\in B}, B\subset A,$ is an \emph{admissible (open) covering} of $\cV$, if it satisfies the following conditions:

\begin{itemize}
\item[(1)] Each $\cC_\alpha\subset \RR_{\geqslant 0}^{\cod_\alpha}$ is a bounded contractible open neighborhood of $0$ satisfying $\cC_\alpha=\ovl{\cC}_\alpha-\partial^\circ \cC_\alpha$, whose boundary $\partial^\circ \cC_\alpha$ is a closed and contractible smooth hypersurface with corners (i.e. $(\partial^\circ \cC_\alpha)\cap \partial \RR_{\geqslant 0}^{\cod_\alpha}=\partial (\partial^\circ\cC_\alpha)$ and the corner structures are inherited from $\RR_{\geqslant 0}^{\cod_\alpha}$), and 
$\cI_\alpha\subset \RR^{d_\alpha}$ is a contractible open region with smooth boundary satisfying $\overline{\nu_{\alpha}^{-1}(0^{\cod_\alpha}\times \cI_\alpha)}\cap S_\beta=\vn, \forall \beta\gneqq\alpha$. The collection $\{\cU_\alpha:=\nu_\alpha^{-1}(\cC_\alpha\times \cI_\alpha)\}_{\alpha\in B}$ should define an open covering of $\cV$, such that $\overline{\cU}_{\alpha}\cap \overline{\cU}_\beta\neq \vn$ if and only if $\alpha>\beta$ or $\beta> \alpha$. 

\item[(2)] Assume $\alpha\gneqq \beta$. 
Then the inclusion $\ovl{\cU}_\alpha\cap \ovl{\cU}_\beta\hookrightarrow \ovl{\cC}_\beta\times\ovl{\cI}_\beta$ under $\nu_\beta$ identifies the former as a product $\ovl{\cC}_\beta\times \ovl{\cI'_{\alpha,\beta}}$, where $\cI'_{\alpha,\beta}\subset \cI_\beta$ is a contractible open subset and $\ovl{\cI'_{\alpha, \beta}}$ is a compact manifold with corners.
Moreover, we require that the natural projection $\pi_{\alpha;\beta}^\nu|_{\overline{\cI'_{\alpha,\beta}}}: \overline{\cI'_{\alpha,\beta}}\to \overline{\cI}_\alpha$ between manifolds with corners is a trivial stratified fiber bundle (cf. \S\ref{sssMWC}) with contractible fibers. 

\end{itemize}
\end{defn}

See Figure \ref{figure: triangle_wtI} for a sketch of an admissible covering of the tropical manifold $\trg$.

\begin{defn}\label{defn: C, I-adapted}
Let $M$ be a tropical manifold of dimension $d$, with stratification $\{S_\alpha\}_{\alpha\in A}$ and structure maps $\nu_\alpha$. 

\begin{itemize}
\item[] We say a closed smooth (cooriented) hypersurface $\cH$ in $\cV$ with boundary contained in $\partial\cV$ is $\cI$-\emph{adapted} (resp. \emph{$\cC$-adapted}) to a given admissible covering  $\{(\cC_\alpha, \cI_\alpha)\}_{\alpha\in B}$,  if for each $\alpha\in B$, $\vn\neq \cH\cap \cU_\alpha=\nu_\alpha^{-1}(H_\alpha\times \cI_\alpha)$, for a closed contractible hypersurface $H_\alpha\subset \cC_\alpha$ with $\partial H_\alpha\subset \partial\cC_\alpha$ that cuts out a contractible open neighborhood $\cU_{H_\a}$ of $0\in \RR_{\geqslant 0}^{\cod_\alpha}$ (resp. $\cH\cap \cU_\a=\nu_\a^{-1}(\cC_\a\times H'_\a)$ for a closed smooth hypersurface $H'_\alpha\subset \cI_\a$ with boundary in $\partial\cI_\a$). 
For any $\cH$ that is $\cC$-adapted to $\{(\cC_\a, \cI_\a)\}_{\alpha\in B}$, we say that the induced covering $\{(\cC_\a, H_\a')\}_{\a\in B}$ is an \emph{induced admissible covering} of $\cH$.

\end{itemize}

\end{defn}

\sss{}\label{ssscUflatsharp} Let $ \cU_\alpha^{\flat; A}:=\cU_\alpha-\bigcup_{\beta\in A_{\gneqq \alpha}}\cU_{\beta}$ and $\cU_\alpha^{\sharp;A}:=\cU_\alpha-\bigcup_{\beta\in A_{\lneqq \alpha}}\cU_{\beta}$ (cf. Figure \ref{figure: triangle_wtI}). We also denote these simply by $\cU_\alpha^{\flat}$ and $\cU_\alpha^\sharp$, respectively, when the poset $A$ is clear from the context. Let $\cU_\alpha^\nu:=\cC_\alpha\times \cI_\alpha$

\begin{remark}\label{remark: C_gamma,alpha}
\begin{itemize}
\item[(o)] Not every open subset $\cV$ of $M$ admits an admissible open covering. It is probably more natural to first define $\{(\cC_\alpha, \cI_\alpha)\}_{\alpha\in B}$ satisfying the conditions (1) and (2) above and then set $\cV=\bigcup_{\alpha\in B}\cU_\alpha$. \\

\item[(i)]
Definition \ref{def: adm cover} condition (2) implies that for any $\alpha\in B$, $\nu_\alpha(\cU^\flat_\alpha)=\cC_\alpha\times \cI^\flat_\alpha$ for some $ \cI^\flat_\alpha\subset \cI_\alpha$. Moreover, it implies that  
we have the commutative diagram (and a similar one by taking the closure of each term)
\begin{equation*}
\begin{tikzcd}[column sep=7em]
\cC_\beta\times \cI'_{\alpha, \beta}\ar[r, hook, "\nu_{\beta, \alpha}|_{\cC_\beta\times \cI'_{\alpha, \beta}}"]\ar[dr, "\pi_{\alpha;\beta}^\nu|_{\cI'_{\alpha,\beta}}\circ\pr_2"']&\cC_{\alpha}\times \cI_\alpha\ar[d, "\pr_2"]\\
\ &\cI_\alpha
\end{tikzcd}.
\end{equation*}
Consequently, given an admissible covering $\frU=\{(\cC_\alpha,\cI_\alpha)\}_{\alpha\in B}$ and any $\alpha\in B$ with $A_{<\alpha}\subset B$, for any tropical manifold structure specified on $\pi_\alpha^{-1}(x)$ as in Remark \ref{rmk: normal,face,tropical}, with $x^\nu\in \ovl{\cI}_\alpha$, there exists an open neighborhood $\cV_{\alpha}$ of $(\nu_\alpha^x)^{-1}(\ovl{\cC}_\alpha)$ in $\pi_\alpha^{-1}(x)$, that naturally inherits an admissible covering $\{(\cC_{\gamma}, \cI^x_{\gamma, \alpha})\}_{\gamma\in B_{<\alpha}=A_{<\alpha}}$, where $\cI^x_{\gamma, \alpha}:=(\pi_{\alpha;\gamma}|_{\cI'_{\alpha,\gamma}})^{-1}(x^\nu)$. In particular, each $\partial^\circ \cC_\alpha$ is a hypersurface that is $\cC$-adapted to the admissible covering $\frU_{\alpha}^x:=\{(\cC_{\gamma}, \cI^x_{\gamma, \alpha})\}_{\gamma\in B_{\lneqq\alpha}}$ of their union. More generally, even when $A_{<\alpha}\not\subset B$, we can still define the admissible covering $\frU_{\alpha}^x:=\{(\cC_{\gamma}, \cI^x_{\gamma, \alpha})\}_{\gamma\in B_{\lneqq\alpha}}$ of their union.
\\

\item[(ii)] We define a partial order on the set of admissible coverings of a fixed open subset $\cV$, denoted by $\cP_\cV$:  $\{(\cC_\alpha, \cI_\alpha)\}_{\alpha\in B}\prec\{(\cC'_\alpha, \cI'_\alpha)\}_{\alpha\in B}$ if and only if $\cC_\alpha\supset \cC'_\alpha$ and $\cI_\alpha\subset \cI'_\alpha$. We say $\{(\cC'_\alpha, \cI'_\alpha)\}_{\alpha\in B}$ is \emph{thinner and wider} than $\{(\cC_\alpha, \cI_\alpha)\}_{\alpha\in B}$.  
More generally, we can define a partial order on the set of admissible coverings of possibly different $\cV$, by the same rule. \\

\item[(iii)] Definition \ref{def: adm cover} condition (2) implies that $\nu_\alpha^{-1}(\partial^\circ\cC_\alpha\times\ovl{\cI}_\alpha)\cap \cV\subset \bigcup_{\beta\in B_{\lneqq \alpha}}\cU_\beta$.

\end{itemize}

\end{remark}

Let $I':=\{1,\cdots, d\}$. Equip $\RR_{\geqslant 0}^{I'}$ with a tropical manifold structure as follows. Let $A=\cP(I')$ and $S_J=\{r_i=0, i\in J\}, J\in A$, be the stratification. For each $J$ fix an identification $\nu_J: M_J\overset{\sim}{\to}\RR_{\geqslant 0}^{J}\times \RR^{I'-J}$ with $\nu_{I'}=id_{\RR_{\geqslant 0}^{I'}}$, such that 
\begin{align*}
\nu_{J, I'}: \RR_{\geqslant 0}^{J}\times \RR^{I'-J}\hookrightarrow \RR^{I'}_{\geqslant 0}
\end{align*}
is given by the composition in the form of (\ref{eq: nu_beta,alpha,compose}). Here $J$ and $I'-J$ have ordering inherited from $I'$.

\begin{defn}\label{def: H0 good}
\item[(i)]
Let $H_0\subset \RR_{\geqslant 0}^{I'}$ be a closed contractible (cooriented) smooth hypersurface with corners lying in $\partial \RR^{I'}_{\geqslant 0}$ such that it cuts out a contractible neighborhood of $0\in \RR^{I'}_{\geqslant 0}$. Let $(r_\beta)_{\beta\in I'}$ denote the standard coordinates on $\RR_{\geqslant 0}^{I'}$. For any collection of positive integers (or more generally positive real numbers) $\bm=(m_\beta)_{\beta\in I'}$, let $\sfZ_{\bm}$ be the vector field on $\RR^{I'}_{\geqslant 0}$ such that $\sfZ_{\bm}(r_\beta)=-m_\beta r_\beta$. 
We say $H_0$ is  \emph{good} with respect to $\sfZ_{\bm}$ and an admissible covering $\frU=\{(\cC_J, \cI_J)\}_{J\in A_{\lneqq I'}}$  of a neighborhood of $H_0$,  if it is $\cC$-adapted to $\frU$ and $\sfZ_{\bm}$ is transverse to $H_0$ everywhere. We simply say  $H_0$ is  \emph{good} when $\sfZ_{\bm}$ is clear from the context and $\frU=\{(\cC_J, \cI_J)\}_{J\in A_{\lneqq I'}}$ exists but is not specified. 
\\

\item[(ii)] We say a function $f$ on $\cV$ is \emph{$\cC$-adapted to} an admissible covering $\frU=\{(\cC_\alpha, \cI_\alpha)\}_{\alpha\in B}$, if for all $\alpha\in B$, $f|_{\cU_\alpha^\sharp}$, is the pullback of a function on $\cI_\alpha$ for the chart defined by $\nu_\alpha$. Similarly, we can define the notion of $\cC$-adaptedness for a function $f$ on a $\cC$-adapted hypersurface $\cH$ in $\cV$ with respect to the induced admissible covering on it. 
\end{defn}

\begin{remark}
Note the condition that $H_0\subset \RR_{\geqslant 0}^{I'}$ is $\cC$-adapted to an admissible covering $\frU$ is independent of the choice of the tropical structure $\nu_J$ on $\RR_{\geqslant 0}^{I'}$, since this is equivalent to that $H_0$ is locally defined by $f_J(r_i, i\not\in J)=0$ in a neighborhood of $S_J$. So is the condition of good on $H_0$. 
\end{remark}

\begin{lemma}\label{lemmaAfrU12sharpflat}
Assume $M$ is compact. For any admissible open coverings $\frU^1\succ \frU^2$, we have
\begin{itemize}
\item[(i)] $(\cU_\alpha^1)^\sharp\subset \bigcup_{\alpha'>\alpha}(\cU_\alpha^2)^{\sharp}$.

\item[(ii)] $(\cU_\alpha^2)^\flat\subset \bigcup_{\beta < \alpha}(\cU_\beta^1)^\flat$. 

\item[(ii)] $\cU^2_\alpha\subset \bigcup_{\beta < \alpha}\cU_\beta^1$. 
\end{itemize}
\end{lemma}

\begin{proof}
This is straightforward to check. 
\end{proof}

\begin{lemma}\label{lemma: admissible, cI contractible}
Let $M$ be a tropical manifold of dimension $d$. Let  $\{(\cC_\alpha, \cI_\alpha)\}_{\alpha\in B}, B\subset A$, be any admissible covering of the union $\cV=\bigcup_{\alpha\in B} \cU_\alpha$.
\begin{itemize}
\item[(i)] Then for any $\vec{\beta}=(\beta^1\gneqq \cdots \gneqq \beta^k)$, with $\beta^i\in B$, the intersection $
\cU_{\vec{\beta}}:=\bigcap_{1\leqsl i\leqsl k}\cU_{\beta^i}$ satisfies that 
\begin{align}\label{eq: I beta1,k}
\nu_{\beta^k}(\cU_{\vec{\beta}})=\cC_{\beta^k}\times \cI_{\vec{\beta}}, 
\end{align}
where $ \cI_{\vec{\beta}}\subset \cI_{\beta^k}$ is a contractible open subset.  

\item[(ii)] Assume $B\subsetneq A$, and let $\alpha\not\in  B$. Assume $B$ contains $A_{\gneqq\alpha}$. Then $\cV\cap \overline{S}_\alpha$ gives a neighborhood of $\partial S_\alpha$ that deformation retracts on to it, and it naturally inherits an admissible covering indexed by $A_{\gneqq\alpha}$. 
\end{itemize}
\end{lemma}

\begin{proof}
(i)
We prove by induction on $k$. First, the case that $k=1$ is trivial. Suppose we have proved the case for $k\leqsl \ell$. For $k=\ell+1$, first by induction we have 
\begin{align*}
\cU^{\nu_{\ell}}_{\vec{\beta}, \leqsl \ell}:=\nu_{\beta_{\ell}}(\cU_{\beta^1,\cdots, \beta^{\ell}})=\cC_{\beta^{\ell}}\times \cI_{\vec{\beta}, \leqsl \ell}
\end{align*}
for some contractible open $\cI_{\vec{\beta}, \leqsl \ell}:=\cI_{\beta^1,\cdots, \beta^{\ell}}\subset \cI_{\beta^{\ell}}$. Let $\nu_{\ell}:=\nu_{\beta^{\ell}}$ and $\cU_{\vec{\beta}, \leqsl \ell}=\nu_{\ell}^{-1}(\cU^{\nu_{\ell}}_{\vec{\beta}, \leqsl \ell})$. 
Using condition (2) in Definition \ref{def: tropicalmfld}, let 
\begin{align*}
\cI_{\vec{\beta}}:=(\pi_{\beta^{\ell+1}, \beta^{\ell}}^\nu|_{\cI_{\beta^{\ell}, \beta^{\ell+1}}})^{-1}(\cI_{\vec{\beta}, \leqsl \ell}).
\end{align*}
Then by Remark \ref{remark: C_gamma,alpha} (i), we directly get 
\begin{align*}
\nu_{\ell+1}(\cU_{\vec{\beta}})=\nu_{\ell+1}\left(\cU_{\vec{\beta}, \leqsl \ell}\cap \cU_{\beta^{\ell+1}}\right)=\cC_{\beta^{\ell+1}}\times \cI_{\vec{\beta}}.  
\end{align*}

(ii) Using Remark \ref{rmk: normal,face,tropical}, it is clear that $\{(\cC_\beta,\cI_\beta), \beta\in A_{\gneqq\alpha}\}$ induces an admissible covering of $\cV\cap \overline{S}_\alpha$ by taking the intersections. Then the statement is a direct consequence from applying (i) to this admissible covering.   

\end{proof}

\begin{lemma}\label{lemma: exist admissible}
For any \emph{compact} tropical manifold $M$ with a choice of $\bm_\alpha, \alpha\in A$, there exist sufficiently thin and wide (notion as in Remark \ref{remark: C_gamma,alpha} (ii)) admissible open coverings $\frU=\{(\cC_\alpha, \cI_\alpha)\}_{\alpha\in A}$ of $M$ such that  
\begin{itemize}
\item[(1)] for all $\alpha\in A$, the hypersurface $\partial^\circ\cC_\alpha\subset \RR_{\geqslant 0}^{\cod_\alpha}$ is good with respect to $\sfZ_{\bm_\alpha}$ and the induced admissible covering $\frU_{\alpha}^x$ for any $x^\nu\in \ovl{\cI}_\alpha$ (cf. Remark \ref{remark: C_gamma,alpha} (i));

\item[(2)] for any $\vec{\beta}=(\beta^1\gneqq \cdots \gneqq \beta^k), \beta^i\in A$, the open contractible subset $\cI_{\vec{\beta}}\subset \cI_{\beta^k}$ from (\ref{eq: I beta1,k}) satisfies that it is cut out by a collection of transverse intersecting (locally closed) smooth hypersurfaces $\{H^{k}_{\partial^\circ\cC, i}: 1\leqsl i\leqsl k-1\}\cup\{H^{k}_{\partial\cI, i}: 1+\delta_{d_{\beta^1}, 0}\leqsl i\leqsl k\}$ in $\RR^{d_{\beta^k}}$, such that there exist small open neighborhoods $\Nb(\overline{\cI}_{\beta^i})\subset \RR^{d_{\beta^i}}, \Nb(\overline{\cC}_{\beta^i})\subset \RR^{\cod_{\beta^i}}_{\geqslant 0}, 1\leqsl i\leqsl k$, with 
\begin{align*}
&\nu_{\beta^k, \beta^i}(H_{\partial^\circ\cC, i}^{k})\subset \left(\partial^\circ\cC_{\beta^i}\times \Nb(\ovl{\cI}_{\beta^i})\right)\cap \left(0^{\cod_{d_{\beta^k}}}\times \RR_{> 0}^{d_{\beta^k}-d_{\beta^i}}\times \RR^{d_{\beta^i}}\right),\\
&\nu_{\beta^k, \beta^i}(H_{\partial\cI, i}^{k})\subset \left(\Nb(\ovl{\cC}_{\beta^i})\times \partial\cI_{\beta^i}\right)\cap \left(0^{\cod_{d_{\beta^k}}}\times \RR_{> 0}^{d_{\beta^k}-d_{\beta^i}}\times \RR^{d_{\beta^i}}\right),
\end{align*}
for all $i\leqsl k$. Here the hypersurfaces are all connected except for $H_{\partial\cI, 1}^{k}$ when $d_{\beta^1}=1$. 
Moreover, the connected components of $\partial \cI_{\vec{\beta}}\cap H^{k}_{\partial^\circ\cC, i}$ and $\partial \cI_{\vec{\beta}}\cap H^{k}_{\partial\cI, i}$ are contractible for all $i$.

\end{itemize}

\end{lemma}

\begin{proof}
First, we assume $M=\RR^{I'}_{\geqslant 0}$, $I'=\{1,\cdots, d\}$ (noncompact for now) with standard coordinates $r_1,\cdots, r_d\geqslant 0$, and $m_\beta=1$ for all $\beta\in I'$.
Let us define $(\cC_{I'}, \cI_{I'}=\{0\})$ as follows. Choose a closed hypersurface $H_{I'}$ that is $C^1$-close to a round sphere of radius $\e>0$ such that for each $J\subsetneq I'$, there exists a neighborhood $\Nb_J(0)$ of $0\in \RR_{\geqslant 0}^J$ such that in $\Nb_J(0)\times \RR^{I'-J}_{> 0}\subset \RR_{\geqslant 0}^{I'}$, $H_{I'}$ is defined by an equation of the form $f_J(r_i, i\not\in J)=0$. By a $C^1$-small perturbation of each $f_J$, we have $H_{I'}$ is good with respect to $\sfZ_{(1,\cdots,1)}$. For a general $\bm$, using that the stratified diffeomorphism 
\begin{align*}
&\RR^{I'}_{\geqsl 0}\overset{\sim}{\to} \RR^{I'}_{\geqsl 0},\quad (r_\beta)_\beta\mapsto (r_\beta^{m_\beta})_{\beta}
\end{align*}
sends $\sfZ_{(1,\cdots, 1)}$ to $\sfZ_{\bm}$ and it preserves admissible open coverings and $\cC$-adaptedness, the image of $H_{I'}$ just defined is good with repsect to $\sfZ_{\bm}$. 
Let $\cC_{I'}$ be the neighborhood of $0$ in $M$ enclosed by $H_{I'}$. 

Now let $M$ be any compact tropical manifold of dimension $d$. We define an admissible covering by induction on $d_\alpha$ (from low to high). First, for $d_\alpha=0$, using $\nu_\alpha: M_\alpha\overset{\sim}{\to} \RR_{\geqslant 0}^{d}$, we are reduced to the above situation and  define $(\cC_\alpha, \cI_\alpha=\{0\})$ as above with $\e$ sufficiently small. Assume we have chosen $\{(\cC_\alpha, \cI_\alpha)\}_{\alpha\in A_{\geqslant d-\ell}}$, where $A_{\geqslant d-\ell}=\{\alpha\in A: \cod_\alpha\geqslant d-\ell\}$, which forms an admissible open covering of the union $\cV_{\geqslant d-\ell}:=\bigcup_{\alpha\in A_{\geqslant d-\ell}}\cU_\alpha$, that satisfies (1) and (2) for $A$ replaced by $A_{\geqslant d-\ell}$. 

Now let $\beta\in A_{d-(\ell+1)}=A_{\geqslant d-(\ell+1)}-A_{\geqslant d-\ell}$, and let $\cU'_{0,\beta}\subset S_\beta$ be a sufficiently large contractible open subset with smooth boundary, such that $\cU'_{0,\beta}\cup (\cV_{\geqslant d-\ell}\cap \overline{S}_\beta)=\overline{S}_\beta$. 
Then on $\overline{S}_\beta$, we have 
\begin{itemize}
\item $\cV_{\geqslant d-\ell}\cap \overline{S}_\beta$ naturally inherits an admissible covering $\{(\cC_{\beta,0;\alpha}, \cI_\alpha)\}_{\alpha\in A_{\gneqq\beta}}$ that satisfies the conditions (1) and (2) for $A$ replaced by $A_{\gneqq \beta}$ (cf. Lemma \ref{lemma: admissible, cI contractible} (ii));\\

\item Let $\wt{h}_\beta: \overline{S}_\beta\to \RR_{\geqslant 0}$ be a subanalytic defining function of $\partial \overline{S}_\beta=\bigcup_{\alpha\in A_{\gneqq \beta}}S_\alpha$, i.e. $\wt{h}_\beta(x)=0$ if and only if $x\in \partial \overline{S}_\beta$. Let $h_\beta: \RR^{d_\beta}\to \RR_{>0}$ be the composition $\wt{h}_\beta|_{S_\beta}\circ \nu_{\beta}^{-1}|_{\{0\}\times  \RR^{d_\beta}}$. Then $\wt{h}_\beta|_{\cU_\alpha\cap S_\beta}$, for $\alpha\in A_{\gneqq \beta}$, is equal to the pullback under $\nu_\alpha$ of the restriction of a  subanalytic function $\eta_{\beta,\alpha}$ on $\RR^{\cod_\alpha}_{\geqslant 0}\times \RR^{d_\alpha}$ that vanishes exactly on $\RR_{\geqslant 0}^{\cod_\beta}\times \partial\RR^{d_\beta-d_\alpha}_{\geqslant 0}\times \RR^{d_\alpha}$ (recall the last sentence in Definition \ref{def: tropicalmfld}). 

\end{itemize}
For $\ep>0$ sufficiently small, let $\cI_\beta=h_\beta^{-1}\left((\ep, \infty)\right)$ and $H_{\partial\cI, \beta}=h_\beta^{-1}(\ep)$, where the latter is smooth. Then $\{(\cC_{\beta,0;\alpha}, \cI_\alpha)\}_{\alpha\in A_{\gneqq\beta}}\cup\{(\{0\}, \cI_\beta)\}$ gives an admissible covering of $\overline{S}_\beta$ satisfying the conditions (1) and (2). 
By choosing a sufficiently small $\cC_{\beta}\subset \RR_{\geqslant 0}^{d-d_\beta}$ with good $\partial^\circ \cC_\beta$, we have $\{(\cC_\alpha, \cI_\alpha)\}_{\alpha\in A_{\geqslant d-\ell}}\cup\{(\cC_\beta, \cI_\beta)\}_{\beta\in A_{d-\ell+1}}$ form an admissible covering of the union $\cV_{\geqslant d-\ell}\cup \cU_\beta$ (note that here we need the compactness of $M$). This finishes the inductive step for defining $\frU$.

The collection of smooth transverse hypersurfaces in (2) also follows from induction. 
We omit the details since the argument is straightforward.  
\end{proof}

A small remark is that Lemma \ref{lemma: exist admissible} (2) does \emph{not} follow directly from the definition of an admissible covering since the hypersurfaces in (2) need \emph{not} all be transverse. It is natural to add Lemma \ref{lemma: exist admissible} (2) into the definition of an admissible covering, but we prefer to make the definition shorter.

\begin{defn}\label{defngoodadm}
For any $\alpha\in A$, fix a $\cod_\alpha$-tuple of positive integers $\bm_\alpha$. Let $\sfZ_{\bm_\alpha}$ be the vector field on $\RR^{\cod_\alpha}_{\geqsl 0}$ defined as above. We say an admissible covering $\{(\cC_\alpha, \cI_\alpha): \alpha\in A\}$ of $M$  is \emph{good} with respect to $\sfZ_{\bm_\alpha}, \alpha\in A$, if it satisfies the conditions in Lemma \ref{lemma: exist admissible}.

More generally, for an admissible covering $\{(\cC_\alpha, \cI_\alpha): \alpha\in B\}, B\subset A$, of $\cV$, we say it is \emph{good} with respect to $\sfZ_{\bm_\alpha}, \alpha\in B$, if it satisfies 
 the same condition in Lemma \ref{lemma: exist admissible} (2) with $A$ replaced by $B$ and the following replacement of (1) \\
 (1') for any $\alpha\in B$ and $x^\nu\in \ovl{\cI}_\alpha$, $\partial^\circ\cC_\alpha\subset \RR_{\geqslant 0}^{\cod_\alpha}$ is good with respect to $\sfZ_{\bm_\alpha}$ and 
 some admissible covering of a neighborhood of $\partial^\circ\cC_\alpha$ that contains the induced admissible covering $\frU_{\alpha}^x$ as the subcollection indexed by $B_{\lneqq\alpha}$ (cf. Remark \ref{remark: C_gamma,alpha} (i)). 
\end{defn}

A direct consequence of Lemma \ref{lemma: exist admissible} and its proof is the following.  
\begin{cor}
For any \emph{compact} tropical manifold $M$ (with a choice of $\bm_\alpha, \alpha\in A$), the poset of all (good) admissible coverings is a filtered poset, and there is a cofinal sequence.  
\end{cor}

\section{Backgrounds on Weinstein sectors and wrapped Fukaya categories}\label{secAppSectors}
For any Liouville manifold/sector $M$, let $\sfZ$ denote the Liouville vector field.  Let $\partial M$ be the finite boundary of $M$ and let $M^\circ=M-\partial M$. For any submanifold $Y\subset M$ with $\partial Y\subset \partial M$ (i.e. $Y\cap M^\circ$ is properly embedded) and $Y$ is transverse to $\partial M$, that is invariant under the Liouville flow near $\partial^\infty M$, we call $Y$ a \emph{cylindrical submanifold of $M$}. For any cylindrical submanifold $Y$, let $\Nb^\sfZ(Y)$ denote a cylindrical open neighborhood of $Y$.

\subsection{Sectorial coverings and descent of wrapped Fukaya categories}
In this section, we recall the definition of a Liouville/Weinstein sectorial covering and the descent property of wrapped Fukaya categories from \cite[\S 12.1]{GPS2}. 
\begin{defn}
Let $(M, \omega=d\vartheta)$ be a Liouville (resp. Weinstein up to deformations) sector. 
\begin{itemize}
\item[(i)]
A finite collection of cylindrical hypersurfaces $H_1,\cdots, H_n\subset M$ is called \emph{sectorial} if (1) their characteristic foliations are $\omega$-orthogonal over their intersections; (2) there exists a function $I_i: \Nb^{\sfZ}(H_i)\to \RR$ that is linear near infinity (with respect to $\sfZ$) such that $dI_i|_{\text{char. fol. of }H_i}\neq 0$ (equivalently, the Hamiltonian vector field $X_{I_i}$ of $I_i$ satisfies that $X_{I_i}\not\in TH_i$), and $dI_i|_{\text{char. fol. of }H_j}=0, \forall j\neq i$ (equivalently, $X_{I_i}\in TH_{j}, \forall j\neq i$). Here $``\text{char. fol. of }H_i"$ is the shorthand of characteristic foliation of $H_i$. 

\item[(ii)] Let $M_i, i=1,\cdots, n,$ be a cover of $M$. If each $M_i$ is a Liouville subsector and the cylindrical hypersurfaces $H_i:=\overline{\partial M_i\cap M^\circ}, i=1,\cdots, n$, together with $\partial M$, are sectorial, then we say $M_i, i=1,\cdots, n$, form a Liouville sectorial covering. 

\end{itemize}
\end{defn}

\begin{remark}
For any sectorial collection of cylindrical hypersurfaces $H_1,\cdots, H_n\subset M$, every finite intersection among them is transverse, and $\bigcap_{j\in J\subset\{1,\cdots,n\}}H_j$ is coisotropic. The latter has a well defined symplectic reduction given by $\{I_j=0, j\in J\}$, which is a Liouville sector with corners (cf. \cite[Lemma 12.8 and 12.10]{GPS2}). 
\end{remark}

For any sectorial covering $M_i, 1\leq i\leq n$, of $M$, stratify $M$ by 
\begin{align*}
M_{A, B, C}:=\bigcap_{a\in A}M_a\cap (\bigcap_{b\in B} \partial M_b)\backslash (\bigcup_{c\in C} M_c)
\end{align*}
for all partitions $A\sqcup B\sqcup C=\{1,\cdots, n\}$. 
Then the closure of each $M_{A, B,C}$ has symplectic reduction a Liouville sector-with-corners. 

\begin{defn}\label{defnWeinsteinsectorcovering}
We say a sectorial covering is \emph{Weinstein} if the convexification of each of the symplectic reduction above is Weinstein (up to  deformations). 
\end{defn}

\begin{exam}\label{exam: Qi, WeinsteinSectorial}
Let $Q$ be a smooth compact manifold-with-boundary. Let $Q_1,\cdots, Q_n\subset Q$ be codimension 0 submanifolds-with-boundary such that $\bigcup_{1\leqsl i\leqsl n}Q_i=Q$, and $\partial^\circ Q_1:=\overline{\partial Q_1\cap Q^\circ}, \cdots, \partial^\circ Q_n:=\overline{\partial Q_n\cap Q^\circ}, \partial Q$ are mutually transverse, then $T^*Q_1, \cdots, T^*Q_n$ is a Weinstein sectorial covering of $T^*Q$. 
\end{exam}

\begin{theorem}[Theorem 1.35 \cite{GPS2}]\label{thm: sectorialcover,cW}
For any Weinstein sectorial covering $M_1,\cdots, M_n$ of a Liouville sector $M$, the induced functor 
\begin{align*}
\colim_{\vn\neq J\subset\{1,\cdots, n\}}\cW(\bigcap_{j\in J}M_j)\overset{\sim}{\to}\cW(M)
\end{align*}
is an equivalence of categories. 
\end{theorem}

\subsection{The contact hypersurface $\cH$ and Weinstein handle(s) in $\cJ_G$}

From this subsection till \S\ref{subsecSectdeformJG}, we review and sketch the main constructions and results in establishing open (generalized) Weinstein sector structure on $\cJ_G$ from \cite{J} and discuss deformations of Liouville sector structures on $\cJ_G$ for later use in \S\ref{subsecSectInclu}. Let $\vartheta_{\cJ_G}$ be the canonical real Liouville 1-form on $\cJ_G$. If $G$ is clear from the context, we simply write $\vartheta_\cJ$ for $\vartheta_{\cJ_G}$. We will assume $G$ is semisimple with a complete set of simple roots $I$, unless otherwise specified. 
All the discussions directly apply to $\cT_G$ by choosing any identification $\cT_G\cong \cJ_G$ (which will be clearly independent of the choices). 

\sss{The critical Weinstein handle(s) in $\cJ_G$}
Consider $\RR^I_{\geqsl 0}$ as a tropical manifold, with stratification $\{0^{J}\times \RR_{>0}^{I-J}, J\subset I\}$ and $\nu_I=id_{\RR_{\geqsl 0}^I}$ and $\nu_{J,I}$ given by  
\begin{align*}
\nu_{J, I}: \RR_{\geqsl 0}^J\times \frz_{J;\RR} &\hookrightarrow \RR_{\geqsl 0}^I\\
(r_j,j\in J; t_{J,\perp})&\mapsto (e^{\lng \omega_j,t_{J,\perp}\rng}r_j, j\in J; e^{\lng \omega_i,t_{J,\perp}\rng}, i\in I-J).
\end{align*}
These determine $\nu_J, J\subset I$, and this tropical manifold structure is the same as the induced structure on $\trg_{\subset I}$ from $\trg$ in \S\ref{sss: set-up, triangle}.

Choose any open $0\in \cC\subset \RR_{\geqsl 0}^I$ such that $H^{sm}:=\partial^\circ C$ is good with respect to $\sfZ_{r, I}$ and some admissible open covering $\frU^{\PdI}:=\{(\cC_J, \cI_J)\}_{J\in \PdI}$ of a neighborhood of $H^{sm}$ (cf. Definition \ref{def: H0 good} (i)). We call such an $H^{sm}$ an \emph{admissible hypersurface} in $\RR_{\geqsl0}^{I}$. It is straightforward to see that the space of admissible hypersurfaces is contractible. 
Then 
\begin{align*}
\cH:=(b_\RR^I)^{-1}(H^{sm}) 
\end{align*}
is a contact hypersurface in $\cJ_G$. Here one defines $b_\RR^I: \cJ_G\to \RR_{\geqsl 0}^{I}$ by choosing any identification $\cJ_G\cong \cT_G$, which is independent of the choice. 
Let $\sfN_\cC: \RR^{I}_{\geqsl 0}\to [0,\infty)$ be the function satisfying $\sfN_\cC^{-1}(1)=H^{sm}$, $\sfN_C$ is homogeneous of weight $-\frac{1}{2}$ with respect to $\sfZ_{r, I}$, and $\sfN_\cC(0)=0$. Let $\sfN_\cJ=\sfN\circ b_\RR^I$.
Let $\|\cdot\|$ be any norm on $\frc$, i.e. $\|\cdot\|:\frc\to [0,\infty)$ is a proper function with a unique zero at $[0]\in \frc$ that is homogeneous of some sufficiently large positive weight with respect to the standard $\CC^\times$-action. For any $\ep>0$, set 
\begin{align}\label{eqfrDepC}
&\frD_{\ep,\ovl{\cC}}:=\frD_{\ep, \ovl{\cC}}^G:=(\|\chi\|\times b_I^{\RR})^{-1}([0,\ep]\times \ovl{\cC})
\end{align}
Let $\sfZ$ be the standard Liouville vector field on $\cJ_G$. Then $\sfZ$ is weakly gradient-like with respect to the function $\|\chi\|^{2k}-(\sfN_\cJ)^k$, for some large $k\in \ZZ_{\geqsl 0}$.

\begin{lemma}\label{lemmafrDepChandle}
Under the above assumptions, for $\ep>0$ sufficiently small, we have each connected component of $\frD_{\ep, \ovl{\cC}}$ (indexed by $z\in Z(G)$) a critical Weinstein handle of index $2r$. Moreover, $\cJ_G$ is the completion of $(b_I^{\RR})^{-1}(\RR_{\geqsl 0}^I-\cC)\cup \frD_{\ep,\ovl{\cC}}$. 
\end{lemma}

\begin{warn}
Although both $b^I_{\RR}$ and $\|\chi\|$ are homogeneous maps, the first part of Lemma \ref{lemmafrDepChandle} is false for large $\ep$ in general. This is due to either of (1) the map $\chi|_{\cH}:\cH\to \frc$ has critical values outside $[0]\in \frc$ since $\cH$ is \emph{not} homotopy equivalent to $T^*S^{2n-1}\times \RR$ unless $n=1$; (2) for $Z(G)\neq 1$, $\frD_{\ep, \ovl{\cC}}$ has $|Z(G)|$ many connected components for small $\ep$, but it is connected for large $\ep$. 
\end{warn}

\subsection{A Weinstein hypersurface $\frF$ in $\cH$}\label{subsec: Hsm}

\sss{Admissible covering of $H^{sm}$}
Given any good admissible covering $\frU_{\PdI}=\{(\cC_J, \cI_J)\}_{J\subsetneq I}$ of an open neighborhood of $H^{sm}$, with respect to which $H^{sm}$ is $\cC$-adapted, let $\frU_{H}:=\{(\cC_J, H^{\cI}_J)\}_{J\subsetneq I}$ denote the induced good admissible covering on $H^{sm}$ as in Definition \ref{defn: C, I-adapted}. Let $\cU_J^{H}=\nu_J^{-1}(\cC_J\times H_J^\cI)=H^{sm}\cap \cU_J$, $H_J^{sm}=H^{sm}\cap (0^J\times \RR_{>0}^{I-J})$ and $H_J^{\frz}=\nu_J(H_J^{sm})\subset \frz_{J;\RR}$, for $J\in \PdI$. 

\sss{Admissible vector fields $\sfZ^{sm}$ on $H^{sm}$}

The following definition is the analog of Definition \ref{defnZadm} in which the datum ($\trg,\frU$) is replaced by $(H^{sm},\frU_{H})$. 

\begin{defn}\label{defnZsmadm}
Let $\delta>0$. Let $\sfZ^{sm}$ be a vector field on $H^{sm}$ (as a manifold with corners). 
Let $\frU_{H}=\{(\cC_J, H_J^{\cI})\}_{J\in \PdI}$ be a good admissible covering. Let $\sfZ^{sm, \frz}_J$ denote $(\nu_J)_*(\sfZ^{sm}|_{H^{sm}_J})$. 
We say $\sfZ_{\flat}$ is \emph{strongly $(\delta,\frU_{H})$-admissible} if the following hold
\begin{itemize}
\item[(ZHAa)]  for each $J\in \cP^\dagg(I)$, 
\begin{align*}
(\nu_J)_*(\sfZ^{sm}|_{\cU_J^H})=\sfZ_{r,J}\times \sfZ_{J}^{sm, \frz},
\end{align*}
with respect to the splitting $\nu_J(\cU_J^H)=\cC_J\times H_J^\cI$, where $\sfZ_{J}^{sm, \frz}$ is convex along $\partial H_J^\cI$ (cf. \S\ref{sssMWC}) and  $\sfZ_{J}^{sm, \frz}$ has no zero outside $H_J^\cI$;

\item[(ZHAb)] The function $\sfI^H_J:=\vartheta_{T^*H_J^\frz}(\sfZ^{sm, \frz}_J)$ satisfies that 
\begin{align*}
&|X_{\sfI^H_J}(\|\bp\|^2)|\leqsl \delta \|\bp\|^2
\end{align*}
with respect to the induced metric on $H_J^{\frz}$ by the Killing form; 

\item[(ZHAc)] $\sfZ_{J}^{sm, \frz}$ is gradient-like with respect to some exhausting Morse function $f_{J, H}^{\frz}: H_J^\frz\to [0,\infty)$. 
\end{itemize}
We say $\sfZ^{sm}$ is \emph{$(\delta,\frU_{H})$-admissible} if it satisfies (ZHAa) and (ZHAb). We also say $\sfZ^{sm}$ is \emph{$\delta$-admissible} (resp. \emph{strongly $\delta$-admissible}) if it is $(\delta,\frU_{H})$-admissible (resp. strongly $(\delta, \frU_{H})$-admissible) for some good admissible covering $\frU_{H}$. 
\end{defn}

We also have the notion of $(\delta, \frU_H')$-semi-admissibility and (weak) $\frU_H$-rigidity on $\sfZ^{sm}$. This is an 
analog of Definition \ref{defnZsemiadm}, in which the projection \eqref{eqfrzJRprojJJprime} in (ZRii) is replaced by 
\begin{align*}
\frz_{J;\RR}^{J'}\times H_{J'}^{\frz}\longrightarrow \frz_{J;\RR}^{J'}.
\end{align*}
Other than that, there is no change for the definition and we omit the details. 
Given a $\delta$-admissible (resp. strongly $\delta$-admissible) $\sfZ^{sm}$ on $\Hsm$, let $\CI_{\sfZ^{sm},\delta}$ (resp. $\CI_{\sfZ^{sm},\delta}^{\sigma}$) be the collection of good admissible coverings $\frU_H$ with respect to which $\sfZ^{sm}$ is $(\delta,\frU_H)$-admissible (resp. strongly $(\delta,\frU_H)$-admissible).

\sss{Existence results about admissible vector fields on $H^{sm}$}

We have the analog of Lemma \ref{lemmaadmZ} and Corollary \ref{corZflatCI} with essentially the same proofs. 

\begin{lemma}\label{lemmaadmZsm}
\begin{itemize}
\item[(i)]
Let $\delta>0$. For any sufficiently thin and wide good admissible open covering $\frU^1_H=\{(\cC_J^1, H_J^{\cI,1})\}_{J\in \PdI}$, there exists a good $\frU_{H}\succ \frU^1_H$ and a strongly $(\delta,\frU_H)$-admissible vector field $\sfZ^{sm}$ on $\Hsm$. 

\item[(ii)] If $\sfZ^{sm}$ is $(\delta,\frU_H)$-admissible (resp. strongly $(\delta,\frU_H)$-admissible), then for any admissible $\frU^1_H\succ \frU_H$, there exists a good admissible $\frU'_H\succ \frU^1_H$ such that $\sfZ^{sm}$ is $(\delta,\frU'_H)$-admissible (resp. strongly $(\delta,\frU'_H)$-admissible). 
\end{itemize}
\end{lemma}

\begin{cor}\label{corZsmCI}
Let $\delta>0$. 
\begin{itemize}

\item[(i)] The space of  $\delta$-admissible (resp. $(\delta,\frU_H)$-rigid)  vector fields on $H^{sm}$ is nonempty and convex, hence contractible. 

\item[(ii)] Given a $\delta$-admissible (resp. strongly $\delta$-admissible) $\sfZ^{sm}$, the inclusion of the subposet $\CI_{\sfZ^{sm}, \delta}$ (resp. $\CI^\sigma_{\sfZ^{sm}, \delta}$) into $\CI_H$ is cofinal. 
\end{itemize}
\end{cor}

\sss{Symplectic hypersurfaces $\frF$ in $\cH$}

Given an admissible hypersurface $H^{sm}$ in $\RR_{\geqsl 0}^I$ and a $\delta$-admissible vector field $\sfZ^{sm}$, define a symplectic hypersurface $\frF^{H^{sm}, Z^{sm}}$ or simply $\frF$ in $\cH=(b_\RR^I)^{-1}(\Hsm)$ as follows. 
Choose any $\frU_H\in \CI_{Z^{sm}, \delta}$. Let 
\begin{align*}
\sfZ^{\tau}:=\sfZ_{r,I}|_{H^{sm}}-Z^{sm}\in \Gamma(H^{sm}, T\RR_{\geqsl 0}^I),
\end{align*}
which is a vector field that is everywhere transverse to $H^{sm}$. Recall the vector field $V_{I,J}^\frz$ on $\frz_{J;\RR}$ from \S\ref{subsubsec: adm,triang}.  Let 
\begin{align*}
\sfZ^{\tau, \frz}_J:=(\nu_J)_*(\sfZ^\tau|_{H^{sm}_J})=V_{I, J}^{\frz}-Z^{sm, \frz}_J.
\end{align*}
Write $\sfZ^{\tau, \cI_J}$ for $\sfZ^{\tau, \frz}_J|_{H_J^\cI}$. Note that on $\cU_J^H$, we have
\begin{align}
(\nu_J)_*(\sfZ^\tau|_{\cU^H_J})=0_{\cC_J}\times \sfZ^{\tau, \cI_J}=0_{\cC_J}\times (V_{I, J}^{\frz}-Z^{sm, \frz}_J)|_{H_J^\cI}
\end{align}

Consider the symplectic hypersurface
\begin{align}\label{eqfrFJflatcI}
\frF_J^{\flat, \cI}:=\{(\bq,\bp):\bq\in H_J^\cI, \lng\bp, \sfZ^{\tau, \cI_J}_\bq\rng=0\}\subset T^*\frz_{J;\RR}|_{H_J^\cI},
\end{align}
which is canonically identified with $T^*H_J^\cI$. Let $\vartheta_{J,\flat}^{\std}$ be the Liouville 1-form on $\frF_J^{\flat, \cI}$ by pulling back the canonical Liouville 1-form on $T^*H_J^\cI$. 
Let $\cJ^{\cC_J}_{L_J^\der}=(b_{\RR,\der}^{J})^{-1}(\cC_J)\subset \cJ_{L_J^\der}$ and define the symplectic hypersurface (under the natural splitting) 
\begin{align}\label{eqfrFJ}
\fF_J:=\cJ^{\cC_J}_{L_J^\der}\times^{Z(L_J^\der)}T^*Z(L_J)_{\cpt}\times \fF_J^{\flat, \cI}\subset (\bRI)^{-1}(\cU^H_J). 
\end{align}
Although $\frF_J$ depends on the choice of $\Hsm$, $Z^{sm}$ and $\cU^H_J$, we will suppress the dependence in the notation when it does not cause any confusion. 

\begin{lemma}
For any fixed choice of $\Hsm$ and $Z^{sm}$, the collection of symplectic hypersurfaces $\frF_J, J\subsetneq I$, glue to be a symplectic hypersurface $\frF$ in $(\bRI)^{-1}(\Hsm)$, which does \emph{not} depend on the choice of $\frU_H\in \CI_{Z^{sm}, \delta}$. 
\end{lemma}

 Recall that the contact Hamiltonian vector field of $\sfN_{\cJ}$ gives the Reeb vector field on $\cH$. 
\begin{prop}\label{propcHfrF}
The symplectic hypersurface $\frF$ in $\cH$ gives a transverse slice to the characteristic foliation in $\cH$, i.e. 
\begin{align*}
\RR\times \frF&\overset{\sim}{\longrightarrow} \cH\\
(s,y)&\mapsto \varphi_{X_{\sfN_\cJ}}^s(y)
\end{align*}
is a diffeomorphism. 
\end{prop}

\sss{The Liouville vector field on $\frF$}

Let $\vartheta_\frF=\vartheta_\cJ|_\frF$. 

\begin{prop}\label{propvarthetafrFJadm}
\begin{itemize}
\item[(i)] For a $(\delta, \frU_H)$-admissible $Z^{sm}$, we have the restriction of the Liouville 1-form
\begin{align*}
\vartheta_{\frF}|_{\frF_J}=\vartheta_{\cJ_{L_J^\der}}|_{\cJ_{L_J^\der}^{\cC_J}}+\vartheta_{T^*Z(L_J)_\cpt}+\left(\vartheta_{T^*H_J^{\cI}}+d\sfI_J^H\right)
\end{align*}
with respect to the splitting \eqref{eqfrFJ} and the canonical identification $\frF_{J}^{\flat, \cI}\cong T^*H_J^\cI$.

\item[(ii)] Under the same conditions as in (i), the Liouville vector field on $\frF$, denoted by $\sfZ_\frF$, satisfies that 
\begin{align*}
\sfZ_{\frF}|_{\frF_J}=\sfZ_{\cJ_{L_J^\der}}+\sfZ_{T^*Z(L_J)_{\cpt}}+\left(\sfZ_{T^*H_J^{\cI}}-X_{\sfI_J^H}\right). 
\end{align*}
\end{itemize}
\end{prop}
We adopt the similar terminology for Liouville 1-forms on $\frF$ as in Remark \ref{remarkvarthetaadmiss} (ii).

\sss{$\frF$ is a (generalized) Weinstein manifold}
We have the obvious analog of \S\ref{sssindWeinsteindomain} and Theorem \ref{thmMGWeinstein} for $\frF$, where $\trg,\ \frU$ and $\sfZ_\flat$ are replaced by $H^{sm},\ \frU_H$ and $\sfZ^{sm}$, respectively.

\subsection{Liouville/Weinstein sector structure(s) and their deformations on $\ovl{\cJ}_G$}\label{subsecSectdeformJG}

\sss{Open Weinstein sector structure(s) on $\cJ_G$ and $\cT_G$ with respect to the standard Liouville 1-form}\label{sssOpenWeinJG}
Let $\cJ_G^{<r}:=\cJ_G-\bigcup_{z\in Z(G)}S_{I,z}$, where $S_{I,z}$ are the Kostant sections. 
Let $\sfH=\frac{1}{\sfN_\cJ}|_{\cJ_G^{<r}}: \cJ_G^{<r}\to \RR_{>0}$, which satisfies $\sfZ\sfH=\frac{1}{2}\sfH$, and let $\sfI: \cJ_G^{<r}\to \RR$ defined by $\sfI|_{\bigcup_{t\in \RR}\varphi_\sfZ^t(\frF)}=0$ and $X_{\sfH}(\sfI)=1$. Then $\sfZ \sfI=\frac{1}{2}\sfI$, and we have an exact symplectomorphism 
\begin{align}\label{eqphicCB}
\phi_{\cC}: \left(\frF\times \CC_{\Real z>0}, \vartheta_{\frF}+\frac{1}{2}(xdy-ydx)\right)&\overset{\sim}{\longrightarrow} \left(\cJ_G^{<r}, \vartheta_\cJ|_{\cJ_G^{<r}}\right)\\
(\frF\ni p, z=x+\sqrt{-1}y)&\mapsto \varphi_{X_H}^y\circ \varphi_\sfZ^{\log x}(p)=\varphi_{X_H}^y\circ\varphi_{-X_{\sfI}}^{x-1}(p).
\end{align}
For each $J\subsetneq I$, under the natural splitting of $\cJ_{L_J}$, the functions $\sfH|_{\cJ_{L_J}}$ and $\sfI|_{\cJ_{L_J}}$ depend only on the factor $T^*Z(L_J)_{0,\RR}\cong T^*\frz_{J;\RR}$, and they are explicitly given by 
\begin{align}
\label{eqsfHsfI}&\sfH(\bq, \bp)=\sfH(\bq, 0);\quad \text{on the zero-section: } \sfH|_{H_J^\frz}=1,\quad V_{I, J}^{\frz}(\sfH)=\frac{1}{2}\sfH;\\
\nonumber &\sfI=2\lng \bp, \sfH^{-1}\cdot V_{I, J}^{\frz}\rng. 
\end{align}

In the following, for any $\alpha\in [0,1]$, let $\vartheta^\alpha=(1-\alpha)xdy-\alpha ydx$ on $\CC_{\Real z\geqsl 0}$ and let $\sfZ^\alpha=(1-\alpha)x\partial_x+\alpha y\partial_y$ be its Liouville vector field. 
Using $\phi_\cC$, define the real partial compactification  
\begin{align*}
\ovl{\cJ}_G=\cJ_G\cup_{(\frF\times \CC_{\Real z>0})} (\frF\times \CC_{\Real z\geqsl 0}),
\end{align*}
which is a Liouville sector (resp. generalized Weinstein sector) for a $\delta$-admissible (resp. strongly $\delta$-admissible) $\vartheta_\frF$. It has the Liouville (resp. generalized Weinstein) completion
\begin{align*}
\wh{\cJ}_G=\cJ_G\cup_{(\frF\times \CC_{\Real z>0})} (\frF\times \CC_{z}). 
\end{align*}
Thus $(\cJ_G, \vartheta_\cJ)$ has an open generalized Weinstein sector structure from attaching $|Z(G)|$ many index $2r$ Weinstein handles to the open generalized Weinstein sector $\frF\times \CC_{\Real z>0}$ (for $\vartheta_\fF$ strongly $\delta$-admissible). We will denote the extension of $\vartheta_\cJ$ on $\ovl{\cJ}_G$ and $\wh{\cJ}_G$ by $\vartheta_{\ovl{\cJ}}$ and $\vartheta_{\wh{\cJ}}$, respectively.

Using the same construction, we get $(\ovl{\cT}_G, \vartheta_{\ovl{\cT}_G})$, and $(\wh{\cT}_G, \vartheta_{\wh{\cT}_G})$. 
Since the function $\sfH$ is invariant under $Z(G)$-action, for any identification $\cT_G\cong \cJ_G$, the identification extends to an identification $\ovl{\cT}_G\cong \ovl{\cJ}_G$ and $\wh{\cT}_G\cong \wh{\cJ}_G$.

\sss{Deformations of Liouville sector structures on $\ovl{\cJ}_G$}\label{sssDeformLiouvillesector}

We will change the standard Liouville 1-form on $\cJ_G$ so that $\cJ_{G}^{\ovl{\cC}}$ is a Liouville/Weinstein subsector with equivalent wrapped Fukaya category as $\ovl{\cJ}_G$. 
Choose a smooth vector field $\sfZ_{\RR_{\geqsl 0}}=\rho(x)x\partial_x$ on $\RR_{\geqsl 0}$, so that $\rho(x)$ is non-decreasing,  $\rho(x)+x\rho'(x)\leqsl \frac{1}{2}-\ep_0$ for some $\ep_0>0$, and 
\begin{align*}
\rho(x)=\begin{cases}-\frac{1}{2}, &x\in [0, N_1]\\
0,&x\in [N_2, \infty)
\end{cases}
\end{align*}
for some $N_2\gg N_1>1$. Let  $\sfI_{\CRez}:=\vartheta^1(\sfZ_{\RR_{\geqsl 0}})=-\rho(x)xy$ and 
let 
\begin{align*}
\vartheta'_{\CRez;t}=\vartheta^{\frac{1}{2}}+t\ d\sfI_{\CRez}, \text{ for }0\leqsl t\leqsl 1.
\end{align*}
Then $(\CRez, \vartheta'_{\CRez; t}), 0\leqsl t\leqsl 1$, is a deformation of Liouville sectors, with corresponding Liouville vector fields given by 
\begin{align*}
\sfZ'_{\CRez;t}=\sfZ^{\frac{1}{2}}-tX_{\sfI_{\CRez}}=\left(\frac{1}{2}+t\rho(x)\right)x\partial_x+\left(\frac{1}{2}-t\rho(x)-tx\rho'(x)\right)y\partial_y.
\end{align*}
Equip $\cJ^{<r}_G\cong \frF\times \CRez$ with the family of Liouville 1-forms $\vartheta_t':=\vartheta_{\frF}+\vartheta'_{\CRez;t}$. Since $\vartheta_t'$ agrees with $\vartheta|_{\cJ^{<r}_G\cap \sfH^{-1}[N_2, \infty)}$, it extends to a Liouville 1-form on $\ovl{\cJ}_G$, denoted by $\vartheta'_t$ as well. Then $(\ovl{\cJ}_G, \vartheta'_t), 0\leqsl t\leqsl 1$, is a deformation of Liouville sectors, with $(\ovl{\cJ}_G, \vartheta'_0)=(\ovl{\cJ}_G, \vartheta_{\ovl{\cJ}})$. Let $\sfZ_t'$ be the Liouville vector field for $\vartheta'_t$. Since the choices of $\rho(x)$ form a nonempty and convex space for any fixed $\ep_0>0$, the space of $\vartheta_t'$ on $\ovl{\cJ}_G$ is contractible.

\sss{Sector structures on $\cJ_G^{\ovl{\cC}}$}\label{sssSectJGovlC}

Consider the Liouville sector $(\ovl{\cJ}_G, \vartheta_1')$. Since $\sfZ'_{\CRez;1}|_{\{x\in [0,N_1]\}}=y\partial_y$, $\cJ_G^{\ovl{\cC}}\hookrightarrow \ovl{\cJ}_G$ is a sectorial inclusion, inducing an equivalence on wrapped Fukaya categories. 
Alternatively, one can define a deformation of Liouville sectors connecting $\ovl{\cJ}_G$ and $\cJ_G^{\ovl{\cC}}$ 
as follows. Let $\sfN_{\ovl{\cJ}}: \ovl{\cJ}_G\to [0,\infty]$ be the continuous function extending $\sfN_{\cJ}$. For $t\in [0,1]$, let $\ovl{\cJ}_G^{[0, 1/t]}=\sfN^{-1}_{\ovl{\cJ}}([0,1/t])=(\frF\times \CC_{\Real z\geqsl t})\cup (b_{\RR}^{I})^{-1}(0)$. Choose a smooth family of diffeomorphisms $\varrho_t: [t, \infty)\to [0,\infty), 0\leqsl t\leqsl 1$, that is the identity on $[R,\infty)$ for some $R>1$ and $\varrho_0=id_{[0,\infty)}$. This defines a family of diffeomorphisms $\wt{\varrho}_t=(\varrho_t, id_y):\CC_{\Real z\geqsl t}\to \CRez$, which results in a family of diffeomorphisms $(id_\frF, \wt{\varrho_t}):\frF\times \CC_{\Real z\geqsl t}\to \frF\times \CC_{\Real z\geqsl 0}$. The latter clearly extends to a family of diffeomorphisms $\ovl{\cJ}_G^{[0, 1/t]}\to \ovl{\cJ}_G$. The pushforward of the Liouville sector structures to $\ovl{\cJ}_G$ gives a deformation of Liouville sectors $(\ovl{\cJ}_G, (\vartheta_1')_t), 0\leqsl t\leqsl 1$ starting with $(\ovl{\cJ}_G, \vartheta_1')$ and ending with  $(\cJ_G^{\ovl{\cC}}, \vartheta_1'|_{\cJ_G^{\ovl{\cC}}})$. Again, the space of $(\varrho_t)_{0\leqsl t\leqsl 1}$ forms a nonempty and convex space, so it is contractible. 

\begin{set}\label{settingcrUPI}
Consider any admissible open coverings $\frU_{\cP(I)}=\{(\cC_J, \cI_J)\}_{J\in \cP(I)}\prec \frU'_{\cP(I)}=\{(\cC'_J, \cI_J')\}_{J\in \cP(I)}$ of open neighborhoods $\cU$ and $\cU'$ of $\ovl{\cC}$, respectively, such that (1) $\partial^\circ \cC\cap \ovl{\cC}_I=\vn$ and 
$H^{sm}=\partial^\circ \cC$ is good with respect to $\sfZ_{r,I}$ and $\frU'_{\cP^\dagg(I)}=\{(\cC_J, \cI_J)\}_{J\in \cP^\dagg(I)}$; (2) $\cI_J'\Supset \cI_J$ and $\cI_J'$ is sufficiently large relative to $\cI_J$. 
\end{set}
We have the notion of $(\delta,\frU_{\cP(I)}')$-semi-admissibility, $\frU_{\cP(I)}$-rigidity and (strong) $(\delta,\frU_{\cP(I)})$-admissibility for a smooth vector field $\sfZ_{\flat;I}$ on $\ovl{\cC}$, defined by the following properties: 
\begin{itemize}
\item $\sfZ_{\flat;I}$ is a smooth vector field on $\ovl{\cC}$ in the sense of \S\ref{sssMWC}---in particular $\sfZ_{\flat;I}$ is tangent to $\partial^\circ \cC$ everywhere;
\item $\sfZ_{\flat;I}$ is the restriction of the respective kinds on $\cU$, analogous to 
Definition \ref{defnZsemiadm} and Definition \ref{defnZadm}, in which $\sfZ_{J,\flat}^{\frz}$ for any $J\subset I$ denotes a smooth extension of the restriction of $\sfZ_{\flat;I}$ to $\nu_J(\cU\cap (\RR_{\geqsl 0}^{J}\times \RR_{>0}^{I-J}))$.
\end{itemize}
Note that if $\sfZ_{\flat;I}$ is $(\delta,\frU_{\cP(I)}',\frU_{\cP(I)})$-rigid on $\cU$, then its restriction on $\ovl{\cC}$ automatically satisfies the first requirement (for example, this follows from Lemma \ref{lemmasemitoadm} (ii) and Remark \ref{rmksplitting} (ii)). We will say such a $\sfZ_{\flat;I}$ is $(\delta,\cP(I))$-rigid, when $\frU_{\cP(I)}'$ and $\frU_{\cP(I)}$ are not specified.

Existence results for admissible $\sfZ_{\flat;I}$ on $\ovl{\cC}$ for sufficiently thin and wide $\frU_{\cP(I)}, \frU'_{\cP(I)}$ follow from the same proofs as in \S\ref{ssexisZflat}. In particular, the space of $\delta$-admissible (resp. $(\delta, \frU)$-rigid) vector field $\sfZ_{\flat, I}$ on $\ovl{\cC}$ is a nonempty and convex space. 
Given a $\delta$-admissible $\sfZ_{\flat;I}$ on $\ovl{\cC}$, we have the canonical Liouville 1-form $\vartheta_{I}$ with its Liouville vector field $\sfZ_{I}$ on $\cJ_G^{\ovl{\cC}}$. 
We will use similar conventions as in Remark \ref{remarkvarthetaadmiss} (ii).

\begin{lemma}\label{lemmadeltaadmCsector}
Assume $\delta>0$ is sufficiently small. 
\begin{itemize}
\item[(i)] Given a $\delta$-admissible $\sfZ_{\flat;I}$ on $\ovl{\cC}$, we have $(\cJ_G^{\ovl{\cC}}, \vartheta_I)$ a Liouville sector. For any two  $\delta$-admissible $\vartheta_I^{(i)}, i=0,1$,  $(\cJ_G^{\ovl{\cC}}, (1-t)\vartheta_I^{(0)}+t\vartheta_I^{(1)}), 0\leqsl t\leqsl 1$, is a (non-compactly supported) deformation of Liouville sectors. Thus, $(\cJ_G^{\ovl{\cC}},\vartheta_I)$, for $\delta$-admissible $\vartheta_I$, form a contractible space of deformation equivalent Liouville sectors. 

\item[(ii)] By choosing $\rho(x)$ so that $|x\rho'(x)|$ has a sufficiently small upper bound (such $\rho(x)$ always exists), the Liouville 1-form $\vartheta_1'|_{\cJ_G^{\ovl{\cC}}}$ defined in \S\ref{sssDeformLiouvillesector} is $\delta$-admissible. 

\item[(iii)] There is a deformation between Liouville sectors $(\ovl{\cJ}_G, \vartheta_0')$ and $(\cJ_G^{\ovl{\cC}},\vartheta_I)$. In particular, this induces a canonical equivalence (up to a contractible space of choices) $\cW(\ovl{\cJ}_G, \vartheta_0'=\vartheta_{\ovl{\cJ}})\simeq \cW(\cJ_G^{\ovl{\cC}}, \vartheta_I)$ for any $\delta$-admissible $\vartheta_I$. 
\end{itemize}
\end{lemma}

\begin{proof}
(i) First, since $H^{sm}$ is $\sfZ_{\flat;I}$-invariant, we have $(b_{\RR}^I)^{-1}(H^{sm})$ is $\sfZ_I$-invariant. Using the setup and notations from \S\ref{sssindWeinsteindomain}, we have for any $R_\bullet=(R_{0}\gg \cdots \gg R_r\gg 1)$, 
$\frD_{R_\bullet}^{\subset I}:=\bigcup_{J\subset I}\frD^J_{R_{|J|}, R_{|J|}}$ a Liouville subdomain in $\cJ_G$. 
Consider $\frD_{R_\bullet}^{\ovl{\cC}}:=\frD^{\subset I}_{R_\bullet}\cap (b_\RR^I)^{-1}(\ovl{\cC})$. 
 We have $\sfZ_I$ is convex along $\partial \frD_{R_\bullet}^{\ovl{\cC}}\cap (b_{\RR}^I)^{-1}(\cC)$, and there is no zero of $\sfZ_{I}$ outside $\frD_{R_\bullet}^{\ovl{\cC}}$ in $\cJ_G^{\ovl{\cC}}$. 
In fact, for a sequence $R_\bullet^{(n)}, n\geqsl 1$ with $R_j^{(n)}\to \infty$, $\frD_{R_\bullet^{(n)}}^{\ovl{\cC}}$ gives an exhausting sequence of compact subdomains with corners in $\cJ_G^{\ovl{\cC}}$.

It follows from the above that 
$\frD_{R_\bullet}^{\cH}:=\partial \frD_{R_\bullet}^{\ovl{\cC}}\cap \cH$ is a compact subdomain in $\cH$ whose boundary has corners such that $\sfZ_{I}$ is convex along its boundary (a reminder: the interior of $\frD_{R_\bullet}^{\cH}\cap (b_\RR^I)^{-1}(\partial(\partial^\circ\cC))$ is \emph{not} part of the boundary of $\frD_{R_\bullet}^{\cH}$, due to the Weinstein handle structure over $\partial(\partial^\circ\cC)$). Again, there is no zero of $\sfZ_{I}$ outside $\frD_{R_\bullet}^{\cH}$ in $\cH$. 
Thus, combining with Proposition \ref{propcHfrF}, we get $(\cJ_G^{\ovl{\cC}}, \vartheta_I)$ is a Liouville sector (by the third equivalent condition in \cite[Definition 2.4]{GPS1}). 

Using $\frD_{R_\bullet}^{\subset I}$, the rest follows from a similar observation as in Remark \ref{remarkMGthetadeform} (although here the deformations of sectors are not known to be conjugate equivalent to compactly supported deformations; see similar situations in \cite{GPS1}). \\

(ii) Choose any sufficiently thin and wide admissible $\frU_{\cP(I)}$ that covers a neighborhood $\cU$ of $\ovl{\cC}$. 
Since $\vartheta_0'=\vartheta_{\ovl{\cJ}}$, we get from \eqref{eqsfHsfI} that for each $J\subset I$, 
\begin{align*}
&\vartheta'_1|_{(b_{\RR}^I)^{-1}(\cU_J)}=\vartheta_{\cJ_{L_J^\der}}+\vartheta_{T^*Z(L_J)_{\cpt}}+\left(\vartheta_{T^*\cI_J}-d\lng \bp, (1+2\rho(\sfH))\cdot V_{I, J}^{\frz}\rng\right)
\end{align*}
under the natural splitting. Let $\sfI_{J, \rho}=-\lng \bp, (1+2\rho(\sfH))\cdot V_{I, J}^{\frz}\rng$.  By a direct calculation of $X_{\sfI_{J, \rho}}$, the statement easily follows. 

(iii) is a direct consequence of the deformation $\vartheta_t'$ defined in \S\ref{sssDeformLiouvillesector} and the result from part (i). 
\end{proof}

\sss{$\delta$-admissible Liouville sector structures for reductive groups}\label{sssdeltaadmReductive}

The definition of admissible vector fields on $\ovl{\cC}\subset \RR^{I}_{\geqsl 0}$ naturally generalizes in the setting of a reductive group $G'$. Let $I'$ be a complete set of simple roots and $\frz_{I';\RR}=\Lie Z(G')_{0,\RR}$, then the tropical manifold in this setting is just the product $\RR_{\geqsl 0}^{I'}\times\frz_{I';\RR}$. Let $\cC\subset \RR_{\geqsl 0}^{I'}$, $\cI\subset \frz_{I';\RR}$ and $\cV:=\cC\times \cI$ as above. Choose any good admissible open covering $\frU_{\cP(I')}=\{(\cC_J, \cI_J)\}_{J\in \cP(I')}\prec \frU'_{\cP(I')}=\{(\cC'_J, \cI_J')\}_{J\in \cP(I')}$ of open neighborhoods $\cU$ and $\cU'$ of $\ovl{\cV}$, respectively, such that (1) $(\partial^\circ\cC\times \ovl{\cI})\cap \ovl{\cU}_{I'}=\vn$ and 
$\partial^\circ \cC\times \Nb(\ovl{\cI})$ is good with respect to $\sfZ_{r,I'}\times 0$ and $\frU_{\cP^\dagg(I')}$; (2) $\cI_J'\Supset \cI_J$ and $\cI_J'$ is sufficiently large relative to $\cI_J$. We have the notion of $(\delta,\frU_{\cP(I')}')$-semi-admissibility, $\frU_{\cP(I')}$-rigidity and (strong) $(\delta,\frU_{\cP(I')})$-admissibility for a smooth vector field $\sfZ_{\flat;I'}$ on $\ovl{\cV}$, defined in the same way as for $\sfZ_{\flat;I}$ in \S\ref{sssSectJGovlC}. The discussion of $\sfZ_{\flat;I}$ in  \S\ref{sssSectJGovlC} generalizes to $\sfZ_{\flat;I'}$ verbatim. In particular, given a $\delta$-admissible $\sfZ_{\flat;I'}$ on $\ovl{\cV}$, we have the canonical Liouville 1-form $\vartheta_{I'}$ for $\sfZ_{\flat;I'}$ on $\cT_{G'}^{\ovl{\cV}}$. We have the analog of Lemma \ref{lemmadeltaadmCsector} (i) with essentially the same proof.

\begin{lemma}\label{lemmreducadeltaadmVsector}
Assume $\delta>0$ is sufficiently small. Given a $\delta$-admissible $\sfZ_{\flat;I'}$ on $\ovl{\cV}$, 
we have $(\cT_{G'}^{\ovl{\cV}}, \vartheta_{I'})$ a Liouville sector. 
\end{lemma}

\begin{proof}
Similarly to the first paragraph of the proof of Lemma \ref{lemmadeltaadmCsector} (i), we define $\frD_{R_\bullet}^{\ovl{\cV}}$ in the same way. 

Let $\cH_{I'}=(b_{\RR}^{I'})^{-1}(\partial^\circ \cV)$. By the same argument as in the second paragraph of that proof, we have $\sfZ_{I'}$ convex along $\frD_{R_\bullet}^{\cH_{I'}}$ and there is no zero of $\sfZ_{I'}$ outside $\frD_{R_\bullet}^{\cH_{I'}}$ in $\cH_{I'}$. It remains to check that  $\cH_{I'}$ (up to smoothing the corners of $\partial^\circ\cV$) has a transverse slice (a symplectic hypersurface) to the characteristic foliation. Note this property is independent of the choice of sector structures.  
Since we can put a product sector structure on $\cT_{G'}^{\ovl{\cV}}\cong\cT_{G'_{\der}}^{\ovl{\cC}}\times^{Z(G'_\der)}T^*Z(G')_{\cpt}\times T^*\ovl{\cI}$, $\cH_{I'}$ as the sector boundary of this particular Liouville sector structure satisfies the property. The proof is complete. 
\end{proof}

Given $\delta>0$, let $\frs_{\delta}(\cT_{G'}^{\ovl{\cV}})$ be the space of $\delta$-admissible $\vartheta_{I'}$ on $\cT_{G'}^{\ovl{\cV}}$, and we also refer to it as the space of $\delta$-admissible Liouville sector structures on $\cT_{G'}^{\ovl{\cV}}$. Clearly, $\frs_{\delta}(\cT_{G'}^{\ovl{\cV}})$ is nonempty and convex, hence it is contractible. 

Given $\delta>0$ and $\frU_{\cP(I')}$, let $\frs_{\delta, \rigid}^{\frU_{\cP(I')}}(\cT_{G'}^{\ovl{\cV}})$ be the space of $(\delta, \frU_{\cP(I')})$-rigid $\vartheta_{I'}$ on $\cT_{G'}^{\ovl{\cV}}$. Since $\frs_{\delta, \rigid}^{\frU_{\cP(I')}}(\cT_{G'}^{\ovl{\cV}})$ is independent of the choice of $\cU_{I'}$  (in particular, we can always assume $\cU_{I'}=\cC\times \Nb(\ovl{\cI})$), which follows from Remark \ref{rmksplitting} (i), we will denote  $\frs_{\delta, \rigid}^{\frU_{\cP(I')}}(\cT_{G'}^{\ovl{\cV}})$ also by  $\frs_{\delta, \rigid}^{\frU_{\cP^\dagg(I')}}(\cT_{G'}^{\ovl{\cV}})$. 
Note that $\frs_{\delta, \rigid}^{\frU_{\cP^\dagg(I')}}(\cT_{G'}^{\ovl{\cV}})$ is a convex sub-space in $\frs_{\delta}(\cT_{G'}^{\ovl{\cV}})$. By Remark \ref{rmksplitting} (i) and Lemma \ref{lemmadeltaadmCsector} (iii), we immediately have the following.

\begin{lemma}\label{lemmarigidproductsector}
For any $\vartheta_{I'}\in \frs_{\delta, \rigid}^{\frU_{\cP^\dagg(I')}}(\cT_{G'}^{\ovl{\cV}})$, we have $(\cT_{G'}^{\ovl{\cV}}, \vartheta_{I'})$ a product Liouville sector 
\begin{align*}
&\cT_{G'_{\der}}^{\ovl{\cC}}\times^{Z(G'_\der)}T^*Z(G')_{\cpt}\times T^*\ovl{\cI},
\end{align*}
where $\cT_{G'_{\der}}^{\ovl{\cC}}$ is equipped with a $(\delta,\cP(I'))$-rigid Liouville sector structure, and both $T^*Z(G')_{\cpt}$ and $ T^*\ovl{\cI}$ are equipped with the standard Liouville (sector) structures. Therefore, there is a canonical equivalence (up to a contractible space of choices) $\cW(\cT_{G'}^{\ovl{\cV}}, \vartheta_{I'})\simeq \cW(\ovl{\cT}_{G'}, \vartheta_{\ovl{\cT}_{G'}})$, for any $\vartheta_{I'}\in \frs_{\delta}(\cT_{G'}^{\ovl{\cV}})$. 
\end{lemma}

\sss{A remark on the Liouville completion of $\cT^{\ovl{\cV}}_{H}$}\label{sssCompletioncTGV}

Let $H$ be a reductive group with maximal torus $T$ and a complete set of simple roots $I_H$. 
let $\hat{\frU}:=\hat{\frU}_{\cP(I_H)}$ be an admissible covering of $\ovl{\cV}$ and choose a $(\delta, \hat{\frU})$-rigid $\vartheta$ on $\cT_H^{\ovl{\cV}}$. The Liouville sector $(\cT_H^{\ovl{\cV}}, \vartheta)$ is a product sector $(\cT_{H^\der}^{\ovl{\cC}}, \vartheta^{\der})\times^{Z(H^\der)}T^*Z(H)_{\cpt}\times T^*\ovl{\cI}$, and it is 
 canonically deformation equivalent to a (generalized) Weinstein sector up to a contractible space of choices. 

Let $\wh{\cT}_{H}^{\ovl{\cV}}$ be the Liouville completion of $\cT^{\ovl{\cV}}_{H}$. 
We now choose an identification between $\wh{\cT}_{H}^{\ovl{\cV}}$ and $\cT_H$, whose restriction to $\cT^{\ovl{\cV}}_{H}$ is the identity map. Let $\frF\subset (b_{I_H,\der}^\RR)^{-1}(\partial^\circ \cC)$ be any Liouville hypersurface that is a transverse slice to the characteristic foliation. Write $\cT_H=\cT_{H^\der}\times^{Z(H^\der)}T^*Z(H)_{\cpt}\times T^*Z(H)_{0,\RR}$. Using a similar consideration as $\phi_\cC$ \eqref{eqphicCB}, 
we have a symplectic embedding
\begin{align*}
\frF\times\CC_{\Real z>0}\hookrightarrow \cT_{H^\der},
\end{align*}
and a sectorial embedding
\begin{align*}
\frF\times \left(T^*[1, 1+\ep],\vartheta_{T^*[1, 1+\ep]}\right)\hookrightarrow (\cT_{H^\der}^{\ovl{\cC}}, \vartheta^\der)
\end{align*}
for small $\ep>0$,  where  $T^*[1, 1+\ep]$ is identified with $\CC_{\Real z\in [1,1+\ep]}$ from the previous line. 
Now write 
\begin{align}\label{eqcTHdercompletion}
\cT_{H^\der}=\cT_{H^\der}^{\ovl{\cC}}\cup_{\frF\times T^*[1, 1+\ep]} (\frF\times \CC_{\Real z\in (0, 1+\ep]}),  
\end{align}
and equip  $\CC_{\Real z\in (0, 1+\ep]}$ with a Liouville 1-form 
\begin{align*}
\vartheta_{\CC_{\Real z\in (0, 1+\ep]}}'=\vartheta^1+d(-\rho(x)y), 
\text{ where }\rho(x)=\begin{cases}
-\frac{1}{2}x, &x\in (0,\ep)\\
0, &x\in [1,1+\ep].
\end{cases}\text{ and }\rho'(x)\leqsl \frac{1}{2}.
\end{align*}
Note that the Liouville vector field on $\CC_{\Real z\in (0, 1+\ep]}$ is complete. 
Then $\cT_{H^\der}$ \eqref{eqcTHdercompletion} with the resulting Liouville 1-form is a Liouville completion of $\cT^{\ovl{\cC}}_{H^\der}$. 

It is standard that one can equip $T^*Z(H)_{0,\RR}$ with a Liouville 1-form so that it is the Liouville completion of $T^*\ovl{\cI}$. Thus we have realized $\cT_H$, equipped with an appropriate Liouville 1-form, as the Liouville completion of $\cT^{\ovl{\cV}}_H$.

%%%%%%%%%%%%%%%%%%%%%%%%%%%%%%%%%%%%%%%%%%%%%%%%%

\subsection{Sectorial inclusions in $\ovl{\cT}_G$ and $\cM_G$}\label{subsecSectInclu}

In this subsection, we will mostly consider $\ovl{\cT}_G$ instead of $\ovl{\cJ}_G$. The reason, as already appeared in \S\ref{secMGHiggs}, is that for any Levi subgroup $L_J$ of $G$, the open embedding $\cT_{L_J}\hookrightarrow \cT_G$ is \emph{canonical}, unlike the open embedding $\cJ_{L_J}\hookrightarrow \cJ_G$ relies on a choice of lifting of Weyl group elements. Similarly for pseudo-Levi subgroups $L_J$, we have the canonical embedding $\cT_{L_J}\hookrightarrow \cT_{L_J'}$ for $J\subset J'\subset_{ft}\wt{I}$, which are compatible under compositions of inclusions.

\sss{Spaces of Liouville sector structures associated with sectorial inclusions on $\cM_G$}\label{sssSpcSectStr}

Fix an admissible open covering $\hat\frU=\{(\hat{\cC}_J,\hat{\cI}_J)\}_{J\in \PfI}$ of $\trg$. To simplify notations, we also denote $\cT_{L_J}^{\nu_J(\ovl{\cU}_J)}$ by $\cT_{L_J}^{\ovl{\cU}_J}$. 
For any contractible open subset with smooth boundary $\cI_J\subset \hat{\cI}_J$, let 
\begin{align}
\label{eqfrsrigdelta}&\frs^{\hat{\frU}}_{\delta,\rigid}(\cT_{L_J}^{\ovl{\cU}_J}):=\frs^{\hat{\frU}_{\subsetneq J}}_{\delta,\rigid}(\cT_{L_J}^{\nu_J(\ovl{\cU}_J)}), \text{ where }\ovl{\cU}_J:=\nu_J^{-1}(\ovl{\hat{\cC}}_J\times \ovl{\cI}_J). 
\end{align}

For any $J^1\supset \cdots\supset J^k$ in $\PfI$, let $(\cP\cI)_{J^1,\cdots, J^k}$ be the poset of $k$-tuples $(\cI_{J^1}, \cI_{J^2},\cdots, \cI_{J^k}), J^1\supset \cdots\supset J^k$, such that 
\begin{align*}
\cI_{J^i}\subset \hat{\cI}_{J^i} \text{ for }1\leqsl i\leqsl k,\quad\text{and for any }i'>i,\  \cU_{J^{i'}}\subset \cU_{J^{i}}. 
\end{align*}
The partial order is defined by $(\cI_{J^i})_{1\leqsl i\leqsl k}\prec (\cI'_{J^i})_{1\leqsl i\leqsl k}$ iff $\cI_{J^i}\subset \cI'_{J^i}$. Let $(\cP\cI)_k$ be the disjoint union of $(\cP\cI)_{J^1,\cdots, J^k}$ over all chains $J^1\supset \cdots\supset J^k$ in $\PfI$. Then $(\cP\cI)_{\bullet}$ is naturally a simplicial ordinary 1-category, where we put $(\cP\cI)_{0}=\pt$.

For any  $(\cI_{J^i})_{1\leqsl i\leqsl k}\in (\cP\cI)_{J^1,\cdots, J^k}$, let $(\frs_{\delta}^{\hat{\frU}})_{(\cI_{J^i})_{1\leqsl i\leqsl k}}$ be the subspace of $\prod_{i=1}^k\frs_\delta(\cT_{L_{J^i}}^{\ovl{\cU}_{J^i}})$ consisting of $(\vartheta_{J^i})_{1\leqsl i\leqsl k}$ satisfying: 
\begin{itemize}
\item[]$(\cT_{L_{J^{i'}}}^{\ovl{\cU}_{J^{i'}}}, \vartheta_{J^{i'}})\hookrightarrow (\cT_{L_{J^{i}}}^{\ovl{\cU}_{J^{i}}}, \vartheta_{J^{i}})$ is a (strict) sectorial inclusion for any $i'>i$, equivalently $\vartheta_{J^{i}}|_{\cT_{L_{J^{i'}}}^{\ovl{\cU}_{J^{i'}}}}=\vartheta_{J^{i'}}$.
\end{itemize}

For any $(\cI_{J^i})_{1\leqsl i\leqsl k}\in (\cP\cI)_{J^1,\cdots, J^k}$, let $(\frs^{\hat{\frU}}_{\delta, \rigid})_{(\cI_{J^i})_{1\leqsl i\leqsl k}}$ be the subspace of  $(\frs_{\delta}^{\hat{\frU}})_{(\cI_{J^i})_{1\leqsl i\leqsl k}}$ consisting of $(\vartheta_{J^i})_{1\leqsl i\leqsl k}$ with $\vartheta_{J^i}\in \frs^{\hat{\frU}}_{\delta,\rigid}(\cT_{L_J}^{\ovl{\cU}_J})$. The $\hat{\frU}_{\subset J^1}$-rigidity condition on $\vartheta_{J^1}$ makes 
the natural map $(\frs^{\hat{\frU}}_{\delta, \rigid})_{(\cI_{J^i})_{1\leqsl i\leqsl k}}\to (\frs^{\hat{\frU}}_{\delta, \rigid})_{\cI_{J^1}}$ an isomorphism, thus $(\frs^{\hat{\frU}}_{\delta, \rigid})_{(\cI_{J^i})_{1\leqsl i\leqsl k}}$ is nonempty and convex. A small (and not important) remark is that for $\cI_{J^1}\subsetneq \cI_{J^1}'$, the natural map from restriction $(\frs^{\hat{\frU}}_{\delta, \rigid})_{\cI'_{J^1}}\to (\frs^{\hat{\frU}}_{\delta, \rigid})_{\cI_{J^1}}$ is always injective but not always an isomorphism of spaces (this can be seen using Remark \ref{rmksplitting} (i)). 
A direct consequence is the following. 

\begin{lemma}\label{lemmafrsdeltaconvex}
For any  $(\cI_{J^i})_{1\leqsl i\leqsl k}\in (\cP\cI)_{J^1,\cdots, J^k}$,  $(\frs_{\delta}^{\hat{\frU}})_{(\cI_{J^i})_{1\leqsl i\leqsl k}}$ is a nonempty and convex space, hence contractible, and it contains $(\frs^{\hat{\frU}}_{\delta, \rigid})_{(\cI_{J^i})_{1\leqsl i\leqsl k}}$ as a convex subspace. 
\end{lemma}

Let $\Spc$ be the $\infty$-category of topological spaces. We have the natural functor
\begin{align*}
(\frs^{\hat{\frU}}_{\delta, \rigid})_k: \left((\cP\cI)_k\right)^{op}&\longrightarrow \Spc\\
(\cI_{J^i})_{1\leqsl i\leqsl k}&\mapsto (\frs^{\hat{\frU}}_{\delta, \rigid})_{(\cI_{J^i})_{1\leqsl i\leqsl k}}(\cong (\frs^{\hat{\frU}}_{\delta, \rigid})_{(\cI_{J^1})}\text{ if }k\geqsl 1). 
\end{align*}
Consider the trivial Cartesian fibration 
\begin{align}\label{eqcEk}
\cE_{(\frs^{\hat{\frU}}_{\delta, \rigid})_k}\to (\cP\cI)_k 
\end{align}
associated with $(\frs^{\hat{\frU}}_{\delta, \rigid})_k$. For any $f: [k]\to [\ell]$ in $\Delta$, let $f^*_{\cPI}: (\cPI)_{\ell}\to  (\cPI)_{k}$. 
Then we have the natural transformation $(\frs^{\hat{\frU}}_{\delta, \rigid})_\ell\to (\frs^{\hat{\frU}}_{\delta, \rigid})_k\circ (f^*_{\cPI})^{op}$, which is certainly compatible with compositions of morphisms in $\Delta$. Therefore, the Cartesion fibrations $\cE_{(\frs^{\hat{\frU}}_{\delta, \rigid})_k}\to (\cP\cI)_k$ for all $k$ form a simplicial object of (trivial) Cartesian fibraions $\cE_{(\frs^{\hat{\frU}}_{\delta, \rigid})_\bullet}\to (\cP\cI)_\bullet$. 

\begin{lemma}\label{lemmacWfrUrigid}
There is a natural simplicial functor (up to a contractible space of choices)
\begin{align*}
\cW^{\cE, \hat{\frU}}_{\delta,\rigid}: \cE_{(\frs^{\hat{\frU}}_{\delta, \rigid})_\bullet}&\longrightarrow \Maps([\bullet], \CatEx)\\
(\cT_{L_{J^{i}}}^{\ovl{\cU}_{J^{i}}}, \vartheta_{J^{i}})_{1\leqsl i\leqsl \bullet}&\mapsto \left((\cW(\cT_{L_{J^{i}}}^{\ovl{\cU}_{J^{i}}}, \vartheta_{J^{i}})\right)_{1\leqsl i\leqsl \bullet},
\end{align*}
where the right-hand-side in the second row represents the corresponding diagram of wrapped Fukaya categories under sectorial inclusions. In particular,  $\cW^{\cE, \hat{\frU}}_{\delta,\rigid}$ descends to a simplicial functor 
\begin{align}\label{eqcWfrUrigid}
\cW^{\hat{\frU}}_{\delta,\rigid}: N_\bullet(\PfI)\simeq |\cPI_\bullet|\longrightarrow \Maps([\bullet], \CatEx).
\end{align}
\end{lemma}

\begin{proof}
The definition of $\cW^{\cE, \hat{\frU}}_{\delta,\rigid}$ is standard: we only need to check that for any Cartesian arrow of \eqref{eqcEk}, $(\cT_{L_{J^{i}}}^{\ovl{\cU}_{J^{i}}}, \vartheta_{J^{i}}=\vartheta'_{J^{i}}|_{\cT_{L_{J^{i}}}^{\ovl{\cU}_{J^{i}}}})_{1\leqsl i\leqsl k}\to (\cT_{L_{J^{i}}}^{\ovl{\cU'}_{J^{i}}}, \vartheta'_{J^{i}})_{1\leqsl i\leqsl  k}$, it induces an equivalence between their respective diagrams of wrapped Fukaya categories. This follows directly from the product sector structures on them (cf. Lemma \ref{lemmarigidproductsector}). Since each $(\cPI)_{J^1,\cdots, J^k}$ has a maximal object $(\cI_{J^i}:=\hat{\cI}_{J^1,\cdots, J^i})_{1\leqsl i\leqsl k}$ (notation as in \eqref{eq: I beta1,k}), 
the rest of the lemma follows immediately. 
\end{proof}

\sss{Review of sectorial inclusions in $\ovl{\cT}_G$ from \cite{J}}\label{sssReviewSectIncl}

Continue on the settings in \S\ref{subsec: Hsm} and \S\ref{sssOpenWeinJG}. In \cite{J}, we constructed Liouville subsectors in $\ovl{\cT}_G$, denoted here by $\ovl{\cY}_J$, for each $J\subsetneq I$. Let us review the constructions below. Note that since in \emph{loc. cit.}, we need the canonical Liouvile 1-form to calculate $\cW(\ovl{\cT}_G)$, the subsectors $\ovl{\cY}_J$ are genuine subsectors in $\ovl{\cT}_G$, without deforming the canonical Liouville 1-form.

For any $\theta_0\in (0,\frac{\pi}{4}]$ (sufficiently small) and $A>0$, let 
$$\ovl{\cP}_{\theta_0, \leqsl A}=\{z=r e^{i\theta}, r\geqsl 0, \theta\in [-\theta_0, \theta_0]\}\cup \{z: 0\leqsl \Real z\leqsl A\}\subset \CC_{\Real z\geqsl 0},$$
which is a Liouville subsector of $(\CC_{\Real z\geqsl 0}, \vartheta^{\frac{1}{2}})$.

Take $\frU'\succ \frU$ as in Setting \ref{settingcrUPI}. They induce admissible covering $\frU_H'\succ \frU_H$ on $H^{sm}$. 
Let $\sfZ^{sm}$ be a $(\delta, \frU_H)$-rigid vector field on $H^{sm}$. 
For any $J\subsetneq I$, $\ovl{\frF}_J$ is a Liouville subsector of $\frF$. This can be easily seen by choosing $\frU'$ such that $\frU'_H\in \cC\cI\frU_{H; Z^{sm}, \delta}$ and use Proposition \ref{propvarthetafrFJadm} for $\frF_J'$ defined using $\frU_H'$. let
\begin{align}\label{eqovlcYJ}
\ovl{\cY}_J:=\ovl{\frF}_J\times \ovl{\cP}_{\theta_0, \leqsl A}
\end{align}
be the product Liouville subsector of $\ovl{\cT}_G$. Technically, one should do a standard ``conification" process to make $\ovl{\cY}_J$ a genuine subsector, but we will ignore this issue. 

For any strictly decreasing chain $\vec{J}=(J^1\supsetneq J^2\supsetneq\cdots\supsetneq J^k)$  in $\PdI$, we have the natural inclusions of Liouville subsectors 
\begin{align}
\nonumber &\ovl{\frF}_{\vec{J}}:=\bigcap_{i=1}^k\ovl{\frF}_{J^i}\hookrightarrow \ovl{\frF}_{J^k},\\
\label{eqcYJ1JkJ1}&\ovl{\cY}_{\vec{J}}:=\bigcap_{i=1}^k\ovl{\cY}_{J^i}\hookrightarrow \ovl{\cY}_{J^k}
\end{align}
in $\frF$ and $\ovl{\cT}_G$, respectively. 
By Lemma \ref{lemma: admissible, cI contractible} (i) (and its analog for the induced admissible covering on $H^{sm}$), these inclusions induce canonical equivalences on the wrapped Fukaya categories. 
In particular, for any $J^1\supsetneq J^2$, the correspondence of sectorial inclusions
\begin{align*}
\ovl{\cY}_{J^2}\hookleftarrow \ovl{\cY}_{J^1,J^2}\hookrightarrow \ovl{\cY}_{J^1}, 
\end{align*}
induces a canonical functor $\cW(\ovl{\cY}_{J^2})\to \cW(\ovl{\cY}_{J^1})$.

\sss{A digression on posets}\label{sssdigposets}

For any poset $(\cP, <)$, let $\left(\sfC(\cP),\prec\right)$ be the poset of (non-empty) finite strictly decreasing chains $C=(c^1\gneqq \cdots \gneqq c^k)$ in $\cP$ with $C_1\prec C_2$ if and only if $C_1$ is contained in $C_2$ (this is also known as the ``barycentric subdivision"). Consider the natural functor 
\begin{align*}
p_{\cP}: \sfC(\cP)^{op}\to \cP,\quad C=(c^1\gneqq \cdots \gneqq c^k)\mapsto c^k. 
\end{align*}

\begin{lemma}\label{lemmapCScofinal}
\begin{itemize}
\item[(i)] 
The functor $p_{\cP}: \sfC(\cP)^{op}\to \cP$ is a Cartesian fibration with weakly contractible fibers. 

\item[(ii)]
The functor $p_{\cP}: \sfC(\cP)^{op}\to \cP$ is cofinal.
\end{itemize}
\end{lemma}
\begin{proof}
(i) For any arrow $c<c'$ in $\cP$ and any chain $C=(c^1\gneqq \cdots\gneqq c^k=c')$ over $c'$, the Cartesian arrow over $c< c'$ and ending at $C$ is $C':=(c^1\gneqq \cdots\gneqq c^k=c'\gneqq c)\succ C$ if $c\lneqq c'$ and is $id_C$ if $c=c'$. Since each fiber contains a final object, the fibers are weakly contractible.

(ii) By \cite[Theorem 4.1.3.1]{Lu1}, it is equivalent to show that $\sfC(\cP)^{op}\times_{\cP}\cP_{c/}$ is weakly contractible, for any fixed  $c\in \cP$. Consider the natural inclusion $\sfC(\cP)^{op}_{\text{end}=c}:=\{(c^1\gneqq\cdots\gneqq c^k): c^k=c\}\subset\sfC(\cP)^{op}\times_{\cP}\cP_{c/}=\{(c^1\gneqq\cdots\gneqq c^k): c^k> c\}$. First, this inclusion admits a right adjoint by adding $c$ at the end to any chain $(c^1\gneqq\cdots\gneqq c^k)$, if $c^k\gneqq c$ and doing nothing if $c^k=c$. Second, $\sfC(\cP)^{op}_{\text{end}=c}$ admits a final object, hence it is weakly contractible. Therefore, $\sfC(\cP)^{op}\times_{\cP}\cP_{c/}$ is weakly contractible as well. 
\end{proof}

Let us put $\ovl{\cY}_I=\ovl{\cT}_G$. Then the sectorial inclusions induce a natural functor 
\begin{align}
\label{eqPhiIscC}\Phi_{I}^{\cY}: \sfC(\cP(I))^{op}&\longrightarrow  \CatEx\\
\nonumber \vec{J}=(J^1\supsetneq \cdots\supsetneq J^k)&\mapsto \cW(\ovl{\cY}_{\vec{J}}, \vartheta_{\ovl{\cT}_G}|_{\ovl{\cY}_{\vec{J}}}). 
\end{align}
This functor $\Phi_{I}^{\cY}$ is actually \emph{canonical} if we consider a cofinal sequence of induced good admissible coverings $\frU_H^{(n)}$ of $H^{sm}$ and the conex nonempty space of $(\delta,\frU_H^{(n)})$-rigid vector fields $Z^{sm}$ (and apply Lemma \ref{lemmasemitoadm} (ii)). 
By Lemma \ref{lemmapCScofinal} (i) and that the sectorial inclusion \eqref{eqcYJ1JkJ1} induces an equivalence on wrapped Fukaya categories, the functor $\Phi_I^{\cY}$ naturally factors as $\Phi^{\cY, \prime}_{I}\circ p_{\cP(I)}$, where 
\begin{align}
\label{eqPhiprimeI}\Phi^{\cY,\prime}_{I}: \cP(I)&\longrightarrow  \CatEx\\
\nonumber J&\mapsto \cW(\ovl{\cY}_J, \vartheta_{\ovl{\cT}_G}|_{\ovl{\cY}_J}).
\end{align}
is the right Kan extension of $\Phi_I^\cY$ along $p_{\cP(I)}$. Let $\Phi_{I}^{\cY,\prime;op}$ be the opposite diagram of $\Phi^{\cY,\prime}_{I}$ in which the target categories are replaced by the corresponding Ind-completion of wrapped Fukaya categories and the target functors are given by restrictions. 

We have the complete analog for a reductive $G'$ as in \S\ref{sssdeltaadmReductive}, where we take the product of the above diagram for $G'_\der$ with $T^*Z(G')_{\cpt}\times T^*\ovl{Z(G')}_{0,\RR}$, and quotient it by $Z(G'_\der)$. 
For any $J\subsetneq I'$, we will use notations $\ovl{\cY}^{G'}_J$, to denote the corresponding subsectors.

\sss{The diagram of wrapped Fukaya categories under sectorial inclusions in $\ovl{\cT}_G$}\label{sssYJcJLJ}

By Lemma \ref{lemmadeltaadmCsector} (iii) for the group $L_J^\der$, we have the canonical equivalences (up to contractible spaces of choices)
\begin{align}\label{eqcYJWcTLJ}
\cW(\ovl{\cY}_J)\simeq \cW(\ovl{\frF}_J)\simeq \cW(\ovl{\cT}_{L_J}), \quad J\subsetneq I,
\end{align}
where $\ovl{\cT}_{L_J}=\ovl{\cT}_{L_J^\der}\times^{Z(L_J^\der)}T^*Z(L_J)_{\cpt}\times T^*\ovl{Z(L_J)}_{0,\RR}$ is equipped with the product sector structure and $\ovl{Z(L_J)_{0,\RR}}$  is a compactification of $Z(L_J)_{0,\RR}$ by a compact ball with smooth boundary. By the main result in \cite{J}, we have $\cW(\ovl{\cY}_J)\simeq \Coh(\cS_{L_J^\vee})$.

For $J\subsetneq I$ and any Kostant section $S_J=S^{1}_{L_J^\der}\times T^*_{z_0}Z(L_J)_{\cpt}\times T^*_{\wt{\bq}_0}Z(L_J)_{0,\RR}$ in $\cT_{L_J}$, where $S^{1}_{L_J^\der}$ is any Kostant section in $\cT_{L_J^\der}$,  $z_0\in Z(L_J)_{\cpt}$ and $\wt{\bq}_0\in Z(L_J)_{0,\RR}$, let 
\begin{align}\label{eqbLJ}
\bL^J= S^{1}_{L_J^\der}\times T^*_{z_0}Z(L_J)_{\cpt}\times (\frF_J^{\flat, \cI})_{\bq_0}\times \bL'_{\CRez}\subset \cY_J,
\end{align}
under the splittings \eqref{eqfrFJ} and \eqref{eqovlcYJ}, where $(\frF_J^{\flat, \cI})_{\bq_0}$ is the cotangent fiber at any point $\bq_0\in \cI_J$ under the canonical isomorphism $\frF_J^{\flat, \cI}\cong T^*H_J^{\cI}$, and $\bL'_{\CRez}$ is any cylindrical Lagrangian contained in $\ovl{\cP}_{\theta_0, \leqsl A}$ that generates $\cW(\ovl{\cP}_{\theta_0, \leqsl A})\simeq \Mod(k)$. 
Then $\bL^J$ is a ``cylindricalization" of $S_J$ in $\cY_J$, in the sense that under the canoincal equivalence \eqref{eqcYJWcTLJ}, 
the image of $\bL^J$ is \emph{canoincally} isomorphic to $S_J$. For $J=I$, set $\bL^I=S_I$. 

For any $I\supset J^1\supset J^2$, let 
\begin{equation}\label{eqAppIndWadjunc}
\begin{tikzcd}
\Ind\cW(\ovl{\cY}_{L_{J^1}})\ar[r, shift right=0.2em, "\res^{J^1}_{J^2}"']&\Ind\cW(\ovl{\cY}_{J^2})\ar[l, shift right=0.2em, "\cores^{J^1}_{J^2}"']
\end{tikzcd}.
\end{equation}
denote the restriction and co-restriction functors induced by the 
sectorial inclusion $\ovl{\cY}_{J^2}\hookrightarrow \ovl{\cY}_{J^1}$.

\begin{prop} [\cite{J}]\label{propcPIcWcY}
Assume $G$ is adjoint.

\begin{itemize}
\item[(i)]
Then for any $J\subset I$, any $\bL_J$  (with the canonical brane structure)  generates $\cW(\ovl{\cY}_J)$ and 
there is a canonical identification
\begin{align*}
\End(\bL^J)= \cO(\cS_{L_J^\vee})=\cO(T^\vee\sslash W_J).
\end{align*}
In particular, we get a natural equivalence 
\begin{align}
\Upsilon_{\bL^J}: \cW(\ovl{\cY}_J)&\overset{\sim}{\longrightarrow} \Coh(\cS_{L_J^\vee})\\
\nonumber \bL^J&\mapsto \cO_{\cS_{L_J^\vee}}. 
\end{align}

\item[(ii)]
For any $I\supset J^1\supsetneq J^2$, 
the restriction functor $\res^{J^1}_{J^2}(\bL^{J^1})$ is isomorphic to any $\bL^{J^2}$, and the morphism
\begin{align*}
\cO(T^\vee\sslash W_{J^1})=\End_{\cW(\cT_{L_{J^1}})}(\bL^{J^1})\lrar \End_{\cW(\cT_{L_{J^2}})}(\res_{J^2}^{J^1}\bL^{J^1})=\End(\bL^{J^2})=\cO(T^\vee\sslash W_{J^2})
\end{align*}
is the natural inclusion of invariant functions.

\end{itemize}
\end{prop}

\sss{A diagram of deformations of Liouville sectors and sectorial inclusions}\label{sssSummdeform}

We have the diagram of Liouvlle sectors, for any given $1<R\leqsl N_1$: 
\begin{equation*}
\begin{tikzcd}
(\ovl{\cT}_G, \vartheta_0')\ar[r, rightsquigarrow, " \vartheta_t' "]&(\ovl{\cT}_G, \vartheta_1')\ar[r,  rightsquigarrow,  "(\vartheta_1')_t"]&(\cT_G^{\ovl{\cC}}, \vartheta_1'|_{\cT_G^{\ovl{\cC}}}=\vartheta_I^1)\\
\ovl{\cY}_J=\ovl{\frF}_J\times \ovl{\cP}_{\theta_0,\leqsl A}\ar[u, hook]&\ovl{\cY}_{J; 0}'=\ovl{\frF}_J\times (\CC_{x\in [0, R]},  \vartheta'_{\CRez; 1})\ar[u, hook]&\ovl{\cY}_{J; 1}'=\ovl{\frF}_J\times (\CC_{x\in [1, R]}, \vartheta'_{\CRez; 1})\ar[u, hook]
\end{tikzcd}. 
\end{equation*}
where 
\begin{itemize}
\item the first row represents deformations of Liouville sector structures on $\ovl{\cT}_G$: $\vartheta_t'$ and $(\vartheta_1')_t$ are from \S\ref{sssDeformLiouvillesector} and \S\ref{sssSectJGovlC} for $\cT_G$ (with trivial change);

\item the second row represents Liouville subsectors for $J\subsetneq I$: cf. \S\ref{sssReviewSectIncl}; by Lemma \ref{lemmadeltaadmCsector} (ii), we have $\vartheta_I^1$ is $\delta$-admissible and $\ovl{\cY}'_{J;1}=\cT^{\ovl{\cU^\dagg_J}}_{L_J}$, for some $\cU_J^\dagg=\nu_J(\cC_J\times\cI_{J}^\dagg)$ with $\cI_{J}^\dagg\subset \cI_J$, is a subsector of $\cT_G^{\ovl{\cC}}$.  
\end{itemize}
Here only ``trival" deformations of $\CRez$ are involved. 
We have similar situations for a reductive $G'$ as in \S\ref{sssdeltaadmReductive}, where we take the product of the above diagram for $G'_\der$ with $T^*Z(G')_{\cpt}\times T^*\ovl{Z(G')}_{0,\RR}$ and quotient it by $Z(G'_\der)$. 
For any $J\subsetneq I'$, we will use  notations $\ovl{\cY}^{G'}_J, \ovl{\cY}^{G',\prime}_{J;i}, i=0,1$, to denote the corresponding subsectors.

\sss{The diagram of wrapped Fukaya categories under sectorial inclusions in $\cM_G$}\label{sssAppFukSecIncl}

Fix a good admissible open covering $\hat{\frU}=\{(\hat{\cC}_J, \hat{\cI}_J)\}_{J\in \PfI}$ of $\trg$. Let $\vartheta$ be a $(\delta, \hat{\frU})$-rigid Liouville 1-form on $\cM_G$. Then the natural functor associated with sectorial inclusions of the sectorial covering $\cT_{L_J}^{\ovl{\hat{\cU}_J}}, J\in \PfI$, naturally factors as (on the respective nerves): 
\begin{equation}\label{eqPhibulletN}
\begin{tikzcd}[column sep=6em]
\Phi_\bullet:N_\bullet(\sfC(\PfI)^{op})\ar[r] &\cE_{(\frs_{\delta,\rigid}^{\hat{\frU}})_\bullet}\ar[r, "(\cW^{\cE, \hat{\frU}}_{\delta,\rigid})_\bullet"]&\Maps([\bullet], \CatEx)\\
\begin{matrix}
\vec{J}_{(1)}=(J^1\supsetneq \cdots \supsetneq J^m)\\
\vee\\
\vdots\\
\vee\\
\vec{J}_{(\bullet)}
\end{matrix}\ar[r, mapsto]& 
(\cT_{L_{p(\vec{J}_{(i)})}}^{\ovl{\cU}_{\vec{J}_{(i)}}}, \vartheta_{\vec{J}_{(i)}}:=\vartheta|_{\cT_{L_{p(\vec{J}_{(i)})}}^{\ovl{\cU}_{\vec{J}_{(i)}}}})_{1\leqsl i\leqsl \bullet}
\ar[r, mapsto] &\cW(\cT_{L_{p(\vec{J}_{(i)})}}^{\ovl{\cU}_{\vec{J}_{(i)}}}, \vartheta_{\vec{J}_{(i)}})_{1\leqsl i\leqsl \bullet}
\end{tikzcd}
\end{equation}
where $p:=p_{\PfI}$ (cf. \S\ref{sssdigposets}) and $\cU_{\vec{J}_{(i)}}=\nu^{-1}_{p(\vec{J}_{(i)})}(\hat{\cC}_{p(\vec{J}_{(i)})}\times\hat{\cI}_{\vec{J}_{(i)}})$. Using Lemma \ref{lemmacWfrUrigid}, this is canonically equivalent to 
\begin{align*}
N_\bullet(\sfC(\PfI)^{op})\longrightarrow N_\bullet(\PfI)\overset{(\cW^{\hat{\frU}}_{\delta, \rigid})_\bullet}{\longrightarrow}\Maps([\bullet], \CatEx).
\end{align*}
By Lemma \ref{lemmapCScofinal}, we have a canononical equivalence
\begin{align}\label{eqPhiWrigU}
\colim_{\sfC(\PfI)^{op}}\Phi\simeq \colim_{\PfI}\cW^{\hat{\frU}}_{\delta, \rigid}. 
\end{align}

\begin{lemma}\label{lemmaAppcWJicY}
For any $(J^1\supset \cdots\supset J^k)$ in $N_k(\PfI)$ and any $\PfI\ni J\supset J^1$, there are canonical equivalence of diagrams
\begin{align*}
(\cW^{\hat{\frU}}_{\delta, \rigid})_k(J^1\supset \cdots\supset J^k)\simeq \left(\cW(\ovl{\cY}^{L_J}_{J^i})\right)_{1\leqsl i\leqsl k}\overset{\eqref{eqcYJWcTLJ}}{\simeq} \left(\cW(\ovl{\cT}_{L_{J^i}}, \vartheta_{\ovl{\cT}_{L_{J^i}}})\right)_{1\leqsl i\leqsl k}
\end{align*}
up to a contractible space of choices, where $\ovl{\cY}^{L_J}_{J^i}$ is the subsector in $(\ovl{\cT}_{L_J}, \vartheta_0')$  defined using any admissible hypersurface $H^{sm}_{J}$ in $\RR_{\geqsl 0}^J$ and any $(\delta, \cP(J))$-rigid vector field on $H^{sm}_{J}$. Moreover, the canonical equivalence are compatible with the morphisms $N_k(\PfI)\to N_\ell(\PfI)$ and  $N_k(\PfI_{\subset J})\to N_\ell(\PfI_{\subset J})$ induced by any $f:[\ell]\to [k]$ in $\Delta$. 
\end{lemma}

\begin{proof}
Take any $\vec{J}_{(1)}>\cdots>\vec{J}_{(k)}$ in $N_\bullet(\sfC(\PfI)^{op})$ such that $p(\vec{J}_{(i)})=J^i$ and each $\vec{J}_{(i)}$ contains $J$. 
Then the subsectors $(\cT_{L_{J^i}}^{\ovl{\cU}_{\vec{J}_{(i)}}}, \vartheta_{\vec{J}_{(i)}}), 1\leqsl i\leqsl k$, are subsectors of $(\cT^{\ovl{\hat{\cU}_J}}_{L_J}, \vartheta_{J})$. Thus we get an object in $\frs:=(\frs^{\hat{\frU}}_\delta)_{\cI^\dagg_{J^1},\cdots, \cI^\dagg_{J^k}}$ (and also an object in $(\frs^{\hat{\frU}}_\delta)_{\hat{\cI}_J, \cI^\dagg_{J^1},\cdots, \cI^\dagg_{J^k}}$), where $\cI^\dagg_{J^i}:=\hat{\cI}_{\vec{J}_{(i)}}$. 
On the other hand, by shrinking $\cI^\dagg_{J^i}$ if necessary, we get another object in $\frs$ by using the right-most column of the diagram in \S\ref{sssSummdeform} for the reductive group $L_J$ in place of $G$ (and $\ovl{\hat{\cU}}_J$ in place of $\ovl{\cC}$). 
The first part of the lemma then follows from the contractibility of $\frs$, the (trivial) deformations in the diagram in \S\ref{sssSummdeform} for the reductive group $L_J$, and the canonical functor $\Phi_{J}^{\cY}$ \eqref{eqPhiIscC} for $L_J$. 
The second part follows from the functoriality of $(\frs_{\delta,\rigid}^{\hat{\frU}})_\bullet$ \eqref{eqcEk} and Lemma \ref{lemmafrsdeltaconvex}. 
\end{proof}

Finally, we remark that $\Phi$ and $\cW^{\hat{\frU}}_{\delta, \rigid}$ for different choices of $\hat{\frU}$ and $(\delta,\hat{\frU})$-rigid $\vartheta$ are canonically equivalent up to a contractible space of choices: one can choose a cofinal sequence of good admissible coverings $\hat{\frU}^{(n)}$ of $\trg$ and use the nonempty and convex space of $(\delta,\hat{\frU}^{(n)})$-rigid $\vartheta$ (and apply Lemma \ref{lemmasemitoadm} (ii)).

\subsection{Constant integer gradings on holomorphic Lagrangians}\label{appconstantgrading}

For any complex symplectic vector space $(V^{2n}, \omega_\CC)$, let $\cH(V, \omega_\CC)$ be the space of quaternionic structures $(I=I_V, J, K)$ on $V$ that is \emph{compatible} with $\omega_{\CC}$ in either of  the following two \emph{equivalent} senses:
\begin{itemize}
\item[(i)] let $g(v_1, v_2)=\Real \omega_\CC(v_1, Jv_2)$, then $(I=I_V, J, K; g)$ is a linear hyperkahler structure on $V$ with hyperkahler metric $g$, satisfying $\omega_\CC=g(-, J-)+ig(-,K-)$;\\
\item[(ii)] $\Real \omega_{\CC}(Jv_1,Jv_2)=\Real\omega_{\CC}(v_1, v_2)$ and $\Imagine \omega_{\CC}(Jv_1,Jv_2)=-\Imagine \omega_{\CC}(v_1, v_2)$, and $g=\Real\omega_{\CC}(-, J-)$ is a metric.   
\end{itemize}

It is easy to check using $\omega_\CC(I_V-,-)=i\omega_\CC(-,-)=\omega_{\CC}(-, I_V-)$ that (ii) implies that $K$ is compatible with $\Imagine \omega_\CC$ and $g(I_V-, I_V-)=\Real\omega_{\CC}(I_V -,JI_V-)=g(-,-)$. 

\begin{lemma}\label{lemmacHVomegaCC}
Given $(V^{2n}, \omega_\CC)$, the space $\cH(V, \omega_\CC)$ is diffeomorphic to $\Sp(2n,\CC)/\Sp(n)$, where $\Sp(n)=\Sp(2n,\CC)\cap U(2n)=\Sp(2n,\CC)\cap O(4n)$ is the compact form of $\Sp(2n,\CC)$. In particular,  $\cH(V, \omega_\CC)$ is contractible. 
\end{lemma}
\begin{proof}
The proof is similar to the proof that the space of compatible linear complex structures $\cJ(E, \omega)$ of a given real symplectic vector space $(E^{2n}, \omega)$ is diffeomorphic to $\Sp(2n,\RR)/U(n)$ (cf.\cite[\S 2.5]{MS}). We omit the details. 
\end{proof}

For a holomorphic symplectic manifold $(M, \omega_{\CC})$, we can talk about \emph{compatible (almost) quaternionic structures} $(\II=\II_M, \JJ, \KK)$, where $\II_M$ is the present complex structure on $M$, $\JJ, \KK$ are almost complex structures, and for each $x\in M$, $(\II=\II_M, \JJ, \KK)_x\in \cH(T_xM, (\omega_\CC)_x)$. 
Denote the space of compatible (almost) quaternionic structures by $\cH(M, \omega_{\CC})$. 

A direct consequence of Lemma \ref{lemmacHVomegaCC} is that 
\begin{cor}
For any holomorphic symplectic manifold $(M, \omega_{\CC})$, the space $\cH(M, \omega_{\CC})$ is contractible. 
\end{cor}

For any $(\II=\II_M, \JJ, \KK)\in \cH(M, \omega_\CC)$, let $g$ be the metric on $M$ defined by $\Real \omega_\CC(-, \JJ -)$ and let $\omega_{\II}=g(-,\II -)$.  Let $\omega_\JJ=\Real  \omega_\CC$ and $\omega_\KK=\Imagine \omega_{\CC}$. Then $\omega_\JJ=g(-, \JJ-)$ and $\omega_\KK=g(-, \KK-)$. 
Using the calculation from \cite[Proposition 5.1]{JholLag}, one gets 

\begin{lemma}\label{lemmaintegergrading}
Let $(\II=\II_M, \JJ, \KK)\in \cH(M^{2n}, \omega_\CC)$, and let $\alpha_M=(\Lambda^{n}(\omega_\KK+i\omega_\II))^{\otimes 2}$ be the trivialization of $\kappa^{\otimes 2}$ ($\kappa$ as in \eqref{eqkappaTMJJ}). Then  for any (connected) holomorphic Lagrangian $\sfL_{\hol}$, $\alpha_M^{\sfL_{\hol}}\equiv 1\in S^1$, and we can choose the \emph{canonical} grading of $\sfL_{\hol}$ to be constantly $0$. 
\end{lemma}

%%%%%%%%%%%%%%%%%%%%%%%%%%%%%%%

\section{Supplementary materials to \S\ref{secProofMirror}}

\subsection{The sheaf $\pi\frZ$ over $\trg$}\label{subsecfrZ}
We use notation from \S\ref{ss:def MG}. Let $\triangle$ be the fundamental alcove attached to the Iwahori subgroup $\bI_0\subset L^+G$, viewed as a poly-simplicial complex in the apartment of $T$. We label the facets of $\triangle$ by subsets $J\sft \wt I$ (subsets that span a root system of finite type): the facet $\cF_J$ corresponding to $J$ is the simultaneous vanishing locus (in $\triangle$) of the affine functions defined by $\{\a_j, j\in J\}$. 

For $J\sft \wt I$, let $\bP_J\subset LG$ be the standard parahoric whose Levi subgroup $L_J\supset T$ has simple roots $\{\a_j; j\in J\}$. 
For $J\subset J'$, there is an inclusion of the centers $Z(L_{J'})\subset Z(L_J)$.

Consider the sheaf $\frZ$ of abelian groups on $\triangle$ whose restriction to $\cF_J$ is the constant sheaf with stalk $\frZ_J:=Z(L_J)$. The cospecialization maps $\frZ_{J'}\to \frZ_{J}$ for $J\subset J'$ are given by the inclusions $Z(L_{J'})\subset Z(L_J)$.

Let $\pi\frZ_J$ be the fundamental groupoid of $\frZ_J$, which is canonically equivalent to the quotient groupoid
\begin{equation*}
\pi\frZ_J\cong \ZZ^J/\xcoch(T)
\end{equation*}
where the action of $\xcoch(T)$ on $\ZZ^J$ is via the map of abelian groups
\begin{equation*}
\xcoch(T)\to \ZZ^J, \quad \l\mt (\j{\a_j,\l})_{j\in J}.
\end{equation*}
As $J$ varies, $J\mt \pi\frZ_J$ gives a sheaf of groupoids over $\triangle$, which we denote by $\pi\frZ$. Its restriction to $\cF_J$ is the constant sheaf with stalk $\pi\frZ_J$, and the cospecialization maps $\pi\frZ_{J'}\to \pi\frZ_{J}$ for $J\subset J'$ are induced by the inclusion $Z(L_{J'})\subset Z(L_J)$.

Let $\partial_j\triangle\subset \triangle$ be the closure of the face $\cF_j$ corresponding to $j\in \wt I$, with inclusion $\io_j: \partial_j\triangle\incl \triangle$. Then $\pi\frZ$ can be presented as a two-step complex
\begin{equation}\label{piZ complex}
    \op_j \j{\a_j,-}:\xcoch(T)\to \op_{j\in \wt I}\io_{j*}\ZZ.
\end{equation}
In other words, $\pi\frZ$ is the quotient of $\op_{j\in \wt I}\io_{j*}\ZZ$ by the translation action of $\xcoch(T)$ via the above map. 

\begin{lemma}\label{l:cohoZ}
We have
\begin{equation*}
    \cohog{i}{\triangle,\pi\frZ}=\begin{cases}
        \ZZ^{\wt I}/\xcoch(T), & i=0,\\
        \xcoch((ZG)^\c) & i=-1,\\
        0 & \textup{otherwise}.
    \end{cases}
\end{equation*}
\end{lemma}
\begin{proof}
Since $\triangle$ and its faces are contractible, the long exact sequence of cohomology associated with \eqref{piZ complex} reads
\begin{equation*}
    0\to \cohog{-1}{\triangle,\pi\frZ}\to \xcoch(T)\xr{\op\j{\a_j,-}} \ZZ^{\wt I}\to  \cohog{0}{\triangle,\pi\frZ}\to 0
\end{equation*}
The lemma follows immediately.
\end{proof}

We will see that the cohomology of $\pi\frZ$ controls the existence and uniqueness of compatible systems of Kostant sections in the open subsets $\cT_{L_J}$ of $\cM_G$.

\subsection{The closure inclusion relations among (the union of) Kostant sections}\label{sssClosureKostant}

For each $J\in \cP(I_H)$, let 
\begin{align*}
\SS^{L_J}&=(b_{\RR}^J)^{-1}(0^J\times Z(L_J)_{0;\RR})\subset \cT_{L_J}\subset \cT_H,\\
&\cong S_{L_J^\ad}\times T^*Z(L_J), 
\end{align*}
which is the union of all Kostant sections in $\cT_{L_J}$. 
Let $\frK_J:=\pi_0(\SS^{L_J})$ (or $\frK_{J;L_J}$ if we want to specify the group $L_J$) be the collection of Kostant sections in $\cT_{L_J}$ modulo $Z(L_J)^\circ$, which is a torsor over $\pi_0(Z(L_J))$. Let $\cX^{L_J}=\chi_{L_J}^{-1}([0]\times \{0\})$, where $\chi_{L_J}: \cT_{L_J}\to \frl_J\sslash L_J=\frl_J^\der\sslash L_J^\der\times \frz_J$ is the characteristic map. There is an obvious canonical identification between $\pi_0(\cX^{L_J})$ and $\frK_J$. For any $\kappa_J\in \frK_J$, let $\SS^{L_J}_{\kappa}$ (resp. $\cX^{L_J}_{\kappa}$) denote the corresponding connected component in $\SS^{L_J}$ (resp. $\cX^{L_J}$).

\begin{lemma}\label{lemmapifrKJpJ}
\begin{itemize}
\item[(i)] For any  $J'\supsetneq J$ and $(\kappa', \kappa)=(\kappa_{J'}, \kappa_J)\in \frK_{J'}\times \frK_J$, the following statements are equivalent: 
\begin{itemize}
\item[(a)] $\SS^{L_{J'}}_{\kappa'}\cap \ovl{\SS^{L_J}_{\kappa}}\neq \vn$; 

\item[(b)] $\SS^{L_{J'}}_{\kappa'}\subset \ovl{\SS^{L_J}_{\kappa}}$; 

\item[(c)] $\cX^{L_{J'}}_{\kappa'}\cap \SS^{L_J}_{\kappa} \neq \vn$ (in this case, the projection of each irreducible component of this intersection to $\PP_J\cong (L_J)_{\ab}$ is a branched covering); 

\item[(d)] $\cX^{L_{J'}}_{\kappa'}\cap S_J^{\kappa}\neq \vn$ for some Kostant section $S_J^{\kappa}\subset \SS^{L_J}_{\kappa}$.
\end{itemize}

\item[(ii)] For any $J'\supsetneq J$ and $\kappa'\in \frK_{J'}$, there exists a unique $\kappa\in \frK_J$ such that $\SS^{L_{J'}}_{\kappa'}\subset \ovl{\SS^{L_J}_{\kappa}}$. In particular, $\ovl{\SS^{L_J}_{\kappa_1}}\cap\ovl{\SS^{L_J}_{\kappa_2}}=\vn$, for $\kappa_1\neq \kappa_2\in \frK_J$. 
\end{itemize}
\end{lemma}

\begin{proof}
To prove the statements, it suffices to reduce to the case that $H=L_{J'}^{\der}$ is semisimple, in which $J'=I_H$. 
Let $\PP_{H}=\Spec\CC[H/N]^{N^-}$, with the natural stratification by $T$-orbits $\PP_J\cong (L_J)_{\ab}$ indexed by $J\in \cP(I_H)$.  
First, since the map $b:=b^{I_H}: \cT_H\to  \PP_H$ has fibers equidimensional of dimension $r_H=|I_H|$ and $b^{-1}(\PP_J)^{\red}$ has irreducible components $\SS^{L_J}_{\kappa}, \kappa\in \frK_J$, we see that 
$$\ovl{\SS^{L_J}_{\kappa}}\cap\SS^{H} =\bigsqcup_{\wt{\kappa}\in \frK_{I_H}^{\kappa}}\SS^{H}_{\wt{\kappa}},$$
for a subset $\frK_{I_H}^{\kappa}\subset \frK_{I_H}$.

(i) Using the above observation, it is clear that (a) $\Leftrightarrow$ (b).  
Let us show (b) $\Leftrightarrow$ (c). The converse is clear (because of $\GG_m$-invariance of the intersection). To see the forward implication, note the following
\begin{itemize}
    \item we have $\dim (\cX^{H}\cap \ovl{\SS_{\kappa}^{L_J}})^{\red}\geqsl r-|J|$ (this holds for each irreducible component in the intersection), since the intersection is non-empty and $\ovl{\SS_{\kappa}^{L_J}}$ is irreducible of dimension $2r_H-|J|$;
    \item $\cX_{\wt{\kappa}}^{H}\to \PP_H$ is a ($\GG_m$-equivariant) finite morphism 
   and 
    $\ovl{\SS_{\kappa}^{L_J}}\subset b^{-1}(\ovl{\PP}_J).$ 
\end{itemize}
Thus by dimension reasons, $\cX_{\wt{\kappa}}^{H}\cap \SS_{\kappa}^{L_J}\neq \vn$, and the projection of each irreducible component of this intersection to $(L_J)_{\ab}$ is a finite and surjective morphism. 
To see that (c) $\Leftrightarrow$ (d), note that different $S_J^{\kappa}\subset \SS_{\kappa}^{L_J}$ are connected by paths in $Z(L_J)_{0,\cpt}$, thus the intersection number (counted with signs and multiplicities) of $\cX^{H}$ and $S_J^{\kappa}$ is independent of the choice of $S_J^{\kappa}$. In fact, the intersection number is 
\begin{align}
\frac{|W|\cdot |\pi_0 Z(L_J)|}{|W_J|\cdot |Z(H)|},
\end{align}
and in the general case (without the reduction to $H=L_{J'}^\der$ at the beginning), we have 
\begin{align}\label{eqcXJpSJ}
\sharp\cX^{L_{J'}}_{\kappa'}\cap S_J^{\kappa}=\frac{|W_{J'}|\cdot |\pi_0Z(L_J)|}{|W_J|\cdot |\pi_0Z(L_{J'})|}. 
\end{align}
This finishes the proof of (i).

(ii) Under the same setting above, let $L_J^\flat$ denote $L_J/Z(H)\subset H^\ad$. It suffices to show that for any Kostant section $S_J^{\kappa}$, the set of components in $Z(H)\cdot S_J^{\kappa}$ that $\cX^H$ intersects forms a single $\pi_1 Z(L_J^\flat)\cong\ker(Z(H)\to \pi_0 Z(L_J))$-orbit.  In $\cT_{H^\ad}$, we have $\sharp\cX^{H^\ad}\cap S_J^\flat=\frac{|W|}{|W_J|}$, for any Kostant section $S_J^\flat$ indexed by $J$. Since $\cX^{H^\ad}$ (resp. $S_J^\flat$) is contractible, its preimage in $\cT_H$ has $|Z(H)|$ many connected components. By the projection formula, we have 
\begin{align}\label{eqsharpcXzetaSJ}
&\sharp \cX^{H}\cap (\bigsqcup_{\zeta\in Z(H)}\zeta\cdot S_J^{\kappa})=\frac{|W|}{|W_J|}. 
\end{align}
Using \eqref{eqcXJpSJ} for $H=L_{J'}^{\der}$, we directly see that $\cX^H$ intersects exactly $|\pi_1 Z(L_J^\flat)|=\frac{|Z(H)|}{|\pi_0Z(L_J)|}$ many components. These components must form a $\pi_1 Z(L_J^\flat)$-orbit, since the Kostant sections in each $\pi_1 Z(L_J^\flat)$-orbit are contained in the same connected component of $\SS^{L_J}$. 
Just a remark that one can also prove (ii) directly for any $J'\supsetneq J$, without the reduction to $H=L_{J'}^\der$ at the beginning. One only needs to replace $Z(H)$ by $\pi_0Z(L_{J'})$, $\pi_1Z(L_J^\flat)$ by $\ker\left(\pi_0 Z(L_{J'})\to \pi_0Z(L_J)\right)$ and $W$ by $W_{J'}$.

\emph{Alternative proof of (ii).} It suffices to show that $\ovl{\SS_{\kappa_1}^{L_J}}\cap \ovl{\SS_{\kappa_2}^{L_J}}=\vn$ for $\kappa_1\neq \kappa_2\in \frK_J$.
It suffices to show that $\ovl{\SS^{L_J}_{\kappa}}$ contains exactly $|\pi_1 Z(L_J^\flat)|$ many distinct Kostant sections of $\cT_H$. By the equivalence condition (d) in part (i), this is equivalent to showing that $S_J^{\kappa}$ intersects exactly $|\pi_1Z(L_J^\flat)|$ many connected components of in $\cX^{H}$.

Similarly to \eqref{eqsharpcXzetaSJ}, by the projection formula, we have 
\begin{align*}
&\sharp \cX^{H}\cap S_J^{\kappa}=\frac{|W|}{|W_J|}.
\end{align*}
Using \eqref{eqcXJpSJ}, we directly see that $S_J^{\kappa}$ intersects exactly $|\pi_1 Z(L_J^\flat)|$ many connected components of in $\cX^{H}$. It is straightforward to see that these components in $\cX^{H}$ must form a single $\pi_1 Z(L_J^\flat)$-orbit. 
This finishes the proof. 
\end{proof}

\begin{cor}\label{corpi1SJZJovl}
\begin{itemize}
\item[(i)]For any $J\subset I_H$ and $\kappa_J\in \frK_J$, the natural morphism $\pi_1(\SS^{L_J}_{\kappa})\to \pi_1(\ovl{\SS^{L_J}_{\kappa}})$ (for any given base point) gives the same image in $\pi_1(\cT_H)$. 

\item[(ii)] For any $J\in\PfI$ and $\kappa\in \frK_J$, the natural morphism $\pi_1(\SS^{L_J}_{\kappa})\to \pi_1(\ovl{\SS^{L_J}_{\kappa}})$ (for any given base point) gives the same image in $\pi_1(\cM_G^\circ)$. 
\end{itemize}
\end{cor}

\begin{proof}
We give the proof of (ii). The proof of (i) is completely similar. 
It follows easily from Lemma \ref{lemmapifrKJpJ} (ii) that the preimage of $\SS^{L_J}_{\kappa}$ in $\cM_{G_{\sc}}$ has the same number of connected components as that  of the closure $\ovl{\SS^{L_J}_{\kappa}}\subset \cM_G^\circ$. Thus the Corollary easily follows. 
\end{proof}

\sss{The $\pi_0\frZ$-torsor $\frK$}\label{ssspi0frZtorsorfrK}
Now consider $J\in \PfI$.
Let $p_{J', J}: \frK_{J'}\to \frK_J, J'\supset J$ be the map $\k\mapsto \k'$ defined in Lemma \ref{lemmapifrKJpJ}. It is clear that the collection of $p_{J', J}$ satisfies the cocycle condition $p_{J'', J'}\circ p_{J', J}=p_{J'', J}$ for any $J\subset J'\subset J''$, and $p_{J', J}$ is equivariant with respect to $\pi_0Z(L_{J'})$ and $\pi_0Z(L_J)$-actions through the natural quotient map $\pi_0Z(L_{J'})\to \pi_0Z(L_J)$. Therefore, the collection $ \frK_J, J\in \PfI^{op}$ is a torsor of the sheaf of abelian groups $\pi_0\frZ$ (cf. \S\ref{subsecfrZ}).

%%%%%%%%%%%%%%%%%%%%%%%%%%%%%%%%%%%%%%%%%%

\subsection{Regular covering maps of Lagrangian skeleta}\label{sssRegCovSke}

 Now we are in the setting of \S\ref{sssCompletioncTGV}. For any $J\in \PdIH$, let $\cU_J=\hat{\cU}_J\cap \cV$, and then $\ovl{\cU}_J=\ovl{\hat{\cU}}_J\cap \ovl{\cV}$. Let $\cU_{I_H}=\cV$. We have subsectors $\cT_{L_J}^{\ovl{\cU}_J}$ of $\cT_H^{\ovl{\cV}}$. 
 Fix any $H^\flat$ with $Z(H^\flat)$ connected. For any $J\subset I_H$, let $L_J^\flat=L_J/\ker(\nu_\flat)$. 
 We have the regular finite covering $\pi_{\cT, J}: \cT_{L_J}^{\ovl{\cU}_J}\to \cT_{L_J^\flat}^{\ovl{\cU}_J}$.

\sss{The Lagrangian skeleton $\Lambda_H$ of $\cT_H^{\ovl{\cV}}$}

Let $\Lambda_H$ be the Lagrangian skeleton of $\cT^{\ovl{\cV}}_{H}$ contained in the Liouville completion $\wh{\cT}_{H}^{\ovl{\cV}}$ of $\cT_{H}^{\ovl{\cV}}$.   As in \S\ref{sssCompletioncTGV}, we identify $\wh{\cT}_H^{\ovl{\cV}}$ with $\cT_H$ (equipped with a non-standard Liouville 1-form). Write $\Lambda_H=\Lambda_{H^\der}\times^{Z(H^\der)}Z(H)_{\cpt}\times Z(H)_{0,\RR}$. For any subset $S\subset \RR^{I_H}_{\geqsl 0}\times \frz_{\frh;\RR}$, let $\Lambda_H^S:=\Lambda_H\cap \cT_H^{S}$. Similar notation applies to $H^\der$ in place of $H$. 
Then we have the identification of Lagrangian skeleta
\begin{align*}
\Lambda_{H^\der}\cong \Lambda_{H^\der}^{\cC}\cup_{\Lambda^{\partial^\circ\cC}_{H^\der}\times (-\ep, \ep)} \left(\Lambda^{\partial^\circ\cC}_{H^\der}\times (-\ep,\infty)\right).
\end{align*}
In the following, we will mostly look at $\Lambda_H^{\cV}=\Lambda_H^\cC\times^{Z(H^\der)}Z(H)_{\cpt}\times \cI$ instead of the entire $\Lambda_H$, since it captures all the information of $\Lambda_H$ as a Lagrangian skeleton and the inclusion $\Lambda_H^{\cV}\hookrightarrow \Lambda_H$ clearly induces an equivalence on $\mu\Sh$-categories. We will say $\Lambda_H^{\cV}$ is the Lagrangian skeleton of $\cT_H^{\ovl{\cV}}$ up to completion.

\sss{The choice of base point and the regular covering map $\pi_{\Lambda, H}$}\label{sssstdbasepointcovering}

Fix the base point $z_\vn^\flat:=1\in T^\flat$ (resp. $z_\vn:=1\in T$) in the zero section of $T^*T^\flat$ (resp. $T^*T$). Without loss of generality, assume $z_\vn\in \Lambda_H^\cV$. 
By Lemma \ref{lemmapi1cMG} and the homotopy equivalence $\Lambda_{H^\flat}^\cV\overset{h.e.}{\hookrightarrow}\cT_{H^\flat}$, the natural map between pointed spaces
\begin{align*}
\pi_{\Lambda, H}: (\Lambda_H, z_\vn)\to (\Lambda_{H^\flat}, z_\vn^\flat)
\end{align*}
corresponds to $K_{H,\flat}:=\ker\nu_\flat$ as a quotient of $\pi_1(\Lambda_{H^\flat}, z_\vn^\flat)$: 
\begin{align}\label{eqKJflattoPpi1cT}
&q_{H, \pi_1}: \pi_1(\Lambda_{H^\flat}, z_\vn^\flat)\cong \pi_1(H^\flat, 1)\twoheadrightarrow K_{H,\flat}=\pi_{\Lambda, H}^{-1}(z_\vn^\flat). 
\end{align}

\sss{The Lagrangian skeleton $U_J$ of $\cT_{L_J}^{\ovl{\cU}_J}$ (up to completion)}\label{sssUJcTLJ}

For any subsector $\scU$ of a Weinstein sector $\ovl{X}$, let 
\begin{align}\label{eqUscU}
U_{\scU}:=(\scU-\partial^{\fin}(\scU))\cap \Lambda 
\end{align}
be the open subset of the Lagrangian skeleton $\Lambda$ of $\ovl{X}$. 
For any $J\in \PdIH$, let $\Lambda_{L_J}^{\cU_J}:=U_{\cT_{L_J}^{\ovl{\cU}_J}}\subset \Lambda_H^\cV$.
Similarly define $\Lambda_{L^\flat_J}^{\cU_J}$. Note that similarly to $\Lambda_H^{\cV}$, $\Lambda_{L_J}^{\cU_J}$ is the Lagrangian skeleton of $\cT_{L_J}^{\ovl{\cU}_J}$ up to completion. 
To simplify notations, we will use $U_J$ (resp. $U_J^\flat$) to denote $\LamJ$ (resp. $\LamJflat$) when $L_J$ and $\cU_J$ (resp. $L_J^\flat$ and $\cU_J$) are clear from the context.

\sss{The set of homotopy classes of paths $\hP_J^\flat$ and $\hP_{H, J}^\flat$}\label{sssPJflatPHJ}

Without loss of generality, assume $z_\vn^\flat\in  U_{\vn}^\flat=(T^\flat)_{\cpt}\times \cU_\vn\subset T^*_{T^\flat}T^\flat$. For each $\vn\neq J\in \PdIH$,  choose a base point $\yvnJ\in U_J$ that is contained in $U_{\vn}\cap T_{\RR}\cap U_J$, and let $\yvnJf=\pi_{\Lambda,H}(\yvnJ)$. Let $y_\vn^{I_H}=z_\vn$. 
By the condition on $\yvnJf$, there is a contractible space of paths connecting from $\zvnf$ to $\yvnJf$ within $U_\vn\cap T_\RR$. Let $\varrho^\flat_{\vn,J,\RR}$ be such a path. 
The restriction of $\pi_{\Lambda, H}$ gives the regular covering map of pointed spaces:
\begin{align*}
\pi_{\Lambda, J}: (U_J, \yvnJ)\to (U_J^\flat, \yvnJf), \quad J\in \cP(I_H).
\end{align*}

Let $U_J^{\flat, \SS}=U_J^\flat\cap \SS^{L_J^\flat}$, and similarly define $U_J^{\SS}$.  By construction of $\Lambda_{H^\flat}$, $(\Lambda_{H^\flat}, U_J^{\flat, \SS})\hookrightarrow (\cT_{H^\flat}, \SS^{L_J^\flat})$ is a homotopy equivalence of pairs. Choose a base point $\yJf$ for each $\SS^{\LJf}$ (note that $\SS^{L_J^\flat}$ is connected), $J\in \cP(I_H)$. 
Let $\hP_{J}^\flat$ (resp. $\hP_{H,J}^{\flat}$) be the set of homotopy classes of paths $\gamma_{J}^\flat$ on $U_{J}^\flat$ (resp. on $\Lambda_{H^\flat}$), starting from $\yvnJf$ (resp. $\zvnf$) and ending at $\yJf$.  
We have the canonical isomorphisms
\begin{align}\label{eqpi1UJyJhPJflat}
\pi_1(U_J, y_J)\backslash \hP_J^\flat&\overset{\sim}{\longrightarrow} \pi_1(\Lambda_H, y_\vn)\backslash \hP_{H, J}^\flat\cong \pi_{\Lambda,J}^{-1}(\yJf)\\
\nonumber [\gamma_J^\flat]&\mapsto [\varrho^\flat_{\vn,J,\RR}\bullet \gamma_J^\flat],
\end{align}
where (1) $\pi_1(U_J, \yvnJ)$ (resp. $\pi_1(\Lambda_H, \zvn)$) acts on $\hP^\flat_{J}$ (resp. $\hP_{H, J}^\flat$) by pre-composing its image based loops at $\yvnJf$ (resp. $\zvnf$) under $\pi_{\Lambda, H}$; (2) the rightmost isomorphism is through the unique lifting of each path in $\hP_{H,J}^\flat$ to a path in $\Lambda_H$ starting from $\zvn$.

\sss{$\frK_J$ as a double coset of homotopy classes of paths on $U_J^\flat$}

Recall the setting of \S\ref{ssspi0frZtorsorfrK}. For any $J^1\supset J^2$, since $\SS^{L_{J^1}^\flat}\subset \ovl{\SS^{L_{J^2}^\flat}}$, there exists a path $\rho^\flat$ connecting from $y_{J^1}^\flat$ to $y_{J^2}^\flat$ within $\ovl{\SS^{L_{J^2}^\flat}}\cap \Lambda_H$ (since $\ovl{\SS^{L_{J^2}^\flat}}$ deformation retracts onto $\ovl{\SS^{L_{J^2}^\flat}}\cap \Lambda_H$ under the inverse of the Liouville flow). 

\begin{lemma}\label{lemmafrKJ1cPJ1}
\begin{itemize}
\item[(i)] For any $J\in \cP(I_H)$, $\frK_{J}$ is canonically bijective to the double cosets 
\begin{align*}
\pi_1(U_J, \yvnJ)\backslash \hP^\flat_{J}/\pi_1(U_J^{\flat, \SS},\yJf)\cong \pi_1(\Lambda_{H}, z_\vn)\backslash \hP^\flat_{H, J}/\pi_1(U_J^{\flat, \SS}, \yJf)\cong \pi_{\Lambda,J}^{-1}(\yJf)/\pi_1(U_J^{\flat, \SS}, \yJf),
\end{align*}
where $\pi_1(U_J^{\flat, \SS}, \yJf)$ acts on $\hP^\flat_{J}$ (resp. $\hP_{H,J}^\flat$) by post-composing based loops at $\yJf$, and it acts on $\pi_{\Lambda,J}^{-1}(\yJf)$ by the standard action from lifting based loops at $\yJf$. 

\item[(ii)] Let $J^1\supset J^2$ and $\gamma^\flat_{J^1}\in \hP_{J^1}^\flat$. 
Under the identification in (i),
the class $p_{J^1,J^2}([\gamma^\flat_{J^1}])$ is 
\begin{align*}
[\varrho^\flat_{\vn,J^1,\RR}\bullet\gamma^\flat_{J^1}\bullet \rho^\flat]\in \pi_1(\Lambda_{H}, z_\vn)\backslash \hP^\flat_{H, J^2}/\pi_1(U_{J^2}^{\flat, \SS}, y_{J^2}^\flat). 
\end{align*}
\end{itemize}
\end{lemma}
\begin{proof}
Both statements are straightforward from the theory of regular covering maps (between pointed spaces). In (ii), note that the class on the last line is indeed independent of the choice of $\rho^\flat$, in view of Corollary \ref{corpi1SJZJovl}. 
\end{proof}

%%%%%%%%%%%%%%%%%%%%%%%%%%%%%%%%%%%%%%%%%%

\subsection{The tensor action of local systems on $\cW(\cT_{\Hf})$}\label{appendsubseclocalsystemaction}

\sss{Geometric composition of Lagrangian correspondences of stopped Liouville manifolds}\label{sssGeometricComposition}
We will consider Lagrangian correspondences of stopped Liouville manifolds coming from exact symplectic morphisms that preserve the stops (and their transpose). Such a Lagrangian correspondence $\bL\subset (X,\frf_X)^-\times (Y, \frf_Y)$ will hit the stop of the product at $\infty$. It is standard to use a small generic negative Hamiltonian perturbation of $\bL$ to avoid the stop, which however does unpleasant things on the geometric composition. 

We comment on how to easily avoid such an issue.  
One can replace $(X, \frf_X)$ (resp. $(Y, \frf_Y)$) by its stabilization $(\wh{X}, \frf)^{\frs}=(\wh{X}, \frf)\times (\CC_{z}, \{\pm\infty\})$ (resp. $(Y, \frf_Y)^{\frs}$) and one replaces 
$\bL$ by $\bL\times \Delta_{(\CC_{z}, \{\pm\infty\})}$. Then one just needs to perturb $\Delta_{(\CC_{z}, \{\pm\infty\})}$ by a small generic negative Hamiltonian. This can be easily realized by replacing $\Delta_{(\CC_{z}, \{\pm\infty\})}$ by the graph $\Gamma_{-\theta}$ of a small clockwise rotation of $\CC_{z}$ by $0<\theta\ll \frac{\pi}{2}$. Then the geometric composition of two Lagrangian correspondences $\bL_1$ and $\bL_2$ is realized by $(\bL_1\circ \bL_2)\times \Gamma_{-2\theta}\cong (\bL_1\circ \bL_2)\times \Gamma_{-\theta}$, where $\bL_1\circ \bL_2$ is the geometric composition of $\bL_1\circ\bL_2$ without perturbation and $\Gamma_{-2\theta}\cong \Gamma_{-\theta}$ is the canonical isomorphism.

\sss{The tensor action of $\Loc(X)^{\bd, \otimes}$ on $\cW(X)$}\label{ssstensorLoc}

Let  $\Loc(X)^{\bd, \otimes}$ be the symmetric monoidal category of local systems on $X$ with bounded cohomology stalks, i.e. the cohomology of each stalk should have bounded degrees and should be perfect modules in each degree. 
For any Weinstein sector $X$, we have the natural tensor action by $\Loc(X)^{\bd, \otimes}$ on $\cW(X)$. If we consider $\Ind\cW(X)$ instead of $\cW(X)$, then we have the tensor action of $\Loc(X)^{\otimes}$.
If we identify $\cW(X)$ with wrapped microlocal sheaves on the Lagrangian skeleton of $X$, then the action by $\Loc(X)^{\bd, \otimes}$ is given by tensoring a microlocal sheaf by a local system $\cL\in \Loc(X)^{\bd, \otimes}$.

On $\cW(X)$, the action is given by the monoidal embedding  
\begin{align*}
\iota: \Loc(X)^{\bd, \otimes}&\hookrightarrow \cW_{\Delta_X}(X^-\times X)\\
\cL&\mapsto (\Delta_X, b, \cL)
\end{align*}
where  (1) $\cW_{\Delta_X}(X^-\times X)\subset \cW(X^-\times X)$ is the monoidal full subcategory generated by $\Delta_X$; (2) $(\Delta_X,b, \cL)$ is the object that equip $\Delta_X$ with the standard brane structure and the local system $\cL$. Here to make $\Delta_X$ a well defined object in $\Ind\cW_{\Delta_X}(X^-\times X)$, we replace the sector $X$ by its corresponding stopped Liouville manifold $(\wh{X}, \frf)$, and 
further its stabilization $(\wh{X}, \frf)^\frs$ as in \S\ref{sssGeometricComposition}. It is trivial to see in this case that geometric compositions involving $(\Delta_X,b, \cL)$ agrees with algebraic compositions.

Alternatively, one can define the action of $\Loc(X)^{\bd, \hs;\otimes}$ on $\cW(X)$ by explicit description of the action on objects, morphisms and $A_\infty$-compositions, and this is sufficient for our purpose below. Let $\cL\in \Loc(X)^{\bd, \hs}$. 
On each object $(\sfL, b, \cE)$ where $b$ represents a brane structure on $\sfL$ and $\cE\in \Loc(\sfL)^{\bd, \hs}$, the action by $\cL$ is just tensoring $\cE$ by $\cL|_{\sfL}$. On the Floer cochain complex of two objects $(\sfL_i, b_i, \cE_i), i=1,2$, the action changes each component $\Hom(\cE_1|_{p}, \cE_2|_{p})[-\deg p]$ for $p\in \sfL_1\cap \sfL_2$ to $\Hom(\cE_1|_{p}\otimes \cL|_{p}, \cE_2|_{p}\otimes \cL|_{p})[-\deg p]$. For the differential on the Floer cochain complex, as well as any $A_\infty$-compositions involving ``counting holomorphic polygons",  the action is the obvious one by taking into account of parallel transport of $\cL$ along the edges of any polygon. 
Moreover, for any Hamiltonian isotopy $\phi^t, \phi^0=id, t\in [0,1]$, on $X$, it defines a canonical isomorphism $(\phi^1)^*\cL\cong \cL$ through parallel transport along the paths $\phi^t(p), p\in X, t\in [0,1]$ (and there are canonical isomorphisms between isomorphisms for homotopies between Hamiltonian isotopies and so on, i.e. we have the canonical $\infty$-groupoid of $\Ham(X)$ acting on $\Loc(X)$). Here $\Ham(X)$ is in the same sense as in \cite{GPS1}. 
Then it is straightforward to see that this gives a well defined action of  $\Loc(X)^{\bd, \hs;\otimes}$ on $\cW(X)$.

\sss{Calculation of the action of tensoring by a rank 1 local system on $\cW(\cT_{\Hf})$}\label{appendssstensorLocTHf}

We have described the canonical tensor action of $\Loc(\cT_{H^\flat})^{\bd, \otimes}$ on $\WTHf$ in \S\ref{ssstensorLoc}, which is canonically compatible with restriction and co-restriction functors.

We will denote the tensor action by 
\begin{align*}
\Loc(\cT_{H^\flat})^{\bd}\otimes \WTHf&\lrar\WTHf\\
(\cL, \cG)&\mapsto \cG\otimes_{\CC}\cL. 
\end{align*}
Now let $\cL$ be a rank $1$ local system on $\cT_{\Hf}$ classified by a character
\begin{align}\label{eqchicLHf}
\chi_\cL: \pi_1(\cT_{\Hf}, \yfv)\cong \pi_1(\Hf, 1)\lrar \CC^\times. 
\end{align}
Then we have the direct analog of $\wt{\psi}_\chi^t$ \eqref{eqwtpsichit} and $\tau_\chi^t$ \eqref{eqvarrhoSyfvres} for $\chi_\cL\in \XX^*( \pi_1(\cT_{\Hf}, \yfv))$ in place of $\chi$. Again since $\pi_1(\Tf, 1)\to \pi_1(\cT_{\Hf}, \yfv)$ is $W$-invariant, $\btau_{\chi_\cL}^t$ preserves $\AHf\subset A^\vn$ and induces an automorphism of $\AHf$.

Let $\cL^J=\cL|_{\cT_{L_J^\flat}}$. 
Using the canonical isomorphism $S_\yf\otimes \cL^{-1}\cong S_\yf\otimes\cL^{-1}|_{\yf}$ by trivializing the local system on the Kostant section through parallel transport, we have the canonical isomorphism \emph{as $\CC$-modules}
\begin{align*}
\mu_{\yf}(\cG\otimes_{\CC} \cL)=\Hom(S_{\yf}, \cG\otimes_{\CC} \cL)\cong \Hom(S_{\yf}\otimes \cL^{-1}, \cG)\cong  \Hom(S_{\yf}\otimes \cL^{-1}|_{\yf}, \cG)\cong \mu_{\yf}(\cG)\otimes_{\CC} \cL|_{\yf}. 
\end{align*} 
Similarly, using $S_{\yfv}\otimes \cL^{-1}\cong S_{\yfv}\otimes\cL^{-1}|_{\yfv}$, we have the canonical isomorphism 
\begin{align*}
\mu_{\yfv}(\cG\otimes_{\CC} \cL)\ovs{\Mod(\CC)}{\cong} \mu_{\yfv}(\cG)\otimes_{\CC} \cL|_{\yfv}. 
\end{align*}

Diagram \eqref{diagUpsilonzJflatyJ} is compatible with tensoring with $\cL$ in \emph{the canonical} way:
\begin{equation}\label{diagtensorcLUJflH} 
\begin{tikzcd}[column sep=3em]
\cG\otimes_{\CC} \cL\ar[d, mapsto ]\ar[r, mapsto]&A_{\Hf}\car\mu_{\yf}(\cG\otimes_{\CC} \cL)\cong \mu_{\yf}(\cG)^{\cL}\otimes_{\CC} \cL|_{\yf}\ar[d, mapsto]\\
\res_\vn(\cG\otimes_{\CC}\cL)\cong \res_\vn(\cG)\otimes_{\CC} \cL^\vn\ar[r, mapsto]\ar[dr, mapsto]& A^{\vn}\car\Hom(\res_\vn S_{\yf}, (\res_\vn\cG)\otimes_\CC\cL^\vn)\cong  \mu_{\yf}(\cG)^\cL\otimes_{\AHf} \Avn\otimes_{\CC} \cL|_{\yf}\ar[d, "\sim"', "\varrho^* "]\\
&A^{\vn}\car\mu_{\yfv}(\cG\otimes_{\CC} \cL)\cong \mu_{\yfv}(\cG)^{\cL^\vn}\otimes_{\CC} \cL|_{\yfv},
\end{tikzcd} 
\end{equation}
in which all isomorphisms are canonical (up to a contractible space of choices) and $ \mu_{\yf}(\cG)^{\cL}$ (resp. $\mu_{\yfv}(\cG)^{\cL^\vn}$) is the $A_{\Hf}$ (resp. $\Avn$)-module with the same underlying $\CC$-module as $ \mu_{\yf}(\cG)$ (resp. $\mu_{\yfv}(\cG)$), but with a possibly twisted  $A_{\Hfv}$ (resp. $\Avn$) module structure determined by $\cL$ (resp. $\cL^\vn$), for which we refer to as \emph{$\cL$-twisted} (resp. \emph{$\cL^\vn$-twisted}) modules.

\begin{remark}\label{remvarrhocL}
Note that for any $\cG\in \Ind\cW(\THf)$, diagram \eqref{diagUpsilonzJflatyJ} gives a natural isomorphism 
\begin{align*}
\Avn\car \mu_{\yf}(\cG)\otimes_{\AHf} \Avn\ovs{\sim}{\lrar} \Avn\car \mu_{\yfv}(\cG)
\end{align*}
determined by $\varrho$ \eqref{eqvarrhoSyfvres}. On the other hand, the bottom two rows in diagram \eqref{diagtensorcLUJflH}  give a natural isomorphism
\begin{align*}
\Avn \car  \mu_{\yf}(\cG)^\cL\otimes_{\AHf} \Avn\otimes_{\CC} \cL|_{\yf}   \ovs{\sim}{\lrar}   \Avn\car \mu_{\yfv}(\cG)^{\cL^\vn}\otimes_{\CC} \cL|_{\yfv}
\end{align*}
determined by the composite isomorphism 
\begin{equation}
\begin{tikzcd}[column sep=4em]
(\varrho_{\cL^{-1}})^{-1}: (\res_\vn S_{\yf})\otimes_\CC \cL^{-1}|_{\yf}&(\res_\vn S_{\yf})\otimes_\CC\cL^{\vn,-1}\ar[l,""', "\sim"]\ar[r, "\varrho^{-1}\otimes \id_{\cL^{\vn, -1}}", "\sim"']&S_{\yfv}\otimes_\CC\cL^{\vn,-1}\cong S_{\yfv}\otimes_\CC \cL^{-1}|_{\yfv}.
\end{tikzcd}
\end{equation}
Fixing an isomorphism of lines $\psi: \cL|_{\yf}\ovs{\sim}{\to}\cL|_{\yfv}$, we have 
\begin{align}\label{eqfvarrhopsidual}
\varrho_{\cL^{-1}}\circ (\varrho^{-1}\otimes_\CC \id_{\cL^{-1}|_{\yfv}})\cong f^\cL_\varrho\otimes_{\CC}\psi^\vee:  (\res_\vn S_{\yf})\otimes_\CC \cL^{-1}|_{\yfv}\ovs{\sim}{\lrar} (\res_\vn S_{\yf})\otimes_\CC \cL^{-1}|_{\yf}, 
\end{align}
for some $f^\cL_\varrho\in (\Avn)^{\times}$ and $\psi^\vee$ is the dual isomorphism between the lines. Clearly, we have $f^\cL_\varrho=f^\cL_{\lambda\cdot\varrho}$, for any $\lambda\in \CC^\times$. 
\end{remark}

\sss{The natural isomorphism $\nu_\varrho^\cL$}\label{sssvucLvarrho}
Fixing any $\varrho$ \eqref{eqvarrhoSyfvres}, the (outer) commutative diagram \eqref{diagUpsilonzJflatyJ} identifies the tensor action of $\cL$ on $\Ind\cW(\THf)$ and $\Ind\cW(\cT_{\Lfvn})$ with a natural diagram 
\begin{equation}\label{diagalphachiJvn}
\begin{tikzcd}[column sep=8em, row sep=3em]
\Mod(A_{\Hf})\ar[r, "\alpha^\cL", "\sim"']\ar[d, "\otimes_{\AHf} \Avn"']&\Mod(A_{\Hf})\ar[d, "\otimes_{\AHf} \Avn"]\ar[dl, Rightarrow, "\nu_\varrho^\cL"', "\sim"]\\
\Mod(A^{\vn})\ar[r, "\alpha^{\cL}_\vn"', "\sim"]&\Mod(A^{\vn})
\end{tikzcd},
\end{equation}
in which (1) $\alpha^\cL\cong (-)^\cL\otimes_\CC \cL|_{\yf}$ and $\alpha^{\cL^\vn}\cong (-)^{\cL^\vn}\otimes_\CC \cL|_{\yfv}$ represent the respective automorphisms determined by the tensor action  of $\cL$ conjugated by $\Upsilon_\yf$ and $\Upsilon^\vn$; (2) the natural 2-isomorphism $\nu_\varrho^\cL$ on each $\mu_{\yf}(\cG)$ is explicitly determined by making the following diagram strictly commutative 
\begin{equation}\label{eqnucLvarrrhocL}
\begin{tikzcd}
(\mu_{\yf}(\cG)^{\cL}\otimes_{\AHf} A^\vn)\otimes_\CC \cL|_{\yf}\ar[r, " \varrho_{\cL^{-1}}^* ", "\sim"']\ar[d, "\nu_\varrho^\cL", "\sim"']&\mu_{\yfv}(\cG)^{\cL^\vn} \otimes_\CC \cL|_{\yfv}\cong \mu_{\yfv}(\cG\otimes_{\CC}\cL^\vn)\\
(\mu_{\yf}(\cG)\otimes_{\AHf} \Avn)^{\cL^\vn}\otimes_\CC\cL|_{\yfv}\ar[ur, "\varrho^{*,\cL^\vn}\otimes id_{\cL_{\yfv}}"', "\sim"]&
\end{tikzcd}, 
\end{equation}
in which $\varrho^{*,\cL^\vn}$ denotes the same isomorphism as $\varrho^*$ (on the underlying $\CC$-module) but viewed between the $\cL^\vn$-twisted module structures.  By Remark \ref{remvarrhocL}, $\nu_\varrho^\cL$ on the underlying $\CC$-modules is given by $f^\cL_\varrho\otimes \psi$ as in \eqref{eqfvarrhopsidual}, where $f_\varrho^\cL$ does multiplication on the factor $\Avn$ and $\psi$ identifies the lines. Note that multiplication by $f^\cL_\varrho$ on $\mu_{\yf}(\cG)\otimes_{\AHf} \Avn$ induces an $\Avn$-module automorphism of $(\mu_{\yf}(\cG)\otimes_{\AHf} \Avn)^{\cL^\vn}$. 
Note that for two local systems $\cL_i, i=1,2$, we have the canonical isomorphism between the composite natural isomorphism $\alpha_\vn^{\cL_2,*}(\nu_{\varrho}^{\cL_1})\circ \alpha^{\cL_1,*}(\nu_\varrho^{\cL_2}): \alpha^{\cL_2}\circ\alpha^{\cL_1}\ovs{\sim}{\Rightarrow}  \alpha_\vn^{\cL_2}\circ\alpha_\vn^{\cL_1}$ and the natural isomorphism $\nu_{\varrho}^{\cL_1\cL_2}: \alpha^{\cL_1 \cL_2}\ovs{\sim}{\Rightarrow}  \alpha_\vn^{\cL_1 \cL_2}$, where we write $\cL_1 \cL_2$ for $\cL_1\otimes\cL_2$. This implies the identity
\begin{align}\label{eqreffvarrhocL1L2}
f_\varrho^{\cL_1} f_{\varrho}^{\cL_2}\otimes \psi_{1}\psi_2= f_\varrho^{\cL_1\cL_2}\otimes \psi_{12}\in \Avn \otimes_\CC (\cL_1 \cL_2)^{-1}|_{\yf}\otimes_\CC (\cL_1\cL_2)|_{\yfv}. 
\end{align}

\begin{remark}\label{remfvarrhocL}
First observe that if $\cL$ is a trivial rank $1$ local system on $\cT_{\Hf}$, then by choosing $\psi^\cL: \cL_{\yf}\ovs{\sim}{\to}\cL_{\yfv}$ to be the parallel transport along any homotopy class of paths from $\yf$ to $\yfv$, we get $f_\varrho^\cL=1$. In other words, $\nu_\varrho^\cL$ on the underlying $\CC$-module is $1\otimes \psi^\cL$. 
 
Now for a general $\cL$, assuming $\chi_\cL$ \eqref{eqchicLHf} factors as $\pi_1(\Hf, 1)\twoheadrightarrow \wt{K}\to \CC^\times$, where $\wt{K}$ is a finite group, 
then $\cL^{n}$ is a trivial local system, for some $n>0$. Let $\psi_\gamma^\cL:  \cL_{\yf}\ovs{\sim}{\to}\cL_{\yfv}$ be the parallel transport along a fixed homotopy class of paths $\gamma$ from $\yf$ to $\yfv$, and let $f_{\varrho,\gamma}^{\cL}$ be the corresponding element in $(\AHf)^\times$. 
By \eqref{eqreffvarrhocL1L2}, we have 
\begin{align*}
&(f^{\cL}_{\varrho, \gamma})^n\otimes (\psi_\gamma^\cL)^n=(f^{\cL}_{\varrho, \gamma})^n\otimes\psi^{\cL^n}_{\gamma}=1\otimes\psi^{\cL^n}\\
\Longrightarrow\ &(f^{\cL}_{\varrho, \gamma})^n=1\\
\Longrightarrow\ &f^{\cL}_{\varrho, \gamma}\in \CC^\times, \text{ and }(f^{\cL}_{\varrho, \gamma})^n=1. 
\end{align*}

\end{remark}

\begin{prop}\label{propfvarrhocLCtimes}
For any rank $1$ local system $\cL$ on $\THf$ and any $\psi^{\cL}: \cL_\yf\ovs{\sim}{\to}\cL_{\yfv}$, we have $f_{\varrho}^\cL\in \CC^\times$.  
\end{prop}

\begin{proof}
First, by Remark \ref{remfvarrhocL}, the statement holds for $\Hf$ semisimple. 

Second, let us check the statement holds for $\Hf=T'$ where $T'$ is a connected torus. In this case, let $z_1=\yf\in T'$ and $z_0=1=\yfv\in T'$, then for any $\varrho: S_{z_0}\ovs{\sim}{\to} S_{z_1}$, we need to find the composition
\begin{align*}
S_{z_1}\otimes\cL^{-1}_{z_0}\ovs{\varrho^{-1}\otimes \id_{\cL^{-1}|_{z_0}}}{\lrar} S_{z_0}\otimes \cL^{-1}_{z_0}\ovs{\sim}{\to} S_{z_0}\otimes \cL^{-1}\ovs{\varrho\otimes \id_{\cL^{-1}}}{\lrar}S_{z_1}\otimes\cL^{-1}\ovs{\sim}{\to}S_{z_1}\otimes\cL^{-1}_{z_1}.
\end{align*}
It is standard that $\varrho=\lambda\cdot \gamma$, where $\gamma$ is an element of the canonical basis of $\Hom(S_{z_0}, S_{z_1})$ represented by a homotopy class of paths $\gamma$ from $z_0$ to $z_1$ (cf. \S\ref{sssconcretepiZH}) and $\lambda\in \CC^\times$. Then it is direct to see that the above composition is just $\id_{S_{z_1}}\otimes \psi_{\gamma^{-1}}^\vee$, where $\psi_{\gamma^{-1}}: \cL_{z_1}\to \cL_{z_0}$ is the parallel transport along the inverse path $\gamma^{-1}$. 

Combining the above two situations, we see the statement holds for $\Hf=H^{\ad}\times H_\ab$. 

For general $\Hf=H^{\flat, \der}\times^{Z(H^{\flat, \der})}Z(\Hf)$, let $\Hff=H^{\ad}\times H_\ab$, and let $\pi: \cT_{\Hf}\to \cT_{\Hff}$ be the projection. It suffices to check that for any $\cL'$ on $\cT_{\Hff}$, the statement holds for $\cL=\pi^*\cL'$ on $\THf$. To see this, we use the compatibility of $F_*: \cW(\THf)\to \cW(\cT_{\Hff})$ with the tensor action of $\cL$ on the top and $\cL'$ on the bottom. Let $\yff$ and $\yffv$ be the image of $\yf$ and $\yfv$, respectively. Then we have the canonical identifications $\cL'_{\yff}\cong \cL_\yf$ and  $\cL'_{\yffv}\cong \cL_{\yfv}$. For any $\psi^\cL:  \cL_{\yf}\ovs{\sim}{\to}\cL_{\yfv}$, let $\psi^{\cL'}:  \cL'_{\yff}\ovs{\sim}{\to}\cL'_{\yffv}$ be the 
induced isomorphism. Let $\varrho'$ is the image of $\varrho$ under $F_*$. 

Using
\begin{align*}
H^0(\Hom_{\cW(\cT_{\Lfvn})}(\res_\vn S_\yf\otimes \cL_{\yf}^{-1}, \res_\vn S_\yf\otimes \cL_{\yfv}^{-1}))\hookrightarrow H^0(\Hom_{\cW(\cT_{\Lffvn})}(\res_\vn S_{\yff}\otimes \cL_{\yff}^{\prime, -1}, \res_\vn S_{\yff}\otimes \cL_{\yffv}^{\prime, -1})),
\end{align*}
and that the map sends $f_{\varrho}^{\cL}\otimes \psi^{\cL}$ to $f_{\varrho'}^{\cL'}\otimes \psi^{\cL'}$. It then follows that $f_{\varrho}^{\cL}\in \CC^\times$. This finishes the proof. 
\end{proof}

For any $\Avn$ (resp. $\AHf$) module $M_\vn$ (resp. $M$), let $M_\vn^{\chi_\cL}$ (resp. $M^{\chi_\cL}$) be the $\chi_\cL$-twisted $\Avn$ (resp. $\AHf$) module defined using $\btau_{\chi_\cL}^t$ in place of $\btau_\chi$ \eqref{eqavarphichia}. 

\begin{remark}\label{remgeneralizetoresHfLf}
All the discussions from \S\ref{appendssstensorLocTHf} till Proposition \ref{propfvarrhocLCtimes} generalize verbatim to restriction functors $\res^{\Hf}_{\Lf}: \IndW(\cT_{\Hf})\lrar \IndW(\cT_{\Lf})$, for any Levi subgroup $\Lf\subset \Hf$ containing $\Tf$, where one replaces $\Tf$ in the aforementioned subsections by $\Lf$. 
In particular, for any standard Levi $L\subset H$ and $\varrho: S_{y_L^{\flat}}\ovs{\sim}{\to}\res^{\Hf}_{\Lf} S_{y_H^\flat}$, we have a natural diagram 
\begin{equation}\label{diagalphachiJvn}
\begin{tikzcd}[column sep=8em, row sep=3em]
\Mod(A_{\Hf})\ar[r, "\alpha^\cL", "\sim"']\ar[d, "\otimes_{\AHf} A_{\Lf}"']&\Mod(A_{\Hf})\ar[d, "\otimes_{\AHf} A_{\Lf}"]\ar[dl, Rightarrow, "\nu_\varrho^\cL"', "\sim"]\\
\Mod(A_{\Lf})\ar[r, "\alpha^{\cL}_{\Lf}"', "\sim"]&\Mod(A_{\Lf})
\end{tikzcd},
\end{equation}
where $\nu_\varrho^\cL\in \cL^{-1}_{\yHf}\otimes\cL_{\yLf}$. This defines a map
\begin{align}\label{eqnuupcL}
\nu^\cL: \Isom_{\cW(\cT_{\Lf})}(S_{\yLf}, \res^{\Hf}_{\Lf} S_{\yHf})/\CC^\times\lrar  \cL^{-1}_{\yHf}\otimes\cL_{\yLf},\quad \varrho\mapsto \nu_\varrho^\cL,
\end{align}
that is multiplicative in $\cL$. Note that $\nu^\cL$ is also functorial in $(y_H^\flat, y_L^\flat; \varrho)$, cf. Remark \ref{remnuLcompatiblevarrho}.
The reason that we specialize $\Lf$ to $\Tf$ is to prove the following lemma. 
\end{remark}

\begin{lemma}\label{lemmaMchicLyf}
Let $\cL$ be a rank $1$ local system on $\cT_\Hf$. 
\begin{itemize}
\item[(i)]
The $A_{\Hf}$-module (resp. $A^{\vn}$-module) structure on $M^{\cL}\otimes_{\CC} \cL_{\yf}$ (resp. $M_\vn^{\cL^\vn}\otimes_{\CC}\cL_{\yfv}$) is canonically identified with $M^{\chi_\cL}\otimes_{\CC}  \cL_{\yf}$ (resp. $M^{\chi_\cL}\otimes_{\CC}\cL_{\yfv}$).

\item[(ii)] The automorphism $\alpha^\cL$ is \emph{canonically} identified with 
\begin{align*}
M\mapsto M^{\chi_\cL}\otimes_\CC \cL_{\yf},
\end{align*}
which is independent of the choice of $\varrho$ \eqref{eqvarrhoSyfvres}.  
\end{itemize}

\end{lemma}

\begin{proof}
(i)  First, on $\THf$, since  (1) $\cL$ is the 1-dimensional module $\cL_{\yfv}$ of $\Avn$ on which acts by $\wt{\psi}_{\chi_\cL}^t$, and (2) tensoring with $\cL$ on $\Ind\cW(\THf)$ is canonically identified with convolution with the corresponding 1-dimensional $\Avn$-module on $\Mod(\Avn)$,  it is clear that we have a canonical identification $\alpha^\cL_\vn\cong (-)^{\chi_\cL}\otimes_\CC \cL_{\yfv}$. 

Second, using the identification of $\nu_\varrho^\cL$ with $f_\varrho\otimes_\CC\psi$ \S\ref{sssvucLvarrho}
we see directly that the $\AHf$-module structure on $\mu_{\yf}(\cG)^{\cL}$ is the $\chi_\cL$-twisted one, which is independent of the choice of $\varrho$. 

(ii) follows directly from (i). 
\end{proof}

\subsection{The canonical $\pi Z(H)$-action on $\cW(\cT_H)$}\label{subsecpiZHactionW}

There are several equivalent models for describing the canonical $\pi Z(H)$-action on $\cW(\cT_H)$. Since $Z(H)$ is a $K(\pi, 1)$-space, it is direct to check that different models match. On microlocal sheaves, the $\pi Z(H)$-action is given by the natural $\pi Z(H)_{\cpt}$-action on the Lagrangian skeleton (here we have an easy situation that $Z(H)_\cpt$ preserves the Lagrangian skeleton).

\sss{Definition(s) of the canonical $\pi Z(H)$-action on $\cW(\cT_H)$}\label{sssDefpiZHaction}
Write $\cT_H=T^*Z(H)\times^{Z(H^\der)}\cT_{H^\der}$, and let 
\begin{align}\label{eqcTHprime}
\cT'_H=T^*Z(H)_{\cpt}\times^{Z(H^\der)}\cT_{H^\der} 
\end{align}
in which 
the sector factor $T^*Z(H)$ is replaced by the Liouville manifold $T^*Z(H)_{\cpt}$. Note $Z(H)_\cpt=Z(H_\cpt)$. 
We have the canonical equivalence $\cW(\cT_H)\simeq \cW(\cT'_H)$ up to a contractible space of choices.  

It is well known that $\cW(T^*Z(H)_{\cpt})\simeq \Loc^w(Z(H)_{\cpt})^{op}$ has a canonical convolution symmetric monoidal structure. 
One defines the canonical $\pi Z(H)_{\cpt}$-action on $\cW(T^*Z(H)_{\cpt})$ by either of the following equivalent ways:
\begin{itemize}
\item[(1)] Since $\Loc(Z(H)_{\cpt})\simeq \Fun(\pi Z(H)_{\cpt}, \Mod(\CC))$, $\pi Z(H)_{\cpt}^{op}$ naturally acts on $\Loc(Z(H)_{\cpt})$ and clearly preserves compact objects;
\item[(2)] There is a canonical  symmetric monoidal functor of $\infty$-groupoids (which are equivalent to usual groupoids in this case)
\begin{align}\label{eqiotaHpiZH}
\iota_H: \pi Z(H)_{\cpt}\lrar \cW(T^*Z(H)_{\cpt};\ZZ)^{\Gpd}. 
\end{align}
Indeed, we have the canonical equivalence $\cW(T^*Z(H)_{\cpt};\ZZ)^{\Gpd}\simeq \wt{\pi Z(H)}_{\cpt}$, where $ \wt{\pi Z(H)}_{\cpt}$ is the groupoid, which has the same objects as $\pi Z(H)_{\cpt}$ and whose morphism set between $x,y\in Z(H)_\cpt$ consists of $([\rho], \varepsilon)$, where $[\rho]$ is a homotopy class of paths from $x$ to $y$ in $Z(H)_\cpt$ and $\varepsilon\in \mu_2$. The composition of morphisms is the natural one (concatenation of paths and multiplication on the $\mu_2$-factors). Thus we have the canonical $\iota$ which is the identity on objects and sends each morphism $[\rho]$ to $([\rho], 1)$ (cf. \S\ref{sssconcretepiZH} for a more concrete description). Then $ \pi Z(H)_{\cpt}$ acts on $\cW(T^*Z(H)_{\cpt};\ZZ)$ and hence on $\cW(T^*Z(H)_{\cpt})$. 
\end{itemize}

Now we have the canonical $\pi Z(H)_{\cpt}$-action on $\wt{\cT}'_H=T^*Z(H)_{\cpt}\times \cT_{H^\der}$. To define the sought-for $\pi Z(H)_{\cpt}$-action on $\Ind\cW(\cT'_{H})$, it suffices to check that the $\pi Z(H)_{\cpt}$-action and the $Z(H^\der)$-action on $\Ind\cW(\wt{\cT}'_H)$ naturally commute, i.e. we have a canonical $\pi Z(H)_\cpt\times Z(H^\der)$-action on  $\Ind\cW(\wt{\cT}'_H)$. But this is obvious from the following
\begin{itemize}
\item[(a)]
the  $Z(H^\der)$-action on $\Ind\cW(\wt{\cT}'_H)\simeq \Ind\cW(T^*Z(H)_{\cpt})\otimes \Ind\cW(\cT_{H^\der})$ is given by the diagonal one from the natural $Z(H^\der)$-action $\Ind\cW(T^*Z(H)_{\cpt})$ and on $\Ind\cW(\cT_{H^\der})$;

\item[(b)] the $\pi Z(H)_{\cpt}$-action $ \Ind\cW(T^*Z(H)_{\cpt})$ is from the left, while the $Z(H^\der)$-action on $\Ind\cW(T^*Z(H)_{\cpt})$ is from the right, thus they naturally commute.  
\end{itemize}

The $\pi Z(H)_\cpt$-action on $\cW(\cT_H')$ is  independent of the choice of finite central covering of $H$ of the form $Z(H)_\cpt \times H_1$, where $H_1$ is (connected) semismiple. Indeed, $H_1$ belongs to the finite collection of finite central covering of $H$ and we have the unique ``maximal" element $H^{sc}$ and the natural adjunction between $\cW(T^*Z(H_{\cpt})\times \cT_{H^{sc}}')$ and $\cW(T^*Z(H_{\cpt})\times \cT_{H_1}')$ are obviously $\pi Z(H)_{\cpt}$-equivariant using the Kunneth formula. 

\sss{A more concrete description using wrapping}\label{sssconcretepiZH}

Fix an Euclidean inner product on $\frz_\cpt=\Lie Z(H)^\circ_\cpt$, and let $\|-\|^2$ be the associated length squared function. 
This defines an invariant flat metric on $Z(H)_\cpt$ and an identification $\frz_\cpt^*\cong \frz_\cpt$. 
Let $c: \RR\to \RR$ be a smooth function, which satisfies 
\begin{align*}
c(x)=\begin{cases}
x,& x\in (-\ep, \ep);\\
\sqrt{x}, &x\geq 1
\end{cases},\quad c'(x)>0\text{ for }x>0,
\end{align*}
where $\ep>0$ is sufficiently small. 
Let $\cH_Z$ be the positive Hamiltonian  on $T^*Z(H)_\cpt$ from pulling back the function $c(\|-\|^2): \frz_\cpt^*\cong \frz_\cpt\lrar \RR_{\geqsl 0}$ under the natural projection $T^*Z(H)_\cpt\to \frz_{\cpt}^*$.

Let $\sfF_{z}$ be the cotangent fiber at $z$. 
For any $z_1, z_2\in Z(H)_\cpt$, using the wrapping by $\cH_Z$, we have $\Hom_{\cW(T^*Z(H)_\cpt)}(\sfF_{z_1}, \sfF_{z_2})$ concentrated in degree $0$, and we 
can canonically identify a $\ZZ$-basis in it as the set of geodesics $\pathz_{1,2}$ (just as oriented curves)
starting from $z_1$ and ending at $z_2$. We will equally view $\pathz_{1,2}$ as the projection of any of the Hamiltonian flow lines $\varphi_{\cH_Z}^{t}(z_1, \xi), \xi\in \frz_{\cpt},  t\in [0, t_0]$, satisfying $\varphi_{\cH_Z}^{t_0}(z_1,\xi)=(z_2, \xi)$ and $\gamma_{1,2}\sim \proj_{Z(H)_\cpt}(\varphi_{\cH_Z}^{t}(z_1, \xi))$, relative to $z_1, z_2$ (clearly such an $\xi$ is determined up to scaling by $\RR_{>0}$). Note the geodesics $\pathz_{1,2}$ are independent of the choice of Euclidean inner products on $\frz_\cpt$ (as they are straight lines in the universal cover $\frz_\cpt$ of each connected component of $\ZHc$). 

Using the natural inclusion $\pr_{H}: U_{H}^{\SS}\hookrightarrow H$, let $U_H^{\SS, \cpt}$ be the subspace of $U_{H}^{\SS}$ consisting of $y\in (\ZHc H^\der)\cap U_{H}^{\SS}=Z(H)_\cpt U_{H^\der}^\SS$. Note that  $U_H^{\SS, \cpt}$ is a $\ZHc$-torsor. 
We have the ``Kostant section" $S'_{y}=\sfF_{z}\times S^\der_{y'}$ in $\cT_H'$, where $y=zy'$ for some $z\in Z(H)_\cpt$, $y'\in U_{H^\der}^{\SS}$, and they generate $\cW(\cT_H')$. Choose a  linear (and strictly positive) near $\infty$ function on $\frc^\der$ that has a non-degenerate global minimum at $[0]$ of value $0$,  and let $\cH_\der$ be the pullback function on $\cT_{H^\der}$ (cf. \cite{J}), then 
the constant Hamiltonian path at $(y', [0])\in S^\der_{y'}\cong \frc^\der$ represents the identity morphism of $S^\der_{y'}$ for all time.

Using the product Hamiltonian function $\cH=\cH_Z\times \cH^\der$, that obviously descends to $\cT'_H$, we have 
\begin{align*}
\Hom_{\WTHp}(S'_{y_1}, S'_{y_2})\cong \colim_{t\to \infty}(\CF^*(\varphi_\cH^t(S'_{y_1}), S'_{y_2}), d).
\end{align*}
Since the Floer cochain $\CF^*(\varphi_\cH^t(S'_{y_1}), S'_{y_2})$ has generators given by Hamiltonian paths connecting $S'_{y_1}$ and $S'_{y_2}$ under $\varphi_\cH^{s}, s\in [0, t]$, we will use a cocyle in $\CF^*(\varphi_\cH^t(S'_{y_1}), S'_{y_2})$ to denote its image morphism between $S'_{y_1}$ and $S'_{y_2}$. In fact, by choosing the Hamiltonian appropriately as in \cite{J}, $\CF^*(\varphi_\cH^t(S'_{y_1}), S'_{y_2})$ are all concentrated in degree $0$. 

The $\pi\ZHc$-action on these generators are explicitly given as follows
\begin{itemize}
\item[(i)] for any $z\in \ZHc$, it sends $S'_{y}$ to $S_{zy}'$;

\item[(ii)] for any $\pathz_{1,2}\in \Hom_{\cW(T^*Z(H)_\cpt)}(\sfF_{z_1}, \sfF_{z_2})$, represented by a Hamiltonian flow line $\varphi_{\cH_Z}^{t}(z_1, \xi), \xi\in \frz_{\cpt},  t\in [0, t_0]$ as above, 
the induced isomorphism $\sigma^H_y(\pathz_{1,2})\in \Hom_{\WTHp}(S'_{z_1y}, S'_{z_2y})$ 
is represented by the product of the Hamiltonian flow line $\varphi_\cH^{t}(zz_1, \xi), \xi\in \frz_{\cpt},  t\in [0, t_0]$ (clearly $\varphi_\cH^{t_0}(zz_1, \xi)=(zz_2, \xi)$) and the constant Hamiltonian flow line at $(y', [0])$. Clearly, $\sigma^H_y(\pathz_{1,2})$ represents a degree $0$ invertible morphism, and this gives an embedding
\begin{align}\label{eqgamma12sigmayH}
\{\pathz_{1,2}\in \Hom_{\cW(T^*Z(H)_\cpt)}(\sfF_{z_1}, \sfF_{z_2})\}&\hookrightarrow H^0(\Hom_{\WTHp}(S'_{z_1y}, S'_{z_2y}))\cong \Hom_{\WTHp}(S'_{z_1y}, S'_{z_2y})\\
\nonumber \gamma_{1,2}&\mapsto \sigma^H_y(\pathz_{1,2}),
\end{align}
which is clearly compatible with compositions. 
\end{itemize}

\begin{remark}\label{remsigmayHindpresentation}
Note that since the Hamiltonian $\cH$ is $\ZHc$-invariant, for any $y_1, y_2\in\UHSc$ and  any two different presentations $y_j=z_j^{(i)}y^{(i)}, 1\leq i,j\leq 2$ (where $z_j^{(i)}\in \ZHc$ and $y^{(i)}\in \UHSc$), the canonical identification 
\begin{align*}
\Hom_{\cW(T^*Z(H)_\cpt)}(\sfF_{z_1^{(1)}}, \sfF_{z_2^{(1)}})&\ovs{\sim}{\lrar} \Hom_{\cW(T^*Z(H)_\cpt)}(\sfF_{z_1^{(2)}}, \sfF_{z_2^{(2)}})\\
 \gaOT^{(1)}&\mapsto \gaOT^{(2)}:= \delta\cdot \gaOT^{(1)}
\end{align*}
by the translation of $\delta=(z_1^{(1)})^{-1}z_1^{(2)}$, identifies the morphism $\sigma_{y^{(1)}}^H(\gaOT^{(1)})$ with $\sigma_{y^{(2)}}^H(\gaOT^{(2)})$ in $\Hom_{\WTHp}(S'_{y_1}, S'_{y_2})$. In particular, the collection of isomorphisms in 
$$H^0(\Hom_{\WTHp}(S'_{y_1}, S'_{y_2}))\cong  \Hom_{\WTHp}(S'_{y_1}, S'_{y_2}), \quad y_i=z_iy$$
that are contained in the image of  \eqref{eqgamma12sigmayH}
are independent of the choice of presentations $y_i=z_iy$; for different presentations, they  are canonically identified.  
\end{remark}

\sss{Compatibility with  co-restriction and restriction functors}\label{ssspiZHcompatible}

For any $J\in \cP(I_H)$, we have the $Z(L_J)_\cpt$-action on the subsector $\cT_{L_J}^{\ovl{\cU}_J}$. Since the sector inclusion $\cT_{L_J}^{\ovl{\cU}_J}\hookrightarrow \cT_{H}^{\ovl{\cV}}$ is $Z(H)_\cpt$-equivariant (where $Z(H)_\cpt$ acts on the former through $Z(H)_\cpt\hookrightarrow Z(L_J)_\cpt$), it is clear from the definition of the $\pi Z(H)_\cpt$-action \S\ref{sssDefpiZHaction} that 
the adjunction of co-restriction and restriction between their ($\Ind$-)wrapped Fukaya categories are both naturally $\pi Z(H)_\cpt$-equivariant.

\bibliographystyle{abbrv}
\bibliography{MG.bib}

\end{document}